\title{Edge-transitive maps}

\author{Gareth A. Jones\\
School of Mathematics\\
University of Southampton\\
Southampton SO17  1BJ, UK\\
{\tt G.A.Jones@maths.soton.ac.uk}\\
}

\documentclass[12pt]{article}
\usepackage[latin1]{inputenc}
\usepackage{latexsym}
\usepackage{amsmath}
\usepackage{amssymb}
\usepackage{amsfonts}
\usepackage{color}
\usepackage{tikz}
\usepackage{enumerate}
\newtheorem{thm}{Theorem}[section]
\newtheorem{lemma}[thm]{Lemma}
\newtheorem{cor}[thm]{Corollary}
\newtheorem{prop}[thm]{Proposition}

\newcommand{\hyp}{\mathbb{H}}

\newcommand{\C}{\mathbb{C}}
\newcommand{\Z}{\mathbb{Z}}
\newcommand{\M}{\mathcal{M}}
\newcommand{\F}{\mathbb{F}}
\newcommand{\G}{\mathcal{G}}
\newcommand{\N}{\mathbb{N}}

\newcommand{\K}{\mathcal{K}}

\date{}
\begin{document} 

%\tableofcontents

\newpage

\maketitle

\begin{abstract}
For each of the $14$ classes of edge-transitive maps described by Graver and Watkins, necessary and sufficient conditions are given for a group to be the automorphism group of a map, or of an orientable map without boundary, in that class. Extending earlier results of \v Sir\'a\v n, Tucker and Watkins, these are used to determine which symmetric groups $S_n$ can arise in this way for each class. Similar results are obtained for all finite simple groups, building on work of Leemans and Liebeck, Nuzhin and others on generating sets for such groups. It is also shown that each edge-transitive class realises finite groups of every sufficiently large nilpotence class or derived length, and also realises uncountably many non-isomorphic infinite groups.
\end{abstract}

\medskip

%30F10, % compact Riemann surfaces
%30F50 % Klein surfaces
 (secondary).
 
 \medskip
 
 
 \newpage
 
 \tableofcontents
 
 \newpage
 
 %%%%%%%%%%%%%%%%%%%
 
 \part{Preliminaries}\label{prelims}

\section{Introduction}\label{intro}

Two of the most important and deeply studied classes of maps on surfaces are those consisting of the regular maps and of the orientably regular chiral maps (here simply called chiral maps for brevity). The regular maps $\M$ are those for which the automorphism group ${\rm Aut}\,\M$ acts transitively on flags,  while the chiral maps are those orientable maps for which ${\rm Aut}\,\M$ is transitive on arcs but not flags, so that $\M$ is not isomorphic to its mirror image.

Such maps are all edge-transitive. In 1997 Graver and Watkins~\cite{GW} partitioned edge-transitive maps $\M$ into $14$ classes, distinguished by the isomorphism class of the quotient map ${\mathcal M}/{\rm Aut}\,{\mathcal M}$ (see \S\ref{14classes} for a summary of this classification); in that year, Wilson gave a similar classification in~\cite{Wil97}. These classes $T$ correspond bijectively to the $14$ isomorphism classes of maps ${\mathcal N}(T)$ with one edge; they include class $1$, consisting of the regular maps, and class $2^P{\rm ex}$, consisting of the chiral maps, together with others such as class~3, consisting of the just-edge-transitive maps, those for which ${\rm Aut}\,\M$ is transitive on edges but not on vertices or faces. The classes are listed in the first column of Table~\ref{GpsRealised} below; here we simply note that the duals of the maps in class $2$ form class $2^*$, while the Petrie duals of the latter form class $2^P$; similar remarks apply to the classes $2\,{\rm ex}, 4$ and $5$, while the classes $1$ and $3$ are invariant under these two operations (see \S\ref{operations} for details).

After Part~\ref{prelims}, consisting of this Introduction and a brief outline of the general method, the basic concepts and techniques required to study these classes are covered in Part~\ref{basics} of this text. As an illustration, the edge-transitive embeddings of complete graphs are classified in Section~\ref{completemaps}, building on earlier work of Biggs, James, Wilson and the author~\cite{Big, Jam83, Jam90, JJ, Wil89}. Part~\ref{autgps} is devoted to studying, for each class $T$, the set $\G(T)$ of groups which can be realised as the automorphism group of a map in $T$. It is widely known, and not hard to prove, that apart from a few small exceptions each finite symmetric or alternating group can be realised as the automorphism group of a regular map, and that the same applies to chiral maps. In 2001 \v Sir\'a\v n, Tucker and Watkins~\cite{STW} showed that for each integer $n\ge 11$ with $n\equiv 3$ or $11$ mod~$(12)$, there are finite, orientable, edge-transitive maps $\mathcal M$ in each of the $14$ classes, with ${\rm Aut}\,{\mathcal M}$ isomorphic to the symmetric group $S_n$. The first aim of this paper is to extend their result by determining, for each class $T$, all integers $n$ such that $S_n\in\G(T)$. This result, together with similar results for the alternating groups $A_n$ and the projective special linear groups $L_2(q)=PSL_2(q)$, is summarised in Theorem~\ref{mainthm} and Table~\ref{GpsRealised}:

\begin{thm}\label{mainthm}
A symmetric group $S_n$, an alternating group $A_n$, or a projective special linear group $L_2(q)$ is isomorphic to the automorphism group of an edge-transitive map in a class $T$ if and only if it satisfies the corresponding condition in Table~\ref{GpsRealised}.
\end{thm}

\begin{table}[ht]
\centering
\begin{tabular}{| p{2.9cm} | p{1.4cm} | p{3.7cm} | p{2.1cm} |}
\hline
Class $T$ & $S_n$ & $A_n$ & $L_2(q)$ \\
\hline\hline
$1$ & $n\ge 1$ & $n=1, 2, 5$ or $n\ge 9$ & $q\ne 3, 7, 9$\\
\hline
$2$, $2^*$, $2^P$ & $n\ge 2$ & $n\ge 5$ & $q\ne 3$ \\
\hline
$2\,{\rm ex}$, $2^*{\rm ex}$, $2^P{\rm ex}$ & $n\ge 6$ & $n\ge 8$ & no $q$ \\
\hline
$3$ & $n\ge 2$ & $n\ge 5$ & $q\ne 3$ \\
\hline
$4$, $4^*$, $4^P$ & $n\ge 2$ & $n\ge 4$ & every $q$ \\
\hline
$5$, $5^*$, $5^P$ & $n\ge 6$ & $n\ge 7$ & no $q$ \\
\hline

\end{tabular}
\caption{Groups $S_n$, $A_n$ and $L_2(q)$ in sets $\G(T)$.}
\label{GpsRealised}
\end{table}

In particular, the groups $S_n$ for $n\ge 6$ and $A_n$ for $n\ge 9$ are represented in all of these classes. In the case $T=1$, the class of regular maps, Sjerve and Cherkassoff dealt with these groups together with $PGL_2(q)$ in~\cite{SC}, while Nuzhin dealt with the alternating groups in~\cite{Nuz92} and the simple groups of Lie type, including $L_2(q)$, in~\cite{Nuz90, Nuz97a, Nuz97}.)  In most cases, all but finitely many of these groups are realised in each class $T$, the exceptions being six of the classes for which no groups $L_2(q)$ arise. The method of proof, both here and for other results stated below, is to follow~\cite{STW} in using necessary and sufficient conditions for a group to be in the various sets $\G(T)$: for instance, when $T=1$ these require the group to have generators $r_i\;(i=0, 1, 2)$ satisfying $r_i^2=(r_0r_2)^2=1$.
In the cases $T\ne 1$ there are similar conditions on generators, together with the requirement that the group should not have certain `forbidden automorphisms' (see \S\ref{parentforbidden} for details). It is then a routine matter to apply these conditions to the groups in Theorem~\ref{mainthm}.

The groups $A_n$ are simple for all $n\ge 5$, as are the groups $L_2(q)$ for all prime powers $q\ge 4$. More generally, it is of interest to determine, for each class $T$, which non-abelian finite simple groups are in $\G(T)$. In fact, it is easier to list those which are not, as in Table~\ref{FSGpsRealised} (where we use ATLAS notation~\cite{ATLAS} for simple groups):

\begin{thm}\label{mainthmsimple}
A non-abelian finite simple group is isomorphic to the automorphism group of an edge-transitive map in a class $T$ if and only if it is not one of the exceptions listed in the corresponding row of Table~\ref{FSGpsRealised}.
\end{thm}

\begin{table}[ht]
\centering
\begin{tabular}{| p{2.6cm} | p{7.9cm}|}
\hline
Class $T$ & Non-abelian finite simple groups $G\not\in\G(T)$   \\
\hline\hline
$1$ & $L_3(q), U_3(q), L_4(2^e), U_4(2^e), U_4(3), U_5(2)$,\\
& $A_6, A_7, M_{11}, M_{22}, M_{23}, McL$  \\
\hline
$2, 2^*, 2^P$  & $U_3(3)$  \\
\hline
$2\,{\rm ex}, 2^*{\rm ex}, 2^P{\rm ex}$  &$L_2(q), L_3(q), U_3(q), A_7$  \\
\hline
$3$ & --  \\
\hline
$4, 4^*,4^P$ & --  \\
\hline
$5, 5^*, 5^P$ & $L_2(q)$  \\
\hline

\end{tabular}
\caption{Non-abelian finite simple groups not in sets  $\G(T)$.}
\label{FSGpsRealised}
\end{table}

% $U_4(3)$ and $U_5(2)$ added to the list for $T=1$.

Since there are no exceptions for class $3$, we immediately have:

\begin{cor}
Every non-abelian finite simple group is isomorphic to the automorphism group of an edge-transitive map.
\end{cor}

Here $M_n$ is the Mathieu group of degree $n$, while $McL$ is the McLaughlin group. In these lists, all integers $e\ge 1$ are allowed, as are all prime powers $q$ provided the corresponding group is simple. Note that the exceptions listed in the first row include the groups $L_2(7)\cong L_3(2)$, $L_2(9)\cong A_6$ and  $A_8\cong L_4(2)$, while those in the third row include $A_5\cong L_2(4)\cong L_2(5)$ and $A_6\cong L_2(9)$. The smallest simple group represented in all classes is the Suzuki group $Sz(8)$, followed by $M_{12}$, the Janko group $J_1$ and $A_9$.

The entry for $T=1$, the class of regular maps, is due to Nuzhin and others, through their answer to Mazurov's (purely algebraic) question~\cite[Problem 7.30]{Kou} asking which non-abelian finite simple groups can be generated by three involutions, two of them commuting (see Section~\ref{gensimple} for details, including the addition of $U_4(3)$ and $U_5(2)$ to previously published lists). The solution for $T=2^P{\rm ex}$, the class of chiral maps, has recently been determined by Leemans and Liebeck~\cite{LL, LL2} in the equivalent context of abstract polyhedra (again, see Section~\ref{gensimple}), and a simple argument using map dualities extends their result to the classes $2\,{\rm ex}$ and $2^*{\rm ex}$. It follows that $M_{11}$ is the smallest simple group in $\G(T)$ for these three classes. The entries for the remaining ten classes are apparently new, and the proofs are given in Section~\ref{gensimple}.

The exceptions for these ten classes are easily explained. The unitary group $U_3(3)$ is not in $\G(T)$ for the classes $T=2$, $2^*$ and $2^P$ since groups realised in such classes must be generated by at most three involutions, and Wagner~\cite{Wag} has shown that this group requires four. By an observation of Singerman~\cite{Sin}, for each generating pair for $L_2(q)$ there is an automorphism inverting both generators; such an automorphism is forbidden for the classes $T=5$, $5^*$ and $5^P$, so $L_2(q)\not\in\G(T)$; the smallest simple group arising for these classes is $A_7$. The exceptions for $T=2\,{\rm ex}$ and $2^*{\rm ex}$ are the same as those found in~\cite{LL} for $T=2^P{\rm ex}$.

The main part of the proofs of these theorems consists of showing that various groups are realised (as an automorphism group) in particular classes. Simple arguments show that if any group is realised in class $1$ or $2^P{\rm ex}$ then it is also realised in various other classes (see Lemma~\ref{regmapslemma} for a precise statement), so this allows one to concentrate on those groups which are not in $\G(1)$ or not in $\G(2^P{\rm ex})$ (in the case of finite simple groups, these are the exceptions in the first and third rows of Table~\ref{FSGpsRealised}). In order to realise such groups, more direct arguments are required, finding specific generators and then showing that these do not admit forbidden automorphisms.

In this text we also consider the set $\G^+(T)$ of groups which can be `evenly realised', that is, as the automorphism group of an orientable map without boundary, in each of the classes $T$. If $T=2^P{\rm ex}$, $5$ or $5^*$ then all maps in that class have these properties, so $\G^+(T)=\G(T)$. For the other eleven classes, each group in $\G^+(T)$ must have a subgroup of index $2$, so in particular no simple groups (other than $C_2$ for $T=1$) can be evenly realised in such a class $T$. For example, it follows from Theorem~\ref{mainthm} that no group $L_2(q)$ is evenly realised by an edge-transitive map. In the case of the groups $S_n$ and $A_n$ we have the following:

\begin{thm}\label{mainthmeven}
A symmetric group $S_n$ or an alternating group $A_n$ is isomorphic to the automorphism group of an edge-transitive orientable map without boundary in a class $T$ if and only if it satisfies the corresponding condition in Table~\ref{GpsEvenlyRealised}.
\end{thm}

\begin{table}[ht]
\centering
\begin{tabular}{| p{1.8cm} | p{2.4cm} | p{2.2cm} |}
\hline
Class $T$ & $S_n\in\G^+(T)$ & $A_n\in\G^+(T)$ \\
\hline\hline
$1$ & $n\ne 1, 5, 6$ & no $n$ \\
\hline
$2$, $2^*$& $n\ne 1, 2, 5, 6$ & no $n$  \\
\hline
$2^P$ & $n\ge 3$ & no $n$ \\
\hline
$2\,{\rm ex}$, $2^*{\rm ex}$ & $n\ge 7$ & no $n$ \\
\hline
$2^P{\rm ex}$ & $n\ge 6$ & $n\ge 8$ \\
\hline
$3$ & $n\ge 3$ & no $n$  \\ 
\hline
$4$, $4^*$, $4^P$ & $n\ge 3$ & no $n$  \\
\hline
$5$, $5^*$ & $n\ge 6$ & $n\ge 7$  \\
\hline
$5^P$ & $n\ge 6$ & no $n$  \\
\hline

\end{tabular}
\caption{Groups $S_n$ and $A_n$ in sets $\G^+(T)$.}
\label{GpsEvenlyRealised}
\end{table}

In particular, $S_n$ is evenly represented in all classes if and only if $n\ge 7$.

With a few small exceptions, the automorphism groups considered above are all non-solvable. 
The following result shows that each class also realises finite groups of every sufficiently large nilpotence class or derived length:

\begin{thm}\label{nilpandsolv}
There is a finite group of nilpotence class $c$, or of derived length $l$, isomorphic to the automorphism group of an edge-transitive map in a class $T$, if and only if $c$ or $l$ satisfy the corresponding condition in Table~\ref{classandlength}.
\end{thm}

\begin{table}[ht]
\centering
\begin{tabular}{| p{2.8cm} | p{3.4cm} | p{3.2cm} |}
\hline
Class $T$ & Nilpotence class $c$ & Derived length $l$ \\
\hline\hline
$2\,{\rm ex}, 2^*{\rm ex}, 2^P{\rm ex}$& $c\ge 5$ & $l\ge 2$  \\
\hline
$5, 5^*, 5^P$ & $c\ge 2$ & $l\ge 2$ \\
\hline
All other $T$ & $c\ge 1$ & $l\ge 1$ \\
\hline
 
\end{tabular}
\caption{Nilpotence class and derived length.}
\label{classandlength}
\end{table}

Although the main emphasis of Part~\ref{autgps} is on compact maps, and hence on finite automorphism groups, non-compact edge-transitive maps, which may have infinite automorphism groups, are considered in Section~\ref{infautogps}. For example, a construction due to B.~H.~Neumann~\cite{Neu} is adapted to prove the following:

\begin{thm}\label{uncountable}
Each of the $14$ classes $T$ of edge-transitive maps contains $2^{\aleph_0}$ maps $\mathcal M$ with empty boundary and with mutually non-isomorphic automorphism groups ${\rm Aut}\,{\mathcal M}$.
\end{thm}

\noindent An embedding theorem of Schupp~\cite{Sch} is used to prove:

\begin{thm}\label{embedctble}
For each of the $14$ edge-transitive classes $T$, every countable group $C$ is isomorphic to a subgroup of ${\rm Aut}\,{\mathcal M}$ for some map $\mathcal M$ in $T$.
\end{thm}

\noindent Questions of growth and decidability are also considered in Section~\ref{infautogps}.

There has been recent interest in arc-transitive maps~\cite{HRS}, motivated in part by the importance of arc-transitive graphs in both algebraic graph theory and permutation group theory. Of course, arc-transitive maps are all edge-transitive, so the results in this paper can be applied to them, simply by restricting attention to the arc-transitive classes $T=1$, $2^*$, $2^P$, $2^*{\rm ex}$ and $2^P{\rm ex}$, those for which ${\mathcal N}(T)$ has a single arc. For instance, it follows from Theorem~1.\ref{FSGpsRealised} that every non-abelian finite simple group, with the single exception of $U_3(3)$, is the automorphism group of an arc-transitive map.

Part~\ref{topology}, which is independent of Part~\ref{autgps}, is devoted to topological properties of edge-transitive maps, in particular their orientability, their Euler characteristic (when compact) and their boundary behaviour. One of the main results is Theorem~\ref{voidtamewild}, which states that of the 14 edge-transitive classes, six are {\em void\/} (meaning that they contain no maps with boundary), four are {\em tame\/} (meaning that they contain such maps, all with dihedral automorphism groups), and four are {\em wild\/} (meaning that they contain such maps, some with non-dihedral automorphism groups). For the tame classes, the maps with boundary are completely classified (see Theorem~\ref{tamemaps}): in each case they consist of at most two infinite families (each containing one infinite map) and at most six sporadic finite examples, all easily described and illustrated (most are on the closed disc). In the case of the wild classes, however, the abundance of examples is such that a classification is not feasible: indeed, in a sense made more precise later, the classification problem in each case is essentially equivalent to that for regular maps or for regular hypermaps. In Section~\ref{boundarycpts} we consider how to count the boundary components of a map, and to determine their combinatorial nature.

If an edge-transitive map contains a free edge, then clearly all edges must be free; this is a very restrictive condition, allowing a classification of such maps in Theorem~\ref{free}. There are two infinite families, each containing one infinite map, together with six sporadic finite examples. All are on the closed disc, apart from one of the infinite families, which are on the sphere, and all have dihedral or trivial automorphism groups. Section~\ref{almostET} revisits medial maps, introduced in Part~\ref{basics}, and considers {\sl almost edge-transitive maps}, those which are not edge-transitive but have an edge-transitive medial map.

\iffalse
Part~III also contains a discussion of the canonical double cover $\M^+$ of a map $\M$, which exists whenever $\M$ is non-orientable or has non-empty boundary. This is an orientable map without boundary, and it most cases it has automorphism group ${\rm Aut}\,\M\times C_2$. However there are examples, including edge-transitive maps, where ${\rm Aut}\,\M^+$ is strictly larger than this; by analogy with terminology for bipartite double covers of graphs, such maps $\M$ are called {\em unstable}. A number of examples are constructed in Section~\ref{stability}.
\fi

%%%%%%%%%%%%%%%%%%%

\section{Outline of the method}\label{method}

The following method is used to prove these results; full details are given in later sections. Maps $\M$ (always assumed to be connected) can be regarded as transitive permutation representations of the group
\[\Gamma=\langle R_0, R_1, R_2\mid R_i^2=(R_0R_2)^2=1\rangle,\]
acting as the monodromy group of $\M$ by permuting its flags. This group $\Gamma$ is the free product of a Klein four group $E=\langle R_0, R_2\rangle\cong V_4$ and a cyclic group $\langle R_1\rangle\cong C_2$. Each map $\M$ has automorphism group ${\rm Aut}\,\M\cong N_{\Gamma}(M)/M$ where $M$ (a {\em map subgroup\/}) is the stabiliser in $\Gamma$ of a flag. Then $\M$ is edge-transitive if and only if $\Gamma=N_{\Gamma}(M)E$. The $14$ classes $T$ of edge-transitive maps correspond to the $14$ conjugacy classes of subgroups $N=N(T)\le \Gamma$ satisfying $\Gamma=NE$: the maps $\M$ in each class are those for which $M$ has normaliser $N$, so that ${\rm Aut}\,\M\cong N/M$ and $\M$ is a regular covering of the single-edge map ${\mathcal N}(T)$ with map subgroup $N$. For instance $N(1)=\Gamma$, so that regular maps correspond to normal subgroups $M$ of $\Gamma$, while $N(2^P{\rm ex})$ is the even subgroup $\Gamma^+$ of index $2$ in $\Gamma$, consisting of the words of even length in the generators $R_i$, and $N(3)$ is the normal closure of $R_1$.

The Reidemeister-Schreier process yields presentations for these subgroups $N=N(T)$ of $\Gamma$ (all of index at most $4$), and these give information about their quotients. Thus the quotients of $\Gamma$ (and hence the automorphism groups of regular maps) are those groups generated by three elements of order dividing $2$, with two of them commuting; for instance, as a result of work of Nuzhin and others~\cite{Maz, Nuz90, Nuz92, Nuz97a, Nuz97}, it is known which finite simple groups have this property (see \S\ref{gensimple} for details). However, if $T\ne 1$ then not all quotients $N/M$ of $N$ are in $\G(T)$: we also require $N/M$ not to have certain `forbiddden automorphisms' (see Table~\ref{forbidden} in \S\ref{parentforbidden}, also~\cite[Condition~3.2]{STW}) which would cause $N_{\Gamma}(M)$ to be strictly larger than $N$. For instance, $N(2^P{\rm ex})=\Gamma^+=\langle X, Y\mid Y^2=1\rangle\cong C_{\infty}*C_2$ where $X=R_1R_2$ and $Y=R_0R_2$, and the corresponding quotients, the groups in $\G(2^P{\rm ex})$, are those with no automorphism inverting the images of $X$ and $Y$; Leemans and Liebeck~\cite{Lee16, LL} have recently characterised the finite simple groups with this property (again, see \S\ref{gensimple} for details).

The task of determining such quotients for each class $T$ is eased by using a group $\Omega$ of map operations, introduced by Wilson~\cite{Wil} and generated by the classical and Petrie dualities for maps:  induced by the outer automorphism group ${\rm Out}\,\Gamma\cong {\rm Aut}\,E\cong S_3$ of $\Gamma$ acting on the conjugacy classes of subgroups $N$, this group permutes the $14$ classes $T$ in six orbits, corresponding to the six rows in Tables~\ref{GpsRealised} and~\ref{FSGpsRealised}; since $\Omega$ preserves the sets of automorphism groups realised by these classes, it is sufficient to consider one representative class from each orbit. For this we will choose the classes $T=1, 2, 3, 4$ and $5$, though in the case of the remaining orbit it is often more convenient to choose the class $2^P{\rm ex}$ of chiral maps, rather than $T=2\,{\rm ex}$.

These groups $N(T)$ have decompositions as free products, namely
\[N(1)=\Gamma\cong V_4*C_2,\quad N(2)\cong C_2*C_2*C_2,
\quad N(2^P{\rm ex})\cong C_{\infty}*C_2,\]
\[N(3)\cong C_2*C_2*C_2*C_2,\quad N(4)\cong C_2*C_2*C_{\infty},
\quad N(5)\cong C_{\infty}*C_{\infty}\cong F_2,\]
leading to further simplifications. There are epimorphisms $N(T)\to\Gamma$ for $T=2, 3$ and $4$, so that every quotient $G$ of $\Gamma$ is also a quotient of $N(T)$; these epimorphisms can be chosen so that in almost all cases $G$ has no forbidden automorphisms, and is therefore realised in classes $T=2$, $3$ and $4$. Similarly, an obvious epimorphism $N(5)\to N(2^P{\rm ex})$ means that any group in $\G(2^P{\rm ex})$ is also in $\G(5)$. Although it does not solve the realisation problem completely, this simplification (stated more precisely in Lemma~\ref{regmapslemma}) means that a large part of the solution can be obtained by concentrating on the two groups $\Gamma$ and $\Gamma^+$ and their quotients, namely the automorphism groups of the regular maps and the orientably regular maps. Fortunately there is a significant amount of useful material already in the literature, or in accessible databases, about these groups. For example, a series of papers by Nuzhin and others characterises the non-abelian finite simple groups in $\G(1)$, while a recent result of Leemans and Liebeck~\cite{LL} does the same for $\G(2^P{\rm ex})$.

Nevertheless, there are cases not covered by these simplifications; for these one has to find explicit generating sets satisfying the relevant conditions for each class $T$, or (in a few exceptional cases) to prove that no such sets exist. In the case of the alternating and symmetric groups, we use classical and more recent results on their generators. For the simple groups of Lie type and the sporadic simple groups, we use results and techniques of Malle, Saxl and Weigel~\cite{MSW}, together with the detailed information in the ATLAS~\cite{ATLAS} concerning automorphisms, maximal subgroups, conjugacy classes and character tables.

%%%%%%%%%%%%%%%%%%%%%%%%%%
%%%%%%%%%%%%%%%%%%%%%%%%%%

\newpage

 \part{Basic ideas}\label{basics}
 
The purpose of this part of the text is to describe the concepts and techniques which will be used in later parts to study the automorphism groups and topological properties of edge-transitive maps. Many of these are already in the literature, but it is useful to bring them together here with a unified treatment.

%%%%%%%%

\section{Algebraic theory of maps}\label{algthymaps}

Here we briefly outline the algebraic theory of maps developed in more detail elsewhere (see~\cite{JT}, for example, and see~\cite{GT} for background in topological graph theory). 

Each map $\mathcal M$ (possibly non-orientable or with non-empty boundary) determines a permutation representation of the group
\[\Gamma=\langle R_0, R_1, R_2\mid R_i^2=(R_0R_2)^2=1\rangle\cong V_4*C_2,\]
on the set $\Phi$ of flags $\phi=(v,e,f)$ of $\mathcal M$, where $v, e$ and $f$ are a mutually incident vertex, edge and face. For each $\phi\in\Phi$ and each $i=0, 1, 2$, there is at most one flag $\phi'\ne \phi$ with the same $j$-dimensional components as $\phi$ for each $j\ne i$ (possibly none if $\phi$ is a boundary flag). Define $r_i$ to be the permutation of $\Phi$ transposing each $\phi$ with $\phi'$ if the latter exists, and fixing $\phi$ otherwise .(See Figures~\ref{flags} and~\ref{fixedflags} for the former and latter cases. In Figure~\ref{fixedflags}, as in all diagrams, the broken line represents part of the boundary of the map.). Since $r_i^2=(r_0r_2)^2=1$ there is a permutation representation
\[\theta:\Gamma\to G:=\langle r_0, r_1, r_2\rangle\le{\rm Sym}\,\Phi\]
of $\Gamma$ on $\Phi$, given by $R_i\mapsto r_i$.

\begin{figure}[h!]
\begin{center}
\begin{tikzpicture}[scale=0.5, inner sep=0.8mm]

\node (c) at (0,0) [shape=circle, fill=black] {};
\node (d) at (8,0) [shape=circle, fill=black] {};
\draw [thick] (c) to (d);
\draw [thick] (c) to (1,-3);
\draw [thick] (c) to (1,3);
\draw [thick] (d) to (7,-3);
\draw [thick] (d) to (7,3);
\draw [thick] (c) to (-2.5,2.5);
\draw [thick] (c) to (-2.5,-2.5);
\draw [thick] (d) to (10.5,2.5);
\draw [thick] (d) to (10.5,-2.5);

\draw (c) to (1,0.5);
\draw (c) to (1,-0.5);
\draw (1,0.5) to (1,-0.5);
\draw (d) to (7,0.5);
\draw (d) to (7,-0.5);
\draw (7,0.5) to (7,-0.5);
\draw (c) to (0.8,0.8);
\draw (0.3,1) to (0.8,0.8);

\node at (-0.8,0) {$v$};
\node at (4,-0.4) {$e$};
\node at (4,2) {$f$};

\node at (1.5,0.5) {$\phi$};
\node at (6.2,0.5) {$\phi r_0$};
\node at (1.4,1.3) {$\phi r_1$};
\node at (1.9,-0.5) {$\phi r_2$};
\node at (5.9,-0.5) {$\phi r_0r_2$};

\end{tikzpicture}

\end{center}
\caption{Generators $r_i$ of $G$ acting on a flag $\phi=(v,e,f)$.} 
\label{flags}
\end{figure}

\begin{figure}[h!]
\begin{center}
\begin{tikzpicture}[scale=0.5, inner sep=0.8mm]

\draw [dashed] (-10,0) to (12,0);

\node (a) at (-8,3) [shape=circle, fill=black] {};
\draw [thick] (a) to (-8,0);
\draw (a) to (-7.5,2) to (-8,2);
\node at (-5.5,2.5) {$\phi r_0=\phi$};

\node (b) at (0,0) [shape=circle, fill=black] {};
\draw [thick] (b) to (0,3);
\draw (b) to (0.5,1) to (0,1);
\node at (2.5,1) {$\phi r_1=\phi$};

\node (c) at (7,0) [shape=circle, fill=black] {};
\draw [thick] (c) to (10,0);
\draw (c) to (8,0.5) to (8,0);
\node at (9,1.2) {$\phi r_2=\phi$};

\end{tikzpicture}

\end{center}
\caption{Flags fixed by $r_0, r_1$ and $r_2$.}
\label{fixedflags}
\end{figure}

Conversely, any permutation representation of $\Gamma$ on a set $\Phi$ determines a map $\mathcal M$ in which the vertices, edges and faces are identified with the orbits on $\Phi$ of the subgroups $\langle R_1, R_2\rangle\cong D_{\infty}$,  $\langle R_0, R_2\rangle\cong V_4$ and  $\langle R_0, R_1\rangle\cong D_{\infty}$, incident when they have non-empty intersection.

The map $\mathcal M$ is connected if and only if $\Gamma$ acts transitively on $\Phi$, as we will always assume. In this case the stabilisers in $\Gamma$ of flags $\phi\in\Phi$ form a conjugacy class of subgroups $M\le\Gamma$, called {\em map subgroups}.

One can represent $\Gamma$ as a group of isometries of the hyperbolic plane $\hyp$, with the generators $R_i$ acting as reflections in the sides of an ideal right-angled triangle $T$. This is shown in Figure~\ref{Riacting} using the Poincar\'e disc model of $\hyp$, with vertices of $T$ at  $0, 1$ and $i$. The images of $T$ form a tessellation of $\hyp$, the barycentric subdivision of the {\em universal map\/} $\M_{\infty}$, a map in which the vertices and edges are the images of $1$ and the unit interval. Then $\M\cong\M_{\infty}/M$, giving a concrete geometric realisation of $\M$.

\begin{figure}[h!]
\begin{center}
\begin{tikzpicture}[scale=0.5, inner sep=0.8mm]

\node (a) at (5,0) [shape=circle, fill=black] {};

\draw [thick, dotted] (5,0) arc (0:360:5);
\draw [ultra thick] (a) to (0,0);
\draw (0,0) to (0,5);
\draw (0,5) arc (180:270:5);

\node at (5.7,0) {$1$};
\node at (-0.5,0) {$0$};
\node at (0,5.6) {$i$};

\draw [<->] (-0.5,1.5) to (0.5,1.5);
\draw [<->] (1.1,1.1) to (1.8,1.8);
\draw [<->] (1.5,0.5) to (1.5,-0.5);

\node at (-1.3,1.5) {$R_0$};
\node at (2.3,2.3) {$R_1$};
\node at (1.6,-1.2) {$R_2$};
\node at (-2,-2) {$\hyp$};

\end{tikzpicture}
\end{center}
\caption{Generators $R_i$ of $\Gamma$ acting on the hyperbolic plane} 
\label{Riacting}
\end{figure}

The map $\mathcal M$ is finite (has finitely many flags) if and only if $M$ has finite index in $\Gamma$, and it has non-empty boundary if and only if $M$ contains a conjugate of some $R_i$, or equivalently some $r_i$ has a fixed point in $\Phi$. In particular, $\M$ is orientable and with empty boundary if and only if $M$ is contained in the even subgroup $\Gamma^+$ of index $2$ in $\Gamma$, consisting of the words of even length in the generators $R_i$, those preserving the orientation of $\hyp$.

The group $G$ is called the {\em monodromy group\/} ${\rm Mon}\,{\mathcal M}$ of $\mathcal M$. The {\em automorphism group\/} $A={\rm Aut}\,{\mathcal M}$ of $\mathcal M$ is the centraliser of $G$ in ${\rm Sym}\,\Phi$. We have $A\cong N/M$ where $N:=N_{\Gamma}(M)$ is the normaliser  of $M$ in $\Gamma$. The map $\mathcal M$ is called {\em regular\/} if $A$ is transitive on $\Phi$, or equivalently $G$ is a regular permutation group, that is, $M$ is normal in $\Gamma$; in this case
\[A\cong G\cong \Gamma/M,\]
and one can identify $\Phi$ with $G$, so that $A$ and $G$ are the left and right regular representations of $G$ on itself. The map $\mathcal M$ is {\em edge-transitive\/} if $A$ acts transitively on its edges; this is equivalent to the condition that $\Gamma=NE$, where $E:=\langle R_0, R_2\rangle\cong V_4$.

The (classical) dual $D({\mathcal M})$ of $\mathcal M$, a map on the same surface formed by transposing the roles of vertices and faces, corresponds to the image of $M$ under the automorphism $\delta$ of $\Gamma$ which fixes $R_1$, and transposes $R_0$ and $R_2$. The Petrie dual $P({\mathcal M})$ embeds the same graph as $\mathcal M$, but the faces are transposed with Petrie polygons, closed zig-zag paths which alternately turn first right and first left at the vertices of $\mathcal M$; this operation corresponds to the automorphism $\pi$ of $\Gamma$ which transposes $R_0$ and $R_0R_2$, and fixes $R_1$ and $R_2$. Both of these operations $D$ and $P$ preserve regularity and automorphism groups, but $P$ may change the underlying surface. For example, the Petrie dual of the tetrahedron $\{3,3\}$, a regular embedding of the complete graph $K_4$ in the sphere, is the antipodal quotient of the cube $\{4,3\}$, a regular embedding of $K_4$ in the real projective plane (see Figures~\ref{tetmap} and \ref{projplKn}). Similarly, an interior free edge, with flags fixed by $r_0r_2$, is transposed by $D$ with a boundary free edge, with flags fixed by $r_0$, so $D$ may cause boundary components to appear or disappear. The group $\Omega=\langle D, P\rangle$ of map operations generated by $D$ and $P$, introduced by Wilson in~\cite{Wil}, is isomorphic to $S_3$, permuting vertices, faces and Petrie polygons; it corresponds to the outer automorphism group ${\rm Out}\,\Gamma\cong {\rm Aut}\, E\cong S_3$ of $\Gamma$ acting on maps by permuting conjugacy classes of map subgroups~\cite{JT}.

%%%%%%%%%%%%%%%

\section{Edge-transitive maps}\label{edgetrans}

Graver and Watkins~\cite{GW} partitioned edge-transitive maps into $14$ classes, mainly for applications to infinite but locally finite edge-transitive planar maps; see their Table~2 for a summary. In this section we will describe and analyse their classification, but from a rather more group-theoretic point of view. This is largely motivated by the work of  \v Sir\' a\v n, Tucker and Watkins in~\cite{STW}, where finite symmetric groups are used as automorphism groups in order to construct finite maps of particular types in all $14$ classes, together with that of Orbani\'c, Pellicer, Pisanski and Tucker in~\cite{OPPT}, where edge-transitive maps of small genus are classified. (In the present paper, `class' refers to the $14$ equivalence classes of edge-transitive maps, whereas `type'  refers to the parameters giving the valencies of the faces and vertices of a map, as in~\cite{CM}.)

\subsection{General properties}

A map $\mathcal M$, corresponding to a conjugacy class in $\Gamma$ of map subgroups $M$, is edge-transitive if and only if $NE=\Gamma$, where $N=N_{\Gamma}(M)$, and $E:=\langle R_0, R_2\rangle\cong V_4$. There are $14$ conjugacy classes of subgroups $N\le\Gamma$ satisfying $NE=\Gamma$: these are the map subgroups for the $14$ isomorphism classes of maps $\mathcal N$ with a single edge, shown in Figure~\ref{basicmaps}. The edge-transitive maps $\mathcal M$ are thus partitioned into $14$ classes $T$, according to the conjugacy class of $N_{\Gamma}(M)$, or equivalently the common isomorphism class of the quotient maps
\[{\mathcal N}={\mathcal N}(T)={\mathcal M}/{\rm Aut}\,{\mathcal M}\]
of the maps $\mathcal M$ in $T$. We call ${\mathcal N}(T)$ the {\em basic map\/} for $T$, and the group $N=N(T)=N_{\Gamma}(M)$ (determined only up to conjugacy in $\Gamma$) the {\em parent group\/} for $T$, since the maps $\mathcal M$ in $T$ are all regular coverings of ${\mathcal N}(T)$ by quotient groups $G\cong N(T)/M$ corresponding to normal subgroups $M$ of $N(T)$. (In essence,  ${\mathcal N}(T)$ is a particular case of the symmetry type graph of a polytope or maniplex, studied by Cunningham, del R\'\i o-Francos, Hubard and Toledo  in~\cite{CDHT}.)

\begin{figure}[h!]
\begin{center}
\begin{tikzpicture}[scale=0.29, inner sep=0.8mm]

\node (a) at (0,53) [shape=circle, fill=black] {};
\draw [very thick] (0,47) arc (-90:90:3);
\draw [dashed] (0,53) arc (90:270:3);
\node at (0,45.5) {$1$ VFP};
% 1

%%%%%%%%%%%%%%%%%%%%%%%

\node (b) at (0,43) [shape=circle, fill=black] {};
\node (c) at (0,37) [shape=circle, fill=black] {};
\draw [very thick] (0,37) arc (-90:90:3);
\draw [dashed] (0,43) arc (90:270:3);
\node at (0,35.5) {$2$ FP};
% 2

%%%%%%%

\node (d) at (16,40) [shape=circle, fill=black] {};
\draw [dashed] (16,40) arc (0:360:3);
\draw [very thick] (d) to (10,40);
\node at (13,35.5) {$2^*$ VP};
% 2*

%%%%%%%

\node (e) at (29,40) [shape=circle, fill=black] {};
\draw [dashed] (29,40) arc (0:360:3);
\draw [very thick] (e) to (26,40);
\node at (26,35.5) {$2^P$ VF};
% 2P

%%%%%%%%%%%%%%%%%%%

\node (f) at (3,30) [shape=circle, fill=black] {};
\draw [very thick] (3,30) arc (0:360:3);
\node at (0,25.5) {$2\,{\rm ex}$ VFP};
% 2ex

%%%%%

\node (g) at (13,30) [shape=circle, fill=black] {};
\draw [dashed] (16,30) arc (0:360:3);
\draw [very thick] (g) to (10,30);
\node at (13,25.5) {$2^*$ex VFP};
% 2*ex

%%%%%%%

\node (h) at (27.5,30) [shape=circle, fill=black] {};
\draw [thin] (29,30) arc (0:360:3);
\draw [very thick] (h) to (24.5,30);
\node at (26,25.5) {$2^P$ex VFP};
% 2Pex

%%%%%%%%%%%%%%%%%%%%

\node (i) at (-3,20) [shape=circle, fill=black] {};
\node (j) at (3,20) [shape=circle, fill=black] {};
\draw [very thick] (i) to (j);
\draw [dashed] (3,20) arc (0:360:3);
\node at (0,15.5) {$3$};
% 3

%%%%%%%%%%%%%%%%%%%%

\node (k) at (0,10) [shape=circle, fill=black] {};
\node (l) at (3,10) [shape=circle, fill=black] {};
\draw [very thick] (k) to (l);
\draw [dashed] (3,10) arc (0:360:3);
\node at (0,5.5) {$4$ FP};
% 4

%%%%%%

\draw [dashed] (16,10) arc (0:360:3);
\node (m) at (10,10) [shape=circle, fill=black] {};
\draw [very thick] (m) to [out=30,in=90] (14,10);
\draw [very thick] (m) to [out=-30,in=-90] (14,10);
\node at (13,5.5) {$4^*$ VP};
% 4*

%%%%%%

\draw [dashed] (23,13) to (29,13);
\draw [dashed] (23,7) to (29,7);
\draw [thick] [dotted] (23,13) to (23,7);
\draw [thick] [dotted] (29,13) to (29,7);
\node (n) at (23,7) [shape=circle, fill=black] {};
\node (o) at (29,13) [shape=circle, fill=black] {};
\draw [very thick] (n) to (o);
\draw [thin] (22.5,9.75) to (23,10.25);
\draw [thin] (23.5,9.75) to (23,10.25);
\draw [thin] (28.5,10.25) to (29,9.75);
\draw [thin] (29.5,10.25) to (29,9.75);
\node at (26,5.5) {$4^P$ VF};
% 4P

%%%%%%%%%%%%%%%%%%%%%

\node (a) at (-1.5,0) [shape=circle, fill=black] {};
\node (b) at (1.5,0) [shape=circle, fill=black] {};
\draw [thin] (3,0) arc (0:360:3);
\draw [very thick] (a) to (b);
\node at (0,-4.5) {$5$ FP};
% 5

%%%%%

\draw [thin] (16,0) arc (0:360:3);
\node (c) at (11,0) [shape=circle, fill=black] {};
\draw [very thick] (c) to [out=30,in=90] (15,0);
\draw [very thick] (c) to [out=-30,in=-90] (15,0);
\node at (13,-4.5) {$5^*$ VP};
% 5*

%%%%%%%

\node (d) at (23,0) [shape=circle, fill=black] {};
\node (e) at (29,0) [shape=circle, fill=black] {};
\draw [very thick] (29,0) arc (0:360:3);
\draw (25.75,3) to (26.25,3.5);
\draw (25.75,3) to (26.25,2.5);
\draw (26.25,-3) to (25.75,-3.5);
\draw (26.25,-3) to (25.75,-2.5);
\node at (26,-4.5) {$5^P$ VF};
% 5P

\end{tikzpicture}

\end{center}
\caption{The basic maps $\mathcal N$ for the $14$ edge-transitive types}
\label{basicmaps}
\end{figure}

For any map $\M$, the quotient map $\M/{\rm Aut}\,\M$ can be found by taking a fundamental region $\mathcal F$ for ${\rm Aut}\,\M$ and identifying equivalent boundary points. When $\M$ is edge-transitive,  $\M/{\rm Aut}\,\M$ must be isomorphic to one of the 14 basic maps ${\mathcal N}(T)$, enabling one to determine the class $T$ to which $\M$ belongs.

\smallskip

\noindent{\bf Example 1.} If $\M$ is a regular map one can take $\mathcal F$ to be any face in its barycentric subdivision. This is a triangle with a vertex of $\M$ at a corner, and half of an edge of $\M$ along one side. In this case ${\rm Aut}\,\M$ is generated by reflections in the three sides of $\mathcal F$, so there are no side-identifications. Thus ${\rm Aut}\,\M\cong{\mathcal F}\cong{\mathcal N}(1)$, so $\M$ is in class~$1$.

\smallskip

\noindent{\bf Example 2.} Suppose that $\M$ is the cuboctahedron, a map on the sphere shown in Figure~\ref{cuboct2*} as a geometric polytope in ${\bf E}^3$. In this case ${\rm Aut}\,\M$ is the common symmetry group $S_4\times C_2$ of the cube and the octahedron. One can take the red triangle $\mathcal F$ as a fundamental region for this group. Again there are no side-identifications, but now we have ${\rm Aut}\,\M\cong{\mathcal F}\cong{\mathcal N}(2^*)$, so $\M$ is in class~$2^*$. The fact that ${\mathcal N}(2^*)$ has one vertex and two faces indicates that each map in this class is vertex- but not face-transitive.

\medskip

\begin{figure}[h!]
\begin{center}
\begin{tikzpicture}[scale=0.3, inner sep=0.5mm]
 
\node (A) at (8,0) [shape=circle, fill=black] {};  
\node (B) at (4,6.92) [shape=circle, fill=black] {};
\node (C) at (-4,6.92) [shape=circle, draw, fill=black] {};
\node (D) at (-8,0) [shape=circle, fill=black] {}; 
\node (E) at (-4,-6.92) [shape=circle, fill=black] {};  
\node (F) at (4,-6.92) [shape=circle, fill=black] {};

\draw [thick] (A) to (B) to (C) to (D) to (E) to (F) to (A);

\node (x) at (0,5) [shape=circle, fill=black] {};
\node (y) at (-4,-2) [shape=circle, fill=black] {};
\draw [thick] (C) to (x) to (y) to (D);

\node (z) at (4,-2) [shape=circle, fill=black] {};
\draw [thick] (B) to (x) to (z) to (A);
\draw [thick] (E) to (y) to  (z) to (F);

\node (a) at (0,0.33) {};  
\node (b) at (0,-4.46) {};  
\draw [thick, red, dashed] (z) to (a) to (b) to (z);

\node at (-11,0) {$\mathcal M$};
\node [red] at (0.8,-1) {$\mathcal F$};

%%%%%

\node (c) at (12,2.75) {};  
\node (d) at (12,7.25) {};  
\node (e) at (16,5) [shape=circle, fill=black] {};
\draw [thick, red, dashed] (e) to (c) to (d) to (e);
\draw [thick] (e) to (12,5) {};
\node at (19,5) {$\mathcal F$};

%%%%%

\node (f) at (16,-5) [shape=circle, fill=black] {};
\draw [dashed] (16,-5) arc (0:360:2.5);
\draw [thick] (f) to (11,-5);
\node at (19.5,-5) {${\mathcal N}(2^*)$};

\end{tikzpicture}

\end{center}
\caption{The cuboctahdron}
\label{cuboct2*}
\end{figure}

\noindent{\bf Example 3.} In Figure~\ref{torusK5}, opposite sides of the outer square are identified to form a torus, giving an orientably regular chiral map $\M$ of type $\{4,4\}$ with ${\rm Aut}\,\M\cong AGL_1(5)\cong C_5\rtimes C_4$. (It is, in fact, a Cayley map for the additive group of the field $\F_5$, with generating set $\F_5^*=\F_5\setminus\{0\}$ taken in the cyclic order $1, 2, 4, 3$, and ${\rm Aut}\,\M$ acting by affine transformations $v\mapsto av+b\;(a\ne 0)$ of the vertices; see \S\ref{completemaps} for a generalisation of this construction.) The red triangle $\mathcal F$ is a fundamental region for ${\rm Aut}\,\M$, but in this case there are side-pairings, induced by rotations of $\M$ of order $2$ and $4$. These identifications yield a spherical map isomorphic to ${\mathcal N}(2^P{\rm ex})$, so $\M$ belongs to the class~$2^P{\rm ex}$. A similar argument applies to every orientably regular chiral map.

\medskip

\begin{figure}

\begin{center}
 \begin{tikzpicture}[scale=0.7, inner sep=0.5mm]

\node (A) at (-6,8) [shape=circle, fill=black] {};
\node (B) at (-6,6) [shape=circle, fill=black] {};
\node (C) at (-4,6) [shape=circle, fill=black] {};
\node (D) at (-2,6) [shape=circle, fill=black] {};
\node (E) at (-8,4) [shape=circle, fill=black] {};
\node (F) at (-6,4) [shape=circle, fill=black] {};
\node (G) at (-4,4) [shape=circle, fill=black] {};
\node (H) at (-4,2) [shape=circle, fill=black] {};

\draw (-7,6) to (D);
\draw (E) to (-3,4);
\draw (A) to (-6,3);
\draw ((-4,7) to (H);

\draw[thick, dotted] (A) to (D) to (H) to (E) to (A);

\draw[thick, red, dashed] (-4,6) to (-5,7) to (-5,5) to (-4,6);
\draw[thick, red] (-4.8, 6.4) to (-5,6.6) to (-5.2,6.4);
\draw[thick, red] (-4.8, 5.6) to (-5,5.4) to (-5.2,5.6);
\draw[thick, red] (-4.3,6.6) to (-4.6, 6.6) to (-4.6,6.3) to (-4.3,6.6);
\draw[thick, red] (-4.3,5.4) to (-4.6, 5.4) to (-4.6,5.7) to (-4.3,5.4);

\node at (-6,2) {${\mathcal M}$};

%%%%%%%

\node (u) at (4.3,7) [shape=circle, fill=black] {};
\draw[thick] (u) to (2.7,7);
\draw[thick, red, dashed] (u) to (2.7,8.6) to (2.7,5.4) to (u);

\draw[thick, red] (3,7.7) to (2.7,8) to (2.4,7.7);
\draw[thick, red] (3,6.3) to (2.7,6) to (2.4,6.3);
\draw[thick, red] (3.3,8) to (3.7,8) to (3.3,7.6) to (3.3,8);
\draw[thick, red] (3.3,6) to (3.7,6) to (3.3,6.4) to (3.3,6);

\node at (6,7) {${\mathcal F}$};

%%%%%%%

\node (z) at (4.3,3) [shape=circle, fill=black] {};
\draw [thin] (5,3) arc (0:360:1.5);
\draw[thick] (z) to (2.7,3);
\node at (7,3) {${\mathcal N}(2^P{\rm ex})$};

\end{tikzpicture}

\end{center}
\caption{An orientably regular chiral map on the torus}
\label{torusK5}
\end{figure}

\smallskip

Note that the map $\M$ corresponding to a normal subgroup $M$ of $N(T)$ need not always be in $T$: if $N_{\Gamma}(M)$ properly contains $N(T)$ then $N_{\Gamma}(M)=N(T')$ for some class $T'\ne T$, and $\mathcal M$ is in class $T'$ (see Section~\ref{order} for details of such pairs $T, T'$). This happens if and only if $\mathcal M$ (or equivalently $N(T)/M$) has certain `forbidden automorphisms' induced by elements of $N(T')\setminus N(T)$, or equivalently by lifting non-identity automorphisms of ${\mathcal N}(T)$ to $\M$. For example, if we take $M$ to be $N(T)$ itself, so that ${\mathcal M}={\mathcal N}(T)$, then Figure~\ref{basicmaps} shows that, except when $T=1$, $\mathcal M$ always has such an automorphism, so that ${\mathcal N}(T)$ is not in class $T$.

%%%%%%%%

\subsection{The 14 classes of edge-transitive maps}\label{14classes}

The $14$ maps $\mathcal N$ with one edge are shown in Figure~\ref{basicmaps}, labelled with the symbol $1$, $2$, $2^*,\ldots$, assigned by Graver and Watkins in~\cite{GW} to the corresponding class $T$ (see also~\cite{Wil97} for a similar classification of edge-transitive maps by Wilson). The underlying surface of each map $\mathcal N$ is the closed unit disc $\overline{\mathbb D}$, with the following exceptions: when $T=2^P{\rm ex}$, $5$ or $5^*$ it is the $2$-sphere $S^2$, when $T=4^P$ it is the M\"obius band, and when $T=5^P$ it is the real projective plane. The letters V, F and P indicate whether $\mathcal N$ has one vertex, face or Petrie polygon, and thus whether the maps $\M$ in class~$T$ are vertex-, face- or Petrie-transitive. The rows correspond to the orbits of the group $\Omega$ of map operations, and the second and third maps in each row (when they exist) are the dual and the Petrie dual of the first map; the same applies to the maps $\M$ in the corresponding classes.

To justify the statement about maps with one edge at the beginning of this section, and to show that there really are $14$ classes of edge-transitive maps, we will now prove the following theorem:

\begin{thm}
There are $14$ classes $T$ of edge-transitive maps, corresponding to the $14$ conjugacy classes of subgroups $N$ of $\Gamma$ satisfying $\Gamma=NE$. These are the map subgroups for the maps ${\mathcal N}(T)$ shown in Figure~\ref{basicmaps}.
\end{thm}

\noindent{\sl Proof.} A map $\mathcal M$ with map subgroup $M$ is edge-transitive if and only if $N:=N_{\Gamma}(M)$ acts transitively on the cosets of $E$ in $\Gamma$, that is, $\Gamma=NE$, or equivalently $E$ acts transitively on the cosets of $N$. Since $|E|=4$, this condition implies that $|\Gamma:N|$ divides $4$, so the different possibilities for $N$ (up to conjugacy) can be found by determining the permutation representations of $\Gamma$ of degrees $1, 2$ and $4$ in which $E$ acts transitively. Since $\Gamma=E*\langle R_1\rangle$, one can consider the actions of $E$ and of $R_1$ independently.

When $|\Gamma:N|=1$ we have the trivial action of $\Gamma$, in which each $R_i$ induces the identity permutation. Equivalently $M$ is a normal subgroup of $\Gamma$, so the corresponding maps $\mathcal M$ are the regular maps, those in class~1 in the Graver-Watkins terminology. The basic map $\mathcal N$ has a single flag, and it embeds a free edge along part of the boundary of the closed disc $\overline{\mathbb D}$.

There are seven subgroups $N$ of $\Gamma$ of index $|\Gamma:N|=2$, namely the kernels of the different epimorphisms $\theta:\Gamma\to C_2$, corresponding to the different choices of which generators $R_i$ to map to the identity (any combination, except all of them). If $R_1, R_2\in N$ or $R_0, R_1\in N$ then $E$ is transitive, so these subgroups $N$ are the parent groups for two classes of edge-transitive maps, called $2$ and $2^*$ respectively in~\cite{GW} (the asterisk indicating that the maps in class $2^*$ are the duals of those in class~$2$); however,  if $R_0, R_2\in N$ then $E$ is intransitive, so there is no corresponding class of edge-transitive maps. If just one generator $R_i$ is in $N$, then $E$ is transitive, and we obtain classes $2^*$ex, $2^P$ or $2\,{\rm ex}$ as $i=0, 1$ or $2$ (with $P$ indicating Petrie duality). Finally, if no $R_i$ is in $N$, so that $N=\Gamma^+$, then again $E$ is transitive and we obtain class $2^P$ex. The corresponding basic maps $\mathcal N$ are all on the closed disc $\overline{\mathbb D}$, apart from that for class~$2^P$ex, which is on the sphere $S^2$. 

If $|\Gamma:N|=4$ then  $E$ acts transitively if and only if it acts regularly, in which case one can number the four flags of $\mathcal N$ so that $R_0$ and $R_2$ induce the permutations $(12)(34)$ and $(14)(23)$. If $R_1$ induces the identity permutation we have class~$3$. If $R_1$ induces $(14)$ or $(23)$, giving isomorphic permutation representations of $\Gamma$, we have class~$4$. If $R_1$ induces $(12)$ or $(34)$ we have class~$4^*$, and if it induces $(13)$ or $(24)$ we have class~$4^P$. Finally, if $R_1$ induces $(14)(23)$, $(12)(34)$ or $(13)(24)$ we have classes~$5$, $5^*$ and $5^P$ respectively. The basic maps $\mathcal N$ are on $\overline{\mathbb D}$ for types $3$, $4$ and $4^*$, but on the closed M\"obius band for type $4^P$. For types $5$, $5^*$ and $5^P$ the underlying surfaces are $S^2$, $S^2$ and the real projective plane.
\hfill$\square$

\medskip

Note that if $T=2^P{\rm ex}$, $5$ or $5^*$ then ${\mathcal N}(T)$ is orientable and without boundary, so the same applies to all maps in these classes, since they are coverings of ${\mathcal N}(T)$. When $T=2^P{\rm ex}$ these maps are the orientably regular chiral maps without boundary; for simplicity we will simply call these maps `chiral'. The maps in the class $5^P$ are also without boundary, but they can be orientable or non-orientable.

It is easily seen that $11$ of these $14$ maps $\mathcal N$ are regular, with elementary abelian automorphism and monodromy groups ${\rm Aut}\,{\mathcal N}\cong\Gamma/N$ where $N$ is normal in $\Gamma$; the exceptions are the maps $\mathcal N$ for $T=4$, $4^*$ and $4^P$, which have automorphism groups isomorphic to $C_2$ and monodromy groups  isomorphic to the dihedral group $D_4$ of order $8$.

Some of these edge-transitive classes have appeared, often with different notation, in other recent studies of highly symmetric maps. For instance $2$-orbit maps are those for which the automorphism group has two orbits on flags, or equivalently $|\Gamma:N|=2$ (see~\cite{Hub}, for example, where Hubard has characterized the automorphism groups of $2$-orbit polyhedra). These maps consist of the six edge-transitive classes $2$, $2^*$, $2^P$, $2\,{\rm ex}$, $2^*{\rm ex}$ and $2^P{\rm ex}$, denoted in~\cite{Hub} by $2_{\{1,2\}}$, $2_{\{0,1\}}$, $2_1$, $2_2$, $2_0$ and $2$, together with the non-edge-transitive class $2_{\{0,2\}}$ consisting of those maps $\M$ for which $\M/{\rm Aut}\,\M$ is isomorphic to the map $\mathcal N(2_{\{0,2\}})$ with one vertex and two free edges, all on the boundary of the closed disc (see Figure~\ref{mapN2{02}}).\footnote{In~\cite{Hub} the subscripts $i$ indicate those $i$ such that flags $\phi$ and $\phi r_i$ are in the same orbit; Wilson~\cite{Wil97} denotes these classes by $\overline 012$, $01\overline 2$, $\overline 01\overline 2$, $\overline 0\overline 12$, $0\overline 1\overline 2$, $\overline 0\overline 1\overline 2$ and $0\overline 12$, with $\overline i$ indicating that $r_i$ transposes the two orbits on flags. In the wider theory of hypermaps (see~\cite{BCD} for example), where $\Gamma$ is replaced with the group $\Delta\cong C_2*C_2*C_2$ obtained by omitting the relation $(R_0R_2)^2=1$, the corresponding classes of $2$-orbit hypermaps are associated with the seven subgroups $\Delta^i$, $\Delta^{\hat i}$ and $\Delta^+$ of index $2$ in $\Delta$, where a superscript $i$, $\hat i$ or $+$ indicates that only $r_i$ preserves or transposes the orbits, or that all three generators do.} Although its maps are not edge-transitive, this class $T=2_{\{0,2\}}$ and its basic map will arise several times in this text (see \S\ref{medialmaps} and \S\ref{class1bdy}, for example). The corresponding subgroup $N=N(2_{\{0,2\}})$ of index $2$ in $\Gamma$ is the normal closure of $E=\langle R_0, R_2\rangle$.

\begin{figure}[h!]
\begin{center}
\begin{tikzpicture}[scale=0.5, inner sep=0.8mm]

\node (E) at (13,0) [shape=circle, fill=black] {};
\draw [dashed] (13,0) arc (0:360:3);
\draw [thick] (13,0) arc (0:90:3);
\draw [thick] (13,0) arc (0:-90:3);

\end{tikzpicture}

\end{center}
\caption{The map $\mathcal N(2_{\{0,2\}})$ on the closed disc} 
\label{mapN2{02}}
\end{figure}

As another example, if a map is arc-transitive then it must be edge-transitive, and conversely an edge-transitive map is arc-transitive if and only if it lies in a class $T$ such that ${\mathcal N}(T)$ has a single arc, that is, its single edge is free. By inspection of Figure~\ref{basicmaps} we therefore have the following:

\begin{thm}\label{arctrans}
A map is arc-transitive if and only if it lies in one of the five edge-transitive classes $T=1$, $2^*$, $2^P$, $2^*{\rm ex}$ or $2^P{\rm ex}$.
\end{thm}

Thus all our results for edge-transitive maps can also be applied to arc-transitive maps by restricting them to these five classes: for instance, it follows from Theorem~\ref{mainthmsimple} that every non-abelian finite simple group except $U_3(3)$ is the automorphism group of an arc-transitive map. The four non-regular arc-transitive classes  have been studied by Hubard, Ramos Rivera and \v Sparl in~\cite{HRS}, where the emphasis is on maps in class $2^*$ which embed rose window graphs.

If $T$ is any of the $14$ classes of edge-transitive maps, let us say that a group $G$ is {\sl realised in\/} $T$, or is {\sl evenly realised in\/} $T$, if $T$ contains a map $\M$, or an orientable map $\M$ with empty boundary, such that $G\cong{\rm Aut}\,\M$. Let us define the following sets of groups:
\[\G(T)=\{G\mid G\;\hbox{is realised in}\; T\},\]
\[\G^+(T)=\{G\mid G\;\hbox{is evenly realised in}\; T\}.\]
The aim of this paper is to determine which members of various families of groups, such as the finite alternating and symmetric groups, are in these sets.

%%%%%%%%%%%%%%

\subsection{Operations on edge-transitive maps}\label{operations}

The operations in $\Omega$ preserve edge-transitivity. They permute the $14$ basic maps $\mathcal N$, and they preserve coverings, so they permute the $14$ edge-transitive classes. In particular, the dual pairs are
\[\{2, 2^*\},\; \{2\,{\rm ex},\; 2^*{\rm ex}\},\;\{4, 4^*\},\; \{5, 5^*\},\]
while the remaining six classes are self-dual. (Of course, this does not mean that the maps themselves are self-dual, only that their duals are in the same class.) Similarly, the Petrie dual pairs are
\[\{2^*, 2^P\},\; \{2^*{\rm ex},\; 2^P{\rm ex}\},\;\{4^*, 4^P\},\; \{5^*, 5^P\},\]
with the remaining six classes are Petrie self-dual. Thus the orbits of $\Omega$ are
\[\{1\},\; \{2, 2^*, 2^P\},\; \{2\,{\rm ex}, 2^*{\rm ex}, 2^P{\rm ex}\},\;
 \{3\},\; \{4, 4^*, 4^P\},\; \{5, 5^*, 5^P\},\]
one for each row in Figure~\ref{basicmaps}, where $D$ transposes the first and second terms in each orbit of length $3$, while $P$ transposes the second and third. Here, to avoid possible confusion, we emphasise that despite the notation, the maps in class $2^P$ are the Petrie duals of those in class $2^*$, not those in class $2$; a similar remark applies to the other three orbits of length $3$.) Since $\Omega$ preserves automorphism groups of maps, the sets $\G(T)$ are invariant under its action on the classes $T$. However, the sets $\G^+(T)$ are only $\langle D\rangle$-invariant.

Since $\Omega$ is generated by the operations $D$ and $P$, in order to understand its effect on maps it is sufficient to understand how these two operations act. The effect of $D$ is well-known, transposing vertices and faces while preserving the underlying surface of a map. By contrast, $P$ preserves the embedded graph, but can (and often does) change the topological properties of the surface, including its boundary behaviour, its orientability and its Euler characteristic (though not its compactness).

Duality $D$ preserves possession of a non-empty boundary, and so does Petrie duality $P$ with the following exception: flags fixed by $R_0$, corresponding to boundary free edges, are transposed by $P$ with flags fixed by $R_0R_2$, corresponding to internal free edges, so that $P$ may cause boundary components to appear or disappear. Among edge-transitive maps, this exceptional behaviour applies only to the Petrie dual pairs $\{2^*, 2^P\}$ and $\{2^*{\rm ex}, 2^P{\rm ex}\}$, containing the only classes allowing free edges. In all other cases, $\Omega$ preserves the properties of having empty or non-empty boundary.

The operations $D$ and $P$ transform a map $\M$ of type $\{m,n\}$ into one of type $\{n,m\}$ or $\{l,n\}$ respectively, where $l$ is the the {\sl Petrie length\/} of $\M$, that is, the order of $r_0r_1r_2$, equal to the least common multiple of the valencies of the Petrie polygons of $\M$. In particular, if $\M$ is regular there are $|\Phi|/2l$ of these, all of valency $l$. Since $P$ preserves the vertices and edges of a map, knowing the number of cycles of $r_0r_1r_2$ allows one to compute the Euler characteristic of $P(\M)$.

A map $\M$ is orientable if one can choose local orientations around the vertices which are preserved by continuation along each edge. It is said to be {\sl pseudo-orientable}~\cite{Wil78} or {\sl bi-orientable}~\cite{BCS} if they can be chosen so that they are reversed by each such continuation (for an orientable map, this is equivalent to the graph being bipartite). If $\M$ has both or neither of these orientation properties, then $P(\M)$ has both or neither, respectively. If $\M$ has only one of them, then $P(\M)$ has only the other. See~\cite[\S 2]{BCS} for details.

One can use the operations in $\Omega$ to create new edge-transitive maps from known ones, with the same automorphism group but in a possibly different class and with a possibly different embedded graph or underlying surface.

\medskip

\noindent{\bf Example} As an example, let us consider the regular and chiral maps on the torus, classified in~\cite[\S 8.3 and \S 8.4]{CM} by Coxeter and Moser, whose notation for these maps we will follow. 

In their notation, $\{4,4\}$ denotes the universal map of that type in the complex (or euclidean) plane: it has the set $\Z[i]=\{z=m+ni\mid m, n\in\Z\}$ of Gaussian integers as its vertex-set, with each vertex $z$ joined by straight edges to $z\pm 1$ and $z\pm i$. Its automorphism group is the extended triangle group $\Delta[4,2,4]$, a semidirect product of a normal translation subgroup, isomorphic to $\Z^2$ and acting regularly on the vertices, by the stabiliser of $0$, a dihedral group $D_4$ generated by the rotation $z\mapsto iz$ and the reflection $z\mapsto\overline z$.

For any integers $b, c\ge 0$ (not both $0$), the map $\M=\{4,4\}_{b, c}$ is the quotient of $\{4,4\}$ by the subgroup of ${\rm Aut}\{4,4\}$ generated by the translations $z\mapsto z+b+ci$ and $z\mapsto z-c+bi$, a subgroup which can be identified with the sublattice $\Lambda$ of $\Z^2$ with basis $(b,c)$ and $(-c,b)$. As shown in~\cite[\S8.3]{CM}, this is a map of type $\{4,4\}$ on the torus $\C/\Lambda$, with $b^2+c^2$ vertices and faces, and $2(b^2+c^2)$ edges. It is bipartite if and only if $b\equiv c$ mod~$(2)$. If $b=0$, $c=0$ or $b=c$ then $\M$ is regular, with automorphism group an extension of a normal translation group $T\cong\Z^2/\Lambda$, an abelian group of order $b^2+c^2$, by $D_4$. Otherwise $\M$ is chiral, with mirror image $\{4,4\}_{c,b}$ and automorphism group an extension of $T$ by $C_4$; for instance, $\{4,4\}_{1,2}$ is shown on the left of Figure~\ref{K5K7}. Up to isomorphisms, induced by powers of the rotation $z\mapsto zi$, we may assume that $(b,c)$ is in the first quadrant, with $b\ge 1$ and $c\ge 0$.

\begin{figure}[h!]

\begin{center}
 \begin{tikzpicture}[scale=0.7, inner sep=1mm]

\node (A) at (-6,3) [shape=circle, fill=black] {};
\node (B) at (-6,1) [shape=circle, fill=black] {};
\node (C) at (-4,1) [shape=circle, fill=black] {};
\node (D) at (-2,1) [shape=circle, fill=black] {};
\node (E) at (-8,-1) [shape=circle, fill=black] {};
\node (F) at (-6,-1) [shape=circle, fill=black] {};
\node (G) at (-4,-1) [shape=circle, fill=black] {};
\node (H) at (-4,-3) [shape=circle, fill=black] {};

\draw [thick] (-7,1) to (D);
\draw [thick] (E) to (-3,-1);
\draw [thick] (A) to (-6,-2);
\draw [thick] ((-4,2) to (H);

\draw[thick, dotted] (A) to (D);
\draw[thick, dotted] (D) to (H);
\draw[thick, dotted] (H) to (E);
\draw[thick, dotted] (E) to (A);

%%%%%%%

\node (a) at (3,0) [shape=circle, draw, fill=black] {};
\node (b) at (5,0) [shape=circle, draw, fill=black] {};
\node (c) at (2,1.73) [shape=circle, fill=black] {}; 
\node (d) at (4,1.73) [shape=circle, fill=black] {}; 
\node (e) at (2,-1.73) [shape=circle, fill=black] {}; 
\node (f) at (4,-1.73) [shape=circle, fill=black] {}; 
\node (g) at (1,0) [shape=circle, draw, fill=black] {};

\draw[thick, dotted] (6,0.7) to (4,2.88) to (1,2.28) to (0,-0.7) to (2,-2.88) to (5,-2.28) to (6,0.7); 

\node (h) at (5.88,0) {};
\node (H) at (0.12,0) {};
\draw [thick] (h) to (H);

\node (i) at (4.48,2.57) {};
\node (I) at (1.51,-2.57){};
\draw [thick] (i) to (I);

\node (j) at (1.5,2.52) {};
\node (J) at (4.5,-2.52) {};
\draw [thick] (j) to (J);

\node (k) at (0.7,1.73) {};
\node (K) at (5.2,1.73) {};
\draw [thick] (k) to (K);

\node (l) at (5.3,-1.73) {};
\node (L) at (0.8,-1.73) {};
\draw [thick] (l) to (L);

\node (m) at (3.38,2.9) {}; 
\node (M) at (5.55,-0.95) {}; 
\draw [thick] (m) to (M);

\node (n) at (2.62,-2.9) {}; 
\node (N) at (0.45,0.95) {}; 
\draw [thick] (n) to (N);

\node (o) at (5.77,1.25) {}; 
\node (O) at (3.47,-2.7) {}; 
\draw [thick] (o) to (O);

\node (p) at (0.23,-1.25) {}; 
\node (P) at (2.53,2.7) {}; 
\draw [thick] (p) to (P);

\end{tikzpicture}

\end{center}
\caption{The maps $\{4,4\}_{1,2}$ and $\{3,6\}_{1,2}$}
\label{K5K7}
\end{figure}

There are $b^2+c^2$ vertices and faces in $\M$, and $2(b^2+c^2)$ edges. If $\M$ is a regular map $\{4,4\}_{b,0}$, then it is easily seen that $\M$ has Petrie length $l=2b$, so $P(\M)$ is a regular map of type $\{2b, 4\}$; having $b^2$ vertices, $2b^2$ edges and $8b^2/2l=2b$ faces, it has characteristic $\chi=b(2-b)$. It is orientable if and only if $\M$ is bipartite, that is, $b$ is even. Similarly, if $\M$ is a regular torus map $\{4,4\}_{b,b}$, with $2b^2$ vertices and $4b^2$ edges, then the Petrie length is again $l=2b$, so $P(\M)$ has type $\{2b,4\}$ as before, but it now has characteristic $\chi=2b(2-b)$. In this case, $\M$ is bipartite for all $b$, so $P(\M)$ is orientable, of genus $b(b-2)+1=(b-1)^2$.

However, finding the Petrie length of a chiral torus map $\{4,4\}_{b,c}$, where $bc(b-c)\ne 0$, is a little more complicated. The element $g=(r_0r_1r_2)^2$ acts as a translation $z\mapsto z+1+i$ of $\M$, contained in the translation subgroup $T$ of order $b^2+c^2$ and of index $4$ in ${\rm Aut}\,\M$. This element has order $|T|/|\overline T|$, where $\overline T=T/\langle g\rangle$ is the quotient of $T$ obtained by adding the relation $g=1$ to its defining relations. As a quotient of $\Z^2$, $\overline T$ thus has defining relations $(b,c)=(-c,b)=(1,1)=0$, so $|\overline T|$ is the greatest common divisor of the three $2\times 2$ minors of the $3\times 2$ relation matrix, that is, ${\rm gcd}(b^2+c^2, b-c, -c-b)=d\,{\rm gcd}(b,c)$ where $d=2$ or $1$ as $b\equiv c$ mod~$(2)$ or not. It follows that $g$ has order $(b^2+c^2)/d\,{\rm gcd}(b,c)$, so the Petrie length of $\M$, the (necessarily even) order of $r_0r_1r_2$, is twice this, namely $l=2(b^2+c^2)/d\,{\rm gcd}(b,c)$. Thus $\M$ has $4(b^2+c^2)/2l=d\,{\rm gcd}(b,c)$ Petrie polygons, so $P(\M)$ has Euler characteristic $d\,{\rm gcd}(b,c)-b^2-c^2$. It is orientable if and only if $b\equiv c$ mod~$(2)$. It is edge-transitive, and is in class~$2^*{\rm ex}$ since $\M$ is in class $2^P{\rm ex}$.

In~\cite[Theorem~3.1]{MU}, Meleko\u glu and Ulusan give a different proof of an equivalent formula for $l$ when $\M$ is chiral.

In~\cite[\S8.4]{CM}, Coxeter and Moser defined analogous torus maps $\{6,3\}_{b,c}$, together with their duals $\{3,6\}_{b,c}$, as the quotients of the hexagonal and triangular euclidean maps $\{6,3\}$ and $\{3,6\}$ by a lattice $\Lambda$. In this case $\Lambda$ is generated by translations $X^{b+c}Y^c$ and $X^{-c}Y^b$, where $X$ and $Y$ are translations, inclined at angle $2\pi/3$, generating the full translation group of these euclidean maps. Then $\{6,3\}_{b,c}$ has $2t$ vertices, $3t$ edges and $t$ faces, while $\{3,6\}_{b,c}$ has $t$, $3t$ and $2t$, where $t=b^2+bc+c^2$. As before, these maps are regular if $bc(b-c)=0$, and otherwise they are chiral; for instance, $\{3,6\}_{1,2}$ is shown on the right of Figure~\ref{K5K7}.

Arguments similar to those given above show that $\{6,3\}_{b,c}$ and $\{3,6\}_{b,c}$ have Petrie length $l=2b$, $6b$ or $2t/{\rm gcd}(b,c)$ as $c=0$, $b=c$ or $bc(b-c)\ne 0$. The maps $\{6,3\}_{b,c}$ are bipartite whereas the maps $\{3,6\}_{b,c}$ are not, so their Petrie duals are respectively orientable or non-orientable.

The maps $\M=\{m,n\}_{b,c}$ described here are the only regular or chiral maps on the torus. These results allow one to determine the type, characteristic and orientability of $P(\M)$ in all cases.

%%%%%%%%%%%%%%

\subsection{Ordering of classes}\label{order}

In the Graver-Watkins notation, each class $T$ is denoted by $t^{\sigma}$ where $t$ is an integer $1,\ldots, 5$ and $\sigma$ is the symbol $*$, $P$ or $\emptyset$ (denoting the absence of a symbol); in the case $t=2$ the two orbits of $\Omega$ are distinguished by the absence or presence of  ``ex" after $2^{\sigma}$. Then $N(T)$ has index
\[n=n(T)=|\Gamma:N(T)|=1, 2\;\;{\rm or}\; \;4\]
in $\Gamma$ as $t=1$, $t=2$ or $t\ge 3$, so that ${\mathcal N}(T)$ has $n(T)$ flags, corresponding to the orbits of ${\rm Aut}\,{\mathcal M}$ on the flags of $\mathcal M$ for each ${\mathcal M}\in T$.

Let us define a relation on the set of $14$ classes by writing $T\to T'$ if each subgroup $N(T)$ is contained in a subgroup $N(T')$, or equivalently there is a covering of maps ${\mathcal N}(T)\to{\mathcal N}(T')$, which will then have degree equal to the index $|N(T'):N(T)|=n(T)/n(T')$. In this case we will say that $T$ {\em covers\/} $T'$. Clearly this relation is invariant under $\Omega$; it is reflexive, antisymmetric and transitive, so it is a partial order relation.

It is useful to know which classes $T'$ are covered by each class $T$, for the following reason. If $M$ is a normal subgroup of $N(T)$ then $N_{\Gamma}(M)\ge N(T)$, so the corresponding map $\mathcal M$ has type $T'$ for some $T'$ covered by $T$; if we want to ensure that $\mathcal M$ is in class $T$ then we must choose $M$ so that it is not normal in any subgroup $N(T')$ properly containing $N(T)$, and these correspond to the classes $T'\ne T$ covered by $T$. Equivalently, we must ensure that $\mathcal M$ has no automorphisms other than those induced by $N(T)/M$.

Apart from the cases where $T=4^{\sigma}$ for some $\sigma$ and $T'=1$, all inclusions $N(T)\le N(T')$ are normal, so that the covering ${\mathcal N}(T)\to{\mathcal N}(T')$ is regular, induced by a subgroup of order $n(T)/n(T')$ in ${\rm Aut}\,{\mathcal N}(T)$. Using this, it is easy to determine those classes $T'$ properly covered by each $T$:

\begin{lemma}\label{covering}
If $T$ is one of the $14$ classes of edge-transitive maps, then the classes $T'$ properly covered by $T$ are as follows:
\begin{enumerate}
\item if $T=1$ there are none;
\item if $T=2^{\sigma}$ or $2^{\sigma}{\rm ex}$ for some $\sigma$ then $T'=1$;
\item if $T=3$ then $T'=1$ or $2^{\sigma}$ for some $\sigma$;
\item if $T=4^{\sigma}$ for some $\sigma$ then $T'=1$ or $2^{\sigma}$;
\item if $T=5^{\sigma}$ for some $\sigma$ then $T'=1$, $2^{\sigma}$ or $2^{\tau}{\rm ex}$, where $\tau\ne\sigma$. \hfill$\square$
\end{enumerate}
\end{lemma}

The resulting lattice of coverings is shown in Figure~\ref{orderclass}, with $T$ below each of the classes $T'$ which it covers. We see from this that every proper covering $T\to T'$ is a composition of (at most two) coverings of degree $2$. In order to show that a normal subgroup $M$ of $N(T)$ corresponds to a map $\mathcal M$ of type $T$, and not of some type $T'$ properly covered by $T$, one may therefore restrict attention to those $T'$ which $T$ covers with degree $2$. For each such $T'$ it is then sufficient to choose an arbitrary element $g\in N(T')\setminus N(T)$ and show that $M^g\ne M$, or equivalently that $g$ does not induce an automorphism of $\mathcal M$.

\begin{figure}

\begin{center}
 \begin{tikzpicture}[scale=0.29, inner sep=0.8mm]

\node (1) at (0,10) [shape=circle, fill=black] {};
\node (2) at (-5,5) [shape=circle, fill=black] {};
\node (2ex) at (5,5) [shape=circle, fill=black] {};
\node (3) at (-5,0) [shape=circle, fill=black] {};
\node (4) at (0,0) [shape=circle, fill=black] {};
\node (5) at (5,0) [shape=circle, fill=black] {};

\draw (1) to (2) to (3);
\draw (2) to (4);
\draw (2) to (5) to (2ex) to (1);

\node at (2,10) {$1$};
\node at (-7,5) {$2^{\sigma}$};
\node at (7.5,5) {$2^{\tau}{\rm ex}$};
\node at (-7,0) {$3$};
\node at (2,0) {$4^{\sigma}$};
\node at (7,0) {$5^{\sigma}$};

\node at (20, 6) {$\sigma, \tau\in\{\emptyset, *, P\}$};

\node at (20,4) {$\sigma \ne \tau$};
\end{tikzpicture}

\end{center}
\caption{The ordering of edge-transitive classes}
\label{orderclass}
\end{figure}

Figure~\ref{orderclass} can be seen as a condensed version of the lattice of subgroups $N(T)$ of $\Gamma$, with the vertices labelled $2^{\sigma}$, $2^{\tau}{\rm ex}$ and $5^{\sigma}$ each representing three such subgroups for $\sigma, \tau=\emptyset, *$ and $P$, and the vertex $4^{\sigma}$ representing three conjugacy classes of two such subgroups. Apart from the subgroups $N(4^{\sigma})$, all subgroups $N(T)$ are normal with abelian quotients, so they contain the commutator subgroup $\Gamma'$ of $\Gamma$. This subgroup has index $8$, so it is the intersection of any two of the subgroups $N(3)$ and $N(5^{\sigma})$. It corresponds to the spherical map $\{2,2\}$, which covers each of these basic maps ${\mathcal N}(T)$. For each $\sigma$, the two conjugate subgroups $N(T)$, $T=4^{\sigma}$, intersect in a normal subgroup of index $8$ in $\Gamma$ with quotient group $D_4$, the monodromy group of ${\mathcal N}(T)$ and the automorphism group of its minimal regular cover; for instance, when $T=4$ this is the antipodal quotient of the spherical map $\{4,2\}$.

%%%%%%%%%%%

\subsection{Parent groups and forbidden automorphisms}\label{parentforbidden}

Using the Reidemeister-Schreier process, described in~\cite[\S II.4]{LS} for example, one can obtain the following presentations for representatives of the orbits of $\Omega$ on the $14$ parent groups $N(T)$:

\begin{prop}\label{parents}
We have the following presentations:
\[N(1)=\Gamma=\langle R_0, R_1, R_2\mid R_i^2=(R_0R_2)^2=1\rangle,\]
\[N(2)=\langle S_1=R_1, S_2=R_1^{R_0}, S_3=R_2\mid S_1^2=S_2^2=S_3^2=1\rangle,\]
\[N(2\,{\rm ex})=\langle S_1=R_2, S=R_0R_1\mid S_1^2=1\rangle,\]
\[N(3)=\langle S_0=R_1,\, S_1=R_1^{R_0},\, S_2=R_1^{R_2},\, S_3=R_1^{R_0R_2}\mid S_i^2=1\rangle,\]
\[N(4)=\langle S_1=R_1, S_2=R_1^{R_2}, S=(R_1R_2)^{R_0}\mid S_1^2=S_2^2=1\rangle,\]
\[N(5)=\langle S=R_1R_2, S'=S^{R_0}\mid - \rangle.\]
In each case, $S$ or $S'$ is an element of infinite order in $\Gamma^+$, while $R_i$ or $S_i$ is a reflection.  \hfill$\square$
\end{prop}

By applying elements of $\Omega$, that is, by permuting $R_0, R_2$ and $R_0R_2$, one can then obtain presentations for the other parent groups; of course, groups in the same orbit are isomorphic to each other. In particular it is often more convenient to choose
\begin{equation}
N(2^P{\rm ex})=\Gamma^+=\langle X=R_1R_2, Y=R_0R_2 \mid Y^2=1\rangle
\end{equation}
as a representative of the third orbit, since it is the parent group of the chiral maps, the subject of a great deal of information in the literature.

One can use these presentations to describe the sets $\G(T)$ of automorphism groups ${\rm Aut}\,{\mathcal M}$ of maps $\mathcal M$ in each class $T$. These all arise as quotients $A=N(T)/M$ of $N(T)$. When $T=1$, so that $N(T)=\Gamma$, all such quotients are realised as automorphism groups. However, when $T\ne 1$ one has to use Lemma~\ref{covering} to exclude those $M$ with $N_{\Gamma}(M)=N(T')>N(T)$ for some class $T'$ properly covered by $T$, corresponding to quotients $A$ with `forbidden automorphisms' induced by elements of $N(T')\setminus N(T)$. As usual, it is sufficient to consider one representative $T$ from each orbit of $\Omega$ on classes.

If $T=1$ (the set of regular maps) then the quotients $A$ of $N(T)=\Gamma$ are those groups generated by elements $r_0, r_1$ and $r_2$ of order dividing $2$ (the images of $R_0, R_1$ and $R_2$), with $r_0$ and $r_2$ commuting. 

If $T=2$ then $N(T)\cong C_2*C_2*C_2$, and the quotients are the groups $A$ generated by elements $s_1, s_2$ and $s_3$ of order dividing $2$, with no requirement that two of them should commute. The corresponding map $\mathcal M$ is in class $T=2$ (rather than $T'=1$) if and only if $M$ is not normalised by the element $R_0\in\Gamma\setminus N(2)$. Since $R_0$, acting by conjugation, transposes $S_1$ and $S_2$ and fixes $S_3$, this is equivalent to $A$ not having an automorphism transposing $s_1$ and $s_2$ and fixing $s_3$.  A similar restriction applies when $T=2^*$ or $2^P$.

The arguments in the other cases are similar. To summarise the results (see also~\cite[Condition~3.2]{STW}), the $14$ classes $T$ are listed in Table~\ref{forbidden}, with the second column giving the isomorphism type of each parent group $N(T)$, and the third column indicating any forbidden automorphisms; these are given in terms of their effect on generators $s_i, s, s'$ of quotients $A$, the images of generators $S_i, S, S'$ of cyclic free factors of $N(T)$ defined in Proposition~\ref{parents}.

\begin{table}[htb]
\centering
\begin{tabular}{|c|c|c|}
\hline
Class $T$&$N(T)$&forbidden automorphisms\\
\hline
$1$&$V_4*C_2$&none\\
$2^{\sigma}$&$C_2*C_2*C_2$&$s_1$ and $s_2$ transposed, $s_3$ fixed\\
$2^{\sigma}{\rm ex}$&$C_2*C_{\infty}$&$s_1$ (or $y$)  fixed, $s$ (or $x$) inverted\\
$3$&$C_2*C_2*C_2*C_2$&double transpositions of generators $s_i$\\
$4^{\sigma}$&$C_2*C_2*C_{\infty}$&$s_1$ and $s_2$ transposed, $s$ inverted\\
$5^{\sigma}$&$C_{\infty}*C_{\infty}=F_2$&$s$ and $s'$ inverted, transposed or both\\ \hline
\end{tabular}
\caption{Parent groups and forbidden automorphisms}
\label{forbidden}
\end{table}

If, instead of $T=2\,{\rm ex}$, we use $T=2^P{\rm ex}$ as a representative of the third orbit of $\Omega$, then $N(T)=\Gamma^+$, as in (1), and the forbidden automorphism is that which respectively inverts and fixes the images $x$ and $y$ of the generators $X$ and $Y$ of $\Gamma^+$ (or equivalently, since $y^2=1$, inverts them both).

%%%%%%%%%%%%%%

\subsection{Some useful lemmas}

Here we note some useful consequences of the information in the preceding section on parent groups and forbidden automorphisms, either excluding various groups from membership of sets $\G(T)$ or finding relationships between different sets $\G(T)$. First we deal with abelian groups (see~\cite[Theorem 6.1]{STW} for a similar result concerning $\G^+(T)$):

\begin{lemma}\label{abelian}
An abelian group $A$ is in $\G(T)$ if and only if either
\begin{itemize}
\item $T=1$ and $A\cong C_2^r$ with $r=0, 1, 2$ or $3$, or
\item $T=2^{\sigma}$ and $A\cong C_2^r$ with $r=1$ or $2$, or
\item $T=3$ and  $A\cong C_2^r$ with $r=1, 2$ or $3$, or
\item $T=4^{\sigma}$ and $A\cong C_n$ or $C_2\times C_n$ where $n$ is even.
\end{itemize}
\end{lemma}

\noindent{\sl Proof.} For instance, if $T=3$ one can realise $A=C_2^r$ for $r=1, 2$ or $3$ by mapping $S_0$ to the identity element, and $S_1, S_2$ and $S_3$ to a generating set of non-identity elements. The rest of the proof is equally straightforward. \hfill$\square$

\begin{lemma}\label{odd}
If a non-trivial group $A$ of odd order is in $\G(T)$ then $T=5^{\sigma}$ for some $\sigma$, and $A$ is a $2$-generator non-abelian solvable group.
\end{lemma}

\noindent{\sl Proof.} If $T=1$, $2^{\sigma}$ or $3$ then $N(T)$ is generated by involutions, so it has no non-trivial quotients of odd order. If $T=2^{\sigma}{\rm ex}$ or $4^{\sigma}$ then the only odd-order quotients of $N(T)$ are cyclic, and these admit forbidden automorphisms. It follows that $T=5$, so $N(T)\cong F_2$ and $A$ is a $2$-generator group, which must be non-abelian to avoid admitting a forbidden automorphism inverting both generators. By the Feit-Thompson Theorem~\cite{FT}, $A$ is solvable.  \hfill$\square$

\medskip

We will show in Theorem~\ref{solv} that each set $\G(T)$ contains finite solvable groups of unbounded derived length. In fact, a small modification of the proof given there shows that when $T=5^{\sigma}$ we can take these to have odd order.

\medskip

The following lemma will prove useful later in realising various groups in certain classes, by showing that in many cases it is sufficient to consider just two classes, namely $T=1$ and $2^P{\rm ex}$.

\begin{lemma}\label{regmapslemma}
{\rm(a)} If $A$ is a non-abelian group in $\G(1)$, then $A\in\G(T)$ for each class $T=2^{\sigma}$, $3$ or $4^{\sigma}$.
\vskip2pt
{\rm(b)} If $A$ is any group in $\G(2^P{\rm ex})$, then $A\in\G(T)$ for each class $T=2^{\sigma}{\rm ex}$, $4^{\sigma}$ or $5^{\sigma}$.
\end{lemma}

\noindent{\sl Proof.} (a) The basic idea is to compose an epimorphism $\Gamma=N(1)\to A$, realising $A$ as an element of $\G(1)$, with an epimorphism $N(T)\to \Gamma$, chosen so that $A$, as a quotient of $N(T)$, does not have any of the forbidden automorphisms associated with the class $T$. We have $A={\rm Aut}\,\M$ for some regular map $\M$, corresponding to an epimorphism $\Gamma\to A,\, R_i\mapsto r_i$, so $A$ has generators $r_0, r_1$ and $r_2$ satisfying $r_i^2=(r_0r_2)^2=1$. Since $A$ is non-abelian, replacing $\M$ with its dual if necessary we can assume that $r_1r_2$ has order  $n>2$.

It is sufficient to consider one representative class $T$ from each orbit of $\Omega$. First let $T=2$. The presentations of $\Gamma$ and $N(2)$ show that there is an epimorphism $N(2)\to\Gamma$ sending $S_i$ to $R_{i-1}$ for $i=1, 2, 3$. Composing this with the epimorphism $\Gamma\to A$, we obtain an epimorphism $\theta: N(2)\to A$, with $S_i\mapsto s_i:=r_{i-1}$ for $i=1, 2, 3$. An automorphism of $A$ transposing $s_1$ and $s_2$, and fixing $s_3$, would send $s_1s_3=r_0r_2$, which has order dividing $2$, to $s_2s_3=r_1r_2$, which has order $n>2$. This is impossible, so $\ker\theta$ is not normal in $\Gamma$, and therefore corresponds to an edge-transitive map $\M'$ in class $2$ with ${\rm Aut}\,\M'\cong A$. 

A similar argument applies to class $4$, with $S_i\mapsto r_{i-1}$ as before for $i=1, 2$, and $S\mapsto s:=r_2$. Again no automorphism of $A$ inverting (equivalently, fixing) $s$ transposes $s_1$ and $s_2$, so the corresponding map is in class $4$. 

When $T=3$ there are several epimorphisms $\theta:N(3)\to A$ one could use. The following is convenient for combinatorial interpretation (see the examples later in this section). In this case we may assume that $r_0r_1$ has order $m>2$. Define $\theta$ by $S_i\mapsto s_i:=r_1, r_0, r_2, r_0$ for $i=0,\ldots, 3$. Any automorphism of $A$ inducing a double transposition on $s_0,\ldots, s_3$ would imply that $(r_0r_1)^2=1$, so ${\rm ker}\,\theta$ corresponds to a map $\M'$ in class $3$ with automorphism group $A$.

%Define an epimorphism $N(3)\to A$, with $S_i\mapsto s_i:=r_{i-1}$ for $i=1,\ldots, 3$ as before, and provided $s_1\ne s_2$ define $S_0\mapsto s_0:=r_2$. Then $s_0=s_3$, so the only double transposition of the generators $s_i$ is that transposing $s_1$ with $s_2$ and $s_0$ with $s_3$. As before, this cannot extend to an automorphism of $A$, so the resulting map is in class~$3$. If $s_1=s_2$, define $S_0\mapsto s_0:=r_0$ instead, so $s_0=s_1=s_2\ne s_3$ and there is no double transposition of the generators $s_i$.

{\rm(b)} The same basic idea as in (a) applies here, except that we now replace $N(1)=\Gamma$ with $N(2^P{\rm ex})=\Gamma^+$. We have $A={\rm Aut}\,\M$ for some chiral map $\M$, corresponding to an epimorphism $N(2^P{\rm ex})=\Gamma^+\to A,\, X\mapsto x, Y\mapsto y$, so that $A$ has generators $x$ and $y$ with $y^2=1$. Composing this with the epimorphism $N(5)\to N(2^P{\rm ex}),\, S\mapsto X, S'\mapsto Y$ we obtain an epimorphism $\theta: N(5)\to A,\, S\mapsto x, S'\mapsto y$. Since $\M$ is chiral, $x$ must have order $n>2$, so an automorphism of $A$ cannot transpose $y$ with $x^{\pm 1}$, and nor can it invert $x$ and $y$. A similar argument applies to class $4$, where we define $\theta: N(4)\to A$ by $S\mapsto x$ and $S_1, S_2\mapsto y$.
\hfill$\square$

\medskip

\noindent{\bf Remarks 1.} Lemma~\ref{abelian} shows why part~(a) of this lemma does not apply to abelian groups.

\medskip

{\bf 2.} Part~(a) does not extend to the remaining classes $T=2^{\sigma}\,{\rm ex}$ or $5^{\sigma}$: for instance, if $\M$ corresponds to a map subgroup $M$ contained in the commutator group $ \Gamma'$ of $\Gamma$ then $A$ requires three generators and hence cannot be a quotient of a $2$-generator group $N(2^{\sigma}{\rm ex})$ or $N(5^{\sigma})$. Similarly, part~(b) does not extend to the classes~$T=1$, $2^{\sigma}$, or $3$, since $\Gamma^+$ has non-trivial quotients of odd order, whereas $N(T)$ does not for these classes $T$.

\medskip

{\bf 3.} We will give a more constructive interpretation of Lemma~\ref{regmapslemma} in the next section, where we will use it to create specific examples of maps in various classes $T$ from known maps in classes $1$ or $2^P{\rm ex}$.

\medskip

{\bf 4.} The smallest non-abelian group $A\in\G(1)$ is the symmetric group $S_3$, generated by involutions $r_i$ with $r_0=r_2\ne r_1$, realised as ${\rm Aut}\,\M$ where

\medskip

Part (c) of the following lemma shows that class~$2$ can also play a similar role to that played by classes~$1$ and $2^P{\rm ex}$ in Lemma~\ref{regmapslemma}. An element $x$ of a group $A$ is {\em strongly real\/} if $x^i=x^{-1}$ for some involution $i\in A$, or equivalently, if $x$ is a product of at most two involutions in $A$.

\begin{lemma}\label{class2lemma}
{\rm (a)} If a group $A$ is generated by involutions $a, b$ and $c$, where $ab$, $ac$ and $bc$ do not all have the same order, then $A\in\G(2)$.
\vskip2pt
{\rm(b)} If a group $A$ is generated by elements $a, b$ and $c$ satisfying $abc=1$, where $a$ is an involution, $b$ is a product of two involutions, and no automorphism of $A$ inverts $c$, then $A\in\G(2)$.
\vskip2pt
{\rm(c)} If $A\in\G(2)$ then $A\in\G(T)$ for each class $T=2^{\sigma}$, $3$ or $4^{\sigma}$.
\end{lemma}

\noindent{\sl Proof.} (a) Without loss of generality we may assume that $ac$ and $bc$ have different orders. There is an epimorphism $N(2)\to A$ given by $S_i\mapsto s_i:=a, b$ or $c$ for $i=1, 2$ or $3$, and no automorphism of $A$ can transpose $s_1$ and $s_2$ while fixing $s_3$.

(b) We can write $b=s_1s_2$ for involutions $s_1, s_2\in A$, and define $s_3=a$, giving involutions $s_1, s_2, s_3$ generating $A$. If an automorphism transposes $s_1$ and $s_2$ while fixing $s_3$, it inverts $a$ and $b$, so composing it with conjugation by $a$ gives an automorphism inverting $c$. However, no such automorphism exists, so $A\in\G(2)$.

(c) The group $N(4)$ can be obtained from $N(2)$ by taking $S=S_3$ and omitting the relation $S_3^2=1$, giving an epimorphism $N(4)\to N(2)$. Thus any quotient $A=\langle s_1, s_2, s_3\rangle$ of $N(2)$ is also a quotient of $N(4)$, with $s:=s_3$; the forbidden automorphisms are the same in each case, since $s$ has order dividing $2$, so $\G(2)\subseteq\G(4)$. Using $\Omega$ then allows $A$ to be realised in all classes $2^{\sigma}$ and $4^{\sigma}$. For class~$3$ one can extend an epimorphism $N(2)\to A$ to $N(3)\to A$ by mapping $S_0$ to $s_3$; again, this introduces no further forbidden automorphisms.  \hfill$\square$

\medskip

The following observation is easily proved:

\begin{lemma}\label{class3lemma}
Suppose that $A$ is generated by involutions $s_0, s_1, s_2, s_3$, and that for at least two of the three partitions $ij\mid kl$ of $\{0, 1, 2, 3\}$ the products $s_is_j$ and $s_ks_l$ have different orders. Then $A\in\G(3)$.  \hfill$\square$
\end{lemma}

\begin{lemma}
Let $A={\rm Aut}\,\M=\langle x, y\mid y^2=1, \ldots\rangle\in\G(2^P{\rm ex})$ for some chiral map $\M$, where the image $x$ of $X$ is strongly real. Then $A\in\G(T)$ for each edge-transitive class $T\ne 1$.
\end{lemma}

\noindent{\sl Proof.} An involution $a\in A$ inverts $x$, so $A=\langle s_1, s_2, s_3\rangle$ where $s_1=a$, $s_2=ax$ and $s_3=y$ all satisfy $s_i^2=1$, and hence $A$ is a quotient of $N(2)$. If an automorphism of $A$ transposes $s_1$ and $s_2$ and fixes $s_3$, it inverts $x$ and fixes $y$, contradicting the chirality of $\M$. Thus no such automorphism exists, so $A\in\G(2)$. Then $A\in\G(T)$ for each $T=2^{\alpha}$, $3$ or $4^{\alpha}$ by Lemma~\ref{class2lemma}(c), and also for $T=2^{\alpha}{\rm ex}$ or $5^{\alpha}$ by Lemma~\ref{regmapslemma}(b), so $A\in\G(T)$ for each $T\ne 1$. \hfill$\square$

%%%%%%%%%%%%%%%

\subsection{Constructing new maps from old ones}

%{\color{blue}Mention functors between permutational categories, Lynne's work?]}

%\medskip

Although the main purpose of Lemma~\ref{regmapslemma} was to simplify the task of proving that various groups $A$ are in $\G(T)$ for certain classes $T$, its proof can also be used to produce specific examples of maps $\M'$ in those classes with automorphism group $A$. Specifically, in Lemma~\ref{regmapslemma}(a) one replaces the flags of a regular map $\M$, permuted regularly by $A:={\rm Aut}\,\M$, with copies of ${\mathcal N}(T)$, regarded as fundamental regions for $A$ and also permuted regularly by $A$; these copies of ${\mathcal N}(T)$ are joined together by certain generators of $A$, namely the images under the epimorphism $\theta:N(T)\to A$ of the generators of $N(T)$ given in Proposition~\ref{parents}, to form the required map $\M'$ in class $T$. Similarly in Lemma~\ref{regmapslemma}(b) it is the directed edges of a chiral map $\M$, again permuted regularly by $A$, which are replaced with copies of ${\mathcal N}(T)$ to form the map $\M'$. We will give examples with $T=2, 3, 4$ and $5$; further examples, in other classes, can be obtained from these by applying suitable operations in $\Omega$.

\begin{figure}[h!]
\begin{center}
\begin{tikzpicture}[scale=0.5, inner sep=0.8mm]

\node (c) at (0,0) [shape=circle, fill=black] {};
\node (d) at (8,0) [shape=circle, fill=black] {};
\draw [thick] (c) to (d);
\draw [thick] (c) to (1,-3);
\draw [thick] (c) to (1,3);
\draw [thick] (d) to (7,-3);
\draw [thick] (d) to (7,3);
\draw [thick] (c) to (-2.5,2.5);
\draw [thick] (c) to (-2.5,-2.5);
\draw [thick] (d) to (10.5,2.5);
\draw [thick] (d) to (10.5,-2.5);

\draw (c) to (1,0.5);
\draw (c) to (1,-0.5);
\draw (1,0.5) to (1,-0.5);
\draw (d) to (7,0.5);
\draw (d) to (7,-0.5);
\draw (7,0.5) to (7,-0.5);
%\draw (c) to (0.8,0.8);
%\draw (0.3,1) to (0.8,0.8);

\node at (-0.8,0) {$v$};
\node at (4,-0.4) {$e$};
\node at (4,2) {$f$};

\node at (1.5,0.5) {$\phi$};

\node at (4,-3) {$\M$};

%%%%%%%%

\node at (12,0) {$\longrightarrow$};

\node (C) at (15,0) [shape=circle, fill=black] {};
\node (D) at (23,0) [shape=circle, fill=black] {};
\node (E) at (19,0) [shape=circle, draw] {};

\draw [thick] (C) to (E) to (D);
\draw [thick, dotted] (C) to (19,2) to (D) to (19,-2) to (C);
\draw [thick, dotted] (19,2) to (E) to (19,-2);

\node at (19,-3) {$\M'$};

%%%%%%%%%%%%%

\node (A) at (7,-11) [shape=circle, fill=black] {};
\node (B) at (17,-11) [shape=circle, draw] {};
\draw [thick] (A) to (B);
\draw [thick, dotted] (A) to (17,-5) to (B);

\draw [<->] (13,-10) to (13,-12);
\draw [<->] (16,-8) to (18,-8);
\draw [<->] (11.8,-6.9) to (13,-8.5);

\node at (10,-7) {$s_3=r_1$};
\node at (20,-8) {$s_1=r_0$};
\node at (13,-13) {$s_3=r_2$};

%%%%%%%%%%%

\node (b) at (0,-5) [shape=circle, fill=black] {};
\node (c) at (0,-11) [shape=circle, fill=black] {};
\draw [very thick] (0,-11) arc (-90:90:3);
\draw [dashed] (0,-5) arc (90:270:3);
\node at (0,-13) {${\mathcal N}(2)$};

\end{tikzpicture}

\end{center}
\caption{Flags $\phi$ of $\M$ replaced with triangular copies of ${\mathcal N}(2)$} 
\label{Class1toClass2}
\end{figure}

\medskip

\noindent{\bf Example 1.} Let $A={\rm Aut}\,\M$ for some regular map $\M$ of valency $n>2$. By composing the epimorphism $N(2)\to\Gamma,\, S_i\mapsto R_{i-1}$ used in the proof of Lemma~\ref{regmapslemma}(a) with the epimorphism $\Gamma\to A$, we obtain an edge-transitive map $\M'$ in class $2$ with ${\rm Aut}\,\M'\cong A$. This is the Walsh bipartite map $W(\M)$ of $\M$, where the latter is now regarded as a hypermap: it is formed by placing a new vertex of valency $2$ at the midpoint of each edge of $\M$ or of valency $1$ at the free end of a free edge of $\M$. (More generally, the Walsh map $W({\mathcal H})$ of a hypermap $\mathcal H$ has vertices, traditionally coloured black and white, at the hypervertices and hyperedges of $\mathcal H$, with edges between them indicating incidence: see~\cite{Wal} for details.) Since $A$ acts transitively on the directed edges and the faces of $\M$, it also acts transitively on the edges and faces of $\M'$, whereas it has two orbits on the vertices and on the flags of $\M'$. In this construction, each flag of $\M$ is replaced with a right-angled triangle (a homeomorphic copy of ${\mathcal N}(2)$), with a new vertex at the right angle and an old one at another corner (see Figure~\ref{Class1toClass2}); the permutations $r_i$ used to join flags of $\M$ together are then replaced with the permutations $s_{i-1}$ joining these triangles along paired sides.

\begin{figure}[h!]
\begin{center}
\begin{tikzpicture}[scale=0.4, inner sep=0.8mm]

\node (a) at (-10,0) [shape=circle, fill=black] {};
\draw [thick] (-5,0) to (a) to (-15,0);
\draw [thick] (-10,5) to (a) to (-10,-5);

\node at (-7.5,-2.5) {$\M$};

%%%%%%

\node (b) at (10,0) [shape=circle, fill=black] {};
\node (c) at (5,0) [shape=circle, draw] {};
\node (d) at (15,0) [shape=circle, draw] {};
\node (e) at (10,5) [shape=circle, draw] {};
\node (f) at (10,-5) [shape=circle, draw] {};
\draw [thick] (c) to (b) to (d);
\draw [thick] (e) to (b) to (f);

\node at (12.5,-2.5) {$\M'$};

\end{tikzpicture}

\end{center}
\caption{Maps $\M$ and $\M'$, where $\M$ has four free edges.} 
\label{Free}
\end{figure}

\medskip

If $\M$ has a free edge, then by regularity all its edges are free, and $\M$ is a map on the sphere with one vertex, one face, and $n$ free edges; in this case $A\cong D_n$, generated by involutions $r_0=r_2$ and $r_1$, and $\M'$ is formed by completing each free edge with a white vertex, as in Figure~\ref{Free} where $n=4$. Otherwise, $\M$ has no free edges, and if it has type $\{m,n\}$ then $\M'$ has old and new vertices of valency $n$ and $2$, while its faces are all $2m$-gons. In either case, the underlying surfaces of these two maps are the same.

\medskip

\noindent{\bf Example 2.} If $\M$ is a chiral map, in class $2^P{\rm ex}$, then $W(\M)$ is an edge-transitive map $\M'$ in class~$5$, corresponding to the composition of epimorphisms $N(5)\to\Gamma^+\to A$ used in proving Lemma~\ref{regmapslemma}(b). If $\M$ has type $\{m,n\}$ then $\M'$ has old and  new vertices of valencies $n$ and $2$, and its faces are all $2m$-gons; the underlying surface is that of $\M$. For instance, if we add a new vertex at the midpoint of each edge of the map $\M$ in Figure~\ref{torusK5} to form $\M'=W(\M)$, we can use the same fundamental region $\mathcal F$ for $A$ to see that $\M'/{\rm Aut}\,\M'\cong{\mathcal N}(5)$. Similarly, when $T=4$ in Lemma~\ref{regmapslemma}(b), the half-turn $y$ about the mid-point of an edge of a chiral map $\M$ induces a reflection of the corresponding edge of the Petrie dual $P(\M)$, so $WP(\M)$ is an edge-transitive map $\M'$ in class~$4$. The embedded graph is the same as that in the case $T=5$, but the underlying surface is that of $P(\M)$, and the faces have valency twice those of $P(\M)$.

\medskip

\noindent{\bf Example 3.} In Lemma~\ref{regmapslemma}(a), when $T=3$ we place two vertices either side of the mid-point of each edge of the regular map $\M$, one in each incident face, join each of them to the two incident vertices, and finally delete all edges of $\M$. This corresponds to replacing each flag of $\M$ with a triangular copy of ${\mathcal N}(3)$, as in Figure~\ref{Class1toClass3}. If $\M$ has type $\{m,n\}$ then $\M'$ has old and new vertices of valencies $2n$ and $2$, and old and new faces of valencies $2m$ and $4$. The underlying surfaces are the same.

%%%%%%%%%%%%%%%%%%%%%

\begin{figure}[h!]
\begin{center}
\begin{tikzpicture}[scale=0.5, inner sep=0.8mm]

\node (c) at (0,0) [shape=circle, fill=black] {};
\node (d) at (8,0) [shape=circle, fill=black] {};
\draw [thick] (c) to (d);
\draw [thick] (c) to (1,-3);
\draw [thick] (c) to (1,3);
\draw [thick] (d) to (7,-3);
\draw [thick] (d) to (7,3);
\draw [thick] (c) to (-2.5,2.5);
\draw [thick] (c) to (-2.5,-2.5);
\draw [thick] (d) to (10.5,2.5);
\draw [thick] (d) to (10.5,-2.5);

\draw (c) to (1,0.5);
\draw (c) to (1,-0.5);
\draw (1,0.5) to (1,-0.5);
\draw (d) to (7,0.5);
\draw (d) to (7,-0.5);
\draw (7,0.5) to (7,-0.5);
%\draw (c) to (0.8,0.8);
%\draw (0.3,1) to (0.8,0.8);

\node at (-0.8,0) {$v$};
\node at (4,-0.4) {$e$};
\node at (4,2) {$f$};

\node at (1.5,0.5) {$\phi$};

\node at (4,-3) {$\M$};

%%%%%%%%

\node at (12,0) {$\longrightarrow$};

\node (C) at (15,0) [shape=circle, fill=black] {};
\node (D) at (23,0) [shape=circle, fill=black] {};
\node (E) at (19,1) [shape=circle, draw] {};
\node (F) at (19,-1) [shape=circle, draw] {};

\draw [thick] (C) to (E) to (D) to (F) to (C);
\draw [thick, dotted] (C) to (19,2) to (D) to (19,-2) to (C);
\draw [thick, dotted] (C) to (D);
\draw [thick, dotted] (19,2) to (E) to (F) to (19,-2);

\node at (19,-3) {$\M'$};

%%%%%%%%%%%%%

\node (A) at (7,-11) [shape=circle, fill=black] {};
\node (B) at (17,-8) [shape=circle, draw] {};
\draw [thick] (A) to (B);
\draw [thick, dotted] (A) to (17,-5) to (B) to (17,-11) to (A);

\draw [<->] (13,-10) to (13,-12);
\draw [<->] (16,-6.5) to (18,-6.5);
\draw [<->] (16,-9.5) to (18,-9.5);
\draw [<->] (11.8,-6.9) to (13,-8.5);

\node at (10,-7) {$s_0=r_1$};
\node at (20,-6.5) {$s_1=r_0$};
\node at (20,-9.5) {$s_3=r_0$};
\node at (13,-13) {$s_2=r_2$};

%%%%%%%%

\node (i) at (-3,-9) [shape=circle, fill=black] {};
\node (j) at (3,-9) [shape=circle, fill=black] {};
\draw [very thick] (i) to (j);
\draw [dashed] (3,-9) arc (0:360:3);
\node at (0,-13) {${\mathcal N}(3)$};

\end{tikzpicture}

\end{center}
\caption{Flags $\phi$ of $\M$ replaced with triangular copies of ${\mathcal N}(3)$} 
\label{Class1toClass3}
\end{figure}

%%%%%%%%

\bigskip

\begin{figure}[h!]
\begin{center}
\begin{tikzpicture}[scale=0.4, inner sep=0.8mm]

\node (c) at (0,0) [shape=circle, fill=black] {};
\node (d) at (8,0) [shape=circle, fill=black] {};
\draw [thick] (c) to (d);
\draw [thick] (c) to (1,-3);
\draw [thick] (c) to (1,3);
\draw [thick] (d) to (7,-3);
\draw [thick] (d) to (7,3);
\draw [thick] (c) to (-2.5,2.5);
\draw [thick] (c) to (-2.5,-2.5);
\draw [thick] (d) to (10.5,2.5);
\draw [thick] (d) to (10.5,-2.5);

\draw (c) to (1,0.5);
\draw (c) to (1,-0.5);
\draw (1,0.5) to (1,-0.5);
\draw (d) to (7,0.5);
\draw (d) to (7,-0.5);
\draw (7,0.5) to (7,-0.5);
%\draw (c) to (0.8,0.8);
%\draw (0.3,1) to (0.8,0.8);

\node at (-0.8,0) {$v$};
\node at (4,-0.4) {$e$};
\node at (4,2) {$f$};
\node at (4,-2) {$f'$};
\node at (-0.5,3) {$f''$};
\node at (1.5,0.5) {$\phi$};

\node at (4,-4) {$\M$};

%%%%%%%%

\node at (12,0) {$\longrightarrow$};

%\node (C) at (15,0) [shape=circle, fill=black] {};
%\node (D) at (23,0) [shape=circle, fill=black] {};
\node (E) at (19,2) [shape=circle, draw] {};
\node (F) at (19,-2) [shape=circle, draw] {};
\node (G) at (15,3) [shape=circle, draw] {};

%\draw [thick] (C) to (E) to (D) to (F) to (C);
%\draw [thick, dotted] (C) to (19,2) to (D) to (19,-2) to (C);
\draw [thick] (F) to [out=120, in=240] (E);
\draw [thick] (F) to [out=60, in=300] (E);
%\draw [thick, dotted] (15,0) to (23,0);
\draw [thick] (G) to [out=20, in=140] (E);
\draw [thick] (G) to [out=320, in=190] (E);
%\draw [thick, dotted] (19,2) to (E) to (F) to (19,-2);
\node at (20,1.9) {$v_f$};
\node at (20.1,-2.1) {$v_{f'}$};
\node at (13.8,3) {$v_{f''}$};
\node at (19,-4) {$2D(\M)$};

%%%%%%%%

\node at (23,0) {$\longrightarrow$};

%\node (C) at (15,0) [shape=circle, fill=black] {};
%\node (D) at (23,0) [shape=circle, fill=black] {};
\node (E) at (29,2) [shape=circle, draw] {};
\node (F) at (29,-2) [shape=circle, draw] {};
\node (G) at (25,3) [shape=circle, draw] {};
\node (H) at (29.7,0) [shape=circle, fill=black] {};
\node (H) at (28.3,0) [shape=circle, fill=black] {};
\node (I) at (27.2,3.1) [shape=circle, fill=black] {};
\node (J) at (26.9,1.95) [shape=circle, fill=black] {};

%\draw [thick] (C) to (E) to (D) to (F) to (C);
%\draw [thick, dotted] (C) to (19,2) to (D) to (19,-2) to (C);
\draw [thick] (F) to [out=120, in=240] (E);
\draw [thick] (F) to [out=60, in=300] (E);
%\draw [thick, dotted] (15,0) to (23,0);
\draw [thick] (G) to [out=20, in=140] (E);
\draw [thick] (G) to [out=320, in=190] (E);
%\draw [thick, dotted] (19,2) to (E) to (F) to (19,-2);
\node at (30,1.9) {$v_f$};
\node at (30.1,-2.1) {$v_{f'}$};
\node at (23.8,3) {$v_{f''}$};
\node at (29,-4) {$\M'$};

%%%%%%%%%%%%%

\node (A) at (25,-11) [shape=circle, fill=black] {};
\node (B) at (27,-8) [shape=circle, draw] {};
\draw [thick] (A) to (B);
\draw [thick, dotted] (B) to (23,-11) to (A) to (27,-11) to (B);

%\draw [<->] (24,-11.5) to (26,-12);
%\draw [<->] (24,-12) to (26,-11.5);
\draw [<->] (26.5,-9.5) to (27.5,-9.5);
\draw [<->] (24.6,-9.1) to (25.4,-9.9);
\draw  (26,-11.2) arc (0:-180:1);
\draw (26.2,-11.4) to (26,-11.2) to (25.8,-11.4);
\draw (24.2,-11.4) to (24,-11.2) to (23.8,-11.4);

%\node at (10,-7) {$s_0=r_1$};
\node at (29.5,-9.5) {$s_1=r_0$};
\node at (25,-13) {$s=r_2$};
\node at (22.5,-9) {$s_2=r_1$};

%%%%%%%%%%%%

\node (k) at (4,-10) [shape=circle, fill=black] {};
\node (l) at (7,-10) [shape=circle, fill=black] {};
\draw [very thick] (k) to (l);
\draw [dashed] (7,-10) arc (0:360:3);
\node at (9,-13) {${\mathcal N}(4)$};

\end{tikzpicture}

\end{center}
\caption{Construction of $\M'$ of class $4$ from a regular map $\M$} 
\label{Class1toClass4}
\end{figure}

\medskip

\noindent{\bf Example 4.} When $T=4$ in Lemma~\ref{regmapslemma}(a), we first take the dual map $D(\M)$, and form its edge-double $2D(\M)$, replacing each edge of $D(\M)$ with a pair of edges between the same two vertices, enclosing a digonal face. We now take the Petrie dual $P(2D(\M))$ of this map (so that the reflection $r_2$ of $\M$ is replaced with a half-turn $s$), and finally form the Walsh map $\M'=W(P(2D(\M)))$ by placing a new vertex of valency $2$ at the midpoint of each edge of $P(2D(\M))$, as shown in Figure~\ref{Class1toClass4}. If $\M$ has type $\{m,n\}$ then $\M'$ has a vertex $v_f$ of valency $2m$ for each face $f$ of $\M$, and two vertices of valency $2$ for each edge of $\M$. The faces of $\M'$, corresponding to the Petrie polygons of $2D(\M)$, are $4m$-gons: for each face $f$ of $\M$ there is a face of $\M'$ bounded by a circuit which successively follows, in cyclic order around $v_f$, the $m$ $4$-gons along the pairs of doubled edges in $2D(\M)$ between $v_f$ and its neighbouris $v_{f'}$.
%{\color{blue}[Underlying surface? Specific example?]}

\medskip

\noindent{\bf Example 5.} As in the case of Lemma~\ref{regmapslemma}, there is also a combinatorial interpretation of Lemma~\ref{class2lemma}(c). For instance, the epimorphism $N(4)\to N(2)$ corresponds to replacing a map $\M$ in class~$2$ with the map $\M'=WPD(\M)$ in class~$4$, while the epimorphism $N(3)\to N(2)$ corresponds to  taking $\M'=WD(\M)$.
%{\color{blue}[Check both. Diagrams?]}

%%%%%%%%%%%%%%%%%%

\subsection{The medial map construction}\label{medialmaps}

If $\M$ is any map, then its {\sl medial map\/} $\M^{\rm med}$ is a map constructed from $\M$, on the same surface as $\M$ (see, for instance,~\cite{ASS}). This operation is relevant here as another way of constructing edge-transitive maps: for instance, if $\M$ is in class $1$ or $2^{\sigma}{\rm ex}$ then  $\M^{\rm med}$ is edge-transitive.

For any map $\M$, the vertices of $\M^{\rm med}$ are the midpoints $v_e$ of the edges $e$ of $\M$ (or their free ends in the case of free edges), with $v_e$ and $v_{e'}$ joined by an edge within $f$ whenever $e$ and $e'$ are consecutive edges round a face $f$ of $\M$. See Figure~\ref{medial} for this construction near a particular edge of $\M$. 

\begin{figure}[h!]
\begin{center}
\begin{tikzpicture}[scale=0.4, inner sep=0.8mm]

\node (a) at (-4,0) [shape=circle, fill=black] {};
\node (b) at (4,0) [shape=circle, fill=black] {};
\node (A) at (0,0) [shape=circle, draw] {};
\node (B) at (-3,3) [shape=circle, draw] {};
\node (C) at (3,3) [shape=circle, draw] {};
\node (D) at (3,-3) [shape=circle, draw] {};
\node (E) at (-3,-3) [shape=circle, draw] {};
\node (F) at (-7,2) [shape=circle, draw] {};
\node (G) at (7,2) [shape=circle, draw] {};
\node (H) at (7,-2) [shape=circle, draw] {};
\node (I) at (-7,-2) [shape=circle, draw] {};

\draw [thick] (a) to (A) to (b);
\draw [thick] (a) to (B) to (-2.5,4.5);
\draw [thick] (a) to (E) to (-2.5,-4.5);
\draw [thick] (b) to (C) to (2.5,4.5);
\draw [thick] (b) to (D) to (2.5,-4.5);
\draw [thick, dashed] (F) to (B) to (A) to (E) to (I) to (F);
\draw [thick, dashed] (D) to (A) to (C) to (G) to (H) to (D);

\draw [thick] (-8.5,3) to (F) to (a) to (I) to (-8.5,-3);
\draw [thick] (8.5,3) to (G) to (b) to (H) to (8.5,-3);

\node (c) at (15,2) [shape=circle, fill=black] {};
\node (d) at (15,-2) [shape=circle, draw] {};
\draw [thick] (13,1) to (c) to (17,3);
\draw [thick, dashed] (13,-3) to (d) to (17,-1);

\node at (16,1) {$\M$};
\node at (16,-3) {$\M^{\rm med}$};

\end{tikzpicture}

\end{center}
\caption{The medial map $\M^{\rm med}$ of a map $\M$} 
\label{medial}
\end{figure}

If $\M$ has empty boundary, then all vertices of $\M^{\rm med}$ have valency $4$. The map $\M^{\rm med}$ is $2$-face colourable, with its faces corresponding to the vertices and faces of $\M$ and inheriting their valencies. A map and its dual have isomorphic medial maps. For example, the medial map of the tetrahedron is the octahedron, that of the cube (see Figure~\ref{cuboctahedron}) and the octahedron is the cuboctahedron, and that of the icosahedron and dodecahedron is the icosidodecahedron, with $20$ triangular and $12$ pentagonal faces.

\begin{figure}[h!]
\begin{center}
\begin{tikzpicture}[scale=0.3, inner sep=0.8mm]

\node (a) at (0,0) [shape=circle, fill=black] {};
\node (b) at (10,0) [shape=circle, fill=black] {};
\node (c) at (10,10) [shape=circle, fill=black] {};
\node (d) at (0,10) [shape=circle, fill=black] {};
\node (e) at (13,4) [shape=circle, fill=black] {};
\node (f) at (13,14) [shape=circle, fill=black] {};
\node (g) at (3,14) [shape=circle, fill=black] {};

\node (A) at (5,0) [shape=circle, draw] {};
\node (B) at (10,5) [shape=circle, draw] {};
\node (C) at (5,10) [shape=circle, draw] {};
\node (D) at (0,5) [shape=circle, draw] {};
\node (E) at (11.5,2) [shape=circle, draw] {};
\node (F) at (13,9) [shape=circle, draw] {};
\node (G) at (8,14) [shape=circle, draw] {};
\node (H) at (1.5,12) [shape=circle, draw] {};
\node (I) at (11.5,12) [shape=circle, draw] {};

\draw [thick] (a) to (A) to (b) to (B) to (c) to (C) to (d) to (D) to (a);
\draw [thick] (b) to (E) to (e) to (F) to (f) to (G) to (g) to (H) to (d);
\draw [thick] (c) to (I) to (f);
\draw [thick, dashed] (A) to (B) to (C) to (D) to (A);
\draw [thick, dashed] (I) to (G) to (H) to (C) to (I);
\draw [thick, dashed] (I) to (B) to (E) to (F) to (I);
\draw [thick, dashed] (A) to (E);
\draw [thick, dashed] (D) to (H);

\node (u) at (25,10) [shape=circle, fill=black] {};
\node (v) at (25,6) [shape=circle, draw] {};
\draw [thick] (23,9) to (u) to (27,11);
\draw [thick, dashed] (23,5) to (v) to (27,7);

\node at (31,10) {$\M=$ cube};
\node at (35,6) {$\M^{\rm med}=$ cuboctahedron};

\end{tikzpicture}

\end{center}
\caption{Solid views of the cube and cuboctahedron} 
\label{cuboctahedron}
\end{figure}

The automorphism groups of $\M^{\rm med}$ and $\M$ are the same, that is, $N_{\tilde\Gamma}(M)=N_{\Gamma}(M)$, unless $\M$ is self-dual, in which case the duality of $\M$ induces an extra automorphism of $\M^{\rm med}$, transposing its two classes of faces, so that ${\rm Aut}\M^{\rm med}$ contains ${\rm Aut}\,\M$ with index $2$. For example, the cuboctahedron is the medial map of the cube and of its dual, the octahedron, and has the same automorphism group as them, isomorphic to $S_4\times C_2$, whereas the octahedron is the medial map of the self-dual tetrahedron, which has a smaller automorphism group $S_4$.

As shown by \v Sir\'a\v n, Tucker and Watkins in~\cite[Lemma~2.2]{STW}, if $\M$ is edge-transitive, then $\M^{\rm med}$ is edge-transitive if and only if $\M$ is in class~$1$ or $2^{\sigma}{\rm ex}$ for some $\sigma=\emptyset, *$ or $P$. This can be seen by constructing the medial maps ${\mathcal N}(T)^{\rm med}$ of the $14$ basic maps ${\mathcal N}(T)$, and noticing that only those for $T=1$ and $2^{\sigma}{\rm ex}$ have a single edge: thus ${\mathcal N}(1)^{\rm med}\cong{\mathcal N}(2^*)$, ${\mathcal N}(2\,{\rm ex})^{\rm med}\cong {\mathcal N}(2^*{\rm ex})^{\rm med}\cong {\mathcal N}(4^*)$ and ${\mathcal N}(2^P{\rm ex})^{\rm med}\cong {\mathcal N}(5^*)$, all with a single edge (see Figure~\ref{basicmaps}), whereas, for example, ${\mathcal N}(3)^{\rm med}$ is a map on the disc with one vertex in the interior and four free edges, shown on the {\color{red}right} in Figure~\ref{mapsA4N0} (this map reappears as the map ${\mathcal A}_4$ in Figure~\ref{mapA4}).

\begin{figure}[h!]
\begin{center}
\begin{tikzpicture}[scale=0.5, inner sep=0.8mm]

\node (a) at (0,3.2) {};
\node (b) at (3,0) [shape=circle, fill=black] {};
\node (c) at (-3,0) [shape=circle, fill=black] {};
\node (d) at (0,-3.2) {};
\draw[thick]  (b) to (c);
\draw [dashed] (3,0) arc (0:360:3);

\node at (-4,-3) {${\mathcal N}(3)$};

%%%%%

\node (E) at (10,0) [shape=circle, draw] {};
\draw [dashed] (13,0) arc (0:360:3);
\draw [very thick, dashed] (8,2) to (E);
\draw [very thick, dashed] (8,-2) to (E);
\draw [very thick, dashed] (12,2) to (E);
\draw [very thick, dashed] (12,-2) to (E);

\node at (6,-3) {${\mathcal N}(3)^{\rm med}$};

\end{tikzpicture}

\end{center}
\caption{The maps ${\mathcal N}(3)$ and ${\mathcal N}(3)^{\rm med}$ on the closed disc} 
\label{mapsA4N0}
\end{figure}

In Section~\ref{almostET} we will show that it is possible for $\M^{\rm med}$ to be edge-transitive even when $\M$ is not. We will also show that the medial map construction is induced by the inclusion of $\Gamma$, regarded as the extended triangle group $\Delta[\infty,2,\infty]$, as a subgroup of index $2$ in the extended triangle group $\Delta[4,2,\infty]$.

%%%%%%%%%%%%%%

\section{Edge-transitive embeddings of the complete graphs}\label{completemaps}

Although it is not the aim of this paper to give a systematic account of edge-transitive embeddings of specific families of graphs, we will do this here for illustrative purposes in the case of the complete graphs, partly because most of the required results have already appeared in the literature.

In 1971 Biggs~\cite{Big} proved:

\begin{thm}
The complete graph $K_n$ has an orientably  regular embedding if and only if $n$ is a prime power.  \hfill$\square$
\end{thm}

The maps Biggs constructed to prove that this condition is sufficient are Cayley maps $\M_n(c)$ for the additive groups of finite fields $\F_n$; the generating set is the multiplicative group $\F_n^*=\F_n\setminus\{0\}$, taken in the cyclic order $1, c, c^2, \ldots, c^{n-2}$ where $c$ is a primitive element for $\F_n$ (that is, a generator for the cyclic group $\F_n^*$).  

\medskip

\noindent{\bf Example} For $n=4$ we have $F_4=\{0, 1, c, c^2=c+1\}$, and the corresponding Cayley map $\M_4(c)$ is the tetrahedral map shown in Figure~\ref{tetmap}.

\begin{figure}[h!]

\begin{center}
 \begin{tikzpicture}[scale=0.2, inner sep=1mm]

\node (A) at (0,0) [shape=circle, fill=black] {};
\node (B) at (0,10) [shape=circle, fill=black] {};
\node (C) at (8.5,-5) [shape=circle, fill=black] {};
\node (D) at (-8.5,-5) [shape=circle, fill=black] {};

\draw [thick] (A) to (B) to (C) to (D) to (A) to (C);
\draw [thick] (B) to (D);

\node at (0,-2) {$0$};
\node at (2,10) {$1$};
\node at (-10.5,-5) {$c$};
\node at (11,-4.5) {$c^2$};

\end{tikzpicture}

\end{center}
\caption{Spherical embedding $\M_4(c)$ of $K_4$ (drawn in the plane)} \label{tetmap}
\end{figure}

\medskip

In 1985 James and the author proved that the maps $\M_n(c)$ constructed by Biggs  are the only orientably regular embeddings of complete graphs:

\begin{thm}\label{orKn}
If $\M$ is an orientably regular embedding of $K_n$ then $\M\cong\M_n(c)$ for some primitive element $c$ of $\F_n$. Moreover, $\M_n(c)$ and $\M_n(c')$ are isomorphic (as oriented maps) if and only if $c$ and $c'$ are equivalent under an automorphism of $\F_n$.  \hfill$\square$
\end{thm}

If follows that if $n=p^e$ for some prime $p$ then there are, up to isomorphism, precisely $\phi(n-1)/e$ orientably regular embeddings of $K_n$, one for each orbit of the Galois group of $\F_n$ (isomorphic to $C_e$, generated by the Frobenius automorphism $t\mapsto t^p$) on the $\phi(n-1)$ primitive elements of the field. For each $n=2, 3$ or $4$ we obtain a single map on the sphere, which is regular (in class~$1$) and of type $\{2,1\}$, $\{3,2\}$ or $\{3,3\}$ respectively; for instance, the last of these is the tetrahedral map, with automorphism group $S_4$.

These maps have orientation-preserving automorphism group $AGL_1(\F_n)$, acting naturally on the vertices by affine transformations $v\mapsto av+b$ where $a, b\in\F_n$ and $a\ne 0$. This group is a semidirect product $\F_n\rtimes \F_n^*$ of the additive group of the field, inducing translations $v\mapsto v+b$, by the multiplicative group, acting by multiplication $v\mapsto av$. For $n>4$ these maps occur in chiral pairs $\M_n(c^{\pm 1})$ and hence lie in class $2^P{\rm ex}$ with full automorphism group $AGL_1(\F_n)$. (Note that in this case $\M_n(c)$ and $\M_n(c^{-1})$ are isomorphic as {\sl unoriented\/} maps, but not as {\sl oriented\/} maps.) For $n=2$, $3$ and $4$ these maps are regular (in class $1$) with automorphism groups $AGL_1(\F_2)\times C_2\cong V_4$, $AGL_1(\F_3)\times C_2\cong S_3\times C_2\cong D_6$ and $AGL_1(\F_4)\rtimes C_2\cong A_4\rtimes C_2\cong S_4$.

For $n\ge 4$ the maps $\M_n(c)$ have type $\{m,n-1\}$ where $m=(n-1)/2$ or $n-1$ as $n\equiv 3$ mod~$(4)$ or not, and have genus $(n^2-7n+4)/4$ or $(n-1)(n-4)/4$ respectively. They have Petrie length $2p$. For example, if $n=5$ or $n=7$ there is a single chiral pair of maps; one of each pair, namely $\M_5(2)$ and $\M_7(3)$, is illustrated in Figure~\ref{K5K7}, where opposite sides of the outer square or hexagon are identified to form a torus carrying the map $\{4,4\}_{1,2}$ or $\{3,6\}_{1,2}$; in each case the mirror image is the map $\M_n(c^{-1})$, where $c^{-1}=3$ or $4$ for $n=5, 7$ respectively.

In~\cite{Jam83} James classified the non-orientable regular embeddings of complete graphs (see also~\cite{Wil89} for an alternative proof due to Wilson):

\begin{thm}\label{nonorKn}
The only non-orientable regular embeddings of a complete graph $K_n$ are the maps $\{6, 2\}_3$, $\{4, 3\}_3$, $\{3, 5\}_5$ and $\{5, 5\}_3$ of characteristic $-1, -1, -1$ and $-3$ for $n=3, 4, 6$ and $6$.  \hfill$\square$
\end{thm}

These maps are all in class~$1$. The first three, on the projective plane, are the antipodal quotients of a hexagon, a cube and an icosahedron on the sphere (see Figure~\ref{projplKn}, where antipodal boundary points of the discs are identified), while the fourth is the Petrie dual of the third. The first two are the Petrie duals of the regular embeddings of $K_3$ and $K_4$ on the sphere, with automorphism groups isomorphic to $D_6\cong S_3\times C_2$ and $S_4$,, while the last two have the icosahedral group $A_5$ as their common automorphism group.

\begin{figure}[h!]
\begin{center}
\begin{tikzpicture}[scale=0.5, inner sep=0.8mm]

\draw [thick, dotted] (-10,0) arc (0:360:3);
\node (A) at (-15,0) [shape=circle, fill=black] {};
\node (B) at (-13,0) [shape=circle, fill=black] {};
\node (C) at (-11,0) [shape=circle, fill=black] {};
\draw [thick] (-16,0) to (-10,0);
\draw (-12.7,3.3) to (-13,3) to (-12.7,2.7);
\draw (-13.3,-3.3) to (-13,-3) to (-13.3,-2.7);

%%%%%%%

\draw [thick, dotted] (0,0) arc (0:360:3);
\node (a) at (-1.6,1.4) [shape=circle, fill=black] {};
\node (b) at (-4.4,1.4) [shape=circle, fill=black] {};
\node (c) at (-4.4,-1.4) [shape=circle, fill=black] {};
\node (d) at (-1.6,-1.4) [shape=circle, fill=black] {};
\draw [thick] (a) to (b) to (c) to (d) to (a);
\draw [thick] (a) to (-0.9,2.1);
\draw [thick] (b) to (-5.1,2.1);
\draw [thick] (c) to (-5.1,-2.1);
\draw [thick] (d) to (-0.9,-2.1);
\draw (-2.7,3.3) to (-3,3) to (-2.7,2.7);
\draw (-3.3,-3.3) to (-3,-3) to (-3.3,-2.7);

%%%%%%%

\draw [thick, dotted] (10,0) arc (0:360:3);
\node (a1) at (7,2) [shape=circle, fill=black] {};
\node (b1) at (8.9,0.618) [shape=circle, fill=black] {};
\node (c1) at (8.17,-1.616) [shape=circle, fill=black] {};
\node (d1) at (5.83,-1.616) [shape=circle, fill=black] {};
\node (e1) at (5.1,0.618) [shape=circle, fill=black] {};
\node (f1) at (7,0) [shape=circle, fill=black] {};
\draw [thick](a1) to (b1) to (c1) to (d1) to (e1) to (a1);
\draw [thick](a1) to (f1) to (b1);
\draw [thick](c1) to (f1) to (d1);
\draw [thick](e1) to (f1);

\draw [thick] (7.927,2.853) to (a1) to (6.033,2.853);
\draw [thick] (4.533,1.763) to (e1) to (4,0);
\draw [thick] (4.573,-1.763) to (d1) to (6.073,-2.853);
\draw [thick] (7.927,-2.853) to (c1) to (9.427,-1.763);
\draw [thick] (10,0) to (b1) to (9.467,1.763);

\draw (7.3,3.3) to (7,3) to (7.3,2.7);
\draw (6.7,-3.3) to (7,-3) to (6.7,-2.7);

\end{tikzpicture}

\end{center}
\caption{Regular embeddings of $K_n$ on the projective plane, $n=3, 4, 6$}
\label{projplKn}
\end{figure}

For $n>4$ the Petrie duals of the Biggs maps $\M_n(c)$ are non-orientable edge-transitive embeddings of $K_n$ in class $2^*{\rm ex}$. They have type $\{2p, n-1\}$ and automorphism group $AGL_1(\F_n)$.

In~\cite{Jam90} James completed the classification of the orientable edge-transitive embeddings of $K_n$. If $3<n=p^e\equiv 3$ mod~$(4)$ and $c$ is a primitive element for $\F_n$, let $\M_n(c,j)$ be the Cayley map for $\F_n$ with generating set $\F_n^*$, where now the cyclic ordering is
\[1, c^j, c^2, c^{j+2}, c^4, c^{j+4}, \ldots, c^{n-3}, c^{j+n-3}\]
for some odd element $j\in\Z_{\color{red}n-1}\setminus\{1\}$. (Taking $j=1$ gives the orientably regular Biggs map $\M_n(c)$.) 

\begin{thm}\label{etransorKn}
If $\M$ is an orientable edge-transitive embedding of $K_n$ which is not orientably regular, then $\M\cong\M_n(c,j)$ for some $n$, $c$ and $j$ as above. As oriented maps, $\M_n(c, j)$ and $\M_n(c', j')$ are isomorphic if and only if $c$ and $c'$ are equivalent under ${\rm Gal}\,\F_n$ and $j'\equiv j$ or $2-j$ mod~$(n-1)$.  \hfill$\square$
\end{thm}

The number of such maps for a given $n$ is therefore $(n-3)\phi(n-1)/4e$. They come in chiral pairs $\M_n(c,j)$ and $\M_n(c^{-1},2-j)$. For $j\ne 1$ the automorphism group of $\M_n(c,j)$ is the unique subgroup $AHL_1(\F_n)$ of index $2$ in $AGL_1(\F_n)$, consisting of the affine transformations $v\mapsto av+b$ for which $a$ is a non-zero square in $\F_n$. This group has two orbits on faces. These are all $p$-gons unless $j^*:=j$ or $2-j\equiv (n-1)/2$ mod~$(n-1)$, in which case one orbit consists of $p$-gons and the other of $(n-1)/2(n-1,j^*)$-gons. There is a single orbit on Petrie polygons, all of them having length $2(n-1)/(n-1,2(j-1))$. These maps are all in class $5^*$. Their Petrie duals, which are non-orientable edge-transitive embeddings of $K_n$, are vertex- and Petrie-transitive, and lie in class $5^P$.

\medskip

\noindent{\bf Example.} The James map $\M_7(5,5)$ is shown in Figure~\ref{JamesK7}, taken from~\cite{Jam90}, with two identification of sides of the outer $14$-gon shown by the letters $A$ and $B$; the others can be found by $C_7$ rotational symmetry. The vertices are identified with the elements $0, 1, \ldots, 6$ of $\F_7$, in clockwise order, starting with $0$ at the top. The mirror image of this map is given by the same diagram, but with the side-identifications reflected. Each of these two maps $\M$ can be drawn on Klein's Riemann surface $\K$ of genus $3$ with automorphism group $L_2(7)$ (see~\cite{Lev}, for instance), so that ${\rm Aut}\,\M$ is contained in ${\rm Aut}\,\K$ as the normaliser of a Sylow $7$-subgroup; in each case we obtain a representation of the Fano plane, with the vertices of $\M$ as its points and the seven triangular faces corresponding to its lines.

\begin{figure}[h!]

\begin{center}
 \begin{tikzpicture}[scale=0.3, inner sep=1mm]
 
\draw [thick, dotted] (10,0) to (9.01,4.34) to (6.24, 7.82) to (2.23, 9.75) to (-2.23, 9.75) to (-6.24, 7.82) to (-9.01, 4.34) to (-10,0) to (-9.01,-4.34) to (-6.24,-7.82) to (-2.23, -9.75) to (2.23, -9.75) to (6.24, -7.82) to (9.01, -4.34) to (10,0);

\draw [thick] (1, -9.75) to (9.7, 1.1);
\draw [thick] (8.24, -5.30) to (5.19, 8.27);
\draw [thick] (9.29, 3.14) to (-3.23, 9.22);
\draw [thick] (-9.29, 3.14) to (3.23, 9.22);
\draw [thick] (-8.24, -5.30) to (-5.19, 8.27);
\draw [thick] (-1, -9.75) to (-9.7, 1.1);
\draw [thick] (7, -6.87) to (-6.9, -6.9);

\node (A) at (7.45,-1.8) [shape=circle, fill=black] {};
\node (B) at (6.04,4.7) [shape=circle, fill=black] {};
\node (C) at (0.106,7.67) [shape=circle, fill=black] {};
\node (D) at (-6.04,4.7) [shape=circle, fill=black] {};
\node (E) at (-7.45,-1.8) [shape=circle, fill=black] {};
\node (F) at (-3.24,-6.95) [shape=circle, fill=black] {};
\node (G) at (3.24,-6.95) [shape=circle, fill=black] {};

\draw [thick] (1,9.75) to (C) to (-1,9.75);
\draw [thick] (-7,6.87) to (D) to (-8.25,5.3);
\draw [thick] (-9.74,-1.19) to (E) to (-9.29,-3.14);
\draw [thick] (-5.15,-8.36) to (F) to (-3.34,-9.22);
\draw [thick] (5.15,-8.36) to (G) to (3.34,-9.22);
\draw [thick] (9.74,-1.19) to (A) to (9.29,-3.14);
\draw [thick] (7,6.87) to (B) to (8.25,5.3);

\node at (0,11) {$A$};
\node at (-8.9,-6.9) {$A$};
\node at (4.9,9.7) {$B$};
\node at (4.9,-9.7) {$B$};

\end{tikzpicture}

\end{center}
\caption{James's edge-transitive embedding of $K_7$} \label{JamesK7}
\end{figure}

With the above information, we can now complete the classification of the edge-transitive embeddings of complete graphs, by dealing with the non-regular non-orientable cases:

\begin{thm}
If $\M$ is an edge-transitive embedding of a complete graph $K_n$, then $\M$ is isomorphic to one of the maps listed in Theorem~\ref{orKn}, \ref{nonorKn} or \ref{etransorKn}, or its Petrie dual. In particular, $n$ is either $6$ or a prime power.
\end{thm}

\noindent{\sl Proof.} Suppose that $\M$ is such an embedding, and is in class $T$. The cases $n\le 3$ are easily dealt with, so we may assume that $n\ge 4$. By Theorems~\ref{orKn}, \ref{nonorKn} and \ref{etransorKn} we may assume that $T\ne 1$ and that $\M$ is non-orientable, so that $T\ne 2^P{\rm ex}$, $5$ or $5^*$. If $T=2^*{\rm ex}$ or $5^P$ then $\M$ is the Petrie dual of an edge-transitive embedding of $K_n$ in class $2^P{\rm ex}$ or $5^*$; these are classified in Theorem~\ref{orKn} or \ref{etransorKn}, so we may also exclude $2^*{\rm ex}$ and $5^P$. This leaves only the possibilities $T=2^{\sigma}$, $2\,{\rm ex}$, $3$ and $4^{\sigma}$; we will show that these cases do not arise.

The edges of $K_n$ may be identified with the $2$-element sets of vertices, so edge-transitivity implies that the group $A={\rm Aut}\,\M$ permutes these sets transitively, that is, $A$ is $2$-homogenous on the vertices. It permutes them faithfully since $n\ge 4$, and is transitive on them since $K_n$ is not bipartite; this rules out classes $2$, $3$ and $4$, so that only $T=2^*$, $2^P$, $2\,{\rm ex}$, $4^*$ and $4^P$ remain. It is easily seen that a $2$-homogenous permutation group of finite degree $n\ge 4$ is either $2$-transitive or is a non-cyclic group of odd order. However, the latter case is impossible since the presentations of $N(T)$ in Proposition~\ref{parents} for $T=2^{\sigma}$, $2^{\sigma}{\rm ex}$ and $4^{\sigma}$ show that their only quotients of odd order are cyclic. 

It follows that $A$ must be $2$-transitive on the vertices of $\M$. Thus each edge of $\M$ is reversed by some automorphism (a reflection or half turn), so the single edge of ${\mathcal N}(T)$ must be free, ruling out $2\,{\rm ex}$, $4^*$ and $4^P$. Only the classes  $2^*$ and $2^P$ remain, with $A$ sharply $2$-transitive on the vertices (since $T\ne 1$), and the stabiliser of an edge generated by a reflection or a half-turn respectively, transposing its two incident vertices. Zassenhaus~\cite{Zas} showed that any sharply $2$-transitive group can be identified with the Frobenius group $AGL_1(F)$, acting naturally on $F$, for some near-field $F$. Since $\M$ is a map, the stabiliser $A_0$ of the vertex $0$ must be a cyclic or dihedral group of order $n-1$, acting regularly on the set $F^*$ of neighbours of $0$. However, $A_0$ is a Frobenius complement, and these cannot be dihedral groups since they contain at most one involution (see~\cite[Satz V.18.1(a)]{Hup} or~\cite[Theorem~18.1(iii)]{Pas}), so $A_0$ must be cyclic. Since $T=2^*$ or $2^P$, $N(T)$ is generated by involutions, and hence so are its epimorphic images $A$ and $A_0$, giving $n-1\le 2$, a contradiction. \hfill$\square$

%\medskip

%{\color{blue}[Extend this proof to cover {\sl all} edge-transitive embeddings of complete graphs, including the earlier classifications?]}

%%%%%%%%%%%%%%%%

\section{Monodromy groups}

If a map $\M$ corresponds to a map subgroup $M$ of $\Gamma$, then its monodromy group $G$ can be identified with the permutation group induced by $\Gamma$ on the cosets of $M$; it is also isomorphic to the automorphism group of the minimal regular cover of $\M$. In this section we will consider how $G$ is related to the automorphism group $A={\rm Aut}\,\M$ of $\M$, depending on which edge-transitive class $\M$ is in.

We have $G\cong \Gamma/M_0$, where
\[M_0:=\bigcap_{g\in\Gamma}M^g\]
 is the core of $M$ in $\Gamma$, the kernel of the action of $\Gamma$ on the flags of $\M$. If $\M$ is edge-transitive, in class $T$, then the number of distinct conjugates $M^g$ of $M$ in $\Gamma$ is
 \[|\Gamma:N_{\Gamma}(M)|=|\Gamma:N(T)|=n(T)=1, 2\;\;{\rm or}\;\;4,\]
 depending only on $T$. 
 
 If $T=1$ then $M$ is a normal subgroup of $\Gamma$, in which case $n(T)=1$, $M_0=M$ and $G\cong A$. We will therefore concentrate on the remaining thirteen edge-transitive classes. Since $G$ is invariant under the action of $\Omega$ on maps, there are just five separate cases to consider, corresponding to the orbits of $\Omega$ on these classes.
 
 If $T\ne 1$ we have strict inclusions $\Gamma>N(T)>M$. Now $M$ is the stabiliser in $\Gamma$ of a flag of $\M$, so these inclusions show that the action of $\Gamma$ on flags  is imprimitive: the $|N(T):M|=|A|$ flags fixed by $M$ form a block of imprimitivity, and the $n(T)=|\Gamma:N(T)|$ images of this block form a system of imprimitivity for $\Gamma$. It follows from a standard result about imprimitive groups that $G$ is isomorphic to a subgroup of the wreath product
\[A\wr G_T=A^n\rtimes G_T\le A^n\rtimes S_n,\]
 where $n=n(T)$ and $G_T$ is the monodromy group of ${\mathcal N}(T)$, the permutation group of degree $n$ induced by $G$ on the set of blocks, or equivalently by $\Gamma$ on the cosets of $N(T)$.
 
 If $T=2^{\sigma}$ or $2^{\sigma}{\rm ex}$ for some $\sigma$ then $n(T)=2$, so $N(T)$ is normal in $\Gamma$ and $G_T\cong \Gamma/N(T)\cong S_2$. Then $M_0$ is the intersection of two conjugates of $M$, namely $M$ itself and $M^g$ for some $g\in\Gamma\setminus N(T)$. We have $M, M^g\le N(T)$, with $M_0=M\cap M^g$, so $|G|=2|N(T):M_0|\le 2|A|^2$. There are cases where this upper bound is attained, with $N(T)/M_0\cong A^2$, for instance if $A$ is simple; see Theorem~\ref{mainthm} for those non-abelian simple groups which can arise here. For an example where the bound is not attained, let $T=2^P{\rm ex}$, so $N(T)=\Gamma^+$, and let $\M$ be an orientably regular embedding of the complete graph $K_q$ for some prime power $q\ge 5$ (see Section~\ref{completemaps}); then $|A|=|AGL_1(\F_q)|=q(q-1)$, but $N(T)/M_0\cong (\F_q\times\F_q)\rtimes \F_q^*$, so $|N(T):M_0|=q^2(q-1)$ and hence $|G|=2q^2(q-1)=2|A|^2/(q-1)$.
 
Similar arguments appliy if $T=3$ or $5^{\sigma}$, since $N(T)$ is again normal in $\Gamma$, though in these cases $n(T)=4$, $G_T\cong V_4$ and $M_0$ is the intersection of four conjugates of $M$. We therefore obtain the upper bound $4|A|^4$ for $|G|$, attained whenever $A$ is a non-abelian simple group (again, see Theorem~\ref{mainthm}). 
 
When $T=4^{\sigma}$, however, we have $N_{\Gamma}(N(T))=N(2^{\sigma})$, of index $2$ in $\Gamma$, and the core $N(T)_0$ of $N(T)$ has index $8$ in $\Gamma$, with $G_T\cong \Gamma/N(T)_0\cong D_4$ (see Section~\ref{14classes}). We therefore have $|G|\le 8|A|^4$ in these cases, again attained whenever $A$ is a non-abelian simple group.

\newpage

%%%%%%%%%%%%%%%%%%%%%
%%%%%%%%%%%%%%%%%%%%%

\part{Automorphism groups}\label{autgps}

This part of the text is devoted to the question of which groups can arise as the automorphism group of a map in the various edge-transitive classes. Although the main emphasis is on finite groups, acting on compact maps, infinite groups are considered in the last section.

%%%%%%%%%%%%%

\section{Realising symmetric groups}

We now return to our main theme of realising various groups as automorphism groups of maps in specific edge-transitive classes. It is well known that, with a few small exceptions, each finite symmetric or alternating group can be realised as the automorphism group of both a regular map and a chiral map: for instance, these facts can be proved using the methods employed by Conder in~\cite{Con80}, where the main emphasis was on realising alternating groups as Hurwitz groups, or more recently by Conder, Huc\'\i kov\'a, Nedela and \v Sir\'a\v n in~\cite{CHNS}, where the emphasis was on constructing chiral maps of a given type, rather than with a given automorphism group. 

In~\cite{STW}, \v Sir\'a\v n, Tucker and Watkins showed that for each integer $n\ge 11$ with $n\equiv 3$ or $11$ mod~$(12)$, there are finite, orientable, edge-transitive maps $\mathcal M$ in each of the $14$ edge-transitive classes $T$, with ${\rm Aut}\,{\mathcal M}$ isomorphic to the symmetric group $S_n$. The aim of this section is to extend this result by determining all the pairs $n$ and $T$ such that $S_n\in\G(T)$. (Realisations by orientable maps without boundary will be considered in the next section.) In order to do this, we will need to prove that various sets of permutations generate $S_n$. The following result is elementary and well-known:

\begin{lemma}\label{c,t}
If $t$ is a transposition $(i,j)$ where $i$ and $j$ are adjacent terms in an $n$-cycle $c\in S_n$, then $c$ and $t$ generate $S_n$.  \hfill$\square$
\end{lemma}

In similar contexts, the following theorem (a simple consequence of~\cite[Th\'eor\`eme~I]{Jor71} and~\cite[Th\'eor\`eme~I]{Jor73}, see also~\cite[Theorem~13.9]{Wie}) has often been used:

\begin{thm}[Jordan]\label{Jordan}
If $G$ is a primitive subgroup of $S_n$, containing a cycle of prime length $l\le n-3$, then $G\ge A_n$.
 \hfill$\square$
 \end{thm}

(See Theorem~\ref{J-Jordan} for a generalisation of Jordan's Theorem, omitting the primality condition.) We will also need the following elementary result:

\begin{lemma}\label{smallSn}
If elements $x$ and $y$ generate $S_n$, where $n\le 5$, then they are simultaneously conjugate in $S_n$ to their inverses.
\end{lemma}

\noindent{\sl Proof.} Suppose that $x$ and $y$ generate $S_5$, and are not simultaneously inverted by conjugation. If $x$ is a transposition then in order for $\langle x, y\rangle$ to be transitive, $y$ can have at most two orbits, so it has cycle structure $14$ or $23$ or $5$; in each case a permutation diagram for $x$ and $y$ is symmetric, admitting a reflection which inverts $x$ and $y$, so $x$ and $y$ are inverted by conjugation in $S_5$, against our assumption. (The four possibilities are shown in Figure~\ref{xyinversion}, where broken and unbroken lines represent the actions of $x$ and $y$, with cycles of the latter cyclically ordered anticlockwise.) Thus $x$ cannot be a transposition, and hence, by symmetry, neither can $y$.

\begin{figure}[h!]

\begin{center}
\begin{tikzpicture}[scale=0.5, inner sep=0.8mm]
 
\node (a) at (-12,0)  [shape=circle, draw, fill=black] {};
\node (b) at (-14,2) [shape=circle, draw, fill=black] {};
\node (c) at (-12,4)  [shape=circle, draw, fill=black] {};
\node (d) at (-10,2) [shape=circle, draw, fill=black] {};
\node (e) at (-12,6) [shape=circle, draw, fill=black] {};

\draw (a) to (b) to (c) to (d) to (a);
\draw [thick, dashed] (c) to (e);

%%%%%%%%%

\node (f) at (-8,0)  [shape=circle, draw, fill=black] {};
\node (g) at (-6,2)  [shape=circle, draw, fill=black] {};
\node (h) at (-4,0) [shape=circle, draw, fill=black] {};
\node (i) at (-6,4)  [shape=circle, draw, fill=black] {};
\node (j) at (-6,6) [shape=circle, draw, fill=black] {};

\draw (f) to (g) to (h) to (f);
\draw [thick, dashed] (g) to (i);
\draw (i) to (j);

%%%%%%%%

\node (k) at (-1.2,0)  [shape=circle, draw, fill=black] {};
\node (l) at (-2,2.3)  [shape=circle, draw, fill=black] {};
\node (m) at (0,4)  [shape=circle, draw, fill=black] {};
\node (n) at (2,2.3)  [shape=circle, draw, fill=black] {};
\node (o) at (1.2,0)  [shape=circle, draw, fill=black] {};

\draw (k) to (l) to (m) to (n) to (o) to (k);
\draw [thick, dashed] (l) to (n);

%%%%%%%%

\node (p) at (4.8,0)  [shape=circle, draw, fill=black] {};
\node (q) at (4,2.3)  [shape=circle, draw, fill=black] {};
\node (r) at (6,4)  [shape=circle, draw, fill=black] {};
\node (s) at (8,2.3)  [shape=circle, draw, fill=black] {};
\node (t) at (7.2,0)  [shape=circle, draw, fill=black] {};

\draw (p) to (q) to (r) to (s) to (t);
\draw (4.8,-0.1) to (7.2,-0.1);
\draw [thick, dashed] (4.8,0.1) to (7.2,0.1);
 
 \end{tikzpicture}

\end{center}
\caption{Permutation diagrams for elements of $S_5$} 
\label{xyinversion}
\end{figure}

If $x$ is a double transposition then $y$ must be odd, of cycle structure $14$ or $23$. Similar arguments using permutation diagrams show that the latter always allows inversion of $x$ and $y$, while the former allows either inversion of $x$ and $y$ or $G:=\langle x, y\rangle\cong AGL_1(5)$ with $\langle xy^2\rangle\triangleleft G$. Thus $x$ cannot be a double transposition, and hence neither can $y$.

If $x$ is a $3$-cycle or a $5$-cycle then $y$ is odd, with cycle-structure $23$ or $14$; all the resulting connected diagrams are symmetric, so $x$ and $y$ are inverted. Thus $x$ cannot be a $3$-cycle or a $5$-cycle, so the same applies to $y$.

This shows that neither $x$ nor $y$ can be even, so both are odd. But then $z:=(xy)^{-1}$ is even, and $\langle x,z\rangle=S_5$, so the preceding arguments show that $x$ and $z$ are inverted, and hence so are $x$ and $y$.

The arguments for $n\le 4$ are similar, and even more elementary, so they are omitted. \hfill$\square$

\medskip

The following result proves the statements in Theorem~\ref{mainthm} concerning symmetric groups:

\begin{thm}\label{symgps}
The symmetric group $S_n$ is the automorphism group of a map in an edge-transitive class $T$ if and only if one of the following holds:
\begin{itemize}
\item $T=1$ and $n\ge 1$;
\item $T=2^{\sigma}$, $3$ or $4^{\sigma}$ for some $\sigma$ and $n\ge 2$;
\item $T=2^{\sigma}\,{\rm ex}$ or $5^{\sigma}$ for some $\sigma$ and $n\ge 6$.
\end{itemize}
\end{thm}

\noindent{\sl Proof.}  Lemma~\ref{abelian} deals with the cases $n\le 2$, so we may assume that $n\ge 3$. We can map the standard generators $R_i$ of $N(1)=\Gamma$ to permutations $r_i\in S_n$, where
\[r_1=(1)(2,n)(3, n-1)\ldots,\quad r_2=(1,2)(3,n)(4,n-1)\ldots\]
are involutions in $S_n$ with
\[r_1r_2=(1, 2, \ldots, n),\]
and where
\[r_0=(1,2)\]
is an involution commuting with $r_2$. This defines a homomorphism $\Gamma\to S_n$, and since $\langle r_1r_2, r_0\rangle=S_n$ by Lemma~\ref{c,t} it is an epimorphism. The kernel $M$ is then the map subgroup of a regular map $\mathcal M$ with ${\rm Aut}\,{\mathcal M}\cong S_n$, so $S_n\in\G(1)$. Since $S_n$ is non-abelian, it follows from Lemma~\ref{regmapslemma}(a) that there are also edge-transitive maps in classes $T=2^{\sigma}$, $3$ and $4^{\sigma}$ with automorphism group $S_n$ for each $\sigma$ and each $n\ge 3$.

Now suppose that $T=2^P{\rm ex}$, so that $N(T)=\Gamma^+=\langle X, Y\mid Y^2=1\rangle$. If $n\ge 6$ let us map the generators $X$ and $Y$ to the permutations
\[x=(1, 2, \ldots, n-1)\quad{\rm and}\quad y=(1,3)(2,4)(n-1,n)\]
in $S_n$. The group $G=\langle x, y\rangle$ is clearly transitive, and in fact doubly transitive since it contains $x$, so it is primitive. If $n\ge 7$ then
\[ [x,y]=(1,n-1,n,3,5),\]
a cycle of length $5$, so Jordan's Theorem implies that $G\ge A_n$ provided $n\ge 8$, and since $y$ is an odd permutation it follows that $G=S_n$. In fact this also holds for $n=7$ since the only primitive groups of this degree containing a $5$-cycle are $A_7$ and $S_7$. In each case it is easy to check, for instance by drawing a permutation diagram, that no automorphism of $S_n$ (necessarily inner) inverts both $x$ and $y$. This argument fails for $n=6$, but in this case we can map $X$ and $Y$ to
\[x=(1,\ldots, 6)\quad{\rm and}\quad y=(1,2)(3,5).\]
Then $xy=(2,5,6)(3,4)$ so $(xy)^3=(3,4)$ and hence $\langle x, y\rangle=S_6$; again a permutation diagram shows that no inner automorphism inverts $x$ and $y$, and since the outer automorphisms of $S_6$ send $6$-cycles to products of disjoint $2$- and $3$-cycles, they are also eliminated. By Lemma~\ref{regmapslemma}(b) the result follows for the remaining classes $2^{\sigma}{\rm ex}$ and $5^{\sigma}$ for all $n\ge 6$.

If $n\le 5$ then Lemma~\ref{smallSn} shows that there are no edge-transitive maps in classes $2^{\sigma}\,{\rm ex}$ or $5^{\sigma}$ with automorphism group $S_n$: in each case the two generators of a quotient group are inverted by some automorphism of $S_n$. \hfill$\square$

\medskip

Cayley's Theorem gives the following consequence of Theorem~\ref{symgps}:

\begin{cor}
If $F$ is any finite group then for each edge-transitive class $T$ there is a compact map $\mathcal M$ in $T$ such that $F$ is isomorphic to a subgroup of ${\rm Aut}\,{\mathcal M}$. \hfill$\square$
\end{cor}

In fact, by using an embedding theorem of Schupp~\cite{Sch}, and allowing non-compact maps, this result can be extended from finite groups to all countable groups, as shown in Section~\ref{embed}.

Note that in the cases where $T=2^{\sigma}$, $3$ or $4^{\sigma}$ and $n=2$, the maps used in the proof of Theorem~\ref{symgps} have non-empty boundaries. This cannot be avoided, since in order to exclude forbidden automorphisms some reflection $S_i\in N(T)$ must be mapped to the identity. In the next section we will consider realisations of symmetric groups as the automorphism groups of edge-transitive maps which are orientable and without boundary.

%%%%%%%%%%%%%%%%%%%

\section{Evenly realising symmetric groups}\label{sym}\label{evensym}

We saw in Theorem~\ref{symgps} that for almost all pairs $T$ and $n$ there is an edge-transitive map in class $T$ with automorphism group $S_n$. In many cases the maps constructed there are non-orientable, and in a few cases they have non-empty boundaries. We will now determine, for each edge-transitive class $T$, the values of $n$ for which $S_n$ is the automorphism group of an orientable map without boundary in the class $T$, that is, $S_n\in\G^+(T)$.

For brevity, in this case  we will often simply say that $S_n$ is {\em evenly realised\/} in class $T$, since a map $\M$ is orientable and without boundary if and only if the corresponding map subgroups $M$ are contained in the even subgroup $\Gamma^+$ of index $2$ in $\Gamma$. If $T=2^P{\rm ex}$, $5$ or $5^*$ then $N(T)\le\Gamma^+$ (since ${\mathcal N}(T)$ is a map on the sphere), so in such cases  {\sl all\/} maps in $T$ are orientable and without boundary, and Theorem~\ref{symgps} gives the required result. For the remaining classes $T$, we need to consider whether one can find a quotient $N(T)/M$ of $N(T)$ isomorphic to $S_n$, with no forbidden automorphisms, and with the additional requirement that $M\le\Gamma^+$. In terms of the induced epimorphism $N(T)\to S_n$, this last condition is equivalent to requiring that all orientation-reversing generators of $N(T)$ should be sent to odd permutations.

As before, we will need to show that various sets of permutations generate $S_n$. In addition to Lemma~\ref{c,t} and Theorem~\ref{Jordan} (Jordan's Theorem) we will need the following result, which uses the classification of finite simple groups to remove the primality condition in Jordan's Theorem (see~\cite{Jon14}):

\begin{thm}\label{J-Jordan}
If $G$ is a primitive subgroup of $S_n$, containing a cycle of length $l\le n-3$, then $G\ge A_n$. \hfill$\square$
\end{thm}

We will also use the following simple corollary (see~\cite [Lemma~6.4(i)]{Jon15}):

\begin{cor}\label{2cycles}
If $G$ is a transitive subgroup of $S_n$, containing a permutation consisting of two cycles with mutually coprime lengths $l, m>1$, and no fixed points, then $G\ge A_n$.  \hfill$\square$
\end{cor}

The main result  of this section is the following:

\begin{thm}\label{Sneven}
The symmetric group $S_n$ is the automorphism group of an orientable map without boundary in an edge-transitive class $T$ if and only if one of the following holds:
\begin{itemize}
\item[\rm(a)] $T=1$ and $n\ne 1, 5$ or $6$;
\item[\rm(b)]  $T=2$ or $2^*$ and $n\ne 1, 2, 5$ or $6$, or $T=2^P$ and $n\ge 3$;
\item[\rm(c)] $T=2\,{\rm ex}$ or $2^*{\rm ex}$ and $n\ge 7$, or $T=2^P{\rm ex}$ and $n\ge 6$;
\item[\rm(d)] $T=3$ and $n\ge 3$;
\item[\rm(e)] $T=4^{\sigma}$ for some $\sigma$ and $n\ge 3$;
\item[\rm(f)] $T=5^{\sigma}$ for some $\sigma$ and $n\ge 6$.
\end{itemize}
\end{thm}

The six parts of this theorem correspond to the six orbits of $\Omega$ on the classes $T$, and thus on the isomorphism classes of parent groups $N(T)$. However, in contrast with Theorem~\ref{symgps}, the results in (b) and (c) vary for different classes in the same orbit, since we now have to take into account which generators of $N(T)$ preserve or reverse orientation.
  
%%%%%%%%

\subsection{Proof of (a), with $T=1$}

First let $T=1$, the class of regular maps, so that $N(T)=\Gamma$ and hence $N(T)\cap\Gamma^+=\Gamma^+$. Given any epimorphism $\theta:\Gamma\to S_n, R_i\mapsto r_i$, the kernel $M$ is contained in a unique subgroup $\theta^{-1}(A_n)$ of index $2$ in $\Gamma$, and this is $\Gamma^+$ if and only if it does not contain any $R_i$, that is, $r_i$ is an odd permutation for each $i=0, 1, 2$. We will therefore look for triples of odd involutions $r_i$ generating $S_n$, with $r_0$ and $r_2$ commuting. 

\medskip

\noindent{\bf Case 1} Let $n\equiv 3$ mod~$(4)$, say $n=4k+3\ge 3$. Define
\[r_0=(1)(2,n)(3,n-1)\ldots(2k+2, 2k+3),\]
\[r_1=(1,2)(3,n)(4,n-1)\ldots (2k+2, 2k+4)(2k+3),\]
\[r_2=(2k+2, 2k+3),\]
so that each $r_i$ is an odd involution, and $r_2$ commutes with $r_0$. Then
\[r_0r_1=(1,2,\ldots, n),\quad{\rm so}\quad\langle r_0r_1, r_2\rangle=S_n\]
by Lemma~\ref{c,t}, and we have an epimorphism $\Gamma\to S_n$ with the required properties.

\medskip

\noindent{\bf Case 2} Let $n\equiv 0$ mod~$(4)$, say $n=4k\ge 4$. Define
\[r_0=(1)(2,n-1)(3,n-2)\ldots(2k, 2k+1)(n),\]
\[r_1=(1,2)(3,n-1)(4,n-2)\ldots (2k, 2k+2)(2k+1)(n),\]
\[r_2=(1,n).\]
Again, each $r_i$ is an odd involution, and $r_2$ commutes with $r_0$. Now
\[r_0r_1=(1, 2, \ldots, n-1)(n),\]
so
\[r_0r_1r_2=(1, 2, \ldots, n)\quad{\rm and}\quad \langle r_0r_1r_2, r_2\rangle =S_n,\]
giving the required epimorphism $\Gamma\to S_n$.

\medskip

These two cases rely on $r_0$ and $r_1$ being natural generators of a `large' dihedral subgroup of $S_n$. The remaining cases are less straightforward, since the corresponding dihedral subgroup is not large enough for our purposes.

\medskip

\noindent{\bf Case 3} Let $n\equiv 1$ mod~$(4)$, say $n=4k+1$, and suppose that $n\ge 9$. Define
\[r_0=(1,2)(3,4)\ldots(n-4,n-3)\quad{\rm and}\]
\[r_2=(3,5)(4,6)(7,9)\ldots(n-5,n-3)(n-2,n-1),\]
commuting odd involutions, and define $r_1$ to be the odd involution
\[(1,3)(5,7)(6,8)(9,11)(10,12)\ldots (n-8,n-6)(n-7,n-5)(n-4,n-2)(n-1,n)\]
or
\[(1,3)(5,8)(6,7)(9,11)(10,12)\ldots (n-8,n-6)(n-7,n-5)(n-4,n-2)(n-1,n)\]
as $k$ is even or odd. 

\begin{figure}[h!]

\begin{center}
\begin{tikzpicture}[scale=0.6, inner sep=0.8mm]
 
\node (1) at (-10,2)  [shape=circle, draw, fill=black] {};
\node (2) at (-10,0) [shape=circle, draw, fill=black] {};
\node (3) at (-8,2)  [shape=circle, draw, fill=black] {};
\node (4) at (-8,0) [shape=circle, draw, fill=black] {};
\node (5) at (-6,2)  [shape=circle, draw, fill=black] {};
\node (6) at (-6,0) [shape=circle, draw, fill=black] {};
\node (7) at (-4,2)  [shape=circle, draw, fill=black] {};
\node (8) at (-4,0) [shape=circle, draw, fill=black] {};
\node (9) at (-2,2)  [shape=circle, draw, fill=black] {};
\node (10) at (-2,0) [shape=circle, draw, fill=black] {};
\node (11) at (-0,2)  [shape=circle, draw, fill=black] {};
\node (12) at (0,0) [shape=circle, draw, fill=black] {};
\node (13) at (2,2)  [shape=circle, draw, fill=black] {};
\node (14) at (2,0) [shape=circle, draw, fill=black] {};
\node (n-6) at (6,2)  [shape=circle, draw, fill=black] {};
\node (n-5) at (6,0) [shape=circle, draw, fill=black] {};
\node (n-4) at (8,2)  [shape=circle, draw, fill=black] {};
\node (n-3) at (8,0) [shape=circle, draw, fill=black] {};
\node (n-2) at (10,2) [shape=circle, draw, fill=black] {};
\node (n-1) at (10,0)  [shape=circle, draw, fill=black] {};
\node (n) at (12,0) [shape=circle, draw, fill=black] {};

\draw (1) to (2);
\draw (3) to (4);
\draw (5) to (6);
\draw (7) to (8);
\draw (9) to (10);
\draw (11) to (12);
\draw (13) to (14);
\draw (n-6) to (n-5);
\draw (n-4) to (n-3);

\draw [thick, dashed] (1) to (3);
\draw [thick, dashed] (5) to (7);
\draw [thick, dashed] (6) to (8);
\draw [thick, dashed] (9) to (11);
\draw [thick, dashed] (10) to (12);
\draw [thick, dashed] (n-4) to (n-2);
\draw [thick, dashed] (n-1) to (n);
\draw [thick, dashed] (13) to (3.5,2);
\draw [thick, dashed] (14) to (3.5,0);
\draw [thick, dashed] (4.5,2) to (n-6);
\draw [thick, dashed] (4.5,0) to (n-5);

\draw [thick, dotted] (3) to (5);
\draw [thick, dotted] (4) to (6);
\draw [thick, dotted] (7) to (9);
\draw [thick, dotted] (8) to (10);
\draw [thick, dotted] (11) to (13);
\draw [thick, dotted] (12) to (14);
\draw [thick, dotted] (n-6) to (n-4);
\draw [thick, dotted] (n-5) to (n-3);
\draw [thick, dotted] (n-2) to (n-1);

\node at (-10,3) {$1$};
\node at (-10,-1) {$2$};
\node at (-8,3) {$3$};
\node at (-8,-1) {$4$};
\node at (-6,3) {$5$};
\node at (-6,-1) {$6$};
\node at (-4,3) {$7$};
\node at (-4,-1) {$8$};
\node at (-2,3) {$9$};
\node at (-2,-1) {$10$};
\node at (0,3) {$11$};
\node at (0,-1) {$12$};
\node at (2,3) {$13$};
\node at (2,-1) {$14$};
\node at (5.8,3) {$n-6$};
\node at (5.8,-1) {$n-5$};
\node at (8,3) {$n-4$};
\node at (8.1,-1) {$n-3$};
\node at (10.1,3) {$n-2$};
\node at (10.3,-1) {$n-1$};
\node at (12,-1) {$n$};

\path [fill=lightgray] (-7.95,0.05) rectangle (-6.05,1.95);
\path [fill=lightgray] (-3.95,0.05) rectangle (-2.05,1.95);
\path [fill=lightgray] (0.05,0.05) rectangle (1.95,1.95);
\path [fill=lightgray] (6.05,0.05) rectangle (7.95,1.95);
 
 \end{tikzpicture}

\end{center}
\caption{The permutations $r_i$ for $n=4k+1$, $k$ even} 
\label{case3keven}
\end{figure}

\begin{figure}[h!]

\medskip

\begin{center}
\begin{tikzpicture}[scale=0.6, inner sep=0.8mm]
 
\node (1) at (-10,2)  [shape=circle, draw, fill=black] {};
\node (2) at (-10,0) [shape=circle, draw, fill=black] {};
\node (3) at (-8,2)  [shape=circle, draw, fill=black] {};
\node (4) at (-8,0) [shape=circle, draw, fill=black] {};
\node (5) at (-6,2)  [shape=circle, draw, fill=black] {};
\node (6) at (-6,0) [shape=circle, draw, fill=black] {};
\node (7) at (-4,2)  [shape=circle, draw, fill=black] {};
\node (8) at (-4,0) [shape=circle, draw, fill=black] {};
\node (9) at (-2,2)  [shape=circle, draw, fill=black] {};
\node (10) at (-2,0) [shape=circle, draw, fill=black] {};
\node (11) at (-0,2)  [shape=circle, draw, fill=black] {};
\node (12) at (0,0) [shape=circle, draw, fill=black] {};
\node (13) at (2,2)  [shape=circle, draw, fill=black] {};
\node (14) at (2,0) [shape=circle, draw, fill=black] {};
\node (n-6) at (6,2)  [shape=circle, draw, fill=black] {};
\node (n-5) at (6,0) [shape=circle, draw, fill=black] {};
\node (n-4) at (8,2)  [shape=circle, draw, fill=black] {};
\node (n-3) at (8,0) [shape=circle, draw, fill=black] {};
\node (n-2) at (10,2) [shape=circle, draw, fill=black] {};
\node (n-1) at (10,0)  [shape=circle, draw, fill=black] {};
\node (n) at (12,0) [shape=circle, draw, fill=black] {};

\draw (1) to (2);
\draw (3) to (4);
\draw (5) to (6);
\draw (7) to (8);
\draw (9) to (10);
\draw (11) to (12);
\draw (13) to (14);
\draw (n-6) to (n-5);
\draw (n-4) to (n-3);

\draw [thick, dashed] (1) to (3);
\draw [thick, dashed] (5) to (8);
\draw [thick, dashed] (6) to (7);
\draw [thick, dashed] (9) to (11);
\draw [thick, dashed] (10) to (12);
\draw [thick, dashed] (n-4) to (n-2);
\draw [thick, dashed] (n-1) to (n);
\draw [thick, dashed] (13) to (3.5,2);
\draw [thick, dashed] (14) to (3.5,0);
\draw [thick, dashed] (4.5,2) to (n-6);
\draw [thick, dashed] (4.5,0) to (n-5);

\draw [thick, dotted] (3) to (5);
\draw [thick, dotted] (4) to (6);
\draw [thick, dotted] (7) to (9);
\draw [thick, dotted] (8) to (10);
\draw [thick, dotted] (11) to (13);
\draw [thick, dotted] (12) to (14);
\draw [thick, dotted] (n-6) to (n-4);
\draw [thick, dotted] (n-5) to (n-3);
\draw [thick, dotted] (n-2) to (n-1);

\node at (-10,3) {$1$};
\node at (-10,-1) {$2$};
\node at (-8,3) {$3$};
\node at (-8,-1) {$4$};
\node at (-6,3) {$5$};
\node at (-6,-1) {$6$};
\node at (-4,3) {$7$};
\node at (-4,-1) {$8$};
\node at (-2,3) {$9$};
\node at (-2,-1) {$10$};
\node at (0,3) {$11$};
\node at (0,-1) {$12$};
\node at (2,3) {$13$};
\node at (2,-1) {$14$};
\node at (5.8,3) {$n-6$};
\node at (5.8,-1) {$n-5$};
\node at (8,3) {$n-4$};
\node at (8.1,-1) {$n-3$};
\node at (10.1,3) {$n-2$};
\node at (10.3,-1) {$n-1$};
\node at (12,-1) {$n$};

\path [fill=lightgray] (-7.95,0.05) rectangle (-6.05,1.95);
\path [fill=lightgray] (-3.95,0.05) rectangle (-2.05,1.95);
\path [fill=lightgray] (0.05,0.05) rectangle (1.95,1.95);
\path [fill=lightgray] (6.05,0.05) rectangle (7.95,1.95);
 
 \end{tikzpicture}

\end{center}
\caption{The permutations $r_i$ for $n=4k+1$, $k$ odd} 
\label{case3kodd}
\end{figure}

These permutations are shown in Figures~\ref{case3keven} and \ref{case3kodd} respectively, with unbroken, dashed and dotted lines representing $2$-cycles of $r_0$, $r_1$ and $r_2$. For clarity, the squares representing regular orbits of the subgroup $\langle r_0, r_2\rangle\cong V_4$ are coloured grey.

The permutations $r_i$ generate a transitive group $G\le S_n$. Now $r_0r_1r_2$ consists of two cycles
\[(1,2,5,10,13,18,\ldots, n-7,n-4,n-5,n-10,n-13,\ldots, 7,4)\quad{\rm and}\]
\[(3,6,9,14,17,\ldots,n-3,n-1,n,n-2,n-6,11,8),\]
or
\[(1,2,5,9,14,17,\ldots,n-4,n-5,n-10,\ldots, 11,8,4)\quad{\rm and}\]
\[(3,6,10,13,18,\ldots,n-3,n-1,n,n-2,n-6,n-9,\ldots,12,7),\]
as $k$ is even or odd. (When $n=9$ or $n=13$ the cycles of $r_0r_1r_2$ reduce to
\[(1,2,5,4)(3,6,8,9,7)\quad{\rm or}\quad (1,2,5,9,8,4)(3,6,10,12,13,11,7)\]
respectively.) In each case, the two cycles of $r_0r_1r_2$ have mutually coprime lengths $2k$ and $2k+1$, so Corollary~\ref{2cycles} implies that $G\ge A_n$. Since $G$ contains odd permutations we have $G=S_n$ as required.

This construction fails if $n=5$. Indeed, the odd involutions in $S_5$ are the transpositions, and it is easy to see that  three of these cannot generate a transitive group, so $S_5$ cannot be the automorphism group of an orientable regular map. However, by Theorem~\ref{symgps} there are regular maps with automorphism group $S_5$: the only such maps are the non-orientable map of genus $5$ and type $\{4, 5\}_6$, and its dual, listed as N5.1 in Conder's list of non-orientable regular maps~\cite{Con}.

\medskip

\noindent{\bf Case 4} Let $n\equiv 2$ mod~$(4)$, say $n=4k+2$. First suppose that $n\ge 18$. Define
\[r_0=(1,2)(3,4)\ldots (n-1,n),\]
\[r_1=(2,4)(6,8)(9,11)(10,12)(13,15)(14,16)\ldots\]
\[\ldots(n-13,n-11)(n-12,n-10)(n-9,n-7)(n-5,n-3)(n-2,n),\]
\[r_2=(3,5)(4,6)(7,9)(8,10)\ldots (n-7,n-5)(n-6,n-4)(n-3,n-2).\]

\begin{figure}[h!]

\begin{center}
\begin{tikzpicture}[scale=0.7, inner sep=0.8mm]
 
\node (1) at (-6,2)  [shape=circle, draw, fill=black] {};
\node (2) at (-6,0) [shape=circle, draw, fill=black] {};
\node (3) at (-4,2)  [shape=circle, draw, fill=black] {};
\node (4) at (-4,0) [shape=circle, draw, fill=black] {};
\node (5) at (-2,2)  [shape=circle, draw, fill=black] {};
\node (6) at (-2,0) [shape=circle, draw, fill=black] {};
\node (7) at (0,2)  [shape=circle, draw, fill=black] {};
\node (8) at (0,0) [shape=circle, draw, fill=black] {};
\node (9) at (2,2)  [shape=circle, draw, fill=black] {};
\node (10) at (2,0) [shape=circle, draw, fill=black] {};
\node (11) at (4,2)  [shape=circle, draw, fill=black] {};
\node (12) at (4,0) [shape=circle, draw, fill=black] {};
\node (13) at (6,2)  [shape=circle, draw, fill=black] {};
\node (14) at (6,0) [shape=circle, draw, fill=black] {};

\draw (1) to (2);
\draw (3) to (4);
\draw (5) to (6);
\draw (7) to (8);
\draw (9) to (10);
\draw (11) to (12);
\draw (13) to (14);

\draw [thick, dashed] (2) to (4);
\draw [thick, dashed] (6) to (8);
\draw [thick, dashed] (9) to (11);
\draw [thick, dashed] (10) to (12);
\draw [thick, dashed] (13) to (7,2);
\draw [thick, dashed] (14) to (7,0);

\draw [thick, dotted] (3) to (5);
\draw [thick, dotted] (4) to (6);
\draw [thick, dotted] (7) to (9);
\draw [thick, dotted] (8) to (10);
\draw [thick, dotted] (11) to (13);
\draw [thick, dotted] (12) to (14);

\node at (-6,3) {$1$};
\node at (-6,-1) {$2$};
\node at (-4,3) {$3$};
\node at (-4,-1) {$4$};
\node at (-2,3) {$5$};
\node at (-2,-1) {$6$};
\node at (0,3) {$7$};
\node at (0,-1) {$8$};
\node at (2,3) {$9$};
\node at (2,-1) {$10$};
\node at (4,3) {$11$};
\node at (4,-1) {$12$};
\node at (6,3) {$13$};
\node at (6,-1) {$14$};

\path [fill=lightgray] (-3.95,0.05) rectangle (-2.05,1.95);
\path [fill=lightgray] (0.05,0.05) rectangle (1.95,1.95);
\path [fill=lightgray] (4.05,0.05) rectangle (5.95,1.95);

%%%%%%

\node (n-11) at (-4,-3) [shape=circle, draw, fill=black] {};
\node (n-10) at (-4,-5) [shape=circle, draw, fill=black] {};
\node (n-9) at (-2,-3)  [shape=circle, draw, fill=black] {};
\node (n-8) at (-2,-5) [shape=circle, draw, fill=black] {};
\node (n-7) at (0,-3)  [shape=circle, draw, fill=black] {};
\node (n-6) at (0,-5) [shape=circle, draw, fill=black] {};
\node (n-5) at (2,-3)  [shape=circle, draw, fill=black] {};
\node (n-4) at (2,-5) [shape=circle, draw, fill=black] {};
\node (n-3) at (4,-3)  [shape=circle, draw, fill=black] {};
\node (n-2) at (4,-5) [shape=circle, draw, fill=black] {};
\node (n-1) at (6,-3)  [shape=circle, draw, fill=black] {};
\node (n) at (6,-5) [shape=circle, draw, fill=black] {};

\draw (n-11) to (n-10);
\draw (n-9) to (n-8);
\draw (n-7) to (n-6);
\draw (n-5) to (n-4);
\draw (4.1,-3) to (4.1,-5);
\draw (n-1) to (n);

\draw [thick, dashed] (-5,-3) to (n-11);
\draw [thick, dashed] (-5,-5) to (n-10);
\draw [thick, dashed] (n-9) to (n-7);
\draw [thick, dashed] (n-5) to (n-3);
\draw [thick, dashed] (n-2) to (n);

\draw [thick, dotted] (n-11) to (n-9);
\draw [thick, dotted] (n-10) to (n-8);
\draw [thick, dotted] (n-7) to (n-5);
\draw [thick, dotted] (n-6) to (n-4);
\draw [thick, dotted] (3.9,-3) to (3.9,-5);

\node at (-4,-2) {$n-11$};
\node at (-4,-6) {$n-10$};
\node at (-2,-2) {$n-9$};
\node at (-2,-6) {$n-8$};
\node at (0,-2) {$n-7$};
\node at (0,-6) {$n-6$};
\node at (2,-2) {$n-5$};
\node at (2,-6) {$n-4$};
\node at (4,-2) {$n-3$};
\node at (4,-6) {$n-2$};
\node at (6,-2) {$n-1$};
\node at (6,-6) {$n$};

\path [fill=lightgray] (-3.95,-4.95) rectangle (-2.05,-3.05);
\path [fill=lightgray] (0.05,-4.95) rectangle (1.95,-3.05);

 \end{tikzpicture}

\end{center}
\caption{The permutations $r_i$ for $n=4k+2$} 
\label{case4}
\end{figure}

These permutations, shown in Figure~\ref{case4}, are odd involutions in $S_n$, and $r_0$ commutes with $r_2$. The element $r_1r_2$ fixes $1$ and $n-1$, has two $2$-cycles $(3,5)$ and $(n-6,n-4)$, and has two cycles 
\[(2,6,10,\ldots,n-8,n-10,n-14,\ldots,8,4)\]
and
\[(7,9,13,17,\ldots,n-5,n-2,n,n-3,n-7,\ldots, 11)\]
of lengths $2k-3$ and $2k-1$. 

We need to show that the group $G=\langle r_0, r_1, r_2\rangle$ is primitive, so suppose not. Now $2(2k-3)$ is coprime to $2k-1$, so $(r_1r_2)^{2(2k-3)}$ is a cycle $c$ of length $2k-1$. This is coprime to $n$, so $c$ cannot be a union of blocks of imprimitivity, and hence some block $B$ meets the support $C$ of $c$ and its complement. Since $B$ contains a fixed point of $c$ it is invariant under $c$, and hence contains all of $C$. Thus $|B|>2k-1$, and since $|B|$ divides $n$ it follows that $|B|=2k+1$ and there are just two blocks. Since $r_1$ and $r_2$ have fixed points they preserve the blocks, so $r_0$ must transpose them. This is impossible, since $r_0$ and $r_2$ both transpose $n-3$ and $n-2$. Thus $G$ is primitive, so we can apply Theorem~\ref{J-Jordan} to the cycle $c$ to show that $G\ge A_n$. Since $G$ contains odd permutations we have $G=S_n$, as required.

\begin{figure}[h!]

\begin{center}
\begin{tikzpicture}[scale=0.7, inner sep=0.8mm]
 
\node (1) at (-6,2)  [shape=circle, draw, fill=black] {};
\node (2) at (-6,0) [shape=circle, draw, fill=black] {};
\node (3) at (-4,2)  [shape=circle, draw, fill=black] {};
\node (4) at (-4,0) [shape=circle, draw, fill=black] {};
\node (5) at (-2,2)  [shape=circle, draw, fill=black] {};
\node (6) at (-2,0) [shape=circle, draw, fill=black] {};
\node (7) at (0,2)  [shape=circle, draw, fill=black] {};
\node (8) at (0,0) [shape=circle, draw, fill=black] {};
\node (9) at (2,2)  [shape=circle, draw, fill=black] {};
\node (10) at (2,0) [shape=circle, draw, fill=black] {};
\node (11) at (4,2)  [shape=circle, draw, fill=black] {};
\node (12) at (4,0) [shape=circle, draw, fill=black] {};
\node (13) at (6,2)  [shape=circle, draw, fill=black] {};
\node (14) at (6,0) [shape=circle, draw, fill=black] {};

\draw (1) to (2);
\draw (3) to (4);
\draw (5) to (6);
\draw (7) to (8);
\draw (9) to (10);
\draw (4.1,2) to (4.1,0);
\draw (13) to (14);

\draw [thick, dashed] (2) to (4);
\draw [thick, dashed] (5) to (7);
\draw [thick, dashed] (6) to (8);
\draw [thick, dashed] (9) to (11);
\draw [thick, dashed] (12) to (14);

\draw [thick, dotted] (3) to (5);
\draw [thick, dotted] (4) to (6);
\draw [thick, dotted] (7) to (9);
\draw [thick, dotted] (8) to (10);
\draw [thick, dotted] (3.9,2) to (3.9,0);

\node at (-6,3) {$1$};
\node at (-6,-1) {$2$};
\node at (-4,3) {$3$};
\node at (-4,-1) {$4$};
\node at (-2,3) {$5$};
\node at (-2,-1) {$6$};
\node at (0,3) {$7$};
\node at (0,-1) {$8$};
\node at (2,3) {$9$};
\node at (2,-1) {$10$};
\node at (4,3) {$11$};
\node at (4,-1) {$12$};
\node at (6,3) {$13$};
\node at (6,-1) {$14$};

\path [fill=lightgray] (-3.95,0.05) rectangle (-2.05,1.95);
\path [fill=lightgray] (0.05,0.05) rectangle (1.95,1.95);

 \end{tikzpicture}

\end{center}
\caption{The permutations $r_i$ for $n=14$} 
\label{case4n=14}
\end{figure}

Now suppose that $n=14$. In this case let
\[r_0=(1,2)(3,4)\ldots(13,14)\]
as before, but let
\[r_1=(2,4)(5,7)(6,8)(9,11)(12,14),\]
\[r_2=(3,5)(4,6)(7,9)(8,10)(11,12),\]
as shown in Figure~\ref{case4n=14}. Then
\[r_1r_2=(2,6,10,8,4)(3,5,9,12,14,11,7),\]
with cycle-structure $1^2,5,7$, so that $c:=(r_1r_2)^5$ is a $7$-cycle. Suppose that $G=\langle r_0, r_1, r_2\rangle$ is imprimitive.  Since $c$ has fixed points, the only possible blocks are the support of $c$ and its complement. However, $r_0$ neither preserves nor transposes these two sets, so $G$ is primitive. Now Jordan's Theorem~\ref{Jordan}, applied to $c$, implies that $G=S_{14}$, as required.

When $n=10$ let
\[r_0=(1,2)(3,4)\ldots(9,10)\]
as before, but let
\[r_1=(2,4)(5,7)(8,10),\]
\[r_2=(3,5)(4,6)(7,8),\]
as shown in Fig.~\ref{case4n=10}. In this case
\[r_1r_2=(2,6,4)(3,5,8,10,7),\]
so $(r_1r_2)^3$ is a $5$-cycle, and an argument similar to that used for $n=14$ shows that $\langle r_0, r_1, r_2\rangle=S_{10}$.

\begin{figure}[h!]

\begin{center}
\begin{tikzpicture}[scale=0.7, inner sep=0.8mm]
 
\node (1) at (-6,2)  [shape=circle, draw, fill=black] {};
\node (2) at (-6,0) [shape=circle, draw, fill=black] {};
\node (3) at (-4,2)  [shape=circle, draw, fill=black] {};
\node (4) at (-4,0) [shape=circle, draw, fill=black] {};
\node (5) at (-2,2)  [shape=circle, draw, fill=black] {};
\node (6) at (-2,0) [shape=circle, draw, fill=black] {};
\node (7) at (0,2)  [shape=circle, draw, fill=black] {};
\node (8) at (0,0) [shape=circle, draw, fill=black] {};
\node (9) at (2,2)  [shape=circle, draw, fill=black] {};
\node (10) at (2,0) [shape=circle, draw, fill=black] {};

\draw (1) to (2);
\draw (3) to (4);
\draw (5) to (6);
\draw (0.1,2) to (0.1,0);
\draw (9) to (10);

\draw [thick, dashed] (2) to (4);
\draw [thick, dashed] (5) to (7);
\draw [thick, dashed] (8) to (10);

\draw [thick, dotted] (3) to (5);
\draw [thick, dotted] (4) to (6);
\draw [thick, dotted] (-0.1,2) to (-0.1,0);

\node at (-6,3) {$1$};
\node at (-6,-1) {$2$};
\node at (-4,3) {$3$};
\node at (-4,-1) {$4$};
\node at (-2,3) {$5$};
\node at (-2,-1) {$6$};
\node at (0,3) {$7$};
\node at (0,-1) {$8$};
\node at (2,3) {$9$};
\node at (2,-1) {$10$};

\path [fill=lightgray] (-3.95,0.05) rectangle (-2.05,1.95);

 \end{tikzpicture}

\end{center}
\caption{The permutations $r_i$ for $n=10$} 
\label{case4n=10}
\end{figure}

However, we will now show that $S_6$ cannot be the automorphism group of an orientable regular map. The generators $r_i$ would need to be either transpositions or triple transpositions. Now $S_6$ has an outer automorphism transposing these two conjugacy classes. Clearly $S_6$ cannot be generated by three transpositions, so the same applies to three triple transpositions. If $r_0$ and $r_2$ are both transpositions, they must have disjoint supports in order to commute, so without loss of generality we can take $r_0=(1,3)$ and $r_2=(2,4)$; since $G=\langle r_0, r_1, r_2\rangle$ must be transitive, we may assume without loss of generality that $r_1=(1,6)(2,3)(4,5)$, so that $G$ is imprimitive, preserving equivalence mod~$(2)$. Similarly, applying the outer automorphism excludes the possibility that $r_0$ and $r_2$ are both triple transpositions, so we may assume that $r_0$ is a triple transposition, say $(1,2)(3,4)(5,6)$, and $r_2$ is a single transposition; since they commute, we may assume that $r_2=(1,2)$; again, since $G$ is transitive we may assume that $r_2=(1,6)(2,3)(4,5)$, so that $G$ is imprimitive as before. Thus $S_6$ cannot be evenly realised.

Finally, $S_2$ can be evenly realised by taking $r_0=r_1=r_2=(1,2)$; the resulting regular map, with $M=\Gamma^+$, is the map ${\mathcal N}(2^P{\rm ex})$ consisting of a single vertex and a semi-edge on the sphere. 

%%%%%%%%%%%

\subsection{Proof of (b), with $T=2^{\sigma}$}

If $T=2^{\sigma}$ we have
\[N(T)=\langle S_1, S_2, S_3\mid S_i^2=1\rangle\cong C_2*C_2*C_2.\]
For epimorphisms $N(T)\to G=\langle s_1, s_2, s_3\rangle$, with $S_i\mapsto s_i$, there is one forbidden automorphism, transposing $s_1$ and $s_2$, while fixing $s_3$. If  $T=2$ or $2^*$ then each generator $S_i$ is orientation-reversing, so we require epimorphisms $N(T)\to S_n$, $S_i\mapsto s_i$, with each involution $s_i$ odd. If $T=2^P$ then only $S_1$ and $S_2$ are orientation-reversing, so we require $s_1$ and $s_2$ to be odd, with no restriction on $s_3$.

By the proof of part~(a), if $n\ne 1, 5$ or $6$ then $S_n$ is evenly realised for class~1 by means of an epimorphism $\Gamma\to S_n$, $R_i\mapsto r_i\;(i=0,1,2)$, so we can define $s_i=r_{i-1}$, giving three odd involutions generating $S_n$. If $n\ge 3$ then $s_1$ commutes with $s_3$, whereas $s_2$ does not (since $r_2$ is not central in $S_n$), so there is no automorphism fixing $s_2$ and transposing $s_1$ and $s_2$. This evenly realises $S_n$ for each of the three classes $T=2^{\sigma}$ for all $n\ne 1, 2, 5$ or $6$.

Clearly the symmetric groups $S_1$ and $S_2$ cannot be evenly realised, and neither can $S_5$ for $T=2$ or $2^*$ since it cannot be generated by three transpositions. However, $S_5$ can be evenly realised for $T=2^P$ by taking
\[s_1=(1,2), \quad s_2=(3,4) \quad{\rm and}\quad s_3=(1,3)(4,5).\]
The group $G$ which these generate is transitive and hence primitive, and it contains a transposition $s_1$, so $G=S_5$. The centraliser of $s_3$ in ${\rm Aut}\,G=S_5$ is the Klein four-group fixing $2$, and this cannot transpose $s_1$ and $s_2$.

To evenly realise $S_6$ for $T=2^P$ take
\[s_1=(1,3), \quad s_2=(1,5)(2,3)(4,6) \quad{\rm and}\quad s_3=(1,2)(3,4).\]
These generate a transitive group $G$, and since $s_1s_2s_3=(2,4,6,3,5)$ this group is doubly transitive and hence primitive. SInce $s_1\in G$ we have $G=S_6$. Since $s_1s_3=(1,4,3,2)$ has order $4$, while $s_2s_3=(1,5,2,4,6,3)$ has order $6$, no automorphism can fix $s_3$ while transposing $s_1$ and $s_2$.

However, $S_6$ cannot be evenly realised for $T=2$ or $2^*$. It clearly cannot be generated by three transpositions, and hence (by applying the outer automorphism of $S_6$) the same applies to three triple transpositions. The only possibilities are therefore two transpositions and one triple transposition, or vice versa, and again these are equivalent under the outer automorphism. In the former case we can assume without loss of generality that $s_1, s_2$ and $s_3$ are $(1,2)(3,4)(5,6)$, $(2,4)$ and $(3,5)$, in which case they generate an imprimitive group, preserving the relation of congruence mod~$(2)$. Thus $S_6$ cannot be evenly realised for these two classes.

%%%%%%%%%%%

\subsection{Proof of (e), with $T=4^{\sigma}$}

If $T=4^{\sigma}$ we have
\[N(T)=\langle S_1, S_2, S\mid S_i^2=1\rangle\cong C_2*C_2*C_{\infty},\]
and for epimorphisms $N(T)\to G=\langle s_1, s_2, s\rangle$ there is one forbidden automorphism, transposing $s_1$ and $s_2$, while inverting $s$. If  $T=4$ or $4^*$ then each $S_i$ is orientation-reversing while $S$ is orientation-preserving, but if $T=4^P$ then all three generators of $N(T)$ are orientation-reversing.

To evenly realise $S_n$ for $n\ne 1, 2, 5$ or $6$ we can use the generators $s_1, s_2$ and $s$ used for types $T=2^{\sigma}$.

To evenly realise $S_5$ for $T=4^{\sigma}$, take
\[s_1=(1,2), \quad s_2=(3,4) \quad {\rm and} \quad s=(2,3,4,5).\]
These generate a group which is transitive and hence primitive since its degree is prime; since it contains $s_1$ this group is $S_5$. Since $s_1s=(1,3,5,5,2)$ and $s_2s^{-1}=(2,5,4)$ have different orders, there is no automorphism transposing $s_1$ and $s_2$ and inverting $s$.

To evenly realise $S_6$ for $T=4^{\sigma}$, take
\[s_1=(1,2), \quad s_2=(3,4) \quad{\rm and}\quad s=(1,5,6)(2,3).\]
These generate a transitive group $G$, and since it contains $s_1s=(1,3,2,5,6)$ it is doubly transitive and hence primitive. Since $s_1\in G$ we have $G=S_6$. Since $s_2s^{-1}=(1,6,5)(2,3,4)$ has order $3$, no automorphism can transpose $s_1$ and $s_2$ while inverting $s$.

%%%%%%%%%%%

\subsection{Proof of (d), with $T=3$}

Now let $T=3$, the class of just-edge-transitive maps. We have
\[N(3)=\langle S_0, \ldots, S_3\mid S_i^2=1\rangle \cong C_2*C_2*C_2*C_2,\]
where each $S_i$ is a reflection, so we require epimorphisms $N(3)\to S_n$, $S_i\mapsto s_i$, with each $s_i$ an odd involution. The forbidden automorphisms are those inducing double transpositions of these generators. 

Write $n=m+r$ where $m\equiv 3$ mod~$(4)$ and $r=0, 1, 2$ or $3$. Let
\[s_0=(1,2)(3,4)\ldots. (m-2,m-1),\]
\[s_1=(2,3)(4,5)\ldots(m-1,m),\]
odd involutions such that
\[s_0s_1=(1,3,5,\ldots,m,m-1,m-3,\ldots, 2)\]
is an $m$-cycle, fixing $r$ points. Let $G:=\langle s_0, s_1, s_2, s_3\rangle$, where $s_2$ and $s_3$ are odd involutions to be defined in the various cases below.

If $n\equiv 3$ mod~$(4)$, so that $m=n$ and $r=0$, define
\[s_2=s_3=(1,2).\]
Then $\langle s_0s_1,s_2\rangle =S_n$ by Lemma~\ref{c,t}, so we have an epimorphism $N(3)\to G=S_n$ with kernel $M\le \Gamma^+$. Since $s_0, s_2$ and $s_3$ commute with each other, whereas none of them commutes with $s_1$, no automorphism of $S_n$ can induce a double transposition on these four generators, so the corresponding map is in class~$3$.

Similar arguments apply in the other cases. If $n\equiv 0$ mod~$(4)$, so that $m=n-1$ and $r=1$, define
\[s_2=s_3=(1,n),\]
so that
\[s_0s_1s_2=(1,3,5,\ldots,n-1,n-2,n-4,\ldots, 2,n)\]
is an $n$-cycle and $\langle s_0s_1s_2, s_2\rangle=S_n$.

If $1<n\equiv 1$ mod~$(4)$, so that $m=n-2$ and $r=2$, define
\[s_2=(1,n)\quad {\rm and}\quad s_3=(2,n-1),\]
so that
\[s_0s_1s_2=(1,3,5,\ldots,n-2,n-3,n-5,\ldots, 4,n-1,2,n)\]
is an $n$-cycle and $\langle s_0s_1s_2s_3, s_2\rangle=S_n$. 

If $2<n\equiv 2$ mod~$(4)$, so that $m=n-3$ and $r=3$, define
\[s_2=(1,n)(2,n-1)(3,n-2)\quad {\rm and}\quad s_3=(1,n),\]
so that 
\[s_0s_1s_2=(1,n-2,3,5,\ldots,n-3,n-4,n-6,\ldots, 4,n-1,2,n)\]
is an $n$-cycle and $\langle s_0s_1s_2, s_3\rangle=S_n$. In all three cases, it is easy to check that the cycle-structures of the generators and their commuting relations prevent any forbidden automorphisms from arising.

However, if $n=2$ the only choice for the generators is $s_i=(1,2)$ for all $i$, so there are forbidden automorphisms, while if $n=1$ there are no odd involutions to choose.

\medskip

\noindent{\bf Example} When $n=3$ the construction described above leads to a just-edge-transitive map on the sphere, in which a vertex of valency $6$ is joined by double edges to three vertices of valency $2$; the faces are three digons and one hexagon.

%%%%%%%%%%%

\subsection{Proof of (c), with $T=2^{\sigma}{\rm ex}$}

If $T=2^{\sigma}{\rm ex}$ we have
\[N(T)=\langle S_1, S\mid S_1^2=1\rangle \cong C_2*C_{\infty},\]
with $S$ preserving orientation, but $S_1$ preserving it if and only if $\sigma=P$ (so that $T$ is the class of orientably regular chiral maps). The forbidden automorphism is that which fixes $s_1$ and inverts $s$.

The proof of Theorem~\ref{symgps} gives a generating pair
\[s_1=(1,3)(2,4)(n-1,n) \quad{\rm and} \quad s=(1,2,\ldots,n-1)\]
for $S_n$ for each $n\ge 7$, with $s_1$ an odd involution and with no forbidden automorphism, so this evenly realises such groups $S_n$ for each class $2^{\sigma}{\rm ex}$. When $n=6$ the generators
\[(1,\ldots,6)\quad{\rm and}\quad (1,2)(3,5)\]
used there deal with the class $2^P{\rm ex}$, but for classes $2\,{\rm ex}$ and $2^*{\rm ex}$ we need the involution $s_1$ to be odd, that is, either a transposition or a triple transposition. Suppose first that $s_1$ is a transposition; if $\langle s_1, s\rangle$ is to be transitive then $s$ can have at most two cycles, in which case one easily sees by drawing permutation diagrams  that there is always a forbidden automorphism. Thus $s_1$ cannot be a transposition, and by applying the outer automorphism of $S_6$ we can eliminate the case where $s_1$ is a triple transposition.

By the proof of Theorem~\ref{symgps}, $S_n$ cannot be evenly realised by any class $2^{\sigma}{\rm ex}$ if $n\le 5$.

%%%%%%%%%%%

\subsection{Proof of (f), with $T=5^{\sigma}$}

If $T=5^{\sigma}$ then
\[N(T)=\langle S, S'\mid --\rangle \cong F_2,\]
with $S$ and $S'$ both even if $\sigma=\emptyset$ or $*$, but both odd if $\sigma=P$. The forbidden automorphisms are those inverting, transposing, or inverting and transposing their images $s$ and $s'$.

First suppose that $n\ge 7$. If $T=5$ or $5^*$, or if $T=5^P$ and $n$ is odd, then the generators
\[(1,3)(2,4)(n-1,n) \quad{\rm and} \quad (1,2,\ldots,n-1)\]
used above for the classes $2^{\sigma}{\rm ex}$ can be used as $s$ and $s'$. We saw earlier that they are not simultaneously inverted, and since they have different orders, no automorphism of $S_n$ can transpose them or invert and transpose them. If $T=5^P$ and $n$ is even let
\[s=(1,2)(3,4,5)\quad{\rm and}\quad s'=(1,2,\ldots, n).\]
These generate $S_n$ since $s^3=(1,2)$. They have different orders, and a permutation diagram shows that they are not simultaneously inverted by conjugation, so there are no forbidden automorphisms.

If $n=6$ then for any $T=5^{\sigma}$ let
\[s=(1,2,5,3)\quad{\rm and}\quad s'=(1,2,\ldots, 6).\]
Then
\[s^2s'=(1,6)(2,4,5),\]
so the elements $(s^2s')^3=(1,6)$ and $s'$ generate $S_6$. They have different orders, so no automorphism can transpose them or invert and transpose them. A permutation diagram shows that no inner automorphism can invert them, and outer automorphisms cannot do so since they send $6$-cycles to elements with cycle structure $123$.

Finally, Lemma~\ref{smallSn} shows that $S_n\not\in\G(5^{\sigma})$ if $n\le 5$. \hfill$\square$

%%%%%%

\section{Realising alternating groups}\label{alt}

In order to prove an analogue of Theorem~\ref{symgps} for the alternating groups, we need to have some suitable generators for these groups. 
\begin{lemma}\label{Angens}
The following sets of elements each generate the alternating group $A_n$:
\begin{enumerate}[\rm(a)]
\item $(1,2,3), (2,3,4),\ldots, (n-2,n-1,n)$ for each $n\ge 3$;
\item $(1,2,3), (1,3,4),\ldots, (1,n-1,n)$ for each $n\ge 3$;
\item $a:=(k, k+1, k+2)$, $c:=(1,2,\ldots,n)$ for each odd $n\ge 3$ and any $k$;
\item $a:=(1,k, k+1)$, $c:=(2,3,\ldots,n)$ for each even $n\ge 4$ and $k\ne 1, n$;
\end{enumerate}
\end{lemma}

\noindent(Here we reduce the entries of $a$ mod $(n)$ if they lie outside $\{1, 2, \ldots, n\}$.)

\medskip

\noindent{\sl Proof.} Statements (a) and (b) can easily be proved by induction on $n$, using the fact that $A_n$ is a maximal subgroup of $A_{n+1}$. For (c) and (d), conjugating $a$ by powers of $c$ gives the generators in (a) or (b).\hfill$\square$

\medskip

The following result proves the statements in Theorem~\ref{mainthm} concerning the alternating groups:

\begin{thm}\label{altgps}
The alternating group $A_n$ is the automorphism group of a map in an edge-transitive class $T$ if and only if one of the following holds:
\begin{itemize}
\item $T=1$, and $n=1, 2, 5$ or $n\ge 9$;
\item $T=2^{\sigma}$ for some $\sigma$ or $T=3$, and $n\ge 5$;
\item $T=2^{\sigma}{\rm ex}$ for some $\sigma$, and $n\ge 8$;
\item $T=4^{\sigma}$ for some $\sigma$, and $n\ge 4$;
\item $T=5^{\sigma}$ for some $\sigma$, and $n\ge 7$.
\end{itemize}
\end{thm}

\noindent{\sl Proof.} The cases $n\le 3$ are covered by Lemma~\ref{abelian}, and $A_4$ is easily dealt with, so we may assume that $n\ge 5$. If $n\ne 6$ then all automorphisms of $A_n $ are induced by conjugation in $S_n$, which preserves cycle structure, so for convenience we will first consider larger values of $n$, leaving smaller values until later.

First let $T=1$, so that $N(T)=\Gamma\cong V_4*C_2$. Nuzhin~\cite{Nuz92}, responding to a question of Mazurov (see Section~\ref{gensimple}), showed that $A_n$ is a quotient of $\Gamma$ if and only if $n=5$ or $n\ge 9$. Specifically, he showed that in these cases one can map the standard generators $R_i$ of $\Gamma$ to generators $r_i$ of $A_n$, where
\[r_0=(1,4)(2,3)(5,6)(n-2,n-1)\quad \hbox{if}\quad n=4k+3\ge 11,\]
\[r_0=(1,2)(3,4)\quad{\rm otherwise},\]
and $r_2$ and $r_1$ are respectively 
\[(1,4)(2,3),\; (2,3)(4,5)\quad{\rm if}\quad n=5,\] 
\[(1,2)(3,4)\ldots(n-2,n-1),\; (2,3)(4,5)\ldots(n-1,n)\;{\rm if}\; n=4k+1\ge 9,\]
\[(3,4)(5,6)\ldots(n-1,n),\; (2,3)(4,5)\ldots(n-2,n-1)\;{\rm if}\; n=4k+2\ge 10,\]
\[(1,2)(3,4)\ldots(n-4,n-3),\; (4,5)(6,7)\ldots(n-1,n)\;{\rm if}\; n=4k+3\ge 11,\]
\[(1,2)(3,4)\ldots(n-1,n),\; (2,3)(4,5)(6,7)\ldots(n-2,n-1)(1,n)\;{\rm if}\; n=4k\ge 12.\]
This deals with $T=1$, and hence by Lemma~\ref{regmapslemma}(a) it also realises $A_n$ for $n=5$ and $n\ge 9$, when  $T=2^{\sigma}$, $3$ or $4^{\sigma}$. However, for these classes we will later give direct arguments which realise $A_n$ for a wider range of values of $n$.

We next consider $T=2^P{\rm ex}$, with $N(T)=\Gamma^+=\langle X, Y\mid Y^2=1\rangle$. For even $n\ge 8$ let
\[x=(2, 3, \ldots, n)\quad{\rm and}\quad y=(1,2)(3,4),\]
so $G:=\langle x, y\rangle$ is $2$-transitive and hence primitive. Now
\[[y, x]=y.y^x=(1,2)(3,4).(1,3)(4,5)=(1,2,3,5,4),\]
so by Jordan's Theorem $G=A_n$.  We can therefore define an epimorphism $\Gamma^+\to A_n$ by $X\mapsto x$ and $Y\mapsto y$, so $A_n$ is a quotient of $\Gamma^+$. Since $n\ne 6$ every automorphism of $A_n$ is induced by conjugation in $S_n$; by inspection, no permutation inverts $x$ and centralises $y$, so there are no forbidden automorphisms. For odd $n\ge 9$, let
\[x=(1, 2, \ldots, n)\quad{\rm and}\quad y=(1, 2)(3, 6),\]
and suppose that $G:=\langle x, y\rangle$ is imprimitive. The only non-trivial equivalence relations invariant under $\langle x\rangle$ are those of congruence mod~$(d)$ for some $d$ dividing $n$, with $3\le d\le n/3$ since $n$ is odd. Each equivalence class contains at least three elements, and $y$ must transpose the classes $[1]$ and $[2]$, which is impossible since it moves only four elements. Thus $G$ is primitive. 
Now
\[[y, x^2]=(1   , 2)(3, 6, 4)(5, 8),\]
so $[y, x^2]^2$ is a $3$-cycle and hence $G=A_n$ by Jordan's Theorem. As before, this shows that $A_n$ is a quotient of $\Gamma^+$, and since no permutation inverts $x$ and commutes with $y$, there are no forbidden automorphisms. By Lemma~\ref{regmapslemma}(b) this realises $A_n$ for $n\ge 8$, when  $T=2^{\sigma}{\rm ex}$ or $5^{\sigma}$.

We now consider small values of $n$ not dealt with above.  In the case $T=1$, Nuzhin~\cite{Nuz92} has already shown that there are regular maps with automorphism group isomorphic to $A_n$ for $n=5$ but not for $n=6, 7$ or $8$. (For example, a simple argument using permutation diagrams shows that if a Klein four group and an involution generate a transitive subgroup $G\le A_6$, then either $G$ is imprimitive, with three blocks of size $2$, or $G\cong L_2(5)$; the latter possibility realises $A_5$, which is isomorphic to $L_2(5)$.) This deals with $T=1$, and by Lemma~\ref{regmapslemma}(a) it also realises $A_5$ in the classes $2^{\sigma}$, $3$ and $4^{\sigma}$.

If $T=2^{\sigma}$ then $N(T)\cong C_2*C_2*C_2$. For $n=6$ take
\[s_1=(1,2)(3,4),\quad s_2=(2,6)(4,5)\quad{\rm and}\quad s_3=(2,3)(4,5),\]
so that
\[s_1s_2=(1,6,2)(3,5,4),\quad s_1s_3=(1,3,5,4,2)\quad{\rm and}\quad s_2s_3=(2,6,3).\]
 For $n=7$ take
\[s_1=(1,2)(3,4),\quad s_2=(2,6)(5,7)\quad{\rm and}\quad s_3=(2,3)(4,5),\]
so that
\[s_1s_2=(1,6,2)(3,4)(5,7),\quad s_1s_3=(1,3,5,4,2)\;\;{\rm and}\;\; s_2s_3=(2,6,3)(4,5,7).\]
 For $n=8$ take
\[s_1=(1,2)(3,4)(5,6)(7,8),\quad s_2=(1,3)(4,6)\quad{\rm and}\quad s_3=(3,4)(6,7),\]
so that
\[s_1s_2=(1,2,3,6,5,4)(7,8),\quad s_1s_3=(1,2)(5,7,8,6)\;\;{\rm and}\;\; s_2s_3=(1,4,7,6,3).\]
In each case, $G:=\langle s_1, s_2, s_3\rangle$ is primitive: for instance, when $n=6$ or $8$ it is doubly transitive, in the latter case because  $\langle(s_1s_2)^2, s_3\rangle$ fixes $8$ and is transitive on the remaining points. Jordan's Theorem, applied to a suitable $3$- or $5$-cycle, then implies that $G=A_n$. In each case $s_1s_3$ and $s_2s_3$ have different orders, so there is no automorphism transposing $s_1$ and $s_2$ and inverting (equivalently fixing) $s_3$. Thus $A_6, A_7$ and $A_8$ are realised in these classes, and the same applies to the classes $3$ and $4^{\sigma}$ by Lemma~\ref{class2lemma}(c).

For the classes $T=2^{\sigma}{\rm ex}$ the  lower bound $n\ge 8$ given above cannot be improved: a result of Singerman (see Proposition~\ref{L2qinvgen} in the next section) eliminates the groups $A_4, A_5$ and $A_6$, which are isomorphic to $L_2(q)$ for $q=3$, $4$ (or $5$) and $9$ respectively, while $A_7$ can be eliminated as follows. If there is an epimorphism $\Gamma^+\to A_7$, with $X\mapsto x$ and $Y\mapsto y$, then the involution $y$ has five cycles, namely two transpositions and three fixed points. In order for $x$ and $y$ to generate a transitive group we therefore need $x$ to have a cycle of length at least $5$, so it must be a $5$-cycle or a $7$-cycle. In the first case, it is easy to see (by drawing permutation diagrams, for instance), that any pair generating a transitive group are inverted by some element of $S_7$, giving a forbidden automorphism. If $x$ is a $7$-cycle, then without loss we may assume that it is $(1, 2, \ldots, 7)$. Replacing $x$ with a suitable power we may also assume that $y$ involves a transposition $(i, i+1)$, and then conjugating $x$ and $y$ by a suitable power of $x$ we may assume that $i=1$, so that $y=(1, 2)(j,k)$ for some $j<k$ in $\{3, \ldots, 7\}$. The ten possibilities for $(j,k)$ lead, in the cases $(3,4), (3,7), (4,5), (4,6), (5,6)$ and $(6,7)$, to a pair $x, y$ inverted in $S_7$; the remaining cases $(3,5), (3,6), (4,7)$ and $(5,7)$ give a pair preserving one of the two $\langle x\rangle$-invariant Fano plane geometries on $\F_7$, where the lines are the translates of the quadratic residues $\{1, 2, 4\}$ or the non-residues $\{3, 5, 6\}$, so that $\langle x,y\rangle$ is a proper subgroup isomorphic to $L_3(2)$. Thus $A_7\not\in\G(2^{\sigma}{\rm ex})$.

If $T=5^{\sigma}$, with $N(T)\cong F_2$, we can realise $A_7$ by using generators $(1,2,3,4,5)$ and $(1,6,7)(2,4,5)$: these generate a transitive group and no maximal subgroup of $A_7$ contains elements of orders $3$, $5$ and $7$. By inspection there are no forbidden automorphisms. However, any generating pair for $A_6\cong L_2(9)$, of $A_5\cong L_2(4)$ or of $A_4\cong L_2(3)$ are inverted by an automorphism, so these groups cannot be realised. \hfill $\square$

\medskip

This, together with Theorem~\ref{Sneven}, completes the proof of Theorem~\ref{mainthmeven}.

\section{Realising $2.A_n$}

%{\color{red}[New section]}

%\medskip

For each $n\ge 4$ there is a unique group $G$ with a central involution $z\in G'$ such that $G/\langle z\rangle\cong A_n$. This is the double cover $2.A_n$, which is in fact the universal central extension of $A_n$ if $n\ne 6, 7$. When $n=4$ or $5$ it is isomorphic to  the binary tetrahedral group $SL_2(3)$ or the binary icosahedral group $SL_2(5)$.  One way of constructing $G$ is to embed $A_n$ in the connected Lie group $SO(n)$ as a group of permutation matrices, and then to lift these to the simply connected double cover ${\rm Spin}(n)$ of $SO(n)$. %{\color{blue}[Need to use $SO(n-1)$? $SO(n)$ seems to be OK.]}

%\medskip

%\noindent{\bf Conjecture} For each edge-transitive class $T$ there is an integer $n(T)$ such that $2.A_n\in\G(T)$ for all $n\ge n(T)$.

%\medskip

%As evidence for this, and as an indication of how one might prove the conjecture, we give the following partial verification:

\begin{thm}\label{2Anthm}
{\rm(a)} If $n\ge 9$ then $2.A_n\in\G(T)$ for each edge-transitive class $T\ne 1$.
\vskip2pt
\noindent{\rm(b)} If $n\ge 33$ then  $2.A_n\in\G(1)$. 
\end{thm}

These conditions on $n$ are sufficient,but not necessary: in (a) we can find better lower bounds for some classes $T$, while in (b) there are also lower values of $n$, starting with $9$, $10$ and $11$, such that $2.A_n\in\G(1)$. More precise details will be given later.

In proving Theorem~\ref{2Anthm}, we will first prove part (a), starting with the classes $T=2^P{\rm ex}$ and $2$, and then deducing the result for the other classes $T\ne 1$. We will then prove part (b) by a separate argument. It might seem more efficient to deal with the classes $T=1$ and $2^P{\rm ex}$ first, and then to use Lemma~\ref{regmapslemma} for the other classes. However, the lower bounds we obtain in the hardest case $T=1$ are weaker than those given by applying direct arguments to the other classes, so the approach taken here gives a stronger result.

Before proving this theorem, we need some basic information about the groups $G=2.A_n$ and the epimorphism $G\to A_n$. Each element $g\in A_n$ lifts to two elements $\tilde g$ and $\tilde gz=z\tilde g$ of $G$. If $g$ has odd order $k$ then these have orders $k$ and $2k$. An involution $g\in A_n$ is a product of $2t$ disjoint transpositions for some integer $t$, and $\tilde g$ and $\tilde gz$ have the same order $2$ or $4$ as $t$ is even or odd; in the former case we will say that the permutation $g$ is {\sl doubly even}. Finally, if elements $g_1,\ldots, g_r$ generate $A_n$, then any choice of their lifts $\tilde g_1,\ldots, \tilde g_r$ generate $G$: the subgroup they generate maps onto $A_n$, so its index in $G$ is at most $2$; however, $G$ has no subgroups of index $2$, so it must be $G$.
%{\color{blue}[Lifting elements of even order $k>2$? Needed? Give a reference for these liftings?]}

We will also need the following lemma:
% {\color{blue}[others, e.g. PJC's double transposition result?]}:

\begin{lemma}\label{longcycle}
%If a transitive group $H$ of degree $n$ contains a cycle $c$ of length $m>n/2$ then any system of imprimitivity has blocks of size dividing $m$ and $n$; if $m$ and $n$ are coprime then $H$ is primitive.

If a transitive permutation group $H$ of degree $n$ contains a cycle $c$ of length $m$, then in any system of imprimitivity, either
\begin{itemize}
\item[\rm(a)] the support of $c$ is contained in a single block, or
\item[\rm(b)] the support and the fixed point set of $c$ are both unions of blocks, of size dividing $m$ and $n$.
\end{itemize}
In particular, if $m>n/2$ and $m$ and $n$ are coprime, then $H$ is primitive.
\end{lemma}

\noindent{\sl Proof.} %Suppose that there is a non-trivial $H$-invariant equivalence relation. Its restriction to the support $\Phi$ of $c$ must be congruence mod~$D$ for some subgroup $D$ of $C:=\langle c\rangle$ of order $d$ dividing $m$, with $1<d<m$ since $m>n/2$. The block containing a point $\alpha\not\in\Phi$ must be a union of orbits of $C$; the orbit $\Phi$ is too large to be one of them, so they must all be contained in $\overline\Phi$. Thus $\overline\Phi$ is a union of blocks, and hence so is $\Phi$, so all blocks have size $d$, which divides both $m$ and $n$. The last conclusion follows immediately.
In any system of imprimitivity for $H$, the blocks meeting the support $\Phi$ of $c$ must do so in the orbits on $\Phi$ of some subgroup $D$ of $C:=\langle c\rangle$, of order $d$ dividing $m$. If $d=m$ then $D=C$ and (a) holds, so assume that $d<m$. Then $c$ moves all $m/d$ of these orbits, while fixing all elements of $\overline\Phi$, so no block can meet both $\Phi$ and $\overline\Phi$, and hence (b) holds. The last conclusion follows since the hypotheses exclude the possibilities (a) and (b). \hfill$\square$

\medskip

Note that if $n/3 < m \le n/2$ then either $H$ is primitive, or $n$ is even and (a) holds, with two blocks of size $n/2$. Note also that it follows from Theorem~\ref{J-Jordan} that if $n/2 < m \le n-3$ and $m$ and $n$ are coprime, then $H\ge A_n$.

\medskip

\noindent{\sl Proof of Theorem~\ref{2Anthm}(a).}  By Lemmas~\ref{regmapslemma}(b) and \ref{class2lemma}(c) it is sufficient to prove this result in the cases $T=2^P{\rm ex}$ and $T=2$.

First let $T=2^P{\rm ex}$. If $n$ is even then instead of the generators $x=(2,3,\ldots,n)$ and $y=(1,2)(3,4)$ for $A_n$ used in proving this case of Theorem~\ref{altgps}, let us use
\[x=(2,3,\ldots, n)\quad {\rm and}\quad y=(1,2)(3,4)(5,7)(6,8)\]
(possible since $n\ge 9$), so that either choice of $\tilde y\in G$ is an involution. As before, $\langle x,y\rangle$ is $2$-transitive and hence primitive. Since
\[[y,x]=(1,2,3,5,9,7,4)\]
is a cycle of length $7$, and $n\ge 10$, it follows from Jordan's Theorem that $\langle x, y\rangle=A_n$. Thus $\langle \tilde x, \tilde y\rangle = G$, with $\tilde y^2=1$, so that $G$ is an epimorphic image of $\Gamma^+=N(2^P{\rm ex})$. A forbidden automorphism of $G$, inverting $\tilde x$ and $\tilde y$, would induce an automorphism of its central quotient $A_n$ inverting $x$ and $y$. Since $n\ne 6$ this must be induced by conjugation in $S_n$. However, the chirality of the permutation diagram shows that no such automorphism exists, so $G\in\G(2^P{\rm ex})$.

For odd $n\ge 9$, let
\[x=(1,2,\ldots, n)\quad {\rm and}\quad y=(1,2)(3,7)(4,5)(6,8).\]
The only equivalence relations invariant under $\langle x\rangle$ are congruence mod~$(d)$ for some $d$ dividing $n$. If $\langle x, y\rangle$ is imprimitive then since $n$ is odd there exists such a relation with $3\le d\le n/3$, so each equivalence class has at least three members. The classes containing $1, 2$ and $3$ are all distinct, and are all moved by $y$ (note that $2-1$ and $7-3$ are powers of $2$), so $y$ must move at least nine symbols, a contradiction. Thus $\langle x, y\rangle$ is primitive, and since it contains a $5$-cycle $[y,x]^2=(1,3,9,7,2)$, Jordan's Theorem gives $\langle x, y\rangle=A_n$ and hence $\langle \tilde x, \tilde y\rangle=G$. The forbidden automorphism is eliminated as before, so $G\in\G(2^P{\rm ex})$. 

Now let $T=2$. Suppose that we can realise $A_n$ as an epimorphic image $S_i\mapsto t_i, S\mapsto t$ of $N(2)$, with no forbidden automorphism of $A_n$ transposing $t_1$ and $t_2$ and fixing $t$, where the generators $t_i$ and $t$ of $A_n$ are doubly even, so that they lift to involutions $s_i$ and $s$ which generate $G$. Any forbidden automorphism of $G$ must fix $z$ and hence induce a forbidden automorphism of $A_n$, whereas no such automorphism exists, so $G\in\G(2)$. It is therefore sufficient to find such permutations $t_i$ and $t$.

Write $n=m+r$ where $m\equiv 1$ mod~$(8)$ and $0\le r\le 7$. We first consider the case $r=0$, so that $n=m$. Define doubly even involutions
\[t_1=(1)(2,n)(3,m-1)(4,m-2)\ldots,\]
\[t_2=(1,2)(3,m)(4,m-1)\ldots((m+3)/2),\]
\[t=(2)(1,m)(3,m-1)(4,m-2)\ldots,\]
so that $t_1t_2=(1,2,\ldots, m)$ and $tt_1=(1,2,m)$, and hence $\langle t_1, t_2, t\rangle=A_m=A_n$ by Lemma~??. A forbidden automorphism of $A_n$, transposing $t_1$ and $t_2$ while fixing $t$, must be induced by a permutation in $S_n$ transposing the fixed points $1$ and $(n+3)/2$ of $t_1$ and $t_2$, so it cannot commute with $t$. Thus $G\in\G(2)$. %for all classes $T=2^{\sigma}$. The same applies to the class $T=3$ since we can take the four generating involutions to be $s_1, s_2$, $s$ and $s$.

If $1\le r\le 7$ we define $t_1$ and $t_2$ as above, and vary the definition of $t$, depending on $r$. To ensure that the group $H:=\langle t_1, t_2, t\rangle$ is transitive, we choose the involution $t$ to transpose the elements of $\Psi:=\{m+1,m+2,\ldots, n\}$ with a subset $\Psi'$ of $\Phi:=\{1,2,\ldots,m\}$.

We also need to choose $t$ so that no $g\in S_n$, acting by conjugation, transposes $t_1$ and $t_2$ and fixes $t$. Any such element $g$ inverts $t_1t_2$, so as a permutation $g$ must act on $\Phi$ as one of the $m$ reflections in the dihedral group $\langle t_1, t_2\rangle\cong D_m$. Since $g$ commutes with $t$ it must leave invariant the subset $\overline\Psi'=\Psi t$ of $\Phi$; since these two sets have $r$ and $m$ elements, with $1\le r\le 7$ and $m\ge 9$, we can always choose $\overline\Psi'$ to be invariant under no reflections, thus eliminating forbidden automorphisms, unless $r=1$ or $2$, or $r=7$ and $m=9$. If $r=1$ or $2$ then $t$, being doubly even, must contain at least two more transpositions, each pairing two elements of $\Phi\setminus\Psi'$, and these can be chosen not to be invariant under the reflection preserving $\Psi'$ (specific choices are given later). If $r=7$ and $m=9$ then $t$ contains only one more transposition, and such a choice is impossible; however, in this case $n=16$, $t_1$ and $t_2$ each have eight fixed points while $t$ has none, so we can avoid a forbidden automorphism simply by transposing the roles of $t_2$ and $t$.

As a first step in proving that $H=A_n$, we would like to show that $H$ is primitive. Since $H$ contains a cycle $t_1t_2$ of length $m>n/2$, Lemma~\ref{longcycle} shows that if $H$ is imprimitive, then the blocks have size $d$ dividing $m$ and $n$, and hence dividing $r=n-m$. Now $m$ is odd, and $1\le r\le 7$, so either $d=3$ and $r=3$ or $6$, or $d=r=5$ or $7$. In each of these cases one can ensure that $H$ is primitive be choosing $\Psi'$ not to be a union of blocks in $\Phi$.

If $r\ge 3$ then since $H$ is primitive and contains a cycle $t_1t_2$ of length $m\le n-3$, the extension of Jordan's Theorem in Theorem~\ref{J-Jordan} shows that $H=A_n$.
%{\color{blue}[Proof avoiding CFSG?]}
If $r=1$, so that $n=m+1$, we can choose
\[t=(2,n)(3,m-1)(4,m-2)(5,m-3),\]
 so that $tt_1=(2,n,m)(6,m-4)(7,m-5)\ldots$ and $(tt_1)^2=(2,m,n)$. Note that $\Psi'=\{2\}$, and the reflection $\alpha\mapsto 4-\alpha$ mod~$(m)$ of $\Phi$ fixing $2$ does not commute with $t$, so there is no forbidden automorphism. Since $H$ is primitive and contains a $3$-cycle we have $H=A_n$. If $r=2$, so that $n=m+2$, we can choose 
\[t=(2,n)(m,n-1)(3,m-1)(4,m-2),\]
so that $tt_1=(2,n,m,n-1)(5,m-3)(6,m-4)\ldots$ and $(tt_1)^2=(2,m)(n-1,n)$. Here $\Psi'=\{2,m\}$, and the reflection of $\Phi$ transposing $2$ and $m$ does not commute with $t$, so again a forbidden automorphism is avoided. A primitive group of degree $n>8$ containing a double transposition contains $A_n$ (a well-known folklore theorem), so $H=A_n$. $\hfill\square$

\medskip

\noindent{\bf Remarks 1.} Note that since $A_n\not\in\G(2^{\sigma}{\rm ex})$ for $n\le 7$ by Theorem~\ref{altgps}, we have $G\not\in\G(2^{\sigma}{\rm ex})$ for $n\le 7$. Similarly, $G\not\in\G(5^{\sigma})$ for $n\le 6$.
%{\color{blue}[$n=7, 8$? Other $T$?]}

\medskip

\noindent{\bf 2.} When $n=8$ the elements $x$ and $y$ used in this proof for $T=2^P{\rm ex}$ and even $n$ generate a proper subgroup of $A_8$, isomorphic to $AGL_3(2)$.

%Random trials with GAP suggest that there is no suitable pair in this case, so that $G\not\in\G(2^P{\rm ex})$.  {\color{blue}[Try to confirm this.]}

\medskip

We now consider the remaining and most important case $T=1$. We have $G\in\G(1)$ if and only if $G$ is an epimorphic image of $\Gamma$. In this case so is $A_n$, so that $A_n\in\G(1)$, giving $n=5$ or $n\ge 9$ by Theorem~\ref{altgps}. When $n=5$ we have $G\cong SL_2(5)$; this group has a single involution, $-I$, so it cannot be an epimorphic image of $\Gamma$. We therefore restrict attention to values $n\ge 9$.
 
Nuzhin's generators $r_i\;(i=0, 1, 2)$ for $A_n$, used in proving the case $T=1$ of Theorem~\ref{altgps}, are all doubly even, and thus lift to involutions $\tilde r_i$ in $G$, if and only if $n\equiv 3$ mod~$(8)$; in this case the element
\[r_0r_2=(1,3)(2,4)(7,8)(9,10)\ldots (n-2)(n-1)\]
also lifts to an involution, so that $\tilde r_0$ and $\tilde r_2$ commute. The involutions $\tilde r_i$ generate $G$, so $G\in\G(1)$ provided $11\le n\equiv 3$ mod~$(8)$. However, finding suitable generators when $n\not\equiv 3$ mod~$(8)$ seems to be much more difficult.

Note that in general, if $r_0$ and $r_2$ are doubly even and commute, it does not follow that $r_0r_2$ is doubly even, so the involutions $\tilde r_0$ and $\tilde r_2$ need not commute. For an example, let
\[r_0=(1,2)(3,4)(5,6)(7,8)\quad {\rm and}\quad r_2=(1,3)(2,4) (9,10)(11,12),\]
so that
\[r_0r_2=r_2r_0=(1,3)(2,4)(5,6)(7,8)(9,10)(11,12)\]
lifts to two elements of order $4$. More generally, suppose that $\langle r_0, r_2\rangle\;(\cong V_4)$ has $m$ regular orbits, and $m_i\;(i=0, 2, 02)$ with kernel generated by $r_0, r_2$ or $r_0r_2$; then each of these three elements is even if and only if the number $t_i=2m+\sum_{j\ne i}m_i$ of its transpositions is even (equivalently $\sum_{j\ne i}m_i$ is even), in which case it is doubly even if and only if $t_i\equiv 0$ mod~$(4)$. Thus $r_0$ and $r_2$ lift to commuting involutions if and only if $m_0, m_2$ and $m_{02}$ are all mutually congruent mod~$(4)$, and are congruent to $m$ mod~$(2)$.
%{\color{blue}[Make this a lemma?]}
In the example above we have $m=1$, $m_0=m_2=2$ and $m_{02}=0$, so $t_0=t_2=4$ but $t_{02}=6$, and hence the involutions $\tilde r_0$ and $\tilde r_2$ do not commute.

\medskip

\noindent{\sl Proof of Theorem~\ref{2Anthm}(b).} First let $9\le n\equiv 1$ mod~$(8)$, so $n=4k+1$ where $k$ is even. Define
\[r_0=(1)(2,3)(4,5)(6,7)\ldots (4k-2,4k-1)(4k,4k+1),\]
\[r_1=(1,2)(3,4)(5,6)\ldots (4k-3,4k-2)(4k-1,4k)(4k+1),\]
\[r_2=(1)(2,4)(3,5)(6,8)(7,9)\ldots (4k-2,4k)(4k-1,4k+1),\]
as shown in Figure~\ref{nequiv1(8)}.

\begin{figure}[h!]

\begin{center}
\begin{tikzpicture}[scale=0.6, inner sep=0.8mm]

\node (1) at (-10,0)  [shape=circle, draw, fill=black] {};
\node (2) at (-8,0)  [shape=circle, draw, fill=black] {};
\node (3) at (-7,1)  [shape=circle, draw, fill=black] {};
\node (4) at (-7,-1)  [shape=circle, draw, fill=black] {};
\node (5) at (-6,0)  [shape=circle, draw, fill=black] {};
\node (6) at (-4,0)  [shape=circle, draw, fill=black] {};
\node (7) at (-3,1)  [shape=circle, draw, fill=black] {};
\node (8) at (-3,-1)  [shape=circle, draw, fill=black] {};
\node (9) at (-2,0)  [shape=circle, draw, fill=black] {};
\node (4k-2) at (4,0)  [shape=circle, draw, fill=black] {};
\node (4k-1) at (5,1)  [shape=circle, draw, fill=black] {};
\node (4k) at (5,-1)  [shape=circle, draw, fill=black] {};
\node (4k+1) at (6,0)  [shape=circle, draw, fill=black] {};

\draw (2) to (3);
\draw (4) to (5);
\draw (6) to (7);
\draw (8) to (9);
\draw (4k-2) to (4k-1);
\draw (4k) to (4k+1);

\draw [thick, dotted] (2) to (4);
\draw [thick, dotted] (3) to (5);
\draw [thick, dotted] (6) to (8);
\draw [thick, dotted] (7) to (9);
\draw [thick, dotted] (4k-2) to (4k);
\draw [thick, dotted] (4k-1) to (4k+1);

\draw [thick, dashed] (1) to (2);
\draw [thick, dashed] (3) to (4);
\draw [thick, dashed] (5) to (6);
\draw [thick, dashed] (7) to (8);
\draw [thick, dashed] (9) to (-1,0);
\draw [thick, dashed] (3,0) to (4k-2);
\draw [thick, dashed] (4k-1) to (4k);

\draw [thick, dotted] (0.5,0) to (1.5,0);

\node at (-10,-0.7) {$1$};
\node at (-8,0.7) {$2$};
\node at (-7,1.7) {$3$};
\node at (-7,-1.7) {$4$};
\node at (-6,-0.7) {$5$};
\node at (-4,0.7) {$6$};
\node at (-3,1.7) {$7$};
\node at (-3,-1.7) {$8$};
\node at (-2,-0.7) {$9$};
\node at (3,0.7) {$4k-2$};
\node at (5,1.7) {$4k-1$};
\node at (5,-1.7) {$4k$};
\node at (7,-0.7) {$4k+1$};
 
 \end{tikzpicture}

\end{center}
\caption{The permutations $r_i$ for $n\equiv 1$ mod~$(8)$} 
\label{nequiv1(8)}
\end{figure}

Here $m=k$ and $m_0=m_2=m_{02}=0$, while $r_1$ is a product of $2k$ transpositions, so $r_0$, $r_1$, $r_2$ and $r_0r_2$ all lift to involutions in $G$. The group $H=\langle r_0, r_1, r_2\rangle$ is clearly transitive. Since
\[r_0r_2r_1=(1,2,6,10,\ldots,4k-2,4k+1,4k-3,4k-7,\ldots,9,5)\]
is a cycle of length $2k+1>n/2$ coprime to $n$, it follows from Lemma~\ref{longcycle} that $H$ is primitive, and since $2k+1\le n-3$ it follows from Theorem~\ref{J-Jordan} that $H=A_n$. It therefore follows from our earlier remarks that $G\in\G(1)$ if $9\le n\equiv 1$ mod~$(8)$.

In the remaining seven cases for $n$ mod~$(8)$ we make minor modifications to the above permutations $r_i$: we retain the permutations $r_0$ and $r_2$, simply adding up to seven further symbols $4k+2, \ldots, 4k+8$ all fixed by them, and we redefine $r_1$ so that it still has $2k$ transpositions, including one pairing each of the new symbols with an original symbol. The constructions ensure that the permutations $r_i$ generate a transitive group $H$, and lift to involutions $\tilde r_i\in G$ with $\tilde r_0$ and $\tilde r_2$ commuting. In each case the cycle of length $2k+1$ in $r_0r_2r_1$ is retained, enabling us to prove that $H=A_n$, so that the elements $\tilde r_i$ generate $G$. 

It is convenient to start with the remaining odd values of $n$. First let $11\le n\equiv 3$ mod~$(8)$, so $n=4k+3$ where $k$ is even. Define
\[r_0=(1)(2,3)(4,5)(6,7)\ldots (4k-2,4k-1)(4k,4k+1)(4k+2)(4k+3),\]
\[r_1=(1,2)(3,4k+2)(4,4k+3)(11,12)\ldots (4k-3,4k-2)(4k-1,4k)(4k+1),\]
\[r_2=(1)(2,4)(3,5)(6,8)(7,9)\ldots (4k-2,4k)(4k-1,4k+1)(4k+2)(4k+3),,\]
as shown in Figure~\ref{nequiv3(8)}.

\begin{figure}[h!]

\begin{center}
\begin{tikzpicture}[scale=0.6, inner sep=0.8mm]

\node (1) at (-10,0)  [shape=circle, draw, fill=black] {};
\node (2) at (-8,0)  [shape=circle, draw, fill=black] {};
\node (3) at (-7,1)  [shape=circle, draw, fill=black] {};
\node (4) at (-7,-1)  [shape=circle, draw, fill=black] {};
\node (5) at (-6,0)  [shape=circle, draw, fill=black] {};
\node (6) at (-4,0)  [shape=circle, draw, fill=black] {};
\node (7) at (-3,1)  [shape=circle, draw, fill=black] {};
\node (8) at (-3,-1)  [shape=circle, draw, fill=black] {};
\node (9) at (-2,0)  [shape=circle, draw, fill=black] {};
\node (10) at (0,0)  [shape=circle, draw, fill=black] {};
\node (11) at (1,1)  [shape=circle, draw, fill=black] {};
\node (12) at (1,-1)  [shape=circle, draw, fill=black] {};
\node (13) at (2,0)  [shape=circle, draw, fill=black] {};
\node (4k-2) at (8,0)  [shape=circle, draw, fill=black] {};
\node (4k-1) at (9,1)  [shape=circle, draw, fill=black] {};
\node (4k) at (9,-1)  [shape=circle, draw, fill=black] {};
\node (4k+1) at (10,0)  [shape=circle, draw, fill=black] {};
\node (4k+2) at (-7,2.5)  [shape=circle, draw, fill=black] {};
\node (4k+3) at (-7,-2.5)  [shape=circle, draw, fill=black] {};

\draw (2) to (3);
\draw (4) to (5);
\draw (6) to (7);
\draw (8) to (9);
\draw (10) to (11);
\draw (12) to (13);
\draw (4k-2) to (4k-1);
\draw (4k) to (4k+1);

\draw [thick, dotted] (2) to (4);
\draw [thick, dotted] (3) to (5);
\draw [thick, dotted] (6) to (8);
\draw [thick, dotted] (7) to (9);
\draw [thick, dotted] (10) to (12);
\draw [thick, dotted] (11) to (13);
\draw [thick, dotted] (4k-2) to (4k);
\draw [thick, dotted] (4k-1) to (4k+1);

\draw [thick, dashed] (1) to (2);
\draw [thick, dashed] (3) to (4k+2);
\draw [thick, dashed] (5) to (6);
\draw [thick, dashed] (4) to (4k+3);
\draw [thick, dashed] (9) to (10);
\draw [thick, dashed] (11) to (12);
\draw [thick, dashed] (13) to (3,0);
\draw [thick, dashed] (7,0) to (4k-2);
\draw [thick, dashed] (4k-1) to (4k);

\draw [thick, dotted] (4.5,0) to (5.5,0);

\node at (-10,-0.7) {$1$};
\node at (-8,0.7) {$2$};
\node at (-6.5,1.7) {$3$};
\node at (-6.5,-1.7) {$4$};
\node at (-6,-0.7) {$5$};
\node at (-4,0.7) {$6$};
\node at (-3,1.7) {$7$};
\node at (-3,-1.7) {$8$};
\node at (-2,-0.7) {$9$};
\node at (0,0.7) {$10$};
\node at (1,1.7) {$11$};
\node at (1,-1.7) {$12$};
\node at (2,-0.7) {$13$};
\node at (7,0.7) {$4k-2$};
\node at (9,1.7) {$4k-1$};
\node at (9,-1.7) {$4k$};
\node at (11,-0.7) {$4k+1$};
\node at (-5.7,3) {$4k+2$};
\node at (-5.7,-3) {$4k+3$};
 
 \end{tikzpicture}

\end{center}
\caption{The permutations $r_i$ for $n\equiv 3$ mod~$(8)$} 
\label{nequiv3(8)}
\end{figure}

The only change from the preceding case $n\equiv 1$ mod~$(8)$ is that we have added two symbols $4k+2$ and $4k+3$, fixed by $r_0$ and $r_2$, and have replaced the transpositions $(3,4)$ and $(7,8)$ in $r_1$ with $(3,4k+2)$ and $(4,4k+3)$.  As before, $m=k$ and $m_0=m_2=m_{02}=0$, while $r_1$ is a product of $2k$ transpositions, so $r_0$, $r_1$, $r_2$ and $r_0r_2$ all lift to involutions in $G$. The group $H=\langle r_0, r_1, r_2\rangle$ is again transitive, but now
\[r_0r_2r_1=(1,2,6,\ldots,4k-2,4k+1,4k-3,4k-7,\ldots,5)(3,4k+3,4,4k+2)(7,8)\]
has cycle structure $1^{2k-4},2, 4, 2k+1$, with the same cycle of length $2k+1$ as before. Thus $(r_0r_2r_1)^4$ is a cycle of length $2k+1$. Since $n/3<2k+1<n/2$ with $n$ odd, it again follows from Lemma~\ref{longcycle} and Theorem~\ref{J-Jordan} that $H=A_n$. Thus $G\in\G(1)$ if $11\le n\equiv 3$ mod~$(8)$.

We can repeat this trick in the cases $21\le n\equiv 5$ mod~$(8)$ and $31\le n \equiv 7$ mod~$(8)$. In each case we add two or four further symbols $4k+4$, $4k+5$ and then $4k+6$, $4k+7$, all fixed by $r_0$ and $r_2$, and replace further transpositions $(11,12)$, $(15,16)$ and then $(19,20)$, $(23,24)$ in the original $r_1$ with transpositions pairing these symbols with $7$, $8$ and then $11$, $12$. In each case $r_0r_2r_1$ consists of the same cycle of length $2k+1$, together with various
%{\color{blue}[how many? irrelevant?]}
cycles of length $1$, $2$ or $4$, so $(r_0r_2r_1)^4$ is a $(2k+1)$-cycle and the proof proceeds as before. 

Now let $10\le n\equiv 2$ mod~$(8)$, so $n=4k+2$ where $k$ is even. Define
\[r_0=(1)(2,3)(4,5)(6,7)\ldots (4k-2,4k-1)(4k,4k+1)(4k+2),\]
\[r_1=(1,2)(3,4k+2)(5,6)\ldots (4k-3,4k-2)(4k-1,4k)(4k+1),\]
\[r_2=(1)(2,4)(3,5)(6,8)(7,9)\ldots (4k-2,4k)(4k-1,4k+1)(4k+2),\]
as shown in Figure~\ref{nequiv2(8)}.

\begin{figure}[h!]

\begin{center}
\begin{tikzpicture}[scale=0.6, inner sep=0.8mm]

\node (1) at (-10,0)  [shape=circle, draw, fill=black] {};
\node (2) at (-8,0)  [shape=circle, draw, fill=black] {};
\node (3) at (-7,1)  [shape=circle, draw, fill=black] {};
\node (4) at (-7,-1)  [shape=circle, draw, fill=black] {};
\node (5) at (-6,0)  [shape=circle, draw, fill=black] {};
\node (6) at (-4,0)  [shape=circle, draw, fill=black] {};
\node (7) at (-3,1)  [shape=circle, draw, fill=black] {};
\node (8) at (-3,-1)  [shape=circle, draw, fill=black] {};
\node (9) at (-2,0)  [shape=circle, draw, fill=black] {};
\node (10) at (0,0)  [shape=circle, draw, fill=black] {};
\node (11) at (1,1)  [shape=circle, draw, fill=black] {};
\node (12) at (1,-1)  [shape=circle, draw, fill=black] {};
\node (13) at (2,0)  [shape=circle, draw, fill=black] {};
\node (4k-2) at (8,0)  [shape=circle, draw, fill=black] {};
\node (4k-1) at (9,1)  [shape=circle, draw, fill=black] {};
\node (4k) at (9,-1)  [shape=circle, draw, fill=black] {};
\node (4k+1) at (10,0)  [shape=circle, draw, fill=black] {};
\node (4k+2) at (-7,2.5)  [shape=circle, draw, fill=black] {};
%\node (4k+3) at (-7,-2.5)  [shape=circle, draw, fill=black] {};

\draw (2) to (3);
\draw (4) to (5);
\draw (6) to (7);
\draw (8) to (9);
\draw (10) to (11);
\draw (12) to (13);
\draw (4k-2) to (4k-1);
\draw (4k) to (4k+1);

\draw [thick, dotted] (2) to (4);
\draw [thick, dotted] (3) to (5);
\draw [thick, dotted] (6) to (8);
\draw [thick, dotted] (7) to (9);
\draw [thick, dotted] (10) to (12);
\draw [thick, dotted] (11) to (13);
\draw [thick, dotted] (4k-2) to (4k);
\draw [thick, dotted] (4k-1) to (4k+1);

\draw [thick, dashed] (1) to (2);
\draw [thick, dashed] (3) to (4k+2);
\draw [thick, dashed] (5) to (6);
\draw [thick, dashed] (7) to (8);
\draw [thick, dashed] (9) to (10);
\draw [thick, dashed] (11) to (12);
\draw [thick, dashed] (13) to (3,0);
\draw [thick, dashed] (7,0) to (4k-2);
\draw [thick, dashed] (4k-1) to (4k);

\draw [thick, dotted] (4.5,0) to (5.5,0);

\node at (-10,-0.7) {$1$};
\node at (-8,0.7) {$2$};
\node at (-6.5,1.7) {$3$};
\node at (-7,-1.7) {$4$};
\node at (-6,-0.7) {$5$};
\node at (-4,0.7) {$6$};
\node at (-3,1.7) {$7$};
\node at (-3,-1.7) {$8$};
\node at (-2,-0.7) {$9$};
\node at (0,0.7) {$10$};
\node at (1,1.7) {$11$};
\node at (1,-1.7) {$12$};
\node at (2,-0.7) {$13$};
\node at (7,0.7) {$4k-2$};
\node at (9,1.7) {$4k-1$};
\node at (9,-1.7) {$4k$};
\node at (11,-0.7) {$4k+1$};
\node at (-5.7,3) {$4k+2$};
%\node at (-5.7,-3) {$4k+3$};
 
 \end{tikzpicture}

\end{center}
\caption{The permutations $r_i$ for $n\equiv 2$ mod~$(8)$} 
\label{nequiv2(8)}
\end{figure}

Here we have modified the original permutations $r_i$ for $n\equiv 1$ mod~$(8)$ by adding an extra symbol $4k+2$, fixed by $r_0$ and $r_2$, and replacing the transposition $(3,4)$ in $r_1$ with $(3,4k+2)$. As before, $r_0$, $r_1$, $r_2$ and $r_{02}$ all lift to involutions in $G$, and the group $H=\langle r_0, r_1, r_2\rangle$ is transitive. In this case, $r_0r_2r_1$ has cycle structure $1^{2k-2},3, 2k+1$, where the cycle of length $2k+1$ is as before, and the $3$-cycle is $(3, 4, 4k+2)$. If $k\not\equiv 1$ mod~$(3)$ then $2k+1$ is coprime to $3$ and so $(r_0r_2r_1)^3$ is a $(2k+1)$-cycle. Since $2k+1=n/2$ it follows from Lemma~\ref{longcycle} that either $H$ is primitive, in which case the proof proceeds as before, or the support $\Phi=\{1,2,5,6,9,\ldots, 4k-2, 4k+1\}$ of this cycle and its complement are unions of blocks of size $d>1$ dividing $2k+1$. However, the latter is impossible, since $r_0$ fixes $1\in\Phi$ but moves all other points in $\Phi$ into $\overline\Phi$.

Now suppose that $k\equiv 1$ mod~$(3)$, so that $3$ divides $2k+1$ and we cannot replace $r_0r_2r_1$ with  a $(2k+1)$-cycle simply by taking a suitable power. Instead we can extend the $3$-cycle to a cycle of length $l$ coprime to $2k+1$, without changing the $(2k+1)$-cycle, by replacing a sufficient number of the transpositions $(7,8)$, $(11,12),\ldots$ in $r_1$ with $(4,8)$, $(7,11)$ and so on. Making $t$ such replacements yields a cycle $(3,8,11,\ldots, 7, 4, 4k+2)$ of length $l=3+2t$ (see Figure~\ref{nequiv2(8)kequiv1(3)} for the case $t=2$), so for example taking $t=k-2$ (note that $k\ge 4$ since $k\equiv 4$ mod~$(6)$) gives $l=2k-1$ coprime to $2k+1$. Then $(r_0r_2r_1)^l$ is the required cycle of length $2k+1$, and the proof is as before. Thus $G\in\G(1)$ if $10\le n\equiv 2$ mod~$(8)$. 

\begin{figure}[h!]

\begin{center}
\begin{tikzpicture}[scale=0.6, inner sep=0.8mm]

\node (1) at (-10,0)  [shape=circle, draw, fill=black] {};
\node (2) at (-8,0)  [shape=circle, draw, fill=black] {};
\node (3) at (-7,1)  [shape=circle, draw, fill=black] {};
\node (4) at (-7,-1)  [shape=circle, draw, fill=black] {};
\node (5) at (-6,0)  [shape=circle, draw, fill=black] {};
\node (6) at (-4,0)  [shape=circle, draw, fill=black] {};
\node (7) at (-3,1)  [shape=circle, draw, fill=black] {};
\node (8) at (-3,-1)  [shape=circle, draw, fill=black] {};
\node (9) at (-2,0)  [shape=circle, draw, fill=black] {};
\node (10) at (0,0)  [shape=circle, draw, fill=black] {};
\node (11) at (1,1)  [shape=circle, draw, fill=black] {};
\node (12) at (1,-1)  [shape=circle, draw, fill=black] {};
\node (13) at (2,0)  [shape=circle, draw, fill=black] {};
\node (4k-2) at (8,0)  [shape=circle, draw, fill=black] {};
\node (4k-1) at (9,1)  [shape=circle, draw, fill=black] {};
\node (4k) at (9,-1)  [shape=circle, draw, fill=black] {};
\node (4k+1) at (10,0)  [shape=circle, draw, fill=black] {};
\node (4k+2) at (-7,2.5)  [shape=circle, draw, fill=black] {};
%\node (4k+3) at (-7,-2.5)  [shape=circle, draw, fill=black] {};

\draw (2) to (3);
\draw (4) to (5);
\draw (6) to (7);
\draw (8) to (9);
\draw (10) to (11);
\draw (12) to (13);
\draw (4k-2) to (4k-1);
\draw (4k) to (4k+1);

\draw [thick, dotted] (2) to (4);
\draw [thick, dotted] (3) to (5);
\draw [thick, dotted] (6) to (8);
\draw [thick, dotted] (7) to (9);
\draw [thick, dotted] (10) to (12);
\draw [thick, dotted] (11) to (13);
\draw [thick, dotted] (4k-2) to (4k);
\draw [thick, dotted] (4k-1) to (4k+1);

\draw [thick, dashed] (1) to (2);
\draw [thick, dashed] (3) to (4k+2);
\draw [thick, dashed] (5) to (6);
\draw [thick, dashed] (4) to (8);
\draw [thick, dashed] (9) to (10);
\draw [thick, dashed] (11) to (7);
\draw [thick, dashed] (13) to (3,0);
\draw [thick, dashed] (7,0) to (4k-2);
\draw [thick, dashed] (4k-1) to (4k);

\draw [thick, dotted] (4.5,0) to (5.5,0);

\node at (-10,-0.7) {$1$};
\node at (-8,0.7) {$2$};
\node at (-6.5,1.7) {$3$};
\node at (-7,-1.7) {$4$};
\node at (-6,-0.6) {$5$};
\node at (-4,0.7) {$6$};
\node at (-3,1.7) {$7$};
\node at (-3,-1.7) {$8$};
\node at (-2,-0.7) {$9$};
\node at (0,0.6) {$10$};
\node at (1,1.7) {$11$};
\node at (1,-1.7) {$12$};
\node at (2,-0.7) {$13$};
\node at (7,0.7) {$4k-2$};
\node at (9,1.7) {$4k-1$};
\node at (9,-1.7) {$4k$};
\node at (11,-0.7) {$4k+1$};
\node at (-5.7,3) {$4k+2$};
%\node at (-5.7,-3) {$4k+3$};
 
 \end{tikzpicture}

\end{center}
\caption{The permutations $r_i$ for $n\equiv 2$ mod~$(8)$, $k\equiv 1$ mod~$(3)$} 
\label{nequiv2(8)kequiv1(3)}
\end{figure}

Just as when $n$ is odd, we can deal with the cases $n\equiv 4, 6$ and $0$ mod~$(8)$ by adjoining two, four or six further symbols $4k+3,\ldots, 4k+8$, fixed by $r_0$ and $r_1$, and redefining certain transpositions in $r_1$ to give a transitive group $H$. The resulting graphs are shown in Figures~\ref{nequiv4(8)}, \ref{nequiv6(8)} and \ref{nequiv0(8)}. In the last two cases, the right-hand part of the diagram, involving the vertices $4k-2,\ldots, 4k+1$, is omitted to save space: it is the same as in earlier cases. Similarly, the labels are omitted in Figure~\ref{nequiv0(8)}; the pattern follows that in Figure~\ref{nequiv6(8)}.

In these three cases, the element $r_0r_2r_1$ has one cycle each of lengths $3$ and $2k+1$, together with $2k-6$, $2k-10$ or $2k-14$ fixed points, and one, two or three cycles of lengths $2$ and $4$. If $2k+1$ is coprime to $3$ then $(r_0r_2r_1)^{12}$ is the required $(2k+1)$-cycle, while if $3$ divides $2k+1$ we can extend the $3$-cycle to one of length $l=3+2t$ coprime to $2k+1$ as before, and use $(r_0r_2r_1)^{4l}$. This shows that $G\in\G(1)$ if $20\le n\equiv 4$ mod~$(8)$ or $30\le n\equiv 6$ mod~$(8)$ or $40\le n\equiv 0$ mod~$(8)$. \hfill$\square$

\begin{figure}[h!]

\begin{center}
\begin{tikzpicture}[scale=0.5, inner sep=0.7mm]

\node (1) at (-10,0)  [shape=circle, draw, fill=black] {};
\node (2) at (-8,0)  [shape=circle, draw, fill=black] {};
\node (3) at (-7,1)  [shape=circle, draw, fill=black] {};
\node (4) at (-7,-1)  [shape=circle, draw, fill=black] {};
\node (5) at (-6,0)  [shape=circle, draw, fill=black] {};
\node (6) at (-4,0)  [shape=circle, draw, fill=black] {};
\node (7) at (-3,1)  [shape=circle, draw, fill=black] {};
\node (8) at (-3,-1)  [shape=circle, draw, fill=black] {};
\node (9) at (-2,0)  [shape=circle, draw, fill=black] {};
\node (10) at (0,0)  [shape=circle, draw, fill=black] {};
\node (11) at (1,1)  [shape=circle, draw, fill=black] {};
\node (12) at (1,-1)  [shape=circle, draw, fill=black] {};
\node (13) at (2,0)  [shape=circle, draw, fill=black] {};
\node (14) at (4,0)  [shape=circle, draw, fill=black] {};
\node (15) at (5,1)  [shape=circle, draw, fill=black] {};
\node (16) at (5,-1)  [shape=circle, draw, fill=black] {};
\node (17) at (6,0)  [shape=circle, draw, fill=black] {};
\node (4k-2) at (12,0)  [shape=circle, draw, fill=black] {};
\node (4k-1) at (13,1)  [shape=circle, draw, fill=black] {};
\node (4k) at (13,-1)  [shape=circle, draw, fill=black] {};
\node (4k+1) at (14,0)  [shape=circle, draw, fill=black] {};
\node (4k+2) at (-7,2.5)  [shape=circle, draw, fill=black] {};
\node (4k+3) at (-7,-2.5)  [shape=circle, draw, fill=black] {};
\node (4k+4) at (1,2.5)  [shape=circle, draw, fill=black] {};

\draw (2) to (3);
\draw (4) to (5);
\draw (6) to (7);
\draw (8) to (9);
\draw (10) to (11);
\draw (12) to (13);
\draw (14) to (15);
\draw (16) to (17);
\draw (4k-2) to (4k-1);
\draw (4k) to (4k+1);

\draw [thick, dotted] (2) to (4);
\draw [thick, dotted] (3) to (5);
\draw [thick, dotted] (6) to (8);
\draw [thick, dotted] (7) to (9);
\draw [thick, dotted] (10) to (12);
\draw [thick, dotted] (11) to (13);
\draw [thick, dotted] (14) to (16);
\draw [thick, dotted] (15) to (17);
\draw [thick, dotted] (4k-2) to (4k);
\draw [thick, dotted] (4k-1) to (4k+1);

\draw [thick, dashed] (1) to (2);
\draw [thick, dashed] (3) to (4k+2);
\draw [thick, dashed] (4) to (4k+3);
\draw [thick, dashed] (5) to (6);
%\draw [thick, dashed] (7) to (8);
\draw [thick, dashed] (9) to (10);
\draw [thick, dashed] (11) to (4k+4);
\draw [thick, dashed] (13) to (14);
\draw [thick, dashed] (15) to (16);
\draw [thick, dashed] (17) to (7,0);
\draw [thick, dashed] (11,0) to (4k-2);
\draw [thick, dashed] (4k-1) to (4k);

\draw [thick, dotted] (8.5,0) to (9.5,0);

\node at (-10,-0.7) {$1$};
\node at (-8,0.7) {$2$};
\node at (-6.5,1.7) {$3$};
\node at (-6.5,-1.6) {$4$};
\node at (-6,-0.7) {$5$};
\node at (-4,0.7) {$6$};
\node at (-3,1.7) {$7$};
\node at (-3,-1.7) {$8$};
\node at (-2,-0.7) {$9$};
\node at (0,0.7) {$10$};
\node at (1.5,1.7) {$11$};
\node at (1,-1.7) {$12$};
\node at (2,-0.7) {$13$};
\node at (4,0.7) {$14$};
\node at (5,1.7) {$15$};
\node at (5,-1.7) {$16$};
\node at (6,-0.7) {$17$};
\node at (11,0.7) {$4k-2$};
\node at (13,1.7) {$4k-1$};
\node at (13,-1.7) {$4k$};
\node at (15,-0.7) {$4k+1$};
\node at (-5.7,3) {$4k+2$};
\node at (-5.7,-3) {$4k+3$};
\node at (2.5,3) {$4k+4$};
 
 \end{tikzpicture}

\end{center}
\caption{The permutations $r_i$ for $n\equiv 4$ mod~$(8)$} 
\label{nequiv4(8)}
\end{figure}

%%%%%%%

\begin{figure}[h!]

\begin{center}
\begin{tikzpicture}[scale=0.5, inner sep=0.7mm]

\node (1) at (-10,0)  [shape=circle, draw, fill=black] {};
\node (2) at (-8,0)  [shape=circle, draw, fill=black] {};
\node (3) at (-7,1)  [shape=circle, draw, fill=black] {};
\node (4) at (-7,-1)  [shape=circle, draw, fill=black] {};
\node (5) at (-6,0)  [shape=circle, draw, fill=black] {};
\node (6) at (-4,0)  [shape=circle, draw, fill=black] {};
\node (7) at (-3,1)  [shape=circle, draw, fill=black] {};
\node (8) at (-3,-1)  [shape=circle, draw, fill=black] {};
\node (9) at (-2,0)  [shape=circle, draw, fill=black] {};
\node (10) at (0,0)  [shape=circle, draw, fill=black] {};
\node (11) at (1,1)  [shape=circle, draw, fill=black] {};
\node (12) at (1,-1)  [shape=circle, draw, fill=black] {};
\node (13) at (2,0)  [shape=circle, draw, fill=black] {};
\node (14) at (4,0)  [shape=circle, draw, fill=black] {};
\node (15) at (5,1)  [shape=circle, draw, fill=black] {};
\node (16) at (5,-1)  [shape=circle, draw, fill=black] {};
\node (17) at (6,0)  [shape=circle, draw, fill=black] {};
\node (18) at (8,0)  [shape=circle, draw, fill=black] {};
\node (19) at (9,1)  [shape=circle, draw, fill=black] {};
\node (20) at (9,-1)  [shape=circle, draw, fill=black] {};
\node (21) at (10,0)  [shape=circle, draw, fill=black] {};
\node (22) at (12,0)  [shape=circle, draw, fill=black] {};
\node (23) at (13,1)  [shape=circle, draw, fill=black] {};
\node (24) at (13,-1)  [shape=circle, draw, fill=black] {};
\node (25) at (14,0)  [shape=circle, draw, fill=black] {};
%\node (4k-2) at (12,0)  [shape=circle, draw, fill=black] {};
%\node (4k-1) at (13,1)  [shape=circle, draw, fill=black] {};
%\node (4k) at (13,-1)  [shape=circle, draw, fill=black] {};
%\node (4k+1) at (14,0)  [shape=circle, draw, fill=black] {};
\node (4k+2) at (-7,2.5)  [shape=circle, draw, fill=black] {};
\node (4k+3) at (-7,-2.5)  [shape=circle, draw, fill=black] {};
\node (4k+4) at (-3,2.5)  [shape=circle, draw, fill=black] {};
\node (4k+5) at (-3,-2.5)  [shape=circle, draw, fill=black] {};
\node (4k+6) at (9,2.5)  [shape=circle, draw, fill=black] {};

\draw (2) to (3);
\draw (4) to (5);
\draw (6) to (7);
\draw (8) to (9);
\draw (10) to (11);
\draw (12) to (13);
\draw (14) to (15);
\draw (16) to (17);
\draw (18) to (19);
\draw (20) to (21);
\draw (22) to (23);
\draw (24) to (25);
%\draw (4k-2) to (4k-1);
%\draw (4k) to (4k+1);

\draw [thick, dotted] (2) to (4);
\draw [thick, dotted] (3) to (5);
\draw [thick, dotted] (6) to (8);
\draw [thick, dotted] (7) to (9);
\draw [thick, dotted] (10) to (12);
\draw [thick, dotted] (11) to (13);
\draw [thick, dotted] (14) to (16);
\draw [thick, dotted] (15) to (17);
\draw [thick, dotted] (18) to (20);
\draw [thick, dotted] (19) to (21);
\draw [thick, dotted] (22) to (24);
\draw [thick, dotted] (23) to (25);
%\draw [thick, dotted] (4k-2) to (4k);
%\draw [thick, dotted] (4k-1) to (4k+1);

\draw [thick, dashed] (1) to (2);
\draw [thick, dashed] (3) to (4k+2);
\draw [thick, dashed] (4) to (4k+3);
\draw [thick, dashed] (5) to (6);
\draw [thick, dashed] (8) to (4k+5);
\draw [thick, dashed] (9) to (10);
\draw [thick, dashed] (7) to (4k+4);
\draw [thick, dashed] (13) to (14);
%\draw [thick, dashed] (15) to (16);
\draw [thick, dashed] (17) to (18);
\draw [thick, dashed] (19) to (4k+6);
\draw [thick, dashed] (21) to (22);
\draw [thick, dashed] (23) to (24);
\draw [thick, dashed] (25) to (15.5,0);
%\draw [thick, dashed] (11,0) to (4k-2);
%\draw [thick, dashed] (4k-1) to (4k);

\draw [thick, dotted] (16.5,0) to (17.5,0);

\node at (-10,-0.7) {$1$};
\node at (-8,0.7) {$2$};
\node at (-6.5,1.7) {$3$};
\node at (-6.5,-1.6) {$4$};
\node at (-6,-0.7) {$5$};
\node at (-4,0.7) {$6$};
\node at (-2.5,1.7) {$7$};
\node at (-2.5,-1.7) {$8$};
\node at (-2,-0.7) {$9$};
\node at (-0.1,0.7) {$10$};
\node at (1,1.7) {$11$};
\node at (1,-1.7) {$12$};
\node at (2.1,-0.7) {$13$};
\node at (3.9,0.7) {$14$};
\node at (5,1.7) {$15$};
\node at (5,-1.7) {$16$};
\node at (6.1,-0.7) {$17$};
\node at (7.9,0.7) {$18$};
\node at (9.5,1.7) {$19$};
\node at (9,-1.7) {$20$};
\node at (10.1,-0.7) {$21$};
\node at (11.9,0.7) {$22$};
\node at (13,1.7) {$23$};
\node at (13,-1.7) {$24$};
\node at (14.1,-0.7) {$25$};
%\node at (11,0.7) {$4k-2$};
%\node at (13,1.7) {$4k-1$};
%\node at (13,-1.7) {$4k$};
%\node at (15,-0.7) {$4k+1$};
\node at (-5.7,3) {$4k+2$};
\node at (-5.7,-3) {$4k+3$};
\node at (-1.7,3) {$4k+4$};
\node at (-1.7,-3) {$4k+5$};
\node at (10.3,3) {$4k+6$};
 
 \end{tikzpicture}

\end{center}
\caption{The permutations $r_i$ for $n\equiv 6$ mod~$(8)$} 
\label{nequiv6(8)}
\end{figure}

%%%%%%

\begin{figure}[h!]

\begin{center}
\begin{tikzpicture}[scale=0.37, inner sep=0.5mm]

\node (1) at (-10,0)  [shape=circle, draw, fill=black] {};
\node (2) at (-8,0)  [shape=circle, draw, fill=black] {};
\node (3) at (-7,1)  [shape=circle, draw, fill=black] {};
\node (4) at (-7,-1)  [shape=circle, draw, fill=black] {};
\node (5) at (-6,0)  [shape=circle, draw, fill=black] {};
\node (6) at (-4,0)  [shape=circle, draw, fill=black] {};
\node (7) at (-3,1)  [shape=circle, draw, fill=black] {};
\node (8) at (-3,-1)  [shape=circle, draw, fill=black] {};
\node (9) at (-2,0)  [shape=circle, draw, fill=black] {};
\node (10) at (0,0)  [shape=circle, draw, fill=black] {};
\node (11) at (1,1)  [shape=circle, draw, fill=black] {};
\node (12) at (1,-1)  [shape=circle, draw, fill=black] {};
\node (13) at (2,0)  [shape=circle, draw, fill=black] {};
\node (14) at (4,0)  [shape=circle, draw, fill=black] {};
\node (15) at (5,1)  [shape=circle, draw, fill=black] {};
\node (16) at (5,-1)  [shape=circle, draw, fill=black] {};
\node (17) at (6,0)  [shape=circle, draw, fill=black] {};
\node (18) at (8,0)  [shape=circle, draw, fill=black] {};
\node (19) at (9,1)  [shape=circle, draw, fill=black] {};
\node (20) at (9,-1)  [shape=circle, draw, fill=black] {};
\node (21) at (10,0)  [shape=circle, draw, fill=black] {};
\node (22) at (12,0)  [shape=circle, draw, fill=black] {};
\node (23) at (13,1)  [shape=circle, draw, fill=black] {};
\node (24) at (13,-1)  [shape=circle, draw, fill=black] {};
\node (25) at (14,0)  [shape=circle, draw, fill=black] {};
\node (26) at (16,0)  [shape=circle, draw, fill=black] {};
\node (27) at (17,1)  [shape=circle, draw, fill=black] {};
\node (28) at (17,-1)  [shape=circle, draw, fill=black] {};
\node (29) at (18,0)  [shape=circle, draw, fill=black] {};
\node (30) at (20,0)  [shape=circle, draw, fill=black] {};
\node (31) at (21,1)  [shape=circle, draw, fill=black] {};
\node (32) at (21,-1)  [shape=circle, draw, fill=black] {};
\node (33) at (22,0)  [shape=circle, draw, fill=black] {};
%\node (4k-2) at (12,0)  [shape=circle, draw, fill=black] {};
%\node (4k-1) at (13,1)  [shape=circle, draw, fill=black] {};
%\node (4k) at (13,-1)  [shape=circle, draw, fill=black] {};
%\node (4k+1) at (14,0)  [shape=circle, draw, fill=black] {};
\node (4k+2) at (-7,2.5)  [shape=circle, draw, fill=black] {};
\node (4k+3) at (-7,-2.5)  [shape=circle, draw, fill=black] {};
\node (4k+4) at (-3,2.5)  [shape=circle, draw, fill=black] {};
\node (4k+5) at (-3,-2.5)  [shape=circle, draw, fill=black] {};
\node (4k+6) at (1,2.5)  [shape=circle, draw, fill=black] {};
\node (4k+7) at (1,-2.5)  [shape=circle, draw, fill=black] {};
\node (4k+8) at (17,2.5)  [shape=circle, draw, fill=black] {};

\draw (2) to (3);
\draw (4) to (5);
\draw (6) to (7);
\draw (8) to (9);
\draw (10) to (11);
\draw (12) to (13);
\draw (14) to (15);
\draw (16) to (17);
\draw (18) to (19);
\draw (20) to (21);
\draw (22) to (23);
\draw (24) to (25);
\draw (26) to (27);
\draw (28) to (29);
\draw (30) to (31);
\draw (32) to (33);
%\draw (4k-2) to (4k-1);
%\draw (4k) to (4k+1);

\draw [thick, dotted] (2) to (4);
\draw [thick, dotted] (3) to (5);
\draw [thick, dotted] (6) to (8);
\draw [thick, dotted] (7) to (9);
\draw [thick, dotted] (10) to (12);
\draw [thick, dotted] (11) to (13);
\draw [thick, dotted] (14) to (16);
\draw [thick, dotted] (15) to (17);
\draw [thick, dotted] (18) to (20);
\draw [thick, dotted] (19) to (21);
\draw [thick, dotted] (22) to (24);
\draw [thick, dotted] (23) to (25);
\draw [thick, dotted] (26) to (28);
\draw [thick, dotted] (27) to (29);
\draw [thick, dotted] (30) to (32);
\draw [thick, dotted] (31) to (33);
%\draw [thick, dotted] (4k-2) to (4k);
%\draw [thick, dotted] (4k-1) to (4k+1);

\draw [thick, dashed] (1) to (2);
\draw [thick, dashed] (3) to (4k+2);
\draw [thick, dashed] (4) to (4k+3);
\draw [thick, dashed] (5) to (6);
\draw [thick, dashed] (8) to (4k+5);
\draw [thick, dashed] (9) to (10);
\draw [thick, dashed] (7) to (4k+4);
\draw [thick, dashed] (11) to (4k+6);
\draw [thick, dashed] (12) to (4k+7);
\draw [thick, dashed] (13) to (14);
%\draw [thick, dashed] (15) to (16);
\draw [thick, dashed] (17) to (18);
\draw [thick, dashed] (21) to (22);
%\draw [thick, dashed] (23) to (24);
\draw [thick, dashed] (25) to (26);
\draw [thick, dashed] (27) to (4k+8);
\draw [thick, dashed] (29) to (30);
\draw [thick, dashed] (31) to (32);
\draw [thick, dashed] (33) to (23.7,0);
%\draw [thick, dashed] (11,0) to (4k-2);
%\draw [thick, dashed] (4k-1) to (4k);

\draw [thick, dotted] (25,0) to (26,0);

\iffalse
\node at (-10,-0.7) {$1$};
\node at (-8,0.7) {$2$};
\node at (-6.5,1.7) {$3$};
\node at (-6.5,-1.6) {$4$};
\node at (-6,-0.7) {$5$};
\node at (-4,0.7) {$6$};
\node at (-2.5,1.7) {$7$};
\node at (-2.5,-1.7) {$8$};
\node at (-2,-0.7) {$9$};
\node at (-0.1,0.7) {$10$};
\node at (1,1.7) {$11$};
\node at (1,-1.7) {$12$};
\node at (2.1,-0.7) {$13$};
\node at (3.9,0.7) {$14$};
\node at (5,1.7) {$15$};
\node at (5,-1.7) {$16$};
\node at (6.1,-0.7) {$17$};
\node at (7.9,0.7) {$18$};
\node at (9.5,1.7) {$19$};
\node at (9,-1.7) {$20$};
\node at (10.1,-0.7) {$21$};
\node at (11.9,0.7) {$22$};
\node at (13,1.7) {$23$};
\node at (13,-1.7) {$24$};
\node at (14.1,-0.7) {$25$};
%\node at (11,0.7) {$4k-2$};
%\node at (13,1.7) {$4k-1$};
%\node at (13,-1.7) {$4k$};
%\node at (15,-0.7) {$4k+1$};
\node at (-5.7,3) {$4k+2$};
\node at (-5.7,-3) {$4k+3$};
\node at (-1.7,3) {$4k+4$};
\node at (-1.7,-3) {$4k+5$};
\node at (10.3,3) {$4k+6$};
\fi
 
 \end{tikzpicture}

\end{center}
\caption{The permutations $r_i$ for $n\equiv 0$ mod~$(8)$} 
\label{nequiv0(8)}
\end{figure}

\medskip
(In fact this proof also shows that $2.A_n\in{\G}(1)$ for $n=9, 10, 11$, for $n= 17, \ldots, 21$, and for $n=25, \dots, 31$. No doubt, with sufficient patience the remaining values of $n$ can be dealt with on an individual basis.)

\section{Realising $L_2(q)$}\label{simple}

In this section we complete the proof of Theorem~\ref{mainthm} by considering the groups $L_2(q)$. We need the following fact, a consequence of Macbeath's work~\cite{Macb} on generators of $L_2(q)$, which was observed by Singerman in proving a result~\cite[Theorem~3]{Sin} concerning symmetries of Riemann surfaces:

\begin{prop}\label{L2qinvgen}
If two elements generate $L_2(q)$, there is an automorphism of $L_2(q)$ inverting them both. \hfill$\square$
\end{prop}

The case $T=1$ of the following theorem is already known as a special case of more extensive results of Nuzhin~\cite{Nuz90, Nuz97} about finite simple groups (see \S\ref{gensimple}), but for completeness we give a simple proof here:  

\begin{thm}\label{L2q}
The group $L_2(q)$ is the automorphism group of a map in an edge-transitive class $T$ if and only if one of the following holds:
\begin{itemize}
\item $T=1$, and $q\ne 3$, $7$ or $9$;
\item $T=2^{\sigma}$ for some $\sigma$ or $T=3$, and $q\ne 3$;
\item $T=4^{\sigma}$ for some $\sigma$. 
\end{itemize}
In particular, this situation does not arise if $T=2^{\sigma}{\rm ex}$ or $5^{\sigma}$ for any $\sigma$.
\end{thm}

\noindent{\sl Proof.} Since $L_2(2)\cong S_3$, $L_2(3)\cong A_4$, $L_2(4)\cong L_2(5)\cong A_5$ and $L_2(9)\cong A_6$, the cases $q\le 5$ and $q=9$ are covered by Theorems~\ref{symgps} and \ref{altgps}, so we may assume from now on that $q=7$ or $q\ge 11$. 

The types $T=2^{\sigma}{\rm ex}$ and $5^{\sigma}$ are eliminated by Proposition~\ref{L2qinvgen}, which shows that there are always forbidden automorphisms in these cases.

Now let $T=1$. We need to find involutions $r_0$, $r_1$ and $r_2$ generating $G:=L_2(q)$, with $(r_0r_2)^2=1$, so that $G$ is a quotient of $\Gamma$. Without loss of generality we can take
\[r_1=\pm\left(\,\begin{matrix}0&1\cr -1&0\cr \end{matrix}\,\right).\]
We now look for elements $x=r_0r_1$ and $z=r_1r_2$ of $G$, inverted by $r_1$, such that $\langle x, z, r_1\rangle=G$ and $xz=r_0r_2$ is an involution. The elements of $G$ inverted by $r_1$ (other than $r_1$ itself) are those of the form
\[\pm\left(\,\begin{matrix}a&b\cr b&d\cr \end{matrix}\,\right)\quad{\rm with}\quad ad-b^2=1,\]
so we take two such elements
\[x=\pm\left(\,\begin{matrix}a&b\cr b&d\cr \end{matrix}\,\right)
\quad{\rm and}\quad
z=\pm\left(\,\begin{matrix}a'&b'\cr b'&d'\cr \end{matrix}\,\right),\]
where $ad-b^2=1$ and similarly for $z$. Their product has trace $aa'+2bb'+dd'$, so for $xz$ to be an involution we require
\[aa'+2bb'+dd'=0.\]
To define $x$ we can take $b=0$ and $d=a^{-1}\ne\pm 1$, so that $x$ is elliptic, of order dividing $(q-1)/h$ where $h=\gcd(q-1,2)$. In this case we require
\[aa'+dd'=0,\]
that is,
\[d'=-aa'/d=-a^2a',\]
so for $z$ to exist we need
\[1=a'd'-b'^2=-(aa')^2-b'^2.\]
Thus we need $-1-(aa')^2$ to be a square in ${\mathbb F}_q$, so that we can take $b'$ to be a square root. For a given $a\ne 0$, this will be true for all $a'$ if $q$ is even, or for $(q+1)/2$ values of $a'$ if $q$ is odd.

If we take $a$ to be a primitive element of ${\mathbb F}_q$, then $D:=\langle r_0, r_1\rangle=\langle r_1, x\rangle$ is a dihedral group of order $2(q-1)/h$. It follows from Dickson's description of the subgroups of $L_2(q)$ in~\cite[Chapter XII]{Dic} that $D$ is a maximal subgroup of $G$ provided $q=8$ or $q\ge 13$. It consists of the elements of $G$ of the form
\[\pm\left(\,\begin{matrix}*&0\cr 0&*\cr \end{matrix}\,\right)
\quad{\rm or}\quad
\pm\left(\,\begin{matrix}0&*\cr *&0\cr \end{matrix}\,\right),\]
so if $z$ does not have this form then $\langle r_0, r_1, r_2\rangle=G$. For this it is sufficient that $a', b'\ne 0$, that is, $a'\ne 0$ and $(aa')^2\ne -1$. For a given value of $a$, this excludes one, two or three values of $a'$ as $q\equiv 3$, $0$ or $1$ mod~$(4)$, so we can choose a suitable value of $a'$ provided $(q+1)/2>3$, which is true since $q>5$.

This realises $G$ for all $q=2^e\ge 8$ and all odd $q\ge 13$. In the case $q=11$, although $D$ is not maximal we can ensure that $\langle r_0, r_1, r_2\rangle=G$ by choosing $a'$ so that $z$ has order $6$; for instance, one can take $a=a'=2$, so that $d'=3$ and $b'=\pm 2$, giving the map N46.3 of type $\{5,6\}$ in Conder's list of non-orientable regular maps~\cite{Con}. Thus $G\in\G(1)$ when $q=8$ or $q\ge 11$, and hence $G\in\G(T)$ for $T=2^{\sigma}$, $3$ or $4^{\sigma}$, by Lemma~\ref{regmapslemma}(a).

The groups $L_2(7)$ and $L_2(9)\;(\cong A_6)$ are not members of $\G(1)$, since Nuzhin has shown that neither group is a quotient of $\Gamma$: for instance, a simple counting argument shows that in $L_2(7)$, a Klein four-group $\langle r_0, r_2\rangle$ and a cyclic group $\langle r_1\rangle$ of order $2$ are always contained in one of the $14$ maximal subgroups isomorphic to $S_4$; we outlined an equally simple argument for $L_2(9)\;(\cong A_6)$ in the proof of Theorem~\ref{altgps}.

Now let $T=2$. In $L_2(7)$ let us define the involutions
\[s_1=\pm\left(\,\begin{matrix}0&1\cr -1&0\cr \end{matrix}\,\right),
\quad
s_2=\pm\left(\,\begin{matrix}0&2\cr 3&0\cr \end{matrix}\,\right)
\quad{\rm and}\quad
s'=\pm\left(\,\begin{matrix}1&3\cr -3&-1\cr \end{matrix}\,\right).\]
Then $s_1s_2$ and $s_1s'$ have trace $\pm 1$ and hence have order $3$, while $s_2s'$ has trace $\pm 3$ and hence has order $4$. These involutions generate $L_2(7)$, since no maximal subgroup (isomorphic to $S_4$ or $C_7\rtimes C_3)$ contains two elements of order $3$, such as $s_2s_1$ and $s_1s'$, with product of order $4$.  Lemma~\ref{class2lemma}(c) now implies that $L_2(7)\in\G(T)$ for all $T=2^{\sigma}$, $3$ or $4^{\sigma}$. \hfill$\square$

\medskip

This complete the proof of Theorem~\ref{mainthm}.

%%%%%%%%%%

\section{Realising other finite simple groups}\label{gensimple}

More generally, we will now consider which finite simple groups are members of $\G(T)$ for each of the $14$ edge-transitive classes $T$. By Lemma~\ref{abelian} the only cyclic group of prime order arising is $C_2$, for $T=1$, $2^{\sigma}$, $3$ or $4^{\sigma}$, so we will restrict our attention to non-abelian finite simple groups. Their classification (see~\cite{ATLAS, Wil09}) allows us to deal with them by inspection. The alternating groups have already been dealt with in Theorem~\ref{mainthm}, so we can concentrate mainly on the groups of Lie type and the sporadic simple groups.

\medskip

\noindent{\bf Example}
Let $G$ be a Suzuki group $Sz(2^e)$ or a `small' Ree group $R(3^e)$, for some odd $e>1$.  
It is known from work of Nuzhin~\cite{Nuz90, Nuz97} that $G$ is a quotient of $\Gamma$,  so $G\in\G(T)$ if $T=1$, $2^{\sigma}$, $3$ or $4^{\sigma}$ by Lemma~\ref{regmapslemma}(a). It is also known from results of Suzuki~\cite{Suz} and Ree~\cite{Ree} that their groups are respectively quotients of the triangle groups $\Delta(2, 4, 5)$ and $\Delta(2, 3, 7)$, with no automorphism inverting the second generator, so $G\in\G(T)$ if $T=2^{\sigma}{\rm ex}$ or $5^{\sigma}$ by Lemma~\ref{regmapslemma}(b). Thus $G\in\G(T)$ for all $T$.

\medskip

As a first step towards a more comprehensive approach, we can use results already in the literature to describe those non-abelian finite simple groups which are quotients of each parent group $N(T)$, postponing until later any consideration of forbidden automorphisms. The following statement, depending on the work of many authors, summarises the situation. First let $\mathbb M$ denote the set of finite simple groups isomorphic to
\[L_3(q),\, U_3(q), L\,_4(2^e),\, U_4(2^e),\, U_4(3),\, U_5(2),\, A_6,\, A_7,\, M_{11},\, M_{22},\, M_{23}\;\;{\rm or}\;\; McL\]
for some prime power $q$ or integer $e\ge 1$. Note that $\mathbb M$ includes the groups $L_2(7)\cong L_3(2)$, $L_2(9)\cong A_6$, $A_8\cong L_4(2)$, the symplectic group $S_4(3)\cong U_4(2)$ and the orthogonal group $O_6^-(3)\cong U_4(3)$.

% $U_4(3)$ and $U_5(2)$ have been added to $\mathbb M$, as a result of a computer search by Martin Ma\v caj, confirmed by Matan Ziv-Av also using GAP and Marston Conder using Magma. 

\begin{thm}\label{genfsgps}
Let $T$ be one of the classes of edge-transitive maps, and let $G$ be a non-abelian finite simple group. Then $G$ is a quotient of the parent group $N(T)$ if and only if one of the following conditions holds:
\begin{itemize}
\item $T=1$, and $G\not\in{\mathbb M}$;
\item $T=2^{\sigma}$ for some $\sigma$, and $G\not\cong U_3(3)$;
\item $T=2^{\sigma}{\rm ex}$, $3$, $4^{\sigma}$ or $5^{\sigma}$ for some $\sigma$.
\end{itemize}
\end{thm}

\noindent{\sl Proof.} This theorem is a simple consequence of a number of well-known results, as follows.

When $T=1$ (the set of regular maps), so that $N(T)=\Gamma\cong V_4*C_2$, we are looking for groups generated by three involutions, two of which commute. In 1980 Mazurov asked in the Kourovka Notebook~\cite[Problem 7.30]{Kou} which finite simple groups have this property. This problem was subsequently solved for the alternating groups and the simple groups of Lie type by Nuzhin~\cite{Nuz90, Nuz92, Nuz97a, Nuz97}. In~\cite{Maz} Mazurov has given an elegant, unified and largely computer-free discussion of which sporadic simple groups are quotients of $\Gamma$, together with a summary of how mathematicians such as Ershov and Nevmerzhitskaya, Nuzhin and Timofeenko, Abasheev, and Norton earlier dealt with various individual sporadic groups, mainly by using computers. The groups $U_4(3)$ and $U_5(2)$ have been added to the published lists, as a computer search using GAP by Martin Ma\v caj~\cite{Maca} (confirmed independently by Matan Ziv-Zv also using GAP, and by Marston Conder using Magma) has shown that they are not quotients of $\Gamma$. The result of this work is that the only non-abelian finite simple groups which do not have this property are those in the set $\mathbb M$ defined above. 

When $T=2^{\sigma}$, so that $N(T)\cong C_2*C_2*C_2$, we are looking for groups generated by three involutions, with no requirement that two of them should commute. Here Malle, Saxl and Weigel~\cite{MSW}, building on earlier work of others such as Dalla Volta (see~\cite{DV} for the sporadic groups, for instance), have shown that the only non-abelian finite simple group without such generators is the unitary group $U_3(3)$, shown by Wagner~\cite{Wag} to  require four involutions. Thus with this one exception, every non-abelian finite simple group is a quotient of $N(2^{\sigma})$, whereas they are all quotients of $N(3)\cong C_2*C_2*C_2*C_2$.

When $T=2^{\sigma}{\rm ex}$, so that $N(T)\cong C_{\infty}*C_2$, we are looking for groups generated by two elements, one of them an involution. Here it is known (again, see~\cite{MSW}) that every non-abelian finite simple group has such a generating pair, so it is a quotient of $N(T)$. Clearly the same applies when $T=4^{\sigma}$, so that $N(T)\cong C_2*C_2*C_{\infty}$, and when $T=5^{\sigma}$, so that $N(T)\cong F_2$. \hfill$\square$

\medskip

Theorem~\ref{genfsgps} immediately shows that the non-abelian finite simple groups $G\in\G(1)$ are those not in $\mathbb M$, since for the class $T=1$ there are no forbidden automorphisms to consider. Moreover,  Lemma~\ref{regmapslemma}(a) then gives:

\begin{cor}\label{Mazcor}
If $G$ is a non-abelian finite simple group not in $\mathbb M$, then $G\in\G(T)$ for each edge-transitive class $T=1$, $2^{\sigma}$, $3$ or $4^{\sigma}$. \hfill$\square$
\end{cor}

In order to prove Theorem~\ref{mainthmsimple} we need to determine, for each of the remaining classes class $T\ne 1$, which of the simple quotients $G$ of $N(T)$ have generators admitting no forbidden automorphisms, so that $G\in\G(T)$. 

Leemans and Liebeck~\cite{Lee16} have recently proved the following theorem in the context of chiral polyhedra, which are equivalent to maps in class $2^P{\rm ex}$. First let $\mathbb L$ denote the set of finite simple groups which are isomorphic to
\[L_2(q),\, L_3(q),\, U_3(q)\;\;{\rm or}\;\; A_7\]
for some prime power $q$.

\begin{thm}[Leemans and Liebeck]\label{LeeLie}
A non-abelian finite simple group is generated by elements $x$ and $y$, where $y$ is an involution and no automorphism inverts $x$ and fixes $y$, if and only if it is not a member of $\,\mathbb L$. \hfill$\square$
\end{thm}

The proof that every non-abelian finite simple group $G\not\in{\mathbb L}$ has such a generating pair is given in~\cite{LL}. It is straightforward to show that $A_7$ and $L_2(q)$ do not have such a pair (see Theorems~\ref{altgps} and~\ref{L2q}, for example); the proof for $L_3(q)$ and $U_3(q)$ is (unsurprisingly) much harder, and the authors of~\cite{LL} intend to present this separately in~\cite{LL2}. Breda d'Azevedo and Catalano~\cite{BC19} have given a proof for $L_3(q)$, based on a similar result for $SL_3(q)$.

%{\color{blue}[update?]}
Lemma~\ref{regmapslemma}(b) immediately implies the following:

\begin{cor}\label{LLcor}
If $G$ is a non-abelian finite simple group which is not in $\mathbb L$, then $G\in\G(T)$ for each edge-transitive class $T=2^{\sigma}{\rm ex}$, $4^{\sigma}$ or $5^{\sigma}$.  \hfill$\square$
\end{cor}

Combining this with Corollary~\ref{Mazcor}, we have:

\begin{cor}
If $G$ is a non-abelian finite simple group which is not in ${\mathbb L}\cup{\mathbb M}$, then $G\in\G(T)$ for each edge-transitive class $T$.  \hfill$\square$
\end{cor}

Thus, in order to prove Theorem~\ref{mainthmsimple} we may restrict attention to the groups $G\in{\mathbb L}\cup{\mathbb M}$, in each case determining those classes $T$ such that $G\in\G(T)$. The groups $L_2(q)$ and $A_n$ have been dealt with in Theorem~\ref{mainthm}, as have $L_3(2)\cong L_2(7)$ and $L_4(2)\cong A_8$, so the following groups remain:
\begin{itemize}
\item $L_3(q)$ and $U_3(q)$ for prime powers $q>2$;
\item $L_4(2^e)$ for $e>1$ and $U_4(2^e)$ for $e\ge 1$;
\item $U_4(3)$ and $U_5(2)$; % added
\item $M_{11}$, $M_{22}$, $M_{23}$ and $McL$.
\end{itemize}

% $U_4(3)$ and $U_5(2)$ added to the list.

\medskip
To deal with specific examples of these groups we will use information and notation concerning their conjugacy classes, characters, maximal subgroups, automorphisms, etc in the ATLAS~\cite{ATLAS}, together with the following formula, due to Frobenius~\cite{Fro}:

\begin{prop}\label{frobchi}
Let $\mathcal A$, $\mathcal B$ and $\mathcal C$ be conjugacy classes in a finite group $G$. Then the number of solutions of the equation $abc=1$ in $G$, with $a\in{\mathcal A}$, $b\in{\mathcal B}$ and $c\in{\mathcal C}$, is 
\begin{equation}\label{trianglechi}
\frac{|{\mathcal A}|.|{\mathcal B}|.|{\mathcal C}|}{|G|}\sum_{\chi}
\frac{\chi(a)\chi(b)\chi(c)}{\chi(1)}
\end{equation}
where $\chi$ ranges over the irreducible complex characters of $G$. \hfill$\square$
\end{prop}

Concerning the groups $L_3(q)$ and $U_3(q)$, we have the following corollary to Theorem~\ref{LeeLie}:

\begin{cor}\label{L3U3cor}
If $G=L_3(q)$, or if $G=U_3(q)$ where $q>3$, then $G\in\G(T)$ for each $T=2^{\sigma}$, $3$ or $4^{\sigma}$.
\end{cor}

\noindent{\sl Proof.} Malle, Saxl and Weigel~\cite{MSW} have shown that every non-abelian finite simple group $G\not\cong U_3(3)$ is generated by a strongly real element $x$ (one such that $x^a=x^{-1}$ for some involution $a\in G$) and an involution $y$. Thus $G$ is generated by three involutions $a$, $b:=ax$ and $y$. If an automorphism of $G$ fixes the involution $y$ and transposes $a$ and $b$, then it inverts $x$. However, by Theorem~\ref{LeeLie} the groups $G=L_3(q)$ and $U_3(q)$ have no generating pairs $x, y$ with this property, so by taking $s_1=a$, $s_2=b$ and $s_3=y$ we see that $G\in\G(T)$ for $T=2$, and hence for $T=2^{\sigma}$, $3$ or $4^{\sigma}$ by Lemma~\ref{class2lemma}(c). \hfill$\square$

\medskip

The group $G=U_3(3)$ was excluded from Corollary~\ref{L3U3cor}, so we treat it here as a special case:

\begin{thm}\label{U33}
The group $U_3(3)$ is in $\G(T)$ for each $T=3$, $4^{\sigma}$ or $5^{\sigma}$.
\end{thm}

\noindent{\sl Proof.} First let $T=3$. The group $G=U_3(3)$ has a single conjugacy class of maximal subgroups $H\cong L_2(7)$. Let $s_1$, $s_2$ and $s_3$ be involutions in such a subgroup $H$ corresponding to the elements
\[\pm\left(\,\begin{matrix}-1&1\cr -2&1\cr \end{matrix}\,\right), \quad
\pm\left(\,\begin{matrix}0&1\cr -1&0\cr \end{matrix}\,\right) \quad{\rm and}\quad
\pm\left(\,\begin{matrix}0&3\cr 2&0\cr \end{matrix}\,\right)\]
of $L_2(7)$, so that $s_1s_2$, $s_1s_3$ and $s_2s_3$ have orders $4$, $4$ and $3$ respectively. Since $s:=s_1s_2s_3$ has order $7$, and no proper subgroup of $L_2(7)$ has order divisible by $14$, these involutions generate $H$.

Now $G$ acts as a primitive group of degree $36$ and rank $4$ on the cosets of $H$, with subdegrees $1, 7, 7$ and $21$. The subgroup $\langle s_1, s_2\rangle\cong D_4$ of $H$ thus fixes more than one point, so it lies in a second subgroup $H'\cong L_2(7)$ of $G$. The subgroups of $L_2(7)$ isomorphic to $D_4$ are all conjugate (they are Sylow $2$-subgroups), so as above there is an involution $s_0\in H'$ such that $\langle s_1, s_2, s_0\rangle=H'$ where $s_1s_2$, $s_1s_0$ and $s_2s_0$ have orders $4$, $4$ and $3$. Then $\langle s_0, s_1, s_2, s_3\rangle=\langle H, H'\rangle=G$ and the partitions $02\mid 13$ and $01\mid 23$ satisfy the condition in Lemma~\ref{class3lemma}, so $G\in\G(3)$.

Next, let $T=4$. The involution $s_1$ and the element $s$ of order $7$ defined above generate the maximal subgroup $H$ of $G$. We need to show that there is another involution $s'_2\in G\setminus H$ such that no automorphism of $G$ inverts $s$ and transposes $s_1$ and $s'_2$. There are seven automorphisms inverting $s$, all of them outer automorphisms in the conjugacy class $2B\subset {\rm Aut}\,G=PGU_3(3)$, forming a coset of $C_G(s)=\langle s\rangle$. Now there are $63$ involutions in $G$, of which $21$ are in $H$. We can therefore choose one of the $42$ involution $s'_2\in G\setminus H$, avoiding the seven images of $s_1$ under the automorphisms inverting $s$. Thus $G\in\G(4^{\sigma})$ for all $\sigma$. 

Finally let $T=5$. There are no maximal subgroups of $G$ containing elements $s$ and $s'$ of orders $7$ and $6$, so any such pair generates $G$. Since they have different orders, it is sufficient to show that there exists such a pair which are not simultaneously inverted by any automorphism. Any element $s$ of order $7$ is inverted by seven automorphisms, all outer, in class $2B$.  Each such automorphism inverts $24$ elements $s'\in G$ of order $6$, so there are at most $168$ such elements $s'$ inverted by automorphisms also inverting $s$. Since $G$ has $504$ elements of order $6$, a suitable pair $s, s'$ exists, so $G\in\G(5^{\sigma})$ for all $\sigma$.\hfill$\square$

\medskip

For arguments involving Singer cycles and their eigenvalues, we will need the following result:

\begin{lemma}\label{priminv}
Let $\lambda$ be a primitive root in a finite field $\F$. If $\lambda$ is conjugate to $\lambda^{-1}$ under ${\rm Gal}\,\F$ then $|\F|\le 4$.
\end{lemma}

\noindent{\sl Proof.} Let $|\F|=p^e$ where $p$ is prime, so that  ${\rm Gal}\,\F$ is generated by the Frobenius automorphism $t\mapsto t^p$, which has order $e$. If $\lambda$ is conjugate to its inverse, then $\lambda^{p^f}=\lambda^{-1}$ and hence $p^f+1$ is divisible by the mulitiplicative order $p^e-1$ of $\lambda$, for some $f=0, 1, \ldots, e-1$. Hence $p^e-1\le p^{e-1}+1$, so $p^{e-1}(p-1)\le 2$, and the result follows immediately.  \hfill$\square$

\begin{thm}
If $q>2$ then $L_3(q), U_3(q)\in\G(5^{\sigma})$ for all $\sigma$.
\end{thm}

\noindent{\sl Proof.} Let $G=L_3(q)$, of order $q^3(q^3-1)(q^2-1)/d$ where $d=(3, q-1)$. Then ${\rm Aut}\,G$ is a semidirect product of $P\Gamma L_3(q)$ by a group of order $2$ generated by the graph automorphism (or polarity) $\gamma$ induced by the matrix operation $A\mapsto (A^{-1})^T$. Let $s$ be a Singer cycle in $G$, that is, an element of order $(q^2+q+1)/d$; its centraliser in $G$ is $\langle s\rangle$, while its centraliser in ${\rm Aut}\,G$ is a cyclic group $C$ of order $q^2+q+1$; this is generated by a Singer cycle in $PGL_3(q)$, which contains $G$ with index $d$. It follows by applying results of Bereczky~\cite{Ber} to $SL_3(q)$ that the only maximal subgroup of $G$ containing $s$ is $M:=N_G(\langle s\rangle)$, a semidirect product of $\langle s\rangle$ by $C_3$.

Simple numerical estimates show that if $q>2$ there exist elements $s'\in G$ such that $s'\not\in M$ (so $\langle s, s'\rangle=G$), $s'$ is not a Singer cycle (so no automorphism can transpose it with $s^{\pm 1}$) and none of the $q^2+q+1$ automorphisms inverting $s$ also inverts $s'$. (It follows from Lemma~\ref{priminv} that the automorphisms inverting $s$ are all conjugate in ${\rm Aut}\,G$ to $\gamma$; the elements of $G$ inverted by $\gamma$ are those corresponding to symmetric or skew-symmetric matrices.) Thus $G\in\G(5^{\sigma})$ for all $\sigma$. 

A similar argument can be applied to the groups $G=U_3(q)$; these have order $q^3(q^3+1)(q-1)^2/d$ where $d=(3, q+1)$, and are simple for each prime power $q>2$. In this case ${\rm Aut}\,G=P\Gamma U_3(q)$; this is an extension of $PGU_3(q)$, which contains $G$ with index $d$, by ${\rm Gal}\,\F_{q^2}$. The analogue in $G$ of a Singer cycle in $L_3(q)$ is an element $s$ of order $(q^2-q+1)/d$; its centraliser in $G$ is $\langle s\rangle$, while its centraliser in ${\rm Aut}\,G$ has order $q^2-q+1$. In each case $s$ is contained in a unique maximal subgroup $M$ of $G$: if $q\ne 3$ or $5$ then $M=N_G(\langle s\rangle)$, a semidirect product of $\langle s\rangle$ by $C_3$, whereas if $q=3$ or $5$ then $M\cong L_2(7)$ or $A_7$. There is a unique conjugacy class $\mathcal A$ of involutions $\alpha\in{\rm Aut}\,G$ inverting such elements $s$, namely the class (denoted by $2B$ for the groups in~\cite{ATLAS}) containing the automorphism $\gamma$ induced by the automorphism $t\mapsto t^q$ of $\F_{q^2}$. Each $s$ is inverted by $q^2-q+1$ automorphisms $\alpha\in\mathcal A$. Each $\alpha\in\mathcal A$ inverts an element $s'\in G$ if and only if $s'\alpha\in\mathcal A$, so the number of elements $s'$ inverted by $\alpha$ is $|{\mathcal A}|=q^2(q-1)(q^2-q+1)$. As before, it follows that there exist such pairs $s, s'$ which generate $G$ with no forbidden automorphisms. \hfill$\square$

\medskip

Theorems~\ref{genfsgps} and \ref{LeeLie} show that the groups $L_3(q)$ and $U_3(q)$ are not in $\G(T)$ for $T=1$ or $2^{\sigma}{\rm ex}$, so this completes the proof of Theorem~\ref{mainthmsimple} for these two families of simple groups. Theorems~\ref{genfsgps} and \ref{LeeLie} also show that the groups $L_4(2^e)$ and $U_4(2^e)$ are not in $\G(1)$, but that they are in $\G(T)$ for $T=2^{\sigma}{\rm ex}$, and hence by Lemma~\ref{regmapslemma}(b) also for $T=4^{\sigma}$ and $5^{\sigma}$. To complete the proof of Theorem~\ref{mainthmsimple} for these groups we prove the following:

\begin{thm}
The groups $L_4(2^e)$ and $U_4(2^e)$ are in $\G(2^{\sigma})$ and $\G(3)$ for all $\sigma$ and all $e\ge 1$, as are $U_4(3)$ and $U_5(2)$.
% $U_4(3)$ and $U_5(2)$ added
\end{thm}

\noindent{\sl Proof.} Let $G=L_4(q)$ where $q=2^e$. We may assume that $e>1$, since the group $L_4(2)\cong A_8$ has been dealt with in Theorem~\ref{altgps}. In~\cite[Theorem~2.1]{MSW},  Malle, Saxl and Weigel showed that there are conjugacy classes ${\mathcal C}_1, {\mathcal C}_2\subset G$, containing elements of orders $q^2+1$ and $q^3-1$ respectively, such that the elements of ${\mathcal C}_1$ are strongly real, and for each non-identity conjugacy class ${\mathcal C}_3\subset G$ there are elements $g_i\in{\mathcal C}_i$ ($i=1, 2, 3$) generating $G$ and satisfying $g_1g_2g_3=1$. We can write $g_1=s_2s_3$ for involutions $s_2, s_3\in G$. There are two conjugacy classes of involutions in $G$, of different sizes, depending on the dimensions of the subspaces of the natural module $\F_q^4$ fixed by their elements. The involutions inverting $g_1$, such as $s_2$ and $s_3$, all lie in one of these two classes, so taking ${\mathcal C}_3$ to be the other class of involutions, and putting $s_1=g_3$, we obtain three involutions $s_1, s_2, s_3$ generating $G$, with no automorphism transposing $s_1$ and $s_2$. This shows that $G\in\G(2)$, so Lemma~\ref{class2lemma}(c) implies that $G\in\G(T)$ for $T=2^{\sigma}$ and for $T=3$.

The same argument can be applied to the groups $G=U_4(q)$ for $q=2^e$, using~\cite[Theorem~2.2]{MSW}, the only difference being that the elements $g_2\in{\mathcal C}_2$ now have order $q^3+1$ rather than $q^3-1$. 

%{\color{blue}[Is $q=2$ really an exception here?]}

In the case of $U_5(2)$, classes ${\mathcal C}_1$ and ${\mathcal C}_2$ in~\cite[Theorem~2.2]{MSW} contain elements of orders $(q^5+1)/d(q+1)=11$ and $q^{4/2}+(-1)^{4/2}=5$ (where $d=(5, q+1)=1$), the latter strongly real. There are two conjugacy classes of involutions, and if we choose ${\mathcal C}_3$ to be the class $2A$ of involutions not inverting elements of order $5$, the proof proceeds as before. 

The group $U_4(3)$, which has only one class of involutions, resists this line of argument; however, a computer search by Matan Ziv-Av~\cite{Ziv}, using GAP, has shown that it satisfies the conditions of Lemma~\ref{class2lemma}(a), for example with $ab$, $ac$ and $bc$ having orders $4, 5$ and $6$, so Lemma~\ref{class2lemma}(c) applies. \hfill$\square$ 

\medskip

\noindent{\bf Problem} Find a computer-free proof for $U_4(3)$.
% Extra arguments added for $U_5(2)$ and $U_4(3)$, omitted from Nuzhin's lists.

\medskip

Finally, to complete the proof of Theorem~\ref{mainthmsimple} we deal with the sporadic simple groups in ${\mathbb L}\cup{\mathbb M}$:

\begin{thm}
The Mathieu groups $M_{11}$, $M_{22}$, $M_{23}$ and the McLaughlin group $McL$ are in $\G(T)$ for all edge-transitive classes $T\ne 1$, but they are not in $\G(1)$.
\end{thm}

\noindent{\sl Proof.} Theorem~\ref{genfsgps} shows that these groups are not in $\G(1)$. Theorem~\ref{LeeLie} shows that they are in $\G(2^P{\rm ex})$, so by Lemma~\ref{regmapslemma}(b) they are in $\G(T)$ for all classes $T=2^{\sigma}{\rm ex}$, $4^{\sigma}$ and $5^{\sigma}$.

For $T=2^{\sigma}$ or $3$ we use Lemma~\ref{class2lemma}(b) and (c). The group $G=M_{11}$ has unique conjugacy classes $2A$ and $4A$ of elements of order $2$ and $4$, and two mutually inverse classes $11A$ and $11B$ of elements of order $11$. Proposition~\ref{frobchi} and the character table of $G$ in~\cite{ATLAS} show that there are elements $a, b, c\in G$, in classes $2A$, $4A$ and $11A$, such that $abc=1$. The maximal subgroups of $G$ are listed in~\cite{ATLAS}, and none contains elements of orders $4$ and $11$, so $\langle a, b\rangle=G$. A second application of Proposition~\ref{frobchi}, or equivalently the existence of subgroups isomorphic to $D_4$, shows that there are involutions $s_1, s_2\in G$ with $s_1s_2=b$. Since $G$ has no outer automorphisms, no automorphism can invert $c$, so Lemma~\ref{class2lemma} shows that $G\in\G(T)$ for all classes $T=2^{\sigma}$ and $3$.

Essentially the same argument gives the result for $M_{23}$, but now with $a, b$ and $c$ in classes $2A$, $8A$ and $23A$, and also for $McL$, using classes $2A$, $12A$ and $11A$. (Even though $McL$ has outer automorphisms, they do not transpose the mutually inverse classes $11A$ and $11B$.)

A similar argument can be applied to $G=M_{22}$, but in this case the details are less straightforward. This group has one class $2A$ of involutions, one class $6A$ of elements of order $6$, and two mutually inverse classes $7A$ and $7B$ of elements of order $7$; the latter are not transposed in ${\rm Aut}\,G$, which contains $G$ with index $2$. Proposition~\ref{frobchi} shows that $G$ contains $12|G|$ triples $(a, b, c)$ of elements of orders $2$, $6$ and $7$ with $abc=1$.  By inspection of the list of maximal subgroups of $G$, each of these triples can generate one of the $2|G:H|$ subgroups $H\cong A_7$, or one of the $|G:H|$ subgroups $H\cong AGL_3(2)$, or $G$ itself. The uniqueness of the corresponding permutation diagram for $a$ and $b$ shows that $A_7$ is generated by $|S_7|=2|A_7|$ such triples, giving a total of $4|G|$ triples $(a, b, c)$ generating subgroups $H\cong A_7$. The affine group $AGL_3(2)$ is a semidirect product $V_8\rtimes GL_3(2)$, and triples $(a, b, c)$ generating this group map onto triples $(\overline a, \overline b, \overline c)$ of type $(2, 3, 7)$ generating $GL_3(2)=L_3(2)$. There are $2|L_3(2)|$ such triples, each lifting to $16$ triples $(a, b, c)$, giving a total of $|G:H|.16.2|L_3(2)|=4|G|$ triples $(a, b, c)$ generating subgroups $H\cong AGL_3(2)$. This leaves $4|G|$ triples $(a, b, c)$ generating $G$ (and corresponding to a unique chiral pair of orientably regular maps of type $\{6, 7\}$ with automorphism group $G$). The element $b\in 6A$ in such a triple is strongly real (since $A_7$, and hence $G$, contains subgroups isomorphic to $D_6$), and the outer automorphisms of $G$ do not transposes the classes $7A$ and $7B$ containing $c$ and $c^{-1}$, so Lemma~\ref{class2lemma} applies as before. \hfill$\square$

\section{Realising nilpotent and solvable groups}\label{nilpsolv}

Apart from a few small examples, all the groups we have considered so far as candidates for automorphism groups have been non-solvable. In this section we will consider nilpotent and solvable groups, proving Theorem~\ref{nilpandsolv}.

 First we need some examples of nilpotent groups. For each prime power $n=p^e$ and each integer $f=1, 2,\ldots, e$, let
\begin{equation}\label{Gpef}
G=G_{p,e,f}=\langle g, h\mid g^n=h^n=1, h^g=h^{p^f+1}\rangle,
\end{equation}
a semidirect product of a normal subgroup $\langle h\rangle\cong C_n$ by $\langle g\rangle\cong C_n$. As a finite $p$-group, $G$ is nilpotent. (These groups were studied in connection with embeddings of complete graphs $K_{n,n}$ for $p>2$ in~\cite{JNS}, and for $p=2$ in~\cite{DJKNS}.)

\begin{thm}\label{nilp}
Let $T$ be an edge-transitive class.

\begin{itemize}
\item If $T=1$, $2^{\sigma}$, $3$ or $4^{\sigma}$ for some $\sigma$, then $\G(T)$ contains finite nilpotent groups of class $c$ for each $c\ge 1$.
\item If $T=2^{\sigma}{\rm ex}$ for some $\sigma$, then $\G(T)$ contains finite nilpotent groups of class $c$ if and only if $c\ge 5$.
\item If $T=5^{\sigma}$ for some $\sigma$, then $\G(T)$ contains finite nilpotent groups of class $c$ if and only if $c\ge 2$.
\end{itemize}
\end{thm}

\noindent{\sl Proof.} Lemma~\ref{abelian} deals with abelian automorphism groups, so we may restrict attention to groups of class $c\ge 2$. If $m=2^e$ then the regular map $\{m, 2\}$, embedding a circuit of length $m$ in the sphere, has automorphism group $D_m\times C_2$, which is nilpotent of class $c=e$.  This proves the result for $T=1$, and hence, by Lemma~\ref{regmapslemma}(a), for $T=2^{\sigma}$, $3$ and $4^{\sigma}$.

Next let $T=5$. For any odd prime power $p^e$ put $f=1$ in~(\ref{Gpef}) and define
\[G=G_{p,e,1}=\langle g, h\mid g^{p^e}=h^{p^e}=1, h^g=h^{p+1}\rangle.\]
As shown in~\cite{JNS}, $G$ has centre $Z=\langle g^{p^{e-1}}, h^{p^{e-1}}\rangle\cong C_p\times C_p$, with $G/Z\cong G_{p,e-1,1}$ if $e>1$, so induction on $e$ shows that $G$ has class $c=e$. There is an epimorphism $N(5)=F_2\to G$ given by $S\mapsto s:=g$ and $S'\mapsto s':=h$.  It is straightforward to check that if $e>1$ then $G$ has no forbidden automorphisms, so $G\in\G(5^{\sigma})$ for each $\sigma$.

Finally, let $T=2^P{\rm ex}$. In~\cite{MNS}, Malni\v c, Nedela and \v Skoviera classified the orientably regular maps with nilpotent automorphism groups of class $c=2$, and none of them are chiral. This was extended to the cases $c=3$ and $4$ by Conder, Du, Nedela and \v Skoviera in~\cite{CDNS}, using a computer search.

On the other hand Du, Kwak, Nedela, \v Skoviera and the author showed in~\cite{DJKNS} that if $n=2^e\ge 4$ then for each $f=2, \ldots, e$ there are orientably regular embeddings $\M$ of the complete bipartite graph $K_{n,n}$ such that the subgroup $G:={\rm Aut}^+_0\M$ of ${\rm Aut}\,\M$ preserving orientation and vertex-colours is isomorphic to the group $G_{2,e,f}$ defined in~(\ref{Gpef}); moreover, if $f\le e-2$ (so that $e\ge 4$) these maps are chiral, and thus in class $T$, with
$A:={\rm Aut}\,\M\cong G\rtimes C_2$. These maps include examples where $f=2$,
\[A=\langle g, h, \alpha\mid g^{2^e}=h^{2^e}=\alpha^2=1, h^g=h^5, g^{\alpha}=gh, h^{\alpha}=h^{-1}\rangle,\]
and the epimorphism $N(T)=\Gamma^+\to A$ is given by $X\mapsto x:=g, Y\mapsto y:=\alpha$. As a finite $2$-group, $A$ is nilpotent.  We will compute its nilpotence class $c$ as the length of its upper central series
\[1=Z_0<Z_1<\cdots<Z_c=A,\]
where $Z_i/Z_{i-1}$ is the centre $Z(A/Z_{i-1})$ of $A/Z_{i-1}$ for each $i\ge 1$.

First note that $Z_1:=Z(A)\le Z(G)$, since each element $a\in A\setminus G$ has the form $a=\alpha g^ih^j$, giving $h^a=h^{-5^i}$ with $-5^i\equiv -1$ mod~$(4)$, so that $h^a\ne h$ since $h^2\ne 1$. Thus $Z_1$ is the subgroup of $Z(G)$ commuting with $\alpha$.

Each element of $G$ has the unique form $g^ih^j$ with $i, j\in{\Z}_{2^e}$. Multiplication in $G$ is given by
\[g^ih^j.g^kh^l=g^{i+k}h^{5^kj+l},\]
from which it follows, as shown in~\cite{DJKNS}, that $G$ has centre
\[Z(G)=\langle u:=g^{2^{e-2}}\negthinspace\negthinspace,\, v:=h^{2^{e-2}}\rangle\cong C_4\times C_4.\]

We need to find the action of $\alpha$ on $Z(G)$. Since $h^{\alpha}=h^{-1}$ we have $v^{\alpha}=v^{-1}$. We also need $u^{\alpha}=(g^{2^{e-2}})^{\alpha}=(gh)^{2^{e-2}}$. Induction on $e$ gives
\[(gh)^{2^{e-2}}=g^{2^{e-2}}h^s\quad
{\rm where}\quad
s=5^{2^{e-2}-1}+5^{2^{e-2}-2}+\cdots+5+1.\]

\begin{lemma}
If $e\ge 3$ then
\[5^{2^{e-2}-1}+5^{2^{e-2}-2}+\cdots+5+1\equiv -2^{e-2}\, {\rm mod}\,(2^e).\]
\end{lemma}

\noindent{\sl Proof.} Use induction on $e$, noting that if $s_e$ denotes the left-hand side, then
$s_e=(5^{2^{e-2}}+1)s_{e-1}$ with $5^{2^{e-3}}\equiv 2^{e-1}+1$ mod $(2^e)$. \hfill$\square$

\medskip

Thus $u^{\alpha}=g^{2^{e-2}}h^{-2^{e-2}}=uv^{-1}$.
It follows that $\alpha$ conjugates $u^iv^j$ to $u^iv^{-i-j}$, so it centralises this element if and only if $i\equiv -2j$ mod~$(4)$. Thus
\[Z_1=\langle u^{-2}v\rangle\cong C_4,\]
so the central quotient $\overline A:=A/Z_1$ has a presentation $\Pi_e$ of the form
\[\langle g, h, \alpha
\mid g^{2^e}=\alpha^2=1, h^{2^{e-2}}=g^{2^{e-1}}\negthinspace\negthinspace, h^g=h^5, g^{\alpha}=gh, h^{\alpha}=h^{-1}\rangle.\]

As before, $Z(\overline A)\le Z(\overline G)$ where $\overline G:=G/Z_1$. Each element of $\overline G$ has the unique form $g^ih^j$ where $i\in{\Z}_{2^e}$ and $j=0, 1, \ldots, 2^{e-2}-1$. Multiplication and the action of $\alpha$ are as above for $G$, but now we use $h^{2^{e-2}}=g^{2^{e-1}}$ to reduce elements to standard form. Calculations similar to those for $G$ show that
\[Z(\overline G)=\langle u:=g^{2^{e-3}}, v:=h^{2^{e-3}}\rangle
=\langle u\rangle\times\langle u^{-2}v\rangle \cong C_8\times C_2,\]
where $u$ and $v$ (now redefined) have orders $8$ and $4$, with $u^4=v^2$. As before, we have $v^{\alpha}=v^{-1}$ and $u^{\alpha}=(gh)^{2^{e-3}}=g^{2^{e-3}}h^{-2^{e-3}}=uv^{-1}$, giving
\[Z_2/Z_1=Z(\overline A)=\langle u^4, u^{-2}v\rangle\cong V_4.\]

The central quotient $\overline{\overline A}:=\overline{A}/Z(\overline{A})=A/Z_2$ therefore has a presentation
\[\langle g, h, \alpha
\mid g^{2^{e-1}}=\alpha^2=1, h^{2^{e-3}}=g^{2^{e-2}}\negthinspace\negthinspace, h^g=h^5, g^{\alpha}=gh, h^{\alpha}=h^{-1}\rangle.\]
This is the presentation $\Pi_{e-1}$, obtained from the presentation $\Pi_e$ for $\overline A$ by replacing $e$ with $e-1$. It follows that $Z_3/Z_2=Z(\overline{\overline A})\cong V_4$, and all further quotients $Z_i/Z_{i-1}$ in the upper central series of $A$ are also isomorphic to $V_4$, until we reach a quotient $A/Z_{e-1}$ of $A$ with presentation $\Pi_2$, namely
\[\langle g, h, \alpha
\mid g^4=\alpha^2=1, h=g^2, h^g=h^5, g^{\alpha}=gh, h^{\alpha}=h^{-1}\rangle\cong D_4.\]
This has centre $\langle g^2\rangle\cong C_2$, with central quotient $V_4$. Thus the successive quotients in the upper central series of $A$ are isomorphic to
\[C_4,\, V_4,\, V_4,\, \ldots, V_4\;\,\hbox{($e-2$ times)},\, C_2\;\;\hbox{and}\;\, V_4,\]
so that $A$ has class $c=e+1$.  This shows that the maps $\M$ defined above realise nilpotent groups of all classes $c\ge 5$, as required. \hfill$\square$

\medskip

\noindent{\bf Example} When $T=2^P{\rm ex}$, taking $e=4$ in the above proof gives a chiral pair of maps of type $\{32, 16\}$ and genus $105$, as shown in~\cite{DJKNS}; these are the duals of the maps C105.25 in~\cite{Con}, with $|A|=512$ and $c=5$. (However, according to~\cite{CDNS} the smallest chiral pair realising a nilpotent group are the duals of C25.1, of type $\{16,4\}$ and genus $25$, with $|A|=256$ and $c=6$.)

\medskip

The automorphism groups constructed in the proof of Theorem~\ref{nilp} have bounded derived length, since $D_m$ and $G_{p,e,f}$ are metabelian. By contrast, we have the following result:

\begin{thm}\label{solv}
Let $T$ be an edge-transitive class. Then $\G(T)$ contains a finite solvable group of derived length $l$ if and only if either
\begin{itemize}
\item $T=1$, $2^{\sigma}$, $3$ or $4^{\sigma}$ for some $\sigma$ and $l\ge 1$, or
\item $T=2^{\sigma}{\rm ex}$ or $5^{\sigma}$ for some $\sigma$ and $l\ge 2$.
\end{itemize}
\end{thm}

\noindent{\sl Proof.} Since Lemma~\ref{abelian} deals with abelian groups, we may assume that $l\ge 2$. First let $T=1$, so that $N(T)=\Gamma\cong V_4*C_2$. Define a sequence of normal (in fact, characteristic) subgroups $\Psi_n$ of $\Gamma$ by
\[\Psi_0=\Gamma, \quad{\rm and}\quad \Psi_{n+1}=\Psi_n'\Psi_n^2 \quad{\rm for}\quad n\ge 0.\]
Thus $\Psi_n/\Psi_{n+1}$ is the maximal quotient of $\Psi_n$ which is an elementary abelian $2$-group, so $\Psi_n$ is a free group of rank $r_n$, and $|\Psi_n/\Psi_{n+1}|=2^{r_n}$, where
\[r_1=3 \quad{\rm and}\quad r_{n+1}=2^{r_n}(r_n-1)+1\quad{\rm for}\quad n\ge 1.\]

Each quotient $\Gamma/\Psi_n$ is a finite solvable group of derived length $l=n$: clearly $l\le n$, and one can prove equality by induction on $n$, using the fact that each two-step quotient $\Psi_{n-1}/\Psi_{n+1}$ is non-abelian, a consequence of the existence of finitely generated non-abelian groups of exponent $4$. The maps $\M_n$ corresponding to the subgroups $\Psi_n$ prove the result for the class $T=1$, and the result follows for the classes $T=2^{\sigma}$, $3$ and $4^{\sigma}$ by Lemma~\ref{regmapslemma}(a) since $\Gamma/\Psi_n$ is non-abelian for $n\ge 2$. 

Now let $T=2^P{\rm ex}$, so that $N(T)=\Gamma^+$. We can argue as before, starting instead with $\Psi_0=\Gamma^+$, so that subsequent terms $\Psi_n$ ($n\ge 1$) are identical to those used above. The corresponding maps $\M_n$ cannot be used in this case, since they are regular. Instead, we can replace them with their joins ${\mathcal N}_n:=\M_n\vee {\mathcal S}$, corresponding to normal subgroups $N_n:=\Psi_n\cap M$ of $\Gamma^+$, where $\mathcal S$ is one of the chiral pair of Edmonds embeddings of the complete graph $K_8$, with a metabelian automorphism group $\Gamma^+/M\cong AGL_1(8)\cong V_8\rtimes C_7$ (see~\cite{JJ}). Since $AGL_1(8)$ has no subgroups of index $2$ we have $\Psi_nM=\Gamma^+$, so
\[\Gamma^+/N_n = \Psi_n/N_n\times M/N_n \cong \Gamma^+/M \times \Gamma^+/\Psi_n.\]
Here both direct factors are characteristic subgroups (the first is generated by the elements of odd order, and the second is its centraliser), so a forbidden automorphism of $\Gamma^+/N_n$ would induce one on $\Gamma^+/M$, contradicting the chirality of $\mathcal S$. Thus ${\mathcal N}_n$ is in class $T$, with
${\rm Aut}\,{\mathcal N}_n\cong \Gamma^+/N_n$,
a finite solvable group of derived length $l=\max\{n, 2\}$. This deals with the class $2^P{\rm ex}$ and hence, by Lemma~\ref{regmapslemma}(b), with the remaining classes $2^{\sigma}{\rm ex}$ and $5^{\sigma}$. \hfill$\square$

\medskip

This completes the proof of Theorem~\ref{nilpandsolv}.

%%%%%%%%%%%%%%%%

\section{Coverings of maps and their groups}

We have considered various classes of finite groups as automorphism groups realised in the various edge-transitive classes $T$. We now widen the class of candidates for automorphism groups by considering coverings of maps and of their automorphism groups.

Suppose that a map $\M$ in class $T$ corresponds to a map subgroup $M$ of $\Gamma$, so it has automorphism group $G\cong N(T)/M$. Covers $\tilde\M$ of $\M$ in $T$ correspond to normal subgroups $\tilde M$ of $N(T)$ contained in $M$, with normaliser $N(T)$, so they have automorphism group $\tilde G\cong N(T)/\tilde M$, which covers $G$.

One simple way of finding subgroups $\tilde M$ of $M$ which are normal in $N(T)$ is to take them to be characteristic subgroups of $M$. This ensures that $N_{\Gamma}(\tilde M)\ge N(T)$, but this inclusion could be proper, so that $\tilde M$ would then be in some class $T'$ properly covered by $T$. If this happens, then $\tilde M$ must be normal in some subgroup $N^*$ of $\Gamma$ containing $N(T)$ with index $2$, so that $N(T)$ is normal in $N^*$. Now $M/\tilde M$ is a normal subgroup of $N(T)/\tilde M$ (since $M$ is normal in $N(T)$); if we can choose $\tilde M$ so that $M/\tilde M$ is a {\sl characteristic\/} subgroup of $N(T)/\tilde M$, then it must be a normal subgroup of $N^*/\tilde M$, so that $M$ is normal in $N^*$, contradicting the fact that $N_{\gamma}(M)=N<N^*$. Thus there is no such subgroup $N^*$, so $N_{\Gamma}(\tilde M)=N(T)$, and hence $\tilde M$ is in class $T$. Thus we have proved:

\begin{prop}
If $N_{\Gamma}(M)=N(T)$, so that the corresponding map $\M$ is in the class $T$, and if $\tilde M$ is a characteristic subgroup of $M$ such that $M/\tilde M$ is a characteristic subgroup of $N(T)/\tilde M$, then $N_{\Gamma}(\tilde M)=N(T)$ and $\tilde M$ corresponds to a map $\tilde\M$ in $T$.
\end{prop}\hfill$\square$

\noindent{\bf Remarks. 1.} One way of ensuring that $\tilde M$ is a characteristic subgroup of $M$ is to use the `Macbeath trick', taking $\tilde M=M'M^n$, the group generated by the commutators and $n$th powers of elements in $M$, for some integer $n\ge 2$; if $G$ is  finite, then $M$ is finitely generated, so $M/\tilde M$ and hence $\tilde G$ are finite.

\smallskip

\noindent{\bf 2.}
One way of ensuring that $M/\tilde M$ is a characteristic subgroup of $N(T)/\tilde M$ is to take it to be a Hall subgroup, that is, a subgroup of order coprime to its index, since a normal Hall subgroup must be a characteristic subgroup. For example, when using the Macbeath trick, described above, one could take $n$ to be coprime to $|G|$.

\medskip

%{\color{blue}[Give some examples.]}

\noindent{\bf Example} In the proof of Theorem~\ref{symgps}, with $T=2^P{\rm ex}$, we used an epimorphism $N(T)=\Gamma^+\to S_6$, mapping the generators $X$ and $Y$ of $\Gamma^+$ to $x=(1,\ldots, 6)$ and $y=(1,2)(3,5)$,
to give a map $\M$ in this class with ${\rm Aut}\,\M\cong S_6$. Since $x$ and $xy$ have order $6$, this map has type $\{6,6\}$, so it has $|V|=|S_6|/6=120$ vertices, and the same number $|F|$ of faces. Since there are $360$ edges it has characteristic $-120$ and genus $g=61$. (It is, in fact, the chiral map C61.7 in Conder's list~\cite{Con}.) The corresponding map subgroup $M$ is a Fuchsian group with signature $(61;-;240)$, meaning that it has $2g=122$ hyperbolic generators $a_i, b_i\;(i=1,\ldots, 61)$, no elliptic generators, and $|V|+|F|=240$ parabolic generators $y_j\;(j=1,\ldots, 240)$, satisfying the single defining relation
\[\prod_i[a_i,b_i].\prod_jy_j=1.\]
Thus $M$ is a free group of rank $r=2g+|V|+|F|-1=361$, freely generated by the elements $a_i, b_i$ and $y_j$ with $j\ne 240$.

Applying the Macbeath trick gives a characteristic subgroup $\tilde M=M'M^n$ of index $n^r$ in $M$. This corresponds to an orientably regular map $\tilde\M$ of type $\{6n, 6n\}$ which is an $n^r$-sheeted regular cover of $\M$, branched over the vertices and face-centres. For any choice of $n\ge 2$, the group $M/\tilde M\cong C_n^r$ is the largest abelian normal subgroup of $\Gamma^+/\tilde M$, so it is a characteristic subgroup and hence $\tilde\M$ is in class $2^P{\rm ex}$. 

There is an intermediate $n^{2g}$-sheeted unbranched covering $\hat\M$ of $\M$, also in class $2^P{\rm ex}$ but now of type $\{6,6\}$, corresponding to the map subgroup $\hat M$ generated by $\tilde M$ and the parabolic generators $y_j$ of $M$, so that $M/\hat M\cong H_1(\M;\Z_n)\cong C_n^{2g}$; this is equivalent to replacing the parent group $\Gamma^+$ with the triangle group $\Delta(6,2,6)$, the parent group for edge-transitive maps of class $2^P{\rm ex}$ and type $\{6,6\}$, by adding the relations $X^6=(XY)^6=1$.

One can identify $M$ with the fundamental group of the surface ${\mathcal S}^{\circ}$ obtained by puncturing the underlying surface $\mathcal S$ of $\M$ at its vertices and face-centres. Then $M/\tilde M$ can be regarded as the first homology group $H_1({\mathcal S}^{\circ};\Z_n)$  of ${\mathcal S}^{\circ}$ over $\Z_n$, with the action of $\Gamma^+/M$ on it by conjugation coinciding with the natural action of ${\rm Aut}\,\M$ on the homology of $\mathcal S^{\circ}$. Similarly, $M/\hat M$ can be regarded as the first homology group $H_1(\mathcal S;\Z_n)$  of $\mathcal S$ over $\Z_n$, with the corresponding action of ${\rm Aut}\,\M$.

%%%%%%%%%%%%%%%%

\section{Infinite automorphism groups}\label{infautogps}

A map $\M$ is compact if and only if the corresponding map subgroups have finite index in $\Gamma$. In the case of edge-transitive maps, this is equivalent to ${\rm Aut}\,\M$ being finite, as we have assumed until now. If we remove this restriction, then there is an even greater abundance and variety of maps and automorphism groups in each edge-transitive class $T$. So much so, in fact, that it is impossible to give as comprehensive an analysis as that given earlier for compact maps, so instead we shall just illustrate a few phenomena which distinguish non-compact edge-transitive maps from compact ones.

%%%%%%%%

\subsection{Uncountably many groups and maps}\label{uncount}

Since $\Gamma$ is finitely generated, it has only countably many conjugacy classes of subgroups of finite index, so there are only countably many isomorphism classes of compact maps. On the other hand, there are uncountably many non-compact maps in each edge-transitive class. The following proof of this is adapted from Bernhard Neumann's proof~\cite{Neu} that there are uncountably many non-isomorphic $2$-generator groups (see also~\cite[\S III.B]{dlH}).

\begin{thm}\label{uncountable}
Each of the $14$ classes $T$ of edge-transitive maps contains $2^{\aleph_0}$ maps $\mathcal M$ with empty boundary and with mutually non-isomorphic automorphism sgroups ${\rm Aut}\,{\mathcal M}$.
\end{thm}

\noindent(Of course, it follows that these maps are also mutually non-isomorphic.)

\medskip

\noindent{\sl Proof.} First let $T=1$. For each integer $n\ge 5$ such that $n\equiv 1$ mod~$(4)$ we define an epimorphism $\theta_n:\Gamma\to A_n$, $R_i\mapsto r_i$ as follows. We number the vertices of a regular $n$-gon $1, 2, \ldots, n$ in cyclic order, and let
\[r_1=(1, n) (2, n-1)\ldots (m-1, m+1)(m),\]
\[r_2=(1)(2, n)(3, n-1)\ldots (m, m+1)\]
be the reflections fixing the vertices $m:=(n+1)/2$ and $1$, so that $r_1$ and $r_2$ are involutions in $A_n$, with
\[x:=r_1r_2=(1, 2, \ldots, n).\]
Now let 
\[r_0=(2, 3)(n-1, n),\]
an involution in $A_n$ commuting with $r_2$. Then
\[[r_0, x^{-2}]=(1, 2, 3)(4, 5)(n-1, n),\]
so that
\[a:=[r_0, x^{-2}]^{-2}=(1, 2, 3).\]
Since $x$ and $a$ generate $A_n$ by Lemma~\ref{Angens}(c), so do $r_0, r_1$ and $r_2$, and thus $\theta_n$ is an epimorphism.

Now let $S$ be any strictly increasing sequence $(n_k)$ of integers satisfying $5\le n_k\equiv 1$ mod~$(4)$ for all $k\in\N$. The epimorphisms $\theta_{n_k}:\Gamma\to A_{n_k}$ define an action $\theta$ of $\Gamma$ on the disjoint union $\Omega$ of the sets $\Omega_k:=\{1, 2, \ldots, n_k\}$ for $k\in\N$, with $\Gamma$ acting as the alternating group $G_k:={\rm Alt}\,\Omega_k\cong A_{n_k}$ on each of its orbits $\Omega_k$. The resulting permutation group $G(S)$ induced by $\Gamma$ on $\Omega$ is thus a quotient of $\Gamma$. We will prove the following:

\medskip

\noindent{\bf Claim:} $G(S)$ has a normal subgroup isomorphic to $A_n$ if and only if $n$ is a term $n_k$ in $S$.

\medskip

\noindent It then follows that the $2^{\aleph_0}$ sequences $S$ satisfying the above conditions give $2^{\aleph_0}$ mutually non-isomorphic quotient groups $G(S)$ of $\Gamma$. The corresponding maps $\M$ with ${\rm Aut}\,\M\cong G(S)$ are without boundary, since each $R_i\;(i=0, 1, 2)$ is represented by a non-identity permutation. 

\medskip

\noindent{\sl Proof of the Claim.} By the simplicity of the alternating groups $G_k$, any normal subgroup of $G(S)$ must act on each $\Omega_k$ either trivially or as $G_k$. It follows that any finite normal subgroup $N$ of $G(S)$ must fix all but finitely many orbits $\Omega_k$ (otherwise it contains alternating groups of unbounded degrees), so $N$ is a subdirect product of finitely many alternating groups $G_k$. Since these groups $G_k$ are mutually non-isomorphic simple groups, $N$ is in fact their direct product. It follows that if any alternating group $A_n$ is isomorphic to a normal subgroup of $G(S)$ then $n$ must be a term $n_k$ in $S$ for some $k$.

Now we prove the converse, that for each $n=n_k$ in $S$ there is a normal subgroup $G_k\cong A_n$ in $G(S)$. For any odd $n$, the elements $x=(1, 2, \ldots, n)$ and $a=(1, 2, 3)$ of $A_n$ satisfy
\[x^{-(i-2)}ax^{i-2}=(i-1,i,i+1)\quad\hbox{for all}\quad i\in{\mathbb N},\]
where we replace entries with their remainders mod~$(n)$ if necessary. It follows that
\[[a, x^{-(i-2)}ax^{i-2}]=1\quad {\rm if}\quad n\ge i+1\ge 6,\]
but
\[[a, x^{-(i-2)}ax^{i-2}]\ne 1\quad {\rm if}\quad n= i.\]
The element $X=R_1R_2$ of $\Gamma$ induces $x$ on each orbit $\Omega_k$, so the element $A:=[R_0, X^2]^{-2}$ induces $a=(1,2,3)$.
It follows that for each $k\in\N$ the element
\[C_k:=[A,X^{-(n_k-2)}AX^{n_k-2}]\] 
of 
$\Gamma$ fixes $\Omega_l$ for all $l>k$, but acts non-trivially on $\Omega_k$ (and possibly on some orbits $\Omega_l$ with $l<k$). The same applies to all conjugates of $C_k$ in $\Gamma$, and hence to its normal closure $N$ in $\Gamma$. Acting on $\Omega$, this induces a finite normal subgroup $N_k$ of $G(S)$, specifically the direct product of $G_k$ and possibly some subgroups $G_l$ with $l<k$. Since these direct factors are mutually non-isomorphic simple groups, $G_k$ is a characteristic subgroup of $N_k$, and hence it is a normal subgroup of $G(S)$, isomorphic to $A_{n_k}$, as claimed.

This proves the Theorem for the class $T=1$, and the result follows for the classes $T=2^{\sigma}$, $3$ and $4^{\sigma}$ by Lemma~\ref{regmapslemma}(a).

Now let $T=2^P{\rm ex}$, so that $N(T)=\Gamma^+$. The method of proof is similar to that used for $T=1$, but now we map the generators $X$ and $Y$ of $\Gamma^+$ to the elements
\[x:=(1, 2, \ldots, n) \quad{\rm and}\quad y:=(1, 2)(3, 4)(5, 6)(8, 9)\]
of $A_n$ for each odd $n\ge 11$. Then
\[[x, y]=(1, 2, 4, 6, 7, 5, 3)(8, 9, 10),\]
so
\[[x, y]^7=(8,9,10).\]
Thus $\langle x, y\rangle=A_n$ by Lemma~\ref{Angens}(c), and we have an epimorphism $\Gamma^+\to A_n$.  The proof now proceeds as in the case $T=1$, again using the elements $x$ and
\[a:=x^7[x,y]^7x^{-7}=(1, 2, 3)\]
to prove the Claim, but this time for any strictly increasing sequence $S$ of odd integers $n=n_k\ge 11$. Since $n>9$, no permutation in $S_n$ inverting $x$ centralises $y$, so these quotients $A_n$ of $\Gamma^+$ have no forbidden automorphisms; the same therefore applies to the quotients $G(S)$ of $\Gamma^+$, since any forbidden automorphism of them would induce forbiddden automorphisms of the characteristic subgroups $G_k$. This deals with the class $T=2^P{\rm ex}$, and Lemma~\ref{regmapslemma}(b) then extends the result to the remaining classes. \hfill$\square$

\medskip

\noindent{\bf Remark} There is a shorter but less explicit proof for $T=1$ (and hence, by Lemma~\ref{regmapslemma}(a), for $T=2^{\sigma}$, $3$ and $4^{\sigma}$). A group $H$ is {\em SQ-universal\/} if every countable group can be embedded in a quotient of $H$. Britton and Levin proved that any free product $A*B$ is SQ-universal, provided $|A|\ge 2$ and $|B|\ge 3$ (see~\cite[Theorem V.10.3]{LS}), so in particular the parent groups $N(T)$ all have this property (see Table~\ref{forbidden}). Now any finitely generated SQ-universal group $H$ must have uncountably many non-isomorphic quotients: as we have seen, B.~H.~Neumann~\cite{Neu} proved that there are uncountably many non-isomorphic finitely generated groups $F$; these are all countable, so each can be embedded in some quotient $Q$ of $H$; however, $Q$ is countable, so it has only countably many finitely generated subgroups; each isomorphism class of quotients $Q$ can therefore embed only countably many of the groups $F$, so $H$ must have uncountably many isomorphism classes of quotients. In particular, this applies to the parent groups $H=N(T)$. This deals with the case $T=1$, where there are no forbidden automorphisms to exclude, and the cases $T=2^{\sigma}$, $3$ and $4^{\sigma}$ then follow by Lemma~\ref{regmapslemma}(a). However, in the cases $T=2^{\sigma}{\rm ex}$ and $5^{\sigma}$ this line of proof breaks down at this point, since we need quotients without the corresponding forbidden automorphisms. By Lemma~\ref{regmapslemma}(b) it would be sufficient to deal with the class $2^P{\rm ex}$.

\medskip

\noindent{\bf Problem} Can one extend this method of proof to the case $T=2^P{\rm ex}$?

%%%%%%%%%

\subsection{An embedding theorem}\label{embed}

The next result shows that for each edge-transitive class $T$, any phenomenon, no matter how pathological, which can happen within a countable group can happen within a group $G\in\G(T)$. As in the case of Theorem~\ref{uncount}, SQ-universality provides a quick proof for $T=1$ (and hence also for $T=2^{\sigma}$, $3$ and $4^{\sigma}$), but in order to eliminate forbidden automorphisms a more explicit argument is needed for the remaining classes.

\begin{thm}\label{embedctble}
For each of the $14$ edge-transitive classes $T$, every countable group $C$ is isomorphic to a subgroup of ${\rm Aut}\,{\mathcal M}$ for some map $\mathcal M$ in $T$.
\end{thm}

\noindent{\sl Proof.} Schupp~\cite[Theorem I]{Sch}, generalising an earlier result of Gorju\v skin~\cite{Gor}, proved that if $A, B$ and $C$ are groups with $|A|\ge 2$, $|B|\ge 3$ and $|C|\le|A*B|$, then $C$ can be embedded in a simple group $S$ generated by a pair of subgroups isomorphic to $A$ and $B$. It follows that for each parent group $N(T)$, every countable group $C$ can be embedded in a simple quotient $S=N(T)/K$ of $N(T)$. Without loss of generality we may assume that $S\not\cong A_9$.

By Theorem~\ref{altgps}, for each $T$ there is a normal subgroup $L$ of $N(T)$ such that $N(T)/L$ is isomorphic to $A_9$ and has no forbidden automorphisms. If we take $M:=K\cap L$ then $G:=N(T)/M$ is a quotient of $N(T)$ isomorphic to $S\times A_9$, so it contains a subgroup isomorphic to $C$. As non-isomorphic simple groups, $S$ and $A_9$ are characteristic subgroups of $G$; any forbidden automorphisms of $G$ would therefore induce forbidden automorphisms of the quotient $A_9$, whereas this quotient has none. Thus $M$ corresponds to a map $\M$ in $T$ with automorphism group $G$ containing a copy of $C$. \hfill$\square$

\medskip

Indeed, by arguing as in the Remark following Theorem~\ref{uncount}, one can use Schupp's result to show that if $T=1$, $2^{\sigma}$, $3$ or $4^{\sigma}$ then $\G(T)$ contains uncountably many non-isomorphic simple groups. As before, however, the problem of eliminating forbidden automorphisms means that a separate argument would be needed to prove this for the remaining classes; once again, it would be sufficient by Lemma~\ref{regmapslemma}(b) to deal with $2^P{\rm ex}$.

\medskip

\noindent{\bf Problem} Can one extend this method of proof to the case $T=2^P{\rm ex}$?

\iffalse

{\color{blue}[GAJ \& JMJ~\cite{JonJon} used small cancellation theory to prove that, with a few exceptions having small periods, every hyperbolic triangle group has uncountably many non-isomorphic infinite simple quotient groups; John Wilson (QJM 1999) extended this. Does it apply to the parent groups $N(T)$, giving uncountably many simple quotients with no forbidden automorphisms?]}

\fi

%%%%%%%%%

\subsection{Growth}

One can adapt the proof of Theorem~\ref{embedctble} to produce, for each edge-transitive class $T$, uncountably many examples of automorphism groups and edge-transitive maps of intermediate growth. In~\cite{Gri80} Grigorchuk produced an example of a finitely-generated group which, as he showed later~\cite{Gri84}, has intermediate growth, in the sense that the number of distinct elements represented by words of length at most $l$ in the generators grows faster than any polynomial in $l$, but slower than exponentially. In~\cite{Gri84} he extended this example by constructing an uncountable set of groups with this property. These groups $G$ are all subgroups of the automorphism group of the infinite rooted binary tree, generated by automorphisms $a, b, c$ and $d$ satisfying\[a^2=b^2=c^2=d^2=bcd=1.\]
Thus $\langle b, c, d\rangle\cong V_4$, so each of these groups $G$ is a quotient of $\Gamma$ and is therefore the automorphism group of a regular map. Like $G$, this map has intermediate growth, meaning that the number of vertices, edges or faces at a given graph-theoretic distance from some chosen base point has a rate of growth which is intermediate between polynomial and exponential (see~\cite{Jon11}).

\begin{thm}
For each class $T$ of edge-transitive maps there are uncountably many maps $\mathcal M$ in $T$ such that $\mathcal M$ and ${\rm Aut}\,{\mathcal M}$ have intermediate growth.
\end{thm}

\noindent{\sl Proof.} Having intermediate growth is inherited by subgroups of finite index, so each parent group $N(T)$ has uncountably many quotients $S=N(T)/K$ of intermediate growth. As in the proof of Theorem~\ref{embedctble}, by taking $M=K\cap L$ where $N(T)/L$ is isomorphic to $A_9$ with no forbidden automorphisms, we obtain uncountably many quotients $Q=N(T)/M$ and associated edge-transitive maps $\mathcal M$ in $T$ with intermediate growth. To show that $Q$ has no forbidden automorphisms, note that $S$ is residually nilpotent (since each Grigorchuk group induces a finite $2$-group on vertices up to any given level in the tree), so $N(T)=KL$, and hence $Q\cong S\times A_9$ with both direct factors as characteristic subgroups. \hfill$\square$

%%%%%%%%%

\subsection{Decidability}

It is reasonable to ask whether there are algorithms to answer  `sensible' questions about edge-transitive maps and their automorphism groups. A  further result of Schupp~\cite{Sch} eliminates such a hope. A property $\mathcal P$ of groups is called a {\em Markov property\/} if
\begin{itemize}
\item it is preserved by isomorphisms,
\item there is a finitely presented group with property $\mathcal P$, and
\item there is a finitely presented group which cannot be embedded in any finitely presented group with property $\mathcal P$. 
\end{itemize}
If $\mathcal P$ is hereditary (inherited by subgroups) then the last two conditions are equivalent to
\begin{itemize}
\item there are finitely presented groups with and without property $\mathcal P$.
\end{itemize}
Examples of Markov properties include being trivial, finite, nilpotent, solvable, free, and torsion-free.

The edge-transitive maps in each class $T$ have automorphism groups obtained by adding further relations to the standard presentation of $N(T)$. For example, the automorphism group of the icosahedron, isomorphic to $A_5\times C_2$, is obtained by adding the relations
\[(R_0R_1)^3=(R_1R_2)^5=1\]
to the presentation of $\Gamma=N(1)$. When discussing decision problems, it makes sense to restrict attention to finitely presented quotients, obtained by adding finitely many relations in the standard generators to parent groups $N(T)$.

Schupp proved in~\cite[Theorem~II]{Sch} that if $\mathcal P$ is a Markov property of groups, and $G$ is a group with a finite presentation $\langle X \mid R\rangle$ and a quotient $H*K$ where $|H|\ge 2$ and $|K|\ge 3$, then the problem of determining $\mathcal P$ for groups with finite presentations $\langle X\mid R\cup S\rangle$ is undecidable. The presentations of the groups $N(T)$ in Proposition~\ref{parents} and their free product decompositions in Table~\ref{forbidden} show that they all satisfy the hypotheses of this result, so we have:

% The free product decompositions of these groups $N(T)$ imply that the following is an immediate consequence of Schupp's Theorem~II in~\cite{Sch} {\color{blue}[Give more details, including specific free product decompositions, and a statement of Schupp's theorem]}:

\begin{thm}
For each class $T$ of edge-transitive maps, and each Markov property $\mathcal P$ of groups, it is undecidable whether or not adding a finite set of relations to the standard presentation of $N(T)$ produces a quotient group with property $\mathcal P$. \hfill$\square$
\end{thm}

For example, it is undecidable whether adding a finite set of relations to those of $N(T)$ yields a map which is compact, or even whether it collapses to the basic map ${\mathcal N}(T)$.

%%%%%%%%%%%%%%%%%%%%
%%%%%%%%%%%%%%%%%%%%

\newpage

\part{Topological properties of edge-transitive maps}\label{topology}

%{\color{blue}[All of this part is new material, not in the arXiv version 2.]}

%\medskip

In this part of the text, we consider the topological properties of edge-transitive maps. For simplicity, and on account of their importance in the theories of Riemann surfaces and dessins d'enfants, we first consider maps without boundary. We then consider edge-transitive maps with boundary, showing that in some of the fourteen classes they do not arise, that in some they arise and are easily classified, and that in the remaining classes the abundance of examples makes such a classification impossible. We also classify edge-transitive maps with free edges, and discuss almost edge-transitive maps, those non-edge-transitive maps which have edge-transitive medial maps.

\section{Topology of maps without boundary}

In this section, we consider topological properties of edge-transitive maps $\M$ without boundary, that is, those for which the map subgroups $M$ contain no reflections (conjugates of $R_i$ for $i=0, 1, 2$). This applies to all maps in classes $T=2^P{\rm ex}$, $5$, $5^*$ and $5^P$ (since $N(T)$ contains no reflections), but only to some maps in the other ten classes.

\subsection{Orientability}\label{orientability}

Orientability is easy to determine for maps without boundary: a simple argument yields the following:

\begin{lemma}
If $\M$ is a map without boundary, with a map subgroup $M\le\Gamma$, then the following are equivalent:
\begin{itemize}
\item $\M$ is orientable;
\item the coset diagram for $M$ in $\Gamma$ with respect to the generators $R_i$ is bipartite;
\item $M\le\Gamma^+$;
\item $\M$ covers the spherical map ${\mathcal N}(2^P{\rm ex})$.
\end{itemize}
\end{lemma}

(Note that since $\Gamma^+$ is a normal subgroup of $\Gamma$, these criteria are independent of the choice of a map subgroup $M$ for $\M$.) In particular, this result implies that all maps in classes $T=2^P{\rm ex}$, $5$ and $5^*$ are orientable (since $N(T)\le\Gamma^+$ in these cases), whereas the other eleven classes contain both orientable and non-orientable maps.

Every non-orientable map $\M$ (edge-transitive or not), corresponding to a conjugacy class of map subgroups $M\not\le\Gamma^+$, has an orientable double cover $\M^+$ corresponding to the subgroups $M^+=M\cap\Gamma^+$, and $\M$ is the quotient of $\M^+$ by its automorphism group of order $2$ corresponding to $M/M^+$. When $\M$ has empty boundary the covering $\M^+\to\M$ is unbranched, with $\chi(\M^+)=2\chi(\M)$ if $\M$ is compact, so that if $\M$ has (non-orientable) genus $g$ then $\M^+$ has (orientable) genus $g-1$. For example, the non-orientable regular embedding $\M$ of the complete graph $K_6$ in the real projective plane (see Figure~\ref{projplKn}) lifts to the icosahedral map $\M^+$ on the sphere, of which it is the antipodal quotient. The regular embeddings of $K_4$ and the cube graph $Q_3$ in these two surfaces (see Figures~\ref{tetmap} and \ref{projplKn}) are related in the same way.

Now $M$ and $\Gamma^+$ are both normalised by $N_{\Gamma}(M)$, and hence so is their intersection $M^+$, so $N_{\Gamma}(M^+)\ge N_{\Gamma}(M)$. The group ${\rm Aut}\,\M^+\cong N_{\Gamma}(M^+)/M^+$ has a subgroup
\[N_{\Gamma}(M)^+/M^+\times M/M^+\cong  N_{\Gamma}(M)/M\times M/M^+\cong {\rm Aut}\,\M\times C_2\]
obtained by lifting each automorphism of $\M$ to two automorphisms of $\M^+$.

Since we are assuming that $\M$ is without boundary, the orientation-reversing involution $i\in{\rm Aut}\,\M^+$ generating the direct factor corresponding to $M/M^+$ has no fixed points. Conversely, if $\K$ is any orientable map without boundary, with a fixed-point-free orientation-reversing involution $i\in{\rm Aut}\,\K$, then $\K$ arises in this way as the orientable double cover $\M^+$ of a non-orientable map $\M=\K/\langle i\rangle$ without boundary. Thus non-orientable maps without boundary can all be obtained from orientable maps without boundary by factoring out such an automorphism $i$.

If $\M$ is edge-transitive, in class~$T$, then the fact that $N_{\Gamma}(M^+)\ge N_{\Gamma}(M)$ means that $\M^+$ is also edge-transitive, in a class $T^+$ covered by $T$ (see Lemma~\ref{covering} for such pairs $T, T^+$), with 
\[|{\rm Aut}\,\M^+:{\rm Aut}\,\M|=2|N(T^+):N(T)|=2\frac{n_T}{n_{T^+}}=2, 4 \;{\rm or}\; 8.\]
The most typical value for this index is $2$, corresponding to the cases where $T^+=T$, or equivalently $N_{\Gamma}(M^+)=N_{\Gamma}(M)$, so that ${\rm Aut}\,\M^+ = {\rm Aut}\,\M\times C_2$. For instance, this always happens if $\M$ is regular (that is, $T=1$), as in the examples on the real projective plane mentioned above. However, as we will see in Section~\ref{stability}, there are edge-transitive maps $\M$ (with and without boundary) for which $\M^+$ has extra automorphisms, so that ${\rm Aut}\,\M^+>{\rm Aut}\,\M\times C_2$ and $T^+\ne T$.

%%%%%%%%%%%%%%%%%%%%%%%%%

\subsection{Type}

A map $\M$ with empty boundary has type $\{m,n\}$ if $m$ and $n$ are the least common multiples of the valencies of its faces and vertices (or $\infty$ if these do not exist). They are the orders of the elements $r_0r_1$ and $r_1r_2$ of the monodromy group $G$ of $\M$, induced by the elements $R_0R_1$ and $R_1R_2$ of $\Gamma$ acting on the cosets of a map subgroup $M$. Maps of a given type $\{m,n\}$ therefore correspond to permutation representations, or conjugacy classes of subgroups, of the extended triangle group
\[\Delta=\Delta[m,2,n]=\langle R_0, R_1, R_2\mid (R_0R_1)^m=(R_0R_2)^2=(R_1R_2)^n=1\rangle.\]

This group $\Delta$ is infinite provided $m^{-1}+n^{-1}\le 1/2$. As a finitely generated linear group, by a theorem of Mal'cev~\cite{Mal} it is residually finite (meaning that the intersection of the normal subgroups of finite index is the identity subgroup), so in this case there are (infinitely many) finite regular maps of each such type $\{m,n\}$. If  $m^{-1}+n^{-1} > 1/2$ then provided $m, n\ge 2$ there is a well-known unique regular map $\M$ of this type on the sphere, and in some cases also, where $\M$ is invariant under the antipodal isometry, one on the real projective plane.

For classes $T\ne 1$ there may be some restrictions on the possible types. For example, there are (again infinitely many) finite chiral maps of each type $\{m,n\}$ with $m^{-1}+n^{-1}\le 1/2$, but there are none with $m^{-1}+n^{-1}> 1/2$ (see~\cite{CHNS} or \cite{Jon15}, for example).

For the remaining classes $T\ne 1, 2^P{\rm ex}$, there are parity restrictions on the type of a map without boundary. If $R_0R_1$ is not an element of $N(T)$ (or of its other conjugate when $T=4^{\sigma}$ for some $\sigma=\emptyset, *$ or $P$) then since $N(T)$ has even index in $\Gamma$ it follows that $r_0r_1$ has a cycle of even length, so the parameter $m$ must be even for each map in class $T$; this applies to the classes $T=2, 2^P, 2^*{\rm ex}, 3, 4^{\sigma}, 5$ and $5^P$. Similarly, by considering $r_1r_2$ we see that $n$ must be even for maps in classes $T=2^*, 2^P, 2\,{\rm ex}, 3, 4^{\sigma}, 5^*$ and $5^P$. In particular, $m$ and $n$ are both even for $T=2^P, 3, 4^{\sigma}$ and $5^P$.
%{\color{blue}[Check these!]}

%%%%%%%%%%%%%%%%%%%%%%%%%

\subsection {Euler characteristic}
 
We now calculate the Euler characteristic
\[\chi=V-E+F\]
of a compact edge-transitive map $\M$, where $V, E$ and $F$ denote the numbers of vertices, edges and faces. We will assume that $\M$ has no free edges, for otherwise edge-transitivity implies that all edges are free, so that $\M$ is a star map, a regular map on the sphere with one vertex and one face. (See Section~\ref{free} for a classification of all edge-transitive maps with free edges, including those maps with non-empty boundary.)

If $\Phi$ is the set of flags of $\M$ then
\[|\Phi|=|\Gamma:M|=|\Gamma:N(T)|.|N(T):M|=n_T|G|,\]
where $\M$ is in class $T$. Each edge is incident with four flags,  so
\[E=\frac{|\Phi|}{4}=\frac{n_T|G|}{4}.\]

By inspection of the maps ${\mathcal N}(T)$, $G$ has at most two orbits on the vertices of $\M$. If $\M$ is vertex-transitive (equivalently, $T=1$, $2^*$, $2^P$, $2^{\sigma}{\rm ex}$, $4^*$, $4^P$, $5^*$ or $5^P$), then all vertices have the same valency $n$, the order of the image of $R_1R_2$ in $G$, so each vertex is incident with $2n$ flags and hence
\[V=\frac{|\Phi|}{2n}=\frac{n_T|G|}{2n}.\]
If, on the other had, $\M$ is not vertex-transitive (equivalently, $T=2$, $3$, $4$ or $5$), then each edge is incident with one vertex from each of the two orbits. If these have valencies $n_1$ and $n_2$ (the lengths of the cycles of $R_1R_2$ on $\Phi$), then $|\Phi|/2$ of the flags must be incident with a vertex in one orbit, and the remaining $|\Phi|/2$ with a vertex in the other, so these orbits have lengths
\[V_i=\frac{|\Phi|/2}{2n_i}=\frac{n_T|G|}{4n_i}\]
and hence
\[V=V_1+V_2=\frac{N_T|G|}{4}\left(\frac{1}{n_1}+\frac{1}{n_2}\right).\]
We can cover both cases by defining $n_1=n_2=n$ when $\M$ is vertex-transitive, so that
\[V=\frac{n_T|G|}{2n^*},\]
where
\[n^*=\frac{2}{n_1^{-1}+n_2^{-1}}=\frac{2n_1n_2}{n_1+n_2}\]
is the harmonic mean of $n_1$ and $n_2$.

A similar argument shows that
\[F=\frac{n_T|G|}{2m^*},\]
where
\[m^*=\frac{2}{m_1^{-1}+m_2^{-1}}=\frac{2m_1m_2}{m_1+m_2}\]
and $m_1$ and $m_2$ denote the length(s) of the orbit(s) of $G$ on faces. (There are two orbits if $T=2^*, 3, 4^*$ or $5^*$, and one in all other cases.) Thus $\M$ has type $\{m,n\}$ where
\[m={\rm lcm}(m_1, m_2)\quad{\rm and}\quad n={\rm lcm}(n_1, n_2).\]
More informatively, we can write that $\M$ has extended type $\{m_1|m_2, n_1|n_2\}$, where $m_1$ and $m_2$ (in either order) are the valencies of the faces, and similarly $n_1$ and $n_2$ are the valencies of the vertices.

It then follows that
\begin{equation}\label{chi}
\chi=\frac{n_T|G|}{4}\left(\frac{2}{m^*}+\frac{2}{n^*}-1\right)
=\frac{n_T|G|}{4}\left(\frac{1}{m_1}+\frac{1}{m_2}+\frac{1}{n_1}+\frac{1}{n_2}-1\right)
\end{equation}
(cf.~\cite[Proposition~2.3]{STW}).

\begin{figure}[h!]

\begin{center}
 \begin{tikzpicture}[scale=0.5, inner sep=0.8mm]

\node (a) at (-3,0) [shape=circle, draw, fill=black] {}; 
\node (b) at (3,0)  [shape=circle, draw, fill=black] {};
\node (c) at (3,6) [shape=circle, draw, fill=black] {};   
\node (d) at (-3,6)  [shape=circle, draw, fill=black] {};
\node (e) at (-0.8,8)  [shape=circle, draw, fill=black] {};
\node (f) at (5,8)  [shape=circle, draw, fill=black] {};
\node (g) at (5,2.2)  [shape=circle, draw, fill=black] {};

\draw [dashed] (a) to (b);
\draw [dashed] (b) to (c);
\draw [dashed] (c) to (d);
\draw [dashed] (d) to (a);
\draw [dashed] (d) to (e);
\draw [dashed] (c) to (f);
\draw [dashed] (b) to (g);
\draw [dashed] (e) to (f);
\draw [dashed] (f) to (g);

\node (w) at (0,0) [shape=circle, draw, fill=black] {}; 
\node (x) at (3,3)  [shape=circle, draw, fill=black] {};
\node (y) at (0,6) [shape=circle, draw, fill=black] {};   
\node (z) at (-3,3)  [shape=circle, draw, fill=black] {};

\node (h) at (-1.9,7)  [shape=circle, draw, fill=black] {};
\node (i) at (2.1,8)  [shape=circle, draw, fill=black] {};
\node (j) at (4,7)  [shape=circle, draw, fill=black] {};
\node (k) at (4,1.1)  [shape=circle, draw, fill=black] {};
\node (l) at (5,5.1)  [shape=circle, draw, fill=black] {};

\draw [thick] (a) to (x) to (d) to (w) to(c) to (z) to (b) to (y) to (a);

%%%%

\node (A) at (11,0) [shape=circle, draw, fill=black] {}; 
\node (B) at (11,4)  [shape=circle, draw, fill=black] {};
\node (C) at (8,6) [shape=circle, draw, fill=black] {};   
\node (D) at (14,6)  [shape=circle, draw, fill=black] {};
\node (E) at (8,2.2)  [shape=circle, draw, fill=black] {};
\node (F) at (14,2.2)  [shape=circle, draw, fill=black] {};
\node (G) at (11,7.6)  [shape=circle, draw, fill=black] {};

\draw [dashed] (A) to (B);
\draw [dashed] (B) to (C);
\draw [dashed] (B) to (D);
\draw [dashed] (C) to (G);
\draw [dashed] (D) to (G);
\draw [dashed] (A) to (E);
\draw [dashed] (A) to (F);
\draw [dashed] (E) to (C);
\draw [dashed] (F) to (D);

\node (X) at (11,2) [shape=circle, draw, fill=black] {}; 
\node (Y) at (9.5,5) [shape=circle, draw, fill=black] {};
\node (Z) at (12.5,5) [shape=circle, draw, fill=black] {}; 

\node (H) at (8,4.1)  [shape=circle, draw, fill=black] {};
\node (I) at (9.5,6.8)  [shape=circle, draw, fill=black] {};
\node (I) at (9.5,1.1)  [shape=circle, draw, fill=black] {};
\node (J) at (12.5,6.8)  [shape=circle, draw, fill=black] {};
\node (K) at (12.5,1.1)  [shape=circle, draw, fill=black] {};
\node (L) at (14,4.1)  [shape=circle, draw, fill=black] {};

\draw [thick] (A) to (Z) to (C) to (X) to (D) to (Y) to (A);

\end{tikzpicture}

\end{center}
\caption{An octagonal face and a hexagonal face of $\mathcal M$} 
\label{knight}
\end{figure}

\medskip

\noindent{\bf Example 1} The only edge-transitive class in which the automorphism group is intransitive on both vertices and faces is class $3$. Here we illustrate the calculation of $\chi$ with a map in this class taken from~\cite{Jon83}{\color{red}, where further examples of class~3 maps are constructed.}

Let $\M$ be the map shown in Figure~\ref{knight}. It has $8+12=20$ vertices, namely the vertices of the unit cube $\mathcal C$ (shown by broken lines) and the midpoints of its edges. Two vertices of $\M$ are joined by an edge in $\M$  if and only if they are distance $\sqrt{5}/2$ apart (that is, separated by a `knight's move' in $\mathcal C$), so they have valency $n_i=6$ or $4$ respectively. There are six octagonal faces (one for each face of $\mathcal C$), and eight hexagonal faces (one for each vertex of $\mathcal C$), so the face-valencies are $m_i=8, 6$ respectively. The symmetry group $G\cong S_4\times C_2$ of $\mathcal C$ is also the automorphism group of $\M$, acting transitively on edges but not on vertices or faces, so $\M$ is in class~$3$. Thus $N_T=4$, $|G|=48$, with $n^*=24/5$ and $m^*=48/7$, so equation (\ref{chi}) gives $\chi=-14$. Since $\M$ is orientable, it has genus $8$.

\begin{figure}[h!]

\begin{center}
 \begin{tikzpicture}[scale=0.7, inner sep=0.8mm]

\node (a) at (-3,0) [shape=circle, draw, fill=black] {}; 
\node (b) at (3,0)  [shape=circle, draw, fill=black] {};
\node (c) at (3,6) [shape=circle, draw, fill=black] {};   
\node (d) at (-3,6)  [shape=circle, draw, fill=black] {};
\node (e) at (-0.8,8)  [shape=circle, draw, fill=black] {};
\node (f) at (5,8)  [shape=circle, draw, fill=black] {};
\node (g) at (5,2.2)  [shape=circle, draw, fill=black] {};

\draw [dashed] (a) to (b);
\draw [dashed] (b) to (c);
\draw [red] (c) to (d);
\draw [red] (d) to (a);
\draw [dashed] (d) to (e);
\draw [dashed] (c) to (f);
\draw [dashed] (b) to (g);
\draw [dashed] (e) to (f);
\draw [dashed] (f) to (g);

\node (w) at (0,0) [shape=circle, draw, fill=black] {}; 
\node (x) at (3,3)  [shape=circle, draw, fill=black] {};
\node (y) at (0,6) [shape=circle, draw, fill=black] {};   
\node (z) at (-3,3)  [shape=circle, draw, fill=black] {};

\draw [red] (x) to (z);
\draw [red] (a) to (c);

\node (h) at (-1.9,7)  [shape=circle, draw, fill=black] {};
\node (i) at (2.1,8)  [shape=circle, draw, fill=black] {};
\node (j) at (4,7)  [shape=circle, draw, fill=black] {};
\node (k) at (4,1.1)  [shape=circle, draw, fill=black] {};
\node (l) at (5,5.1)  [shape=circle, draw, fill=black] {};

\draw [thick] (c) to (z);

\node at (-0.5,4.7) {$e$};
\node at (-3.5,4.5) {$l_0$};
\node at (1.5,6.5) {$l_1$};
\node at (-2.1,2.5) {$l_2$};
\node at (2,4.3) {$l_3$};

\end{tikzpicture}

\end{center}
\caption{Four reflections generating $G$} 
\label{redknight}
\end{figure}

Group-theoretically, this example arises from an epimorphism
\[\theta: N(3)\to G={\rm Aut}\,{\mathcal C}\cong S_4\times C_2\]
sending each of the four involutions $S_i\;(i=0,\ldots, 3)$ generating $N(3)\cong C_2*C_2*C_2*C_2$ to the reflection $s_i$ of $\mathcal C$ fixing the line $l_i$ shown in red in Figure~\ref{redknight}: these send the indicated edge $e$ of $\M$ to the four edges sharing a vertex and a face of $\mathcal M$ with $e$. The vertex-valencies $n_i=4$ and $6$ of $\mathcal M$ are twice the orders of the rotations $s_0s_2=(r_1r_2)^2$ and $s_1s_3=r_0(r_1r_2)^2r_0$, while the face-valencies $m_i=6$ and $8$ are twice the orders of $s_1s_0=(r_0r_1)^2$ and $s_3s_2=r_2(r_0r_1)^2r_2$.

\medskip

\noindent{\bf Example 2} More generally, if $\mathcal C$ is any regular map of type $\{p,q\}$ with automorphism group $G$, where $(p,3)=(q,2)=1$, then a similar knight's move construction based on $\mathcal C$ gives an edge-transitive map $\M$ in class $3$, with automorphism group $G$. It has face-valencies $m_i=2p$ and $2q$, and vertex valencies $n_i=4$ and $2q$, so $m^*=4pq/(p+q)$ and $n^*=8q/(q+2)$; it follows that if $G$ is finite then $\M$ has characteristic
\[\chi=|G|\frac{(4p+2q-3pq)}{4pq}.\]
This map $\M$ is orientable if and only if $\mathcal C$ is.

For instance, the regular embedding $\mathcal C$ of Petersen's graph in the real projective plane (the antipodal quotient of the dodecahedral map) is a regular map of characteristic $1$ and type $\{5,3\}$, so applying this construction gives a non-orientable map $\M$ of characteristic $-19$ in class $3$, with automorphism group $A_5$.
Similarly, taking $\mathcal C$ to be Klein's orientable regular map of genus $3$ and type $\{7,3\}$ gives an orientable map $\M$ of characteristic $-116$ and genus $59$ in class $3$, with automorphism group $L_2(7)\times C_2$.

%Section on Hurwitz-type bounds omitted

\section{Maps with non-empty boundary}

\subsection{Void, tame and wild classes}

In this section we will consider maps with non-empty boundary in the 14 edge-transitive classes $T$. This generalises the discussion in~\cite{JonHSymm}, where the emphasis was on regular maps with a brief disgression concerning those in class~$3$.

A map $\M$, corresponding to a conjugacy class of subgroups $M\le\Gamma$, has non-empty boundary if and only if $M$ contains a reflection, that is, a conjugate of $R_i$ for some $i=0, 1$ or $2$. In particular, if $\M$ is edge-transitive, in class $T$, then it has non-empty boundary if and only if $M$ contains a reflection in $N(T)$. In the presentations of the six parent groups $N(T)$ in Proposition~4.2, the generators denoted by $R_i$ or $S_i$ are representatives of the conjugacy classes of reflections in $N(T)$. From this it is easy to find representatives for these classes in the remaining eight cases. As shown in more detail in~\cite[\S5]{Jonrefl}, the information on the number $c={\rm cr}(N(T))$ of conjugacy classes of reflections in each parent group $N(T)$, and on representatives for those classes, is as in Table~\ref{conjclasses}.

\begin{table}[ht]
\centering
\begin{tabular}{| p{1.8cm} | p{0.5cm} | p{8cm} |}
\hline
Class $T$ & $c$ & representatives \\
\hline\hline
$1$ & $3$ & $R_0, R_1, R_2$  \\
\hline
$2$ & $3$ & $S_1=R_1$, $S_2=R_1^{R_0}$, $S_3=R_2$ \\
\hline
$2^*$ & $3$ & $R_0$, $R_1$, $R_1^{R_2}$ \\
\hline
$2^P$ & $2$ & $R_1$, $R_1^{R_0}$ \\
\hline
$2\,{\rm ex}$ & $1$ & $S_1=R_2$ \\
\hline
$2^*{\rm ex}$ & $1$ & $R_0$ \\
\hline
$2^P{\rm ex}$ & $0$ & \\
\hline
$3$ & $4$ & $S_0=R_1$, $S_1=R_1^{R_0}$, $S_2=R_1^{R_2}$, $S_3=R_1^{R_0R_2}$ \\
\hline
$4, 4^*, 4^P$ & $2$ & $S_1=R_1$, and $S_2=R_1^{R_2}$, $R_1^{R_0}$ or $R_1^{R_0R_2}$ resp. \\
\hline
$5, 5^*, 5^P$ & $0$ & \\
\hline
 
\end{tabular}
\caption{Representatives of conjugacy classes of reflections in $N(T)$.}
\label{conjclasses}
\end{table}

% {\color{blue}[Definition of $S_2$ in~\cite{Jonrefl} incorrect for $T=4^{\sigma}$?]}

\iffalse
\begin{itemize}

\item Class 1: $c=3$, with representatives $R_0, R_1$ and $R_2$.

\item  Class 2: $c=3$, with representatives $S_1=R_1, S_2=R_1^{R_0}$ and $S_3=R_2$.

\item  Class $2^*$: $c=3$, with representatives $R_0, R_1$ and $R_1^{R_2}$.

\item Class $2^P$: $c=2$, with representatives $R_1$ and $R_1^{R_0}$.

\item  Class $2\,{\rm ex}$: $c=1$, with representative $S_1=R_2$.

\item Class $2^*{\rm ex}$: $c=1$, with representative $R_0$.

\item Class $2^P{\rm ex}$: $c=0$.

\item  Class $3$: $c=4$; representatives $S_0=R_1, S_1=R_1^{R_0}, S_2=R_1^{R_2}$ and $S_3=R_1^{R_0R_2}$. 

\item Classes $4$, $4^*$ and $4^P$: $c=2$, with representatives $S_1=R_1$ and $S_2=R_1^{R_2}$, $R_1^{R_0}$ or $R_1^{R_0R_2}$ respectively.

\item Classes $5^{\sigma}$: $c=0$.

\end{itemize}

\fi

%\medskip

Since $M$ is a normal subgroup of $N(T)$, whenever it contains a reflection $R$ in $N(T)$, it contains the corresponding conjugacy class $\mathcal R$ of reflections in $N(T)$; thus $\M$ is a quotient of the map corresponding to the subgroup $\langle{\mathcal R}\rangle$ of $N(T)$ generated by $\mathcal R$, and similarly its automorphism group $G$ is a quotient of the group $N(T)/\langle{\mathcal R}\rangle$ formed by adding the relation $R=1$ to the standard presentation for $N(T)$.

We will say that an edge-transitive class $T$ is {\em void\/} if it contains no maps with non-empty boundary, that it is {\em tame} if it does contain such maps, but they all have dihedral automorphism groups (here we include $D_1\cong C_2$, $D_2\cong V_4$, and $D_{\infty}\cong C_2*C_2$), and that it is {\em wild} if it contains maps with non-dihedral automorphism groups. Our aim is to prove the following theorem, determining which classes are of these three types, to classify the maps and automorphism groups in the tame classes, and to justify using the term `wild' for the third type of class.

\begin{thm}\label{voidtamewild}
Of the $14$ classes of edge-transitive maps,
\begin{itemize}
\item six, namely $2\,{\rm ex}$, $2^*${\rm ex}, $2^P${\rm ex}, $5$, $5^*$ and $5^P$, are void,
\item four, namely $1$, $2$, $2^*$ and $2^P$, are tame, and
\item four, namely $3$, $4$, $4^*$ and $4^P$, are wild.
\end{itemize}
\end{thm}

First of all, it is clear that the classes $T=2^P{\rm ex}$ and $5^{\sigma}$ are all void, since in these cases $N(T)$ contains no reflections, so we may set these four classes aside. We will prove the theorem by considering the remaining ten classes in turn in the subsequent sections. Inspection of Figure~\ref{basicmaps} shows that if ${\mathcal N}(T)$ has non-empty then it is either the closed disc or the M\"obius band, each with a single boundary component. Since $\M/{\rm Aut}\,\M\cong {\mathcal N}(T)$ for each map $\M$ in the class $T$, we immediately have the following:

\begin{lemma}\label{bdytrans}
 If $\M$ is any edge-transitive map with non-empty boundary, then ${\rm Aut}\,\M$ acts transitively on the boundary components of $\M$. \hfill$\square$
 \end{lemma}
 
Thus, in order to understand the structure of the boundary components of an edge-transitive map $\M$, it is sufficient to consider just one of them. Moreover, in order to understand how ${\rm Aut}\,\M$ permutes them, it is sufficient to determine the subgroup leaving that component invariant, since the action is isomorphic to the action on the cosets of that subgroup. This strategy is not so important in the case of the tame classes, where it is easy to classify the relevant maps and automorphism groups by direct arguments, but it is very useful in the case of the wild classes, where no such explicit classifications are feasible.

%%%%%%%%

\subsection{Maps in class $1$}\label{class1bdy}

The maps in class~$1$ are the regular maps. Those with non-empty boundary have been classified in~\cite{JonHSymm}. These maps $\M$ correspond to the normal subgroups $M$ of $\Gamma$ containing at least one of the three generating reflections $R_i$. This implies that the group $G={\rm Aut}\,\M\cong\Gamma/M$ is dihedral, generated by the images $r_j$ of the other two generators $R_j$ ($j\ne i$) of $\Gamma$. It is now straightforward to classify the pairs of involutions generating a finite or infinite dihedral group, and to construct the corresponding regular maps. In all cases the underlying surface is $\overline{\mathbb D}$. The finite maps are as follows.

\begin{itemize}

\item an infinite family $\{{\mathcal A}_n\mid n\ge 1\}$, each embedding a semi-star of valency $n$, with $n$ boundary faces and automorphism group $D_n$ (see Figure~\ref{mapA4} for ${\mathcal A}_4$),

\item an infinite family $\{{\mathcal C}_n=D({\mathcal A}_n)\mid n\ge 1\}$, each embedding a circuit of $n$ vertices and $n$ edges around the boundary $S^1$ of $\overline{\mathbb D}$, with one face and automorphism group $D_n$ (see Figure~\ref{mapA4} for ${\mathcal C}_4$),

\item an embedding ${\mathcal B}\cong{\mathcal N}(3)$ of $K_2$ as a diameter of  $\overline{\mathbb D}$, with two boundary faces and automorphism group $V_4$ (see Figure~\ref{basicmaps}),

\item the quotients of $\mathcal B$ by the three subgroups of ${\rm Aut}\,{\mathcal B}$ of order $2$, isomorphic to ${\mathcal N}(2)$, ${\mathcal N}(2^*)$ and ${\mathcal N}(2^P)$, each with automorphism group $D_1$ (see Figure~\ref{basicmaps}), 

\item an embedding ${\mathcal D}\cong{\mathcal N}(2_{\{0,2\}})$ of a vertex and two half-edges along the boundary of $\overline{\mathbb D}$, with automorphism group $D_1$ (see Figures~\ref{mapA4} and \ref{mapN0}),

\item the trivial map ${\mathcal T}\cong{\mathcal N}(1)$, embedding a half-edge along the boundary of  $\overline{\mathbb D}$, with automorphism group $C_1$ (see Figure~\ref{basicmaps}).

\end{itemize}

\begin{figure}[h!]
\begin{center}
\begin{tikzpicture}[scale=0.5, inner sep=0.8mm]

\node (a) at (0,3.2) {};
\node (b) at (3.2,0) {};
\node (c) at (-3.2,0) {};
\node (d) at (0,-3.2) {};
\node (e) at (0,0) [shape=circle, fill=black] {};

\draw [thick] (a) to (e);
\draw[thick]  (b) to (e);
\draw [thick] (c) to (e);
\draw [thick] (d) to (e);

\draw [dashed] (3,0) arc (0:360:3);

%%%%%

\node (A) at (10,3) [shape=circle, fill=black] {};
\node (B) at (13,0) [shape=circle, fill=black] {};
\node (C) at (7,0) [shape=circle, fill=black] {};
\node (D) at (10,-3) [shape=circle, fill=black] {};

\draw [thick] (13,0) arc (0:360:3);

%%%%%

\node (E) at (23,0) [shape=circle, fill=black] {};
\draw [dashed] (23,0) arc (0:360:3);
\draw [thick] (23,0) arc (0:90:3);
\draw [thick] (23,0) arc (0:-90:3);

\end{tikzpicture}

\end{center}
\caption{The maps ${\mathcal A}_4$, ${\mathcal C}_4$ and $\mathcal D$} 
\label{mapA4}
\end{figure}

%%%%%%%

In the infinite case we have the following:

\begin{itemize}

\item an embedding ${\mathcal A}_{\infty}$ of a semi-star of countably infinite valency, with one boundary vertex and infinitely many semi-edges joining it to boundary points which accumulate at the vertex; this map is shown in Figure~\ref{mapAinfty}, where the semi-edges are the images of the geodesics ${\rm Re}\,z=n\in\Z$ under the M\"obius transformation $f:{\mathbb H}\to{\mathbb D}$, $z\mapsto (z-i)/(z+i)$, so the vertex is at $f(\infty)=1$, and the semi-edges meet $S^1$ at the points $f(n)=(n-i)/(n+i)$ for $n\in\Z$;

\item an embedding ${\mathcal C}_{\infty}$ of a doubly-infinite path along the boundary of $\overline{\mathbb D}$, with vertices accumulating at a single boundary point; this embedding, isomorphic to the dual of ${\mathcal A}_{\infty}$, is shown in Figure~\ref{mapCinfty} as the image under $f$ of the path along $\mathbb R$ with vertices at the integers.

\end{itemize}

%%%%%%

\begin{figure}[h!]
\begin{center}
\begin{tikzpicture}[scale=0.5, inner sep=0.8mm]

\node (a) at (5,0) [shape=circle, fill=black] {};
\node (b) at (-5.2,0) {};

\draw [dashed] (5,0) arc (0:360:5);

\draw [thick] (a) to (b);
\draw [thick] (5,0) arc (270:180:5);
\draw [thick] (5,0) arc (90:180:5);
\draw [thick] (5,0) arc (270:140:2.4);
\draw [thick] (5,0) arc (90:220:2.4);
\draw [thick] (5,0) arc (270:125:1.7);
\draw [thick] (5,0) arc (90:235:1.7);
\draw [thick] (5,0) arc (270:118:1.3);
\draw [thick] (5,0) arc (90:242:1.3);

\node at (5.5,0) {$1$};
\node at (-6,0) {$-1$};
\node at (0,6) {$i$};
\node at (0,-6) {$-i$};
\node at (5.3,4.5) {$(3+4i)/5$};
\node at (5.3,-4.5) {$(3-4i)/5$};
\node at (6.2,3.5) {$(4+3i)/5$};
\node at (6.2,-3.5) {$(4-3i)/5$};
\node at (7.1,2.6) {$(15+8i)/17$};
\node at (7.1,-2.6) {$(15-8i)/17$};

\end{tikzpicture}

\end{center}
\caption{The map ${\mathcal A}_{\infty}$} 
\label{mapAinfty}
\end{figure}

%%%%%%%%

\begin{figure}[h!]
\begin{center}
\begin{tikzpicture}[scale=0.5, inner sep=0.8mm]

\node at (-5,0) [shape=circle, fill=black] {};
\node at (0,5) [shape=circle, fill=black] {};
\node at (3,4) [shape=circle, fill=black] {};
\node at (4,3) [shape=circle, fill=black] {};
\node at (4.41,2.36) [shape=circle, fill=black] {};
\node at (0,-5) [shape=circle, fill=black] {};
\node at (3,-4) [shape=circle, fill=black] {};
\node at (4,-3) [shape=circle, fill=black] {};
\node at (4.41,-2.36) [shape=circle, fill=black] {};

\draw [thick] (5,0) arc (0:360:5);

\node at (5.5,0) {$1$};
\node at (-6,0) {$-1$};
\node at (0,6) {$i$};
\node at (0,-6) {$-i$};
\node at (5.3,4.5) {$(3+4i)/5$};
\node at (5.3,-4.5) {$(3-4i)/5$};
\node at (6.2,3.5) {$(4+3i)/5$};
\node at (6.2,-3.5) {$(4-3i)/5$};
\node at (7.1,2.6) {$(15+8i)/17$};
\node at (7.1,-2.6) {$(15-8i)/17$};
\node at (7.1,1.5) {$\vdots$};
\node at (7.1,-1) {$\vdots$};

\end{tikzpicture}

\end{center}
\caption{The embedding ${\mathcal C}_{\infty}$}
\label{mapCinfty} 
\end{figure}

%%%%%%%%%%%%%%

\subsection{Maps in class~$2$}\label{class2bdy}

The parent group $N=N(2)$ for class $T=2$ has a presentation
\[N=\langle S_1=R_1, S_2=R_1^{R_0}, S_3=R_2\mid S_1^2=S_2^2=S_3^2=1\rangle.\]
It is isomorphic to $C_2*C_2*C_2$, and it has three conjugacy classes ${\mathcal R}_i$ of reflections, represented by the elements $R=S_i$ for $i=1, 2, 3$.

Adding the relation $S_1=1$ to the above presentation for $N$ gives a quotient group
\[N_1=\langle S_2, S_3\mid S_2^2=S_3^2=1\rangle\cong C_2*C_2\cong D_{\infty}.\]
The proper normal subgroups of this are the identity subgroup, the normal closure of $(S_2S_3)^n$, with quotient group $A\cong D_n$, for each integer $n\ge 1$, and two more subgroups of index $2$, the normal closures of $S_2$ and $S_3$, all lifting to proper normal subgroups $M$ of $N$. Conjugation by $R_0$ transposes the generators $S_1$ and $S_2$ of $N$, while fixing $S_3$, so the only one of these subgroups $M$ which is invariant under $R_0$, and thus normal in $\Gamma$, is the normal closure of $\{S_1, S_2\}$; all others have normaliser $N$, so they correspond to maps $\mathcal M$ in class~$2$. The maps ${\mathcal M}_{2,n}$ and ${\mathcal M}_{2,\infty}$ with automorphism group $A\cong D_n$ or $D_{\infty}$ are on the closed disc, with a central vertex and either $n$ or infinitely many interior edges to vertices on the boundary. Figure~\ref{mapM2,4} shows ${\mathcal M}_{2,4}$, while ${\M}_{2,\infty}$ can be obtained from $\mathcal C_{\infty}$, shown in Fig~\ref{mapCinfty}, by adding a vertex at $0$ and replacing all the edges with radii from $0$ to the other vertices.

\begin{figure}[h!]
\begin{center}
\begin{tikzpicture}[scale=0.5, inner sep=0.8mm]

\node (a) at (0,3) [shape=circle, fill=black] {};
\node (b) at (0,-3) [shape=circle, fill=black] {};
\node (c) at (-3,0) [shape=circle, fill=black] {};
\node (d) at (3,0) [shape=circle, fill=black] {};
\node (e) at (0,0) [shape=circle, fill=black] {};

\draw [dashed] (3,0) arc (0:360:3);
\draw [thick] (a) to (b);
\draw [thick] (c) to (d);

\end{tikzpicture}

\end{center}
\caption{The map ${\mathcal M}_{2,4}$} 
\label{mapM2,4}
\end{figure}

The normal closure of $\{S_1, S_3\}$ corresponds to a map ${\mathcal M}_{2,0}$ embedding a path of three vertices and two edges along the boundary of the closed disc (see Figure~\ref{mapM2,0}).
 
 \begin{figure}[h!]
\begin{center}
\begin{tikzpicture}[scale=0.5, inner sep=0.8mm]

\node (b) at (0,3) [shape=circle, fill=black] {};
\node (c) at (0,-3) [shape=circle, fill=black] {};
\node (d) at (3,0) [shape=circle, fill=black] {};

\draw [thick] (0,-3) arc (-90:90:3);
\draw [dashed] (0,3) arc (90:270:3);

\end{tikzpicture}

\end{center}
\caption{The map ${\mathcal M}_{2,0}$} 
\label{mapM2,0}
\end{figure}

Since conjugation by $R_0$ transposes $S_1$ and $S_2$, putting $S_2=1$ simply gives conjugate subgroups $M$ to those arising above, so the same maps arise.

If we add the relation $S_3=1$ to the presentation for $N$ we again get a quotient $N_3\cong D_{\infty}$, but now all the induced proper normal subgroups $M$ of $N$ are invariant under $R_0$, and hence normal in $\Gamma$, with the exception of the normal closures of $\{S_1, S_3\}$ and $\{S_2, S_3\}$, corresponding to the map ${\mathcal M}_{2,0}$. This shows that class~$2$ is tame.

%%%%%%%%%%%%%%%%

\subsection{Maps in class~$2^*$}\label{class2*bdy}

The maps in class~$2^*$ are the duals $D({\mathcal M})$ of the maps $\mathcal M$ in class~$2$ already described. Their map subgroups can be obtained from those for class~$2$ by applying the automorphism $\delta$ of $\Gamma$ which fixes $R_1$ and transposes $R_0$ and $R_2$. Thus this class is also tame.

For each $n\ge 1$, $D({\mathcal M}_{2,n})$ is a map ${\mathcal M}_{2^*,n}$ which embeds a circuit of $n$ boundary vertices and $n$ interior edges in the closed disc $\overline{\mathbb D}$ (see Figure~\ref{mapM2*,4} for ${\mathcal M}_{2^*,4}$). Similarly, ${\mathcal M}_{2^*,\infty}$ embeds a doubly infinite path in $\overline{\mathbb D}$ , with interior edges and boundary vertices accumulating at a boundary point; this map can be obtained from ${\mathcal C}_{\infty}$, shown in Figure~\ref{mapCinfty}, by replacing its edges, along the boundary, with straight line segments.

\begin{figure}[h!]
\begin{center}
\begin{tikzpicture}[scale=0.5, inner sep=0.8mm]

\node (a) at (2.1,2.1) [shape=circle, fill=black] {};
\node (b) at (2.1,-2.1) [shape=circle, fill=black] {};
\node (c) at (-2.1,-2.1) [shape=circle, fill=black] {};
\node (d) at (-2.1,2.1) [shape=circle, fill=black] {};

\draw [dashed] (3,0) arc (0:360:3);
\draw [thick] (a) to (b) to (c) to (d) to (a);

\end{tikzpicture}

\end{center}
\caption{The map ${\mathcal M}_{2^*,4}$} 
\label{mapM2*,4}
\end{figure}

The map ${\mathcal M}_{2^*,0}=D({\mathcal M}_{2,0})$ has a single vertex on the boundary of $\overline{\mathbb D}$, and two semi-edges across the interior to two boundary points (see Figure~\ref{mapM2*,0}).

\begin{figure}[h!]
\begin{center}
\begin{tikzpicture}[scale=0.5, inner sep=0.8mm]

\node (a) at (3,0) [shape=circle, fill=black] {};

\draw [dashed] (3,0) arc (0:360:3);
\draw [thick] (a) to (-2.15,2.15);
\draw [thick] (a) to (-2.15,-2.15);

\end{tikzpicture}

\end{center}
\caption{The map ${\mathcal M}_{2^*,0}$} 
\label{mapM2*,0}
\end{figure}

Maps in class~$2^*$ have also been studied in~\cite[\S5]{Jon15}, as part of a classification of arc-transitive maps with non-empty boundary. (A map is {\em arc-transitive\/} if its automorphism group acts transitively on arcs, that is, on directed edges; this is equivalent to the condition that $\Gamma=N\langle R_2\rangle$, so the arc-transitive maps are those edge-transitive maps in class $1$, $2^*$, $2^P$, $2^*{\rm ex}$ or $2^P{\rm ex}$.) The maps ${\mathcal M}_{2^*,0}$,  ${\mathcal M}_{2^*,n}$ and ${\mathcal M}_{2^*,\infty}$ are denoted in~\cite{Jon15} by $\mathcal E$, ${\mathcal F}_n$ and ${\mathcal F}_{\infty}$.

%%%%%%%%%%%%%%%%%%

\subsection{Maps in class~$2^P$}\label{class2Pbdy}

The group $N(2^P)$ has a presentation
\[N=\langle S_1=R_1, S_2=R_1^{R_0}, S=R_0R_2\mid S_1^2=S_2^2=S^2=1\rangle,\]
so that $N\cong C_2*C_2*C_2$.
There are two conjugacy classes of reflections in $N$, represented by $S_1$ and $S_2$. Like the maps in class~$2^*$, those in class~$2^P$ have been discussed in~\cite[\S5]{Jon15}.  The following maps arise, showing that this class is also tame.

There is a map ${\mathcal M}_{2^P,0}$ on the closed disc, denoted in~\cite{Jon15} by $\mathcal G$; this has a single boundary vertex, two interior semi-edges, and a single face. Its automorphism group is generated by a reflection, transposing the semi-edges.

\begin{figure}[h!]
\begin{center}
\begin{tikzpicture}[scale=0.5, inner sep=0.8mm]

\node (a) at (3,0) [shape=circle, fill=black] {};
\node (b) at (0,1) {};
\node (c) at (0,-1) {};

\draw [thick] (a) to (b);
\draw [thick] (a) to (c);

\draw [dashed] (3,0) arc (0:360:3);
\end{tikzpicture}

\end{center}
\caption{The map ${\mathcal M}_{2^P,0}$.} 
\label{mapM2p,0}
\end{figure}

There is an infinite series of maps ${\mathcal M}_{2^P,n}$ ($n\ge 1$), denoted by ${\mathcal H}_n$ in~\cite{Jon15}, with automorphism group $D_n$; these embed a circuit of $n$ boundary vertices and interior edges in an annulus or M\"obius band as $n$ is even or odd. (See Figure~\ref{mapM2P,n}, where the left and right sides of the rectangular strip are identified orientably or non-orientably in these two cases.)  In Figure~\ref{map2P,infty} the pattern is the same, except that the strip now extends infinitely far in both directions; the resulting map ${\mathcal M}_{2^P,\infty}={\mathcal H}_{\infty}$ embeds an infinite path with vertices on alternate boundary components, and edges in the interior. The automorphism group $D_{\infty}$ acts as the frieze group~${\rm p2mg}$: the cyclic subgroup of index $2$ is generated by a glide reflection, and the involutions are either reflections or half-turns.

\begin{figure}[h!]
\begin{center} \begin{tikzpicture}[scale=0.5, inner sep=0.8mm]

\node (a) at (-6,2) [shape=circle, fill=black] {};
\node (b) at (-3,-2) [shape=circle, fill=black] {};
\node (c) at (0,2) [shape=circle, fill=black] {};
\node (d) at (3,-2) [shape=circle, fill=black] {};
\node (e) at (6,2) [shape=circle, fill=black] {};

\node (A) at (-7,0.67) {};
\node (B) at (-8,-0.67) {};
\draw [thick] (-6,2) to (-7,0.67);
\draw [thick] [dotted] (A) to (B);

\node (C) at (7,0.67) {};
\node (D) at (8,-0.67) {};
\draw [thick] (6,2) to (7,0.67);
\draw [thick] [dotted] (C) to (D);

\draw [thick] (a) to (b) to (c) to (d) to (e);
\draw [dashed] (-12,2) to (12,2);
\draw [dashed] (-12,-2) to (12,-2);
\draw [thick] [dotted] (-12,2) to (-12,-2);
\draw [thick] [dotted] (12,2) to (12,-2);

\end{tikzpicture}

\end{center}
\caption{The map ${\mathcal M}_{2^P,n}$ (identify left and right sides).} 
\label{mapM2P,n}
\end{figure}

%%%%%%%

\begin{figure}[h!]
\begin{center} \begin{tikzpicture}[scale=0.5, inner sep=0.8mm]

\node (a) at (-6,2) [shape=circle, fill=black] {};
\node (b) at (-3,-2) [shape=circle, fill=black] {};
\node (c) at (0,2) [shape=circle, fill=black] {};
\node (d) at (3,-2) [shape=circle, fill=black] {};
\node (e) at (6,2) [shape=circle, fill=black] {};

\node (A) at (-7,0.67) {};
\node (B) at (-8,-0.67) {};
\draw [thick] (-6,2) to (-7,0.67);
\draw [thick] [dotted] (A) to (B);

\node (C) at (7,0.67) {};
\node (D) at (8,-0.67) {};
\draw [thick] (6,2) to (7,0.67);
\draw [thick] [dotted] (C) to (D);

\draw [thick] (a) to (b) to (c) to (d) to (e);
\draw [dashed] (-12,2) to (12,2);
\draw [dashed] (-12,-2) to (12,-2);

\end{tikzpicture}

\end{center}
\caption{The map ${\mathcal M}_{2^P,\infty}$.} 
\label{map2P,infty}
\end{figure}

%%%%%%%%%%%%%%%%%%%%%%

\subsection{Maps in classes $2\,{\rm ex}$ and $2^*{\rm ex}$}

The group $N(2\,{\rm ex})$ has a presentation
\[N=\langle S_1=R_2, S=R_0R_1\mid S_1^2=1\rangle\cong C_2*C_{\infty}.\]
The reflections in $N$ form a single conjugacy class, represented by $S_1$, and putting $S_1=1$ gives a quotient $N_1\cong C_{\infty}$. Now $R_0$ acts by conjugation on $N$ by centralising $S_1$ and inverting $S$, so all the map subgroups $M$ obtained from quotients of $N_1$ are normal in $\Gamma$; the corresponding maps $\mathcal M$ are therefore all regular, and not in class~$2\,{\rm ex}$.
(They are, in fact, the maps ${\mathcal C}_n$ and ${\mathcal C}_{\infty}$ described in \S\ref{class1bdy}.) Thus there are no edge-transitive maps with non-empty boundary in class $2\,{\rm ex}$. By the duality of maps, the same applies to the class $2^*{\rm ex}$, so both of these classes are void.

%%%%%%%%%%%%%%%%%%%%%%

\subsection{Summary for void and tame classes}

As we will shortly prove, the four remaining classes $T=3$ and $4^{\sigma}$ are wild, so here we pause briefly to summarise the results obtained so far, listing all the maps with non-empty boundary in the four tame classes (recall that the classes $2^{\sigma}{\rm ex}$ and $5^{\sigma}$ are void, contributing no examples).

\begin{thm}\label{tamemaps}
The maps with non-empty boundary in the tame classes $T=1$ and $2^{\sigma}$ are as follows:
\begin{itemize}
\item $T=1$: ${\mathcal A}_n$ and ${\mathcal C}_n=D({\mathcal A}_n)$ for $n\in\N\cup\{\infty\}$, ${\mathcal B}={\mathcal N}(3)$, ${\mathcal B}/C_2\cong{\mathcal N}(2^{\sigma})$ ($\sigma=\emptyset, *, P$), $\mathcal D$, ${\mathcal T}={\mathcal N}(1)$;
\item $T=2$: ${\mathcal M}_{2,n}$ for $n\in\N\cup\{\infty\}$, and ${\mathcal M}_{2,0}$;
\item $T=2^*$: ${\mathcal M}_{2^*,n}=D({\mathcal M}_{2,n})$ for $n\in\N\cup\{\infty\}$, and ${\mathcal M}_{2^*,0}=D({\mathcal M}_{2,0})$;
\item $T=2^P$: ${\mathcal M}_{2^P,n}=P({\mathcal M}_{2^*,n})$ for $n\in\N\cup\{\infty\}$, and ${\mathcal M}_{2^P,0}=P({\mathcal M}_{2^*,0})$.
\end{itemize}
All are on the closed disc $\overline{\mathbb D}$, apart from ${\mathcal M}_{2^P,n}$ on the annulus or M\"obius band for even or odd $n\ge 1$, and on the doubly infinite strip for $n=\infty$.
\end{thm}

%%%%%%%%%%%%%%%%%%%%%%

\subsection{Maps in class $3$}\label{class3bdy}

The group $N(3)$ has a presentation
\[N=\langle S_0:=R_1,\, S_1:=R_1^{R_0},\, S_2:=R_1^{R_2},\, S_3:=R_1^{R_0R_2}\mid S_i^2=1\rangle,\]
so $N\cong C_2*C_2*C_2*C_2$. There are four conjugacy classes of reflections in $N$, represented by $S_0,\ldots, S_3$, so a map $\M$ in this class has non-empty boundary if and only if $S_i\in M$ for some $i$. Moreover, since $\Gamma$ is a semirect product of $N$ by $E$, it follows that $\M$ is in class $3$ if and only if only the identity element of $E$ normalises $M$.
 
If $S_i\in M$ for every $i=0,\ldots, 3$ then $M=\langle S_0,\ldots, S_3\rangle=N$; thus $\M=\mathcal N(3)$, which is regular, against our assumption, so this case does not arise.
 
If $S_i\in M$ for three values of $i$ then only the identity element of $E$ induces an automorphism of $G$, so $N=N_{\Gamma}(M)$ as required. There are four possible subgroups $M$, each of index $2$ in $N$, depending on which reflections $S_i$ are in $M$; these subgroups form a single conjugacy class in $\Gamma$, corresponding to a single map ${\mathcal M}_{3,2}$ which is a double covering of $\mathcal N(3)$ (see Figure~\ref{mapM3,2}). Clearly ${\rm Aut}\,\M_{3,2}\cong D_1\cong C_2$.
 
  \begin{figure}[h!]
\begin{center}
\begin{tikzpicture}[scale=0.5, inner sep=0.8mm]

\node (b) at (-2,2.3) [shape=circle, fill=black] {};
\node (c) at (-2,-2.3) [shape=circle, fill=black] {};
\node (d) at (3,0) [shape=circle, fill=black] {};
\draw [thick] (b) to (d) to (c);

\draw [dashed] (3,0) arc (0:360:3);

\end{tikzpicture}

\end{center}
\caption{The map ${\mathcal M}_{3,2}$.} 
\label{mapM3,2}
\end{figure}

If $S_i\in M$ for two values of $i$ then the involution in $E$ transposing them (and also transposing the other pair of generators) gives a forbidden automorphism of $G$, so that $\mathcal M$ is not in class~$3$. (By Lemma~\ref{covering}(3) it is in class~$1$ or $2^{\sigma}$ for some $\sigma$.)

Finally suppose that $S_i\in M$ for just one $i$, so that $N_{\Gamma}(M)=N$ as required. The corresponding map subgroups $M$ are the inverse images in $N$ of the normal subgroups of $N_i:=N/K_i$, where $K_i$ is the normal closure of $S_i$ in $N$. Since $N_i\cong C_2*C_2*C_2$, the corresponding automorphism groups $G$ are those groups which can be generated by three involutions.

This shows that class~$3$ is wild since there are (very many) non-dihedral groups $G$ with this property. Indeed, $N_i$ is isomorphic to the parent group $N(2)$ for the classes $2^{\sigma}$, and results such as Theorems 1.1 to 1.6 show how abundant its quotient groups are. Here, without the restriction of avoiding forbidden automorphisms, the sets of possible automorphism groups and maps are even larger, and far too large for us to expect a reasonable classification. For instance, every non-abelian finite simple group except $U_3(3)$ is a quotient of $N_i$ (see Theorem~\ref{mainthmsimple} and~\cite{MSW}).

To give an alternative perspective on the difficulty of this classification problem, the parent group for all hypermaps is also isomorphic to $C_2*C_2*C_2$, so (ignoring ${\mathcal M}_{3,2}$) the problem of classifying maps with boundary in class~$3$ is equivalent that of classifying all regular hypermaps. The following construction (which generalises one described in~\cite[\S6]{JonHSymm} for regular maps) shows how every regular hypermap $\mathcal H$ without boundary can be converted into a map $\M$ in class $3$ with non-empty boundary and with the same automorphism group. Conversely, every such map $\M$, with the single exception of $\M_{3,2}$, arises in this way.

A regular hypermap $\mathcal H$ of type $(l,m,n)$ can be represented as a three-coloured triangular map $\mathcal T$, with black, white and red vertices of valencies $2l$, $2m$ and $2n$ representing the hypervertices, hyperedges and hyperfaces, and edges indicating incidence. One can convert this into an edge-transitive map $\mathcal M$, with the same automorphism group as $\mathcal H$, as follows. First place a new vertex at the midpoint of each edge between white and red vertices, and join it by two new edges to the two black vertices in the triangles on either side of that edge. Next, around each red vertex remove the interior of an $n$-gon spanned by the adjacent new vertices, leaving this $n$-gon as a boundary component; if $\mathcal H$ is without boundary, these will be the only boundary components of $\M$. Finally, delete all the white vertices and all the edges of $\mathcal T$, leaving just the black vertices, the new vertices, and the new edges. This process is illustrated, with $l=m=n=3$, in Figure~\ref{class3map}, where $\mathcal T$ is shown on the left, the new vertices (in blue) and new edges (dotted lines) are added in the centre, and the final map $\M$ is shown on the right, with broken lines for boundary components, the removed regions in grey, and the edges now unbroken lines.

\begin{figure}[h!]

\begin{center}
 \begin{tikzpicture}[scale=0.5, inner sep=0.8mm]

\node (a) at (-1,3.46) [shape=circle, draw, fill=black] {};
\node (b) at (1,3.46) [shape=circle, draw] {};
\node (c) at (-2,1.73) [shape=circle, draw] {};
\node (d) at (0,1.73) [shape=circle, fill=red] {};
\node (e) at (2,1.73) [shape=circle, draw, fill=black] {};
\node (f) at (-3,0) [shape=circle, draw, fill=red] {};
\node (g) at (-1,0) [shape=circle, draw, fill=black] {};
\node (h) at (1,0) [shape=circle, draw] {};
\node (i) at (3,0) [shape=circle, draw, fill=red] {};
\node (j) at (-2,-1.73) [shape=circle, draw] {};
\node (k) at (0,-1.73) [shape=circle, fill=red] {};
\node (l) at (2,-1.73) [shape=circle, draw, fill=black] {};
\node (m) at (-1,-3.46) [shape=circle, draw, fill=black] {};
\node (n) at (1,-3.46) [shape=circle, draw] {};

\draw (a) to (b);
\draw (c) to (d) to (e);
\draw (f) to (h) to (i);
\draw (j) to (k) to (l);
\draw (m) to (n);

\draw (a) to (c) to (f);
\draw (b) to (d) to (j);
\draw (e) to (h) to (k) to (m);
\draw (i) to (n);

\draw (b) to (i);
\draw (a) to (d) to (h) to (l);
\draw (c) to (k) to (n);
\draw (f) to (j) to (m);

\node at (3,-4) {$\mathcal T$};
\node at (5,0) {$\longrightarrow$};

%%%%%%%%%%%%%%%

\node (a1) at (9,3.46) [shape=circle, draw, fill=black] {};
\node (b1) at (11,3.46) [shape=circle, draw] {};
\node (c1) at (8,1.73) [shape=circle, draw] {};
\node (d1) at (10,1.73) [shape=circle, fill=red] {};
\node (e1) at (12,1.73) [shape=circle, draw, fill=black] {};
\node (f1) at (7,0) [shape=circle, draw, fill=red] {};
\node (g1) at (9,0) [shape=circle, draw, fill=black] {};
\node (h1) at (11,0) [shape=circle, draw] {};
\node (i1) at (13,0) [shape=circle, draw, fill=red] {};
\node (j1) at (8,-1.73) [shape=circle, draw] {};
\node (k1) at (10,-1.73) [shape=circle, fill=red] {};
\node (l1) at (12,-1.73) [shape=circle, draw, fill=black] {};
\node (m1) at (9,-3.46) [shape=circle, draw, fill=black] {};
\node (n1) at (11,-3.46) [shape=circle, draw] {};

\draw (a1) to (b1);
\draw (c1) to (d1) to (e1);
\draw (f1) to (h1) to (i1);
\draw (j1) to (k1) to (l1);
\draw (m1) to (n1);

\draw (a1) to (c1) to (f1);
\draw (b1) to (d1) to (j1);
\draw (e1) to (h1) to (k1) to (m1);
\draw (i1) to (n1);

\draw (b1) to (i1);
\draw (a1) to (d1) to (h1) to (l1);
\draw (c1) to (k1) to (n1);
\draw (f1) to (j1) to (m1);

\node (bd1) at (10.5,2.66) [shape=circle, draw, fill=blue] {};
\node (cd1) at (9,1.73) [shape=circle, draw, fill=blue] {};
\node (dh1) at (10.5,0.86) [shape=circle, draw, fill=blue] {};
\node (cf1) at (7.5,0.86) [shape=circle, draw, fill=blue] {};
\node (hi1) at (12,0) [shape=circle, draw, fill=blue] {};
\node (fj1) at (7.5,-0.86) [shape=circle, draw, fill=blue] {};
\node (hk1) at (10.5,-0.86) [shape=circle, draw, fill=blue] {};
\node (jk1) at (9,-1.73) [shape=circle, draw, fill=blue] {};
\node (kn1) at (10.5,-2.66) [shape=circle, draw, fill=blue] {};

\draw [very thick, dotted] (a1) to (bd1) to (e1);
\draw [very thick, dotted] (a1) to (cd1) to (jk1) to (m1);
\draw [very thick, dotted] (fj1) to (dh1) to (e1);
\draw [very thick, dotted] (cf1) to (l1);
\draw [very thick, dotted] (e1) to (hi1) to (l1);
\draw [very thick, dotted] (l1) to (kn1) to (m1);

%%%%%%%%%%%%%%%%%%%%%

\node at (15,0) {$\longrightarrow$};

\node (a2) at (19,3.46) [shape=circle, draw, fill=black] {};
\node (e2) at (22,1.73) [shape=circle, draw, fill=black] {};
\node (g2) at (19,0) [shape=circle, draw, fill=black] {};
\node (l2) at (22,-1.73) [shape=circle, draw, fill=black] {};
\node (m2) at (19,-3.46) [shape=circle, draw, fill=black] {};

\node (bd2) at (20.5,2.66) [shape=circle, draw, fill=blue] {};
\node (cd2) at (19,1.73) [shape=circle, draw, fill=blue] {};
\node (dh2) at (20.5,0.86) [shape=circle, draw, fill=blue] {};
\node (cf2) at (17.5,0.86) [shape=circle, draw, fill=blue] {};
\node (hi2) at (22,0) [shape=circle, draw, fill=blue] {};
\node (fj2) at (17.5,-0.86) [shape=circle, draw, fill=blue] {};
\node (hk2) at (20.5,-0.86) [shape=circle, draw, fill=blue] {};
\node (jk2) at (19,-1.73) [shape=circle, draw, fill=blue] {};
\node (kn2) at (20.5,-2.66) [shape=circle, draw, fill=blue] {};

\draw [thick] (a2) to (bd2) to (e2);
\draw [thick] (a2) to (cd2) to (jk2) to (m2);
\draw [thick] (fj2) to (dh2) to (e2);
\draw [thick] (cf2) to (l2);
\draw [thick] (e2) to (hi2) to (l2);
\draw [thick] (l2) to (kn2) to (m2);

\draw [fill=lightgray, dashed] (20.45,2.61) -- (19.05,1.73) -- (20.45,0.91) -- (20.45,2.61);
\draw [fill=lightgray, dashed] (20.45,-0.91) -- (19,-1.73) -- (20.45,-2.61) -- (20.45,-0.91);
\draw [dashed] (16.5,0.15) to (cf2) to (fj2) to (16.5,-0.15);
\draw [lightgray, fill=lightgray] (16.5,0.15) to (17.45,0.8) to (17.45,-0.8) to (16.5,-0.15) to (16.5,0.15);
\draw [dashed] (23, 0.6) to (hi2) to (23, -0.6);
\draw [lightgray, fill=lightgray] (23, 0.6) to (22.1,0) to (23, -0.6) to (23,0.6);

\node at (22,-4) {$\M$};

\end{tikzpicture}

\end{center}
\caption{The construction of $\M$ from $\mathcal T$}
\label{class3map} 
\end{figure}

The resulting map $\M$, which is orientable if and only if $\mathcal H$ is, has interior black vertices of valency $2l$, and new boundary vertices of valency $2$.  Its edges, one for each triangle in $\mathcal T$, each connect an interior vertex and a boundary vertex, so $\M$ is bipartite. It has interior faces of valency $2m$ (one for each white vertex of $\mathcal T$), together with triangular boundary faces, $l$ of them incident with each black vertex. Now ${\rm Aut}\,\mathcal H$ acts regularly on the triangles in $\mathcal T$, and hence also on the edges of $\M$, so $\M$ is edge-transitive. Clearly ${\rm Aut}\,\M$ is intransitive on the vertices and on the faces of $\M$, so $\M$ is just-edge-transitive, that is, in class $3$, with ${\rm Aut}\,\M={\rm Aut}\,{\mathcal H}$.  Conversely, every class~$3$ map $\M$ with non-empty boundary (other than $\M_{3,2}$) is obtained in this way from a regular hypermap $\mathcal H$ without boundary: the hypervertices, hyperedges and hyperfaces of $\mathcal H$ correspond to the non-boundary vertices, the non-boundary faces and the boundary components of $\M$, with obvious definitions of incidence.

Further properties of $\M$ can be seen from Figure~\ref{class3map}: for instance, there is a single orbit of ${\rm Aut}\,\M$ on boundary components, one for each red vertex of $\mathcal H$, each component having boundary length (number of segments between incident vertices) equal to $n$. There are two orbits of ${\rm Aut}\,\M$ on Petrie polygons of $\M$: one orbit, corresponding bijectively  to the set of red vertices of $\mathcal H$, consists of paths of length $2n$, each following the edges of $\M$ meeting a particular boundary component; the other orbit, corresponding bijectively to the incident vertex/interior face pairs of $\M$, consists of paths of length $4$, each passing once in each direction along the two edges of $\M$ incident with such a pair. One Petrie polygon of each type is shown in red in Figure~\ref{twoPetries} for the map $\M$ in Figure~\ref{class3map}, where $n=3$.

\begin{figure}[h!]

\begin{center}
 \begin{tikzpicture}[scale=0.7, inner sep=0.8mm]

\node (a2) at (19,3.46) [shape=circle, draw, fill=black] {};
\node (e2) at (22,1.73) [shape=circle, draw, fill=black] {};
\node (g2) at (19,0) [shape=circle, draw, fill=black] {};
\node (l2) at (22,-1.73) [shape=circle, draw, fill=black] {};
\node (m2) at (19,-3.46) [shape=circle, draw, fill=black] {};

\node (bd2) at (20.5,2.66) [shape=circle, draw, fill=blue] {};
\node (cd2) at (19,1.73) [shape=circle, draw, fill=blue] {};
\node (dh2) at (20.5,0.86) [shape=circle, draw, fill=blue] {};
\node (cf2) at (17.5,0.86) [shape=circle, draw, fill=blue] {};
\node (hi2) at (22,0) [shape=circle, draw, fill=blue] {};
\node (fj2) at (17.5,-0.86) [shape=circle, draw, fill=blue] {};
\node (hk2) at (20.5,-0.86) [shape=circle, draw, fill=blue] {};
\node (jk2) at (19,-1.73) [shape=circle, draw, fill=blue] {};
\node (kn2) at (20.5,-2.66) [shape=circle, draw, fill=blue] {};

\draw [thick] (a2) to (bd2) to (e2);
\draw [thick] (a2) to (cd2) to (jk2) to (m2);
\draw [thick] (fj2) to (dh2) to (e2);
\draw [thick] (cf2) to (l2);
\draw [thick] (e2) to (hi2) to (l2);
\draw [thick] (l2) to (kn2) to (m2);

\draw [fill=lightgray, dashed] (20.45,2.61) -- (19.05,1.73) -- (20.45,0.91) -- (20.45,2.61);
\draw [fill=lightgray, dashed] (20.45,-0.91) -- (19,-1.73) -- (20.45,-2.61) -- (20.45,-0.91);
\draw [dashed] (16.5,0.15) to (cf2) to (fj2) to (16.5,-0.15);
\draw [lightgray, fill=lightgray] (16.5,0.15) to (17.45,0.8) to (17.45,-0.8) to (16.5,-0.15) to (16.5,0.15);
\draw [dashed] (23, 0.6) to (hi2) to (23, -0.6);
\draw [lightgray, fill=lightgray] (23, 0.6) to (22.1,0) to (23, -0.6) to (23,0.6);

\node at (23,-3.5) {$\M$};

\draw [thick, red, rounded corners] (21.25,2.25) to (21.25,2.5) to (20.05,3.15) to (19.8,2.8) to (19.2, 3.13)
to (19.2,2.5) to (18.7,2.4) to (18.7,1.1) to (19.2,1) to (19.2,0.3) to (19.8,0.7) to (20.2,0.35) to (21.25,1) to (21.1,1.5) to (21.7,1.73) to (21.15,2) to (21.25,2.25);

\draw [thick, red, rounded corners] (22.2,-0.3) to (22.2,-0.8) to (21.8,-0.85) to (21.8,-1.4) to (21.25,-1.1) to (20.95,-1.4) to (20.65,-1.2);

\end{tikzpicture}

\end{center}
\caption{Two Petrie polygons in $\M$}
\label{twoPetries} 
\end{figure}
\medskip

(There is, in fact, an alternative and possibly simpler construction of $\M$: represent $\mathcal H$ as a bipartite map, with black and white vertices; in each $2n$-gonal face join successive black vertices by chords to form an $n$-gon; remove an open $n$-gon spanned by the mid-points of these new edges; finally delete the white vertices and edges of $\M$.)

%%%%%%%%%%%%%%%%%%%%%

\subsection{Maps in classes $4$, $4^*$ and $4^P$}\label{class4sigmabdy}

The group $N=N(4)$ has a presentation
\[N=\langle S_1=R_1, S_2=R_1^{R_2}, S=(R_1R_2)^{R_0}\mid S_1^2=S_2^2=1\rangle,\]
so that $N\cong C_2*C_2*C_{\infty}$. There are two conjugacy classes of reflections in $N$, represented by $S_1$ and $S_2$. Since $R_2$ transposes these by conjugation, it is sufficient to consider the quotient group $N_1$ obtained by putting $S_1=1$. 

We have $N_1\cong C_2*C_{\infty}$, so the resulting quotients $N/M$ are those groups generated by two elements $s_2$ and $s$, with $s_2$ of order $2$. (Putting $S_1=S_2=1$ results in an infinite cyclic quotient, with all subgroups normalised by $R_2$ and thus not corresponding to maps in class~$4$, so we may assume that $s_2\ne 1$.) As in the case of class~$3$, these quotients  include non-dihedral groups, so class~$4$ is wild and there is no realistic hope of classifying such groups or their associated maps. For instance, $N_1$ is isomorphic to the parent group $N(T)$ for each of the classes $T=2^{\sigma}{\rm ex}$, including the class $2^P{\rm ex}$ of all chiral maps, so the classification problem for class~$4$ maps with boundary is at least as hard as that for any of them. Indeed, any chiral map $\mathcal C$ can be converted into a map $\M$ in class~$4$, with non-empty boundary and with the same automorphism group, by the following process.

%removing a small open disc around the mid-point of each edge, and then adding two new vertices where that edge meets the new boundary component. This is illustrated in Figure~\ref{ChiraltoClass4}, which also shows a fundamental region for ${\rm Aut}\,\M'={\rm Aut}\,\M$ in $\M'$, together with the side-pairings induced by its standard generators, the rotations $x$ and $y$.} 

\begin{figure}[h!]
\begin{center}
\begin{tikzpicture}[scale=0.6, inner sep=0.8mm]

\node (c) at (0,0) [shape=circle, fill=black] {};
\node (d) at (8,0) [shape=circle, fill=black] {};
\draw [thick] (c) to (d);
\draw [thick] (c) to (1,-3);
\draw [thick] (c) to (1,3);
\draw [thick] (d) to (7,-3);
\draw [thick] (d) to (7,3);
\draw [thick] (c) to (-2.5,2.5);
\draw [thick] (c) to (-2.5,-2.5);
\draw [thick] (d) to (10.5,2.5);
\draw [thick] (d) to (10.5,-2.5);

\draw [blue, rounded corners] (1.3,3) to (0.3,0.2) to (3.9,0.2) to (4,0) to (4.1,-0.2) to (7.7,-0.2) to (6.7,-3);
\draw [red, rounded corners] (1.3,-3) to (0.3,-0.2) to (3.9,-0.2) to (4,0) to (4.1,0.2) to (7.7,0.2) to (6.7,3);

\draw (14,0) circle [radius=0.5];
\draw (14.2,-0.2) to (14.5,0.2) to (14.75,-0.2);

\node at (4,-0.6) {$e$};
\node at (14.2,-1) {orientation};
\node at (4,-3) {$\mathcal C$};

\end{tikzpicture}

\end{center}
\caption{Right and left Petrie polygons at the edge $e$ of $\mathcal C$} 
\label{RLPetrie}
\end{figure}

Let $\mathcal C$ be a chiral map, with a chosen orientation of its underlying surface. Each Petrie polygon of $\mathcal C$ represents a path turning first right or first left at alternate vertices, where we define `right' and `left' to mean that the rotation from the incoming edge of the path to the outgoing edge either follows the chosen orientation around the vertex or goes against it. Each edge $e$ of $\mathcal C$ lies on two Petrie polygons of $\mathcal C$, which we will call the right and left Petrie  polygons at $e$ as they enter $e$ (in either direction) through a right or left turn. (See Figure~\ref{RLPetrie}, where the right and left Petrie polygons at $e$ are shown in red and blue.) Similarly, the two faces of $P(\mathcal C)$ incident with $e$ can be called the right face $r_e$ or left face $l_e$ at $e$, depending on their bounding polygons. Thus each face of $P(\mathcal C)$ is alternately the right or left face at its successive edges. We next form the Walsh map $WP(\mathcal C)$ of $P(\mathcal C)$, by inserting a new vertex $v_e$ of valency $2$ at the midpoint of each edge $e$ (of $\mathcal C$ or equivalently of $P(\mathcal C)$); note that, being chiral, $\mathcal C$ has empty boundary, and hence so has $WP(\mathcal C)$. We now form $\M$ by removing an open disc from each face of $WP(\mathcal C)$, so that the resulting boundary component of $\M$ meets incident edges $e$ only at those vertices $v_e$ for which that face is the left face $l_e$ at $e$ (see Figure~\ref{ChiraltoClass4}). Now ${\rm Aut}\,\mathcal C$ acts regularly on the directed edges of $\mathcal C$, and hence also on the edges of $WP(\mathcal C)$, so this map $\M$ is edge-transitive. Clearly ${\rm Aut}\,\mathcal C\le{\rm Aut}\,\M$, and the chirality of $\mathcal C$ guarantees that these are the same group $G$. Figure~\ref{ChiraltoClass4} shows a fundamental region $F$ for $G$ on $\M$, with side-identifications induced by its standard generators $x$ and $y$. (Note that the involution $y$, which induces a half-turn on $\mathcal C$, induces a reflection on $P(\mathcal C)$ and hence on $\M$: as in similar examples in \S\ref{parentforbidden}, this is the reason for first replacing $\mathcal C$ with $P(\mathcal C)$ in this construction.) Applying these identifications to $F$ yields a quotient map isomorphic to ${\mathcal N}(4)$, with $v$ and $v_e$ the interior and boundary vertices, so $\M$ is in class~$4$.

The interior vertices of $\M$ are the vertices of $\mathcal C$, of the same valency, while there is a boundary vertex on each edge of $\mathcal C$; these two sets of vertices are orbits of $G$, permuted by $G$ as it permutes the vertices or edges of $\mathcal C$. The edges of $\M$ are the half-edges of $\mathcal C$, both sets being permuted regularly by $G$. Each face of $\M$ is a pentagon, bounded by four edges of $\M$ and a segment $S$ of a boundary component; the faces correspond bijectively to the boundary vertices of $\M$ (opposite to $S$), and thus to the edges of $\mathcal C$, all permuted by $G$ in the same way. Similarly, the set of boundary components of $\M$ is $G$-isomorphic to the set of Petrie polygons of $\mathcal C$, and the boundary length of $\M$ is half the Petrie length of $\mathcal C$ (the latter being even since $\mathcal C$ is orientable). There is a single orbit of $G$ on Petrie polygons of $\M$: each of these has length $8$, following (once in each direction) a path of length $4$ in $\M$ between two vertices of valency $2$ around a face of $W(\mathcal C)$ (see the path indicated in red in Figure~\ref{ChiraltoClass4}). Apart from its initial and final vertices. such a path passes through a unique third vertex $v_e$ of valency $2$, so the set of Petrie polygons of $\M$ is $G$-isomorphic to the set of edges $e$ of $\mathcal C$.

\begin{figure}[h!]
\begin{center}
\begin{tikzpicture}[scale=0.4, inner sep=0.8mm]

\node (c) at (0,0) [shape=circle, fill=black] {};
\node (d) at (8,0) [shape=circle, fill=black] {};
\draw [thick] (c) to (d);
\draw [thick] (c) to (1,-3);
\draw [thick] (c) to (1,3);
\draw [thick] (d) to (7,-3);
\draw [thick] (d) to (7,3);
\draw [thick] (c) to (-2.5,2.5);
\draw [thick] (c) to (-2.5,-2.5);
\draw [thick] (d) to (10.5,2.5);
\draw [thick] (d) to (10.5,-2.5);

\node at (-0.8,0) {$v$};
\node at (4,-0.4) {$e$};
\node at (4,2) {$f$};
\node at (4,-2) {$f'$};
\node at (4,-4) {$\mathcal C$};

\node at (14,0) {$\longrightarrow$};

%%%%%%%%

\node (C) at (20,0) [shape=circle, fill=black] {};
\node (D) at (28,0) [shape=circle, fill=black] {};
\node (E) at (24,0) [shape=circle, draw] {};
\node (F) at (20,3) [shape=circle, draw] {};
\node (G) at (20,-3) [shape=circle, draw] {};
\node (H) at (28,3) [shape=circle, draw] {};
\node (I) at (28,-3) [shape=circle, draw] {};
\node (J) at (20,6) [shape=circle, fill=black] {};
\node (K) at (24,6) [shape=circle, draw] {};
\node (L) at (28,6) [shape=circle, fill=black] {};
\draw [thick] (C) to (E) to (D);
\draw [thick] (C) to (F) to (J);
\draw [thick] (J) to (K) to (L);
\draw [thick] (C) to (G) to (20,-5);
\draw [thick] (D) to (H) to (L);
\draw [thick] (D) to (I) to (28,-5);
\draw [thick] (C) to (17.5,1.5);
\draw [thick] (C) to (17.5,-1.5);
\draw [thick] (D) to (30.5,1.5);
\draw [thick] (D) to (30.5,-1.5);

\node at (19,0) {$v$};
\node at (24,0.8) {$v_e$};
\node at (24,-3) {\color{red}$r_e$};
\node at (24,3) {\color{blue}$l_e$};
\node at (16,-4) {$WP(\mathcal C)$};
\node at (24,-6.5) {$\downarrow$};
\draw (24,-6.5) to (24,-5.5);

%%%%%%%%%%%%%

\node (A) at (2,-15) [shape=circle, fill=black] {};
\node (B) at (6,-15) [shape=circle, draw] {};
\draw [thick] (A) to (B);
\draw [thick, dotted] (4.6,-12.7) to (A) to (4.6,-17.3);
\draw [thick, dotted] (B) to (6,-17);
\node at (1.2,-15) {$v$};

\draw [thick, dashed] (5.7,-14.8) arc [radius=4, start angle=225, end angle=190];
\draw [thick, dashed] (4.6,-17.3) arc [radius=3.1, start angle=120, end angle=90];

\draw [thin] (3.2,-13.3) arc [radius=2, start angle=60, end angle=300];
\draw [thin] (3.2,-16.8) to (2.5, -16.6);
\draw [thin] (3.2,-16.8) to (3.0, -17.4);
\draw [thin] (3.2,-13.3) to (2.5, -13.5);
\draw [thin] (3.2,-13.3) to (3.0, -12.7);
\node at (-0.5,-15) {$x$};

\draw [thin] (5,-16) to (7,-16);
\draw [thin] (5.5,-15.7) to (5,-16) to (5.5, -16.3);
\draw [thin] (6.5,-15.7) to (7,-16) to (6.5, -16.3);
\node at (7.7,-16.1) {$y$};
\node at (7,-15) {$v_e$};

\node at (4,-19) {$F$};

%%%%%%%%

\node (c) at (20,-15) [shape=circle, fill=black] {};
\node (d) at (28,-15) [shape=circle, fill=black] {};
\node (e) at (24,-15) [shape=circle, draw] {};
\node (f) at (20,-12) [shape=circle, draw] {};
\node (g) at (20,-18) [shape=circle, draw] {};
\node (h) at (28,-12) [shape=circle, draw] {};
\node (i) at (28,-18) [shape=circle, draw] {};
\node (j) at (20,-9) [shape=circle, fill=black] {};
\node (k) at (24,-9) [shape=circle, draw] {};
\node (l) at (28,-9) [shape=circle, fill=black] {};
\draw [thick] (c) to (e) to (d);
\draw [thick] (c) to (f) to (j);
\draw [thick] (j) to (k) to (l);
\draw [thick] (c) to (g) to (20,-20);
\draw [thick] (d) to (h) to (l);
\draw [thick] (d) to (i) to (28,-20);
\draw [thick] (c) to (17.5,-13.5);
\draw [thick] (c) to (17.5,-16.5);
\draw [thick] (d) to (30.5,-13.5);
\draw [thick] (d) to (30.5,-16.5);

\draw [thick, dashed] (24,-9.3) arc [radius=4, start angle=135, end angle=225];
\draw [thick, dashed] (24,-9.3) arc [radius=4, start angle=45, end angle=-45];

\draw [thick, dashed] (27.8,-18) arc [radius=5.4, start angle=45, end angle=135];
\draw [thick, dashed] (27.8,-18) arc [radius=5.4, start angle=-45, end angle=-135];

\draw [thick, dashed] (28.2,-12) arc [radius=5.4, start angle=135, end angle=120];
\draw [thick, dashed] (28.2,-12) arc [radius=5.4, start angle=-135, end angle=-120];

\draw [thick, dashed] (19.8,-12) arc [radius=5.4, start angle=45, end angle=60];
\draw [thick, dashed] (19.8,-12) arc [radius=5.4, start angle=-45, end angle=-60];

\draw [thick, red, rounded corners] (19.7,-12.4) to (19.7, -13.5) to (20.3,-13.5) to (20.3,-14.7)
 to (22,-14.7) to (22,-15.4) to (26,-15.4) to (26,-14.7) to (27.7,-14.7) to (27.7,-13.5) to (28.3,-13.5) to (28.3,-12.4);

\node at (19,-15) {$v$};
%\node at (24,0.8) {$v_e$};
%\node at (24,-3) {\color{red}$r_e$};
%\node at (24,3) {\color{blue}$l_e$};
\node at (17,-19) {$\M$};

\end{tikzpicture}

\end{center}
\caption{Construction of $\M$, of class $4$, from a chiral map $\mathcal C$} 
\label{ChiraltoClass4}
\end{figure}

%{\color{blue}[Gray shading for removed discs in Figure~\ref{ChiraltoClass4}?]}

\medskip

\noindent{\bf Remarks 1.} There is nothing special about choosing the left face $l_e$ at each edge, rather than the right face $r_e$: it would be equally valid to remove discs from faces so that the boundary components meet edges $e$ at $v_e$ whenever the face is $r_e$. This is equivalent to replacing $\mathcal C$ with its mirror image, by reversing the chosen orientation.
 
\smallskip

\noindent{\bf 2.} In this construction a more direct approach would be to start with a map in class $2^*{\rm ex}$, the Petrie dual of the class $2^P{\rm ex}$ of chiral maps, and then to remove open discs from its Walsh map. This yields the same set of maps $\M$. However the chiral maps have been more widely studied than those in the class $2^*{\rm ex}$, so it is easier to find examples of them in the literature (as in the examples below), and then to take their Petrie duals. 

\smallskip

\noindent{\bf 3.} The maps $\mathcal C$ in the class $2^*{\rm ex}$ are those for which a map subgroup $M$ has normaliser $N_{\Gamma}(M)=N(2^*{\rm ex})$. More generally, provided $M$ is normal in $N(2^*{\rm ex})$ (this allows those regular maps $\mathcal C$ which cover ${\mathcal N}(2^*{\rm ex})$) one can construct a map $\M$ with boundary, as above, by removing open discs from $W(\mathcal C)$: such maps $\M$ have a $2$-colouring of the sides of edges $e$, invariant under rotation around each vertex, which allows a consistent choice of which face at $e$ from which to remove an open disc.

\smallskip

\iffalse

\noindent{\bf 4.} In this construction, each boundary component of $\M'$ contains just one vertex. However, we can construct examples where they each contain any given number $n\ge 1$ of vertices by forming an $n$--sheeted covering of $\M'$, unbranched over the interior and with the same $n$-cycle as the monodromy permutation at each boundary component, which therefore lifts to a single boundary component in the covering space, containing $n$ vertices.

\medskip

\fi

\noindent{\bf Example 1.} James and the author showed in~\cite{JJ} that for each prime power $q=p^e$ there are, up to isomorphism, $\phi(q-1)/e$ orientably regular embeddings $\mathcal C$ of the complete graph $K_q$, each with ${\rm Aut}\,\mathcal C\cong AGL_1(q)$. (These maps were originally constructed by Biggs~\cite{Big}, generalising Heffter's examples in~\cite{Hef} for $e=1$.) If $q>4$ then each such map $\mathcal C$ is chiral, so by applying the above process we obtain an edge-transitive map $\M$ in class $4$, with non-empty boundary and with the same automorphism group. Since the Petrie polygons of $\mathcal C$ all have length $2p$, there are $q(q-1)/2p$ of them, so $P(\mathcal C)$ and hence $WP(\mathcal C)$ have characteristic
\[\chi=q-\frac{q(q-1)}{2}+\frac{q(q-1)}{2p}.\]
Since $\mathcal C$ is orientable and $K_q$ is not bipartite, $P(\mathcal C)$ and $WP(\mathcal C)$ are non-orientable. The map $\M$ is formed from $WP(\mathcal C)$ by removing $q(q-1)/2p$ open discs, one for each Petrie polygon of $\mathcal C$. In the simplest case, namely $q=5$, $\M$ has $\chi=-3$ with two boundary components of length $5$.

%{\color{blue}[Diagram for the case $q=5$?]}

\medskip

\noindent{\bf Example 2.} For a family of orientable examples, one could apply this process to the torus maps $\mathcal C=\{6,3\}_{b,c}$ for integers $b{\color{red}\,>\,}c>0$ (see~\cite[\S8.4]{CM}). These are orientable and bipartite, so in each case $P(\mathcal C)$ is also orientable, and hence so is $\M$. Now $\mathcal C$ has $2t$ vertices, $3t$ edges and $t$ faces, where $t=b^2+bc+c^2$. As shown by Melekoglu and Ulusan~\cite{MU}, $\mathcal C$ has Petrie length $2t/\gcd(b,c)$, so it has $3\gcd(b,c)$ Petrie polygons, and hence $WP(\mathcal C)$ has characteristic $3\gcd(b,c)-t$. There are $3\gcd(b,c)$ boundary components of $\M$, each of length $t/\gcd(b,c)$. {In the simplest case, where ${\mathcal C}=\{6,3\}_{2,1}$, $\M$ has $\chi=-4$ with three boundary components of length $7$.

% {\color{blue}[Diagram for the case $\{6,3\}_{2,1}$?]}

\medskip

Since $\mathcal N(4)$ has no flags fixed by $R_0, R_2$ or $R_0R_2$, the property of having a non-empty boundary is invariant under the group $\Omega$ of map operations, so the classes $4^*$ and $4^P$ are also wild. The maps with boundary in these classes can be obtained from those in class $4$ by applying duality $D$ and then Petrie duality $P$, so those in all three classes can be obtained from chiral maps $\mathcal C$, using the construction described above. Since $D$ transposes vertices and faces, while leaving Petrie polygons invariant, and $P$ transposes faces and Petrie polygons, while leaving vertices invariant, one can read off the properties of these features of the maps, including the actions of $G$ on them, from those of $\mathcal C$, as described above for maps in class~$4$.

\section{Boundary components}\label{boundarycpts}

It is useful to be able to count and describe the boundary components of an edge-transitive map with boundary. Clearly we can ignore the six void classes $2^{\sigma}{\rm ex}$ and $5^{\sigma}$, while maps in the four tame classes $1$ and $2^{\sigma}$ have all been described in Sections~\ref{class1bdy}, \ref{class2bdy}, \ref{class2*bdy} and \ref{class2Pbdy}: for classes $1$, $2$ and $2^*$ the maps with boundary are all on the closed disc, while in the case of $2^P$ they can also be on the annulus, the M\"obius band or the doubly infinite strip. It is therefore sufficient to restrict attention to the remaining classes, namely the four wild classes $T=3$ and $4^{\sigma}$. In all four cases, the only reflections in $N(T)$ are conjugates of $R_1$, so any boundary component contains vertices but not edges of the map. By Lemma~\ref{bdytrans}, all boundary components are equivalent under automorphisms, so it is sufficient to describe just one of them. For this, we need to determine its length, together with the valencies of the incident boundary vertices.

%%%%%%%%%%%%%%%%

\subsection{Boundary components in class~$3$}

First let $\M$ be in class $T=3$. For convenience of exposition, let us also assume assume that $\M$ is not the exceptional map $\M_{3,2}$ on the disc, shown in Figure~\ref{mapM3,2}; this map has one boundary component, of valency $3$, incident with vertices of valencies $1$, $1$ and $2$.

We saw in Section~\ref{class3bdy} that, apart from this exceptional case, just one of the four conjugacy classes of reflections in $N(3)$ lies in $M$, so for each edge $e$ of $\M$, just one of its four incident flags is on the boundary. It follows that $e$ has two vertices, namely a boundary vertex $v_e$ and an interior vertex $v'_e$. We can choose $\phi$ to be the flag  on $e$ fixed by the generator $S_0=R_1$ of $N(3)$, so $S_2=R_1^{R_2}$ induces a reflection $s_2$ of $\M$ fixing $v_e$, while $S_1=R_1^{R_0}$ and $S_3=R_1^{R_0R_2}$ induce reflections $s_2$ and $s_3$ fixing $v_e'$ (see Figure~\ref{class3edgebdy}, where the square bounded by fixed points of the reflections $s_i$ is a fundamental region for $G$). The existence of $s_2$ in $G={\rm Aut}\,\M$ implies that $v_e$ has valency $2$, while the valency of $v_e'$ is the order of the rotation $s_1s_3$ around it. The reflections $s_1$ and $s_2$ preserve the boundary component $B$ containing $v_e$, and $s_1s_2$ sends $v_e$ to an adjacent vertex on $B$, so the boundary length $b$ of $\M$ is the order of $s_1s_2$. There are $2b$ flags incident with each component $B$ (two at each vertex on $B$), while there are $|\Gamma:M|=4|G|$ flags in $\M$, a quarter of them incident with the boundary, so the number of boundary components is $|G|/2b$.

\begin{figure}[h!]
\begin{center}
\begin{tikzpicture}[scale=0.4, inner sep=0.8mm]

\node (a) at (-10,-0) [shape=circle, draw] {};
\node (b) at (10,0) [shape=circle, draw] {};
\node (c) at (0,10) [shape=circle, fill=black] {};
\draw [thick, dashed] (-15,0) to (a) to (b) to (15,0);
\draw [thick] (-15,5) to (a) to (c) to (b) to (15,5);
\draw [thick] (-2,12) to (c) to (2,12);

\draw [thick, dotted] (b) to (10,10) to (c) to (0,0);
\draw (b) to (7,1) to (8,2);

\draw (-2,5) to (2,5);
\draw (-1.5,5.5) to (-2,5) to (-1.5,4.5);
\draw (1.5,5.5) to (2,5) to (1.5,4.5);

\draw (8,5) to (12,5);
\draw (8.5,5.5) to (8,5) to (8.5,4.5);
\draw (11.5,5.5) to (12,5) to (11.5,4.5);

\draw (5,12) to (5,8);
\draw (4.5,11.5) to (5,12) to (5.5,11.5);
\draw (4.5,8.5) to (5,8) to (5.5,8.5);

\node at (16,0) {$B$};
\node at (2,3.8) {$s_1$};
\node at (8.5,3.8) {$s_2$};
\node at (6.5,8.2) {$s_3$};
\node at (5,4) {$e$};
\node at (10,-1) {$v_e$};
\node at (0,11.5) {$v'_e$};
\node at (6.5,1.5) {$\phi$};

\end{tikzpicture}

\end{center}
\caption{Part of a boundary component in a class $3$ map} 
\label{class3edgebdy}
\end{figure}

\bigskip

%%%%%%%%%%%%%%%

\subsection{Boundary components in class~$4$}

Now let $T=4$, so $N(T)=\langle S_1, S_2, S\mid S_i^2=1\rangle$. As shown earlier, we may assume that $S_1\in M$ but $S_2\not\in M$. The interior vertex of ${\mathcal N}(4)$ must lift to an orbit of $G={\rm Aut}\,\M$ consisting of interior vertices of $\M$. As in the preceding case $T=3$, with reflections $s_1$ and $s_2$ (induced by $S_1$ and $S_2$) now taking the roles played there by $s_0$ and $s_2$, it follows that each edge $e$ of $\M$ joins a boundary vertex $v_e$ of valency $2$ to an interior vertex $v_e'$. In this case, however, $G$ contains a rotation $s$ around $v'_e$, induced by $S=(R_1R_2)^{R_0}$ and with order equal to the valency of $v'_e$, but $G$ contains no reflection fixing this vertex. All faces are equivalent under $G$; in particular, we see in Figure~\ref{class4edgebdy} that the face containing the flag $\phi$ is invariant under a reflection (conjugate to $s_2$), so it is a $5$-gon, bounded by part of the boundary component $B$ containing $v_e$ and by four edges of $\M$. The same therefore applies to all faces. The boundary length $b$ of $\M$ is the order of the commutator $[s_2,s]$, and as in the case of class~$3$ the number of boundary components is $|G|/2b$.

\begin{figure}[h!]
\begin{center}
\begin{tikzpicture}[scale=0.4, inner sep=0.8mm]

\node (a) at (-10,-0) [shape=circle, draw] {};
\node (b) at (10,0) [shape=circle, draw] {};
\node (c) at (0,10) [shape=circle, draw] {};
\node (d) at (-7,7) [shape=circle, fill=black] {};
\node (e) at (7,7) [shape=circle, fill=black] {};
\draw [thick, dashed] (-15,0) to (a) to (b) to (15,0);
\draw [thick] (-12.5,5) to (a) to (d) to (c) to (e) to (b) to (12.5,5);
\draw [thick, dashed] (-3,11) to (c) to (3,11);
\draw [thick] (8,9) to (e) to (9,8);
\draw [thick] (-8,9) to (d) to (-9,8);

\draw [thick, dotted] (b) to (10,11);
\draw [thick, dotted] (c) to (0,0);
\draw (b) to (8.5,1.5) to (9.2,1.9);

\draw (-2,5) to (2,5);
\draw (-1.5,5.5) to (-2,5) to (-1.5,4.5);
\draw (1.5,5.5) to (2,5) to (1.5,4.5);

\draw (8,10) to (12,10);
\draw (8.5,10.5) to (8,10) to (8.5,9.5);
\draw (11.5,10.5) to (12,10) to (11.5,9.5);
\node at (2,4) {$ss_2s^{-1}$};

%\draw (5,12) to (5,8);
%\draw (4.5,11.5) to (5,12) to (5.5,11.5);
%\draw (4.5,8.5) to (5,8) to (5.5,8.5);

\node at (16,0) {$B$};
%\node at (2,3.8) {$s_1$};
\node at (12,9) {$s_2$};
%\node at (6.5,8.2) {$s_3$};
\node at (7.5,4) {$e$};
\node at (10,-1) {$v_e$};
\node at (8,6.8) {$v'_e$};
\node at (7.8,1.5) {$\phi$};

\draw (5.7,7) arc [radius=1, start angle=140, end angle=300];
\draw (4.9,6.4) to (5.6,6.9) to (6.1,6.3 );
\node at (5.2,5.2) {$s$};

\end{tikzpicture}

\end{center}
\caption{Part of a boundary component in a class $4$ map} 
\label{class4edgebdy}
\end{figure}

\medskip

The situation is similar when $T=4^*$, since each map in this class is the dual $D(\M)$ of a map $\M$ in class~$4$. Duality leaves the underlying surface unchanged, but now $G$ has a single orbit on vertices, all of which have valency $4$ and lie on the boundary, and there are two orbits on faces: boundary faces, each bounded by two edges and a segment of a boundary component, and interior faces, all  of valency equal to the order of $s$. As when $T=4$, there are $|G|/2b$ boundary components, all of length equal to the order $b$ of $[s_2,s]$.

The class~$4$ map $\M$ on the left of Figure~\ref{class4,48edgebdy}, with interior vertices of valency~$4$, is shown centred on such a vertex, displaying the cyclic but not dihedral symmetry about that vertex.  Parts of four boundary components are shown. The map on the right is the dual map $D(\M)$, in class~$4^*$.

\begin{figure}[h!]
\begin{center}
\begin{tikzpicture}[scale=0.3, inner sep=0.8mm]

\node (a) at (0,0) [shape=circle, fill=black] {};
\node (b) at (5,0) [shape=circle, draw] {};
\node (c) at (0,5) [shape=circle, draw] {};
\node (d) at (-5,0) [shape=circle, draw] {};
\node (e) at (0,-5) [shape=circle, draw] {};
\node (b') at (6,4) [shape=circle, fill=black] {};
\node (c') at (-4,6) [shape=circle, fill=black] {};
\node (d') at (-6,-4) [shape=circle, fill=black] {};
\node (e') at (4,-6) [shape=circle, fill=black] {};
\node (b") at (2,6) [shape=circle, draw] {};
\node (c") at (-6,2) [shape=circle, draw] {};
\node (d") at (-2,-6) [shape=circle, draw] {};
\node (e") at (6,-2) [shape=circle, draw] {};

\draw [thick] (b") to (b') to (b) to (a) to (d) to (d') to (d");
\draw [thick] (c") to (c') to (c) to (a) to (e) to (e') to (e");

\draw [thick, dashed] (7.5,1) to (6,1) to (b) to (e") to (7.5,-2);
\draw [thick, dashed] (-1,7.5) to (-1,6) to (c) to (b") to (2,7.5);
\draw [thick, dashed](-7.5,-1) to (-6,-1) to (d) to (c") to (-7.5,2);
\draw [thick, dashed] (1,-7.5) to (1,-6) to (e) to (d") to (-2,-7.5);

\draw [thick] (6.5,5.5) to (b') to (7.5,4.5);
\draw [thick] (-5.5,6.5) to (c') to (-4.5,7.5);
\draw [thick] (-6.5,-5.5) to (d') to (-7.5,-4.5);
\draw [thick] (5.5,-6.5) to (e') to (4.5,-7.5);

\node at (8,-6) {$\M$};

%%%%%%%

\node (B) at (25,-2) [shape=circle, fill=black] {};
\node (C) at (22,5) [shape=circle, fill=black] {};
\node (D) at (15,2) [shape=circle, fill=black] {};
\node (E) at (18,-5) [shape=circle, fill=black] {};

\node (B') at (26.5,1) [shape=circle, fill=black] {};
\node (C') at (19,6.5) [shape=circle, fill=black] {};
\node (D') at (13.5,-1) [shape=circle, fill=black] {};
\node (E') at (21,-6.5) [shape=circle, fill=black] {};

\draw [thick] (B') to (C) to (24,6);;
\draw [thick] (C') to (D) to (14,4);;
\draw [thick] (D') to (E) to (16,-6);
\draw [thick] (E') to (B) to (26,-4);

\draw [thick, dashed] (28,1) to (B') to (25,0) to (B) to (27,-2);
\draw [thick, dashed] (19,8) to (C') to (20,5) to (C) to (22,7);
\draw [thick, dashed] (12,-1) to (D') to (15,0) to (D) to (13,2);
\draw [thick, dashed] (21,-8) to (E') to (20,-5) to (E) to (18,-7);

\draw [thick] (B) to (C) to (D) to (E) to (B);
\draw [thick] (27,3) to (B') to (28.3,1.8); \draw [thick] (B') to (27.8,2.5);
\draw [thick] (17,7) to (C') to (18.2,8.3); \draw [thick] (C') to (17.5,7.8);
\draw [thick] (13,-3) to (D') to (11.7,-1.8); \draw [thick] (D') to (12.2,-2.5);
\draw [thick] (23,-7) to (E') to (21.8,-8.3); \draw [thick] (E') to (22.5,-7.8);

\node at (28,-6) {$D(\M)$};

\end{tikzpicture}

\end{center}
\caption{A dual pair of maps with boundary, in classes $4$ and $4^*$.} 
\label{class4,48edgebdy}
\end{figure}

\medskip

The maps in class $4^P$ are the Petrie duals of those in class $4^*$.  In the latter class (as in class $4$), all boundary flags are fixed by $R_1$ rather than $R_0$ or $R_2$, so the same applies to the corresponding flags in the Petrie dual maps. Thus the number of boundary flags (and in particular the property of having non-empty boundary) is preserved by the pairing $P$ between these two classes. The embedded graphs are also preserved, each having vertices of valency $4$, all on the boundary.

%\medskip

%{\color{blue}[Give more information in the case $T=4^P$.]}

%%%%%%%%%%%%%

\section{Edge-transitive maps with free edges}\label{free}

Having, at least to some extent, dealt with edge-transitive maps with non-empty boundary, it is also useful at this point to deal with those with a free edge. Here we can provide a complete and concise classification, allowing one to restrict attention mainly to edge-transitive maps with no free edges.

If an edge-transitive map $\M$ has a free edge, then all its edges must be free, so in the monodromy group
\[G=\langle r_0, r_1, r_2\mid r_i^2= \ldots =1\rangle\]
we have $r_0=1$ or $r_0=r_2$ as the free ends of the edges are in the boundary or the interior of $\M$. Thus $G=\langle r_1, r_2\mid r_i^2= \ldots =1\rangle$ is either the trivial group or a dihedral group $D_n$ for some $n\in\N\cup\{\infty\}$. The first case corresponds to the trivial map $\M={\mathcal N}(1)$, so we will assume from now on that $G\cong D_n$.  Here we include the cases $n=1$ where $G\cong C_2$, and $n=2$ where $G\cong C_2\times C_2$, in addition to the non-abelian cases $n\ge 3$.

By definition, $G$ acts faithfully on the flags of $\M$, so the stabiliser $G_{\phi}$ of a flag $\phi$ has trivial core in $G$. Since every subgroup of $G^+:=\langle r_1r_2\rangle$ is normal in $G$, this implies that $G_{\phi}\cap G^+=1$, so $|G_{\phi}|\le |G:G^+|\le 2$. The edge-transitivity of $\M$ is equivalent to the condition that
\[G=N_G(G_{\phi}).\langle r_0, r_2\rangle=N_G(G_{\phi}).\langle r_2\rangle,\]
which implies that $|G:N_G(G_{\phi})|\le 2$. These inequalities, and the fact that $G_{\phi}$ has trivial core, give us just two possibilities: either
\begin{itemize}
\item $G_{\phi}=1$, so $G$ acts regularly and $\M$ is a regular map, or
\item $G$ is non-abelian with $|G_{\phi}|=|G:N_G(G_{\phi})|=2$.
\end{itemize}
It is an easy exercise in the algebra of dihedral groups to see that the second possibility arises if and only if $n=4$, with $G_{\phi}$ conjugate to $\langle r_1\rangle$. This case corresponds to the maps $\M_{2^*,0}$ and $\M_{2^P,0}$ in classes $2^*$ and $2^P$, described in Sections~\ref{class2*bdy} and \ref{class2Pbdy} and illustrated in Figures~\ref{mapM2*,0} and \ref{mapM2p,0}, as $r_0=1$ or $r_0=r_2$. We may therefore assume from now on that the first possibility arises, and $\M$ is a regular map.

\iffalse
In order to describe the maps arising in this situation, we first introduce some notation. Given any map $\M$ on the closed disc $\mathbb D$, let $\M^{\circ}$ be the map obtained by adding an open disc to extend the underlying surface to a sphere, embedding the same graph in a hemisphere.
\fi

Suppose first that $r_1, r_2\ne 1$ (as must be the case when $n\ne 1$). When $r_0=1$ we obtain the regular map ${\mathcal A}_n$ on the disc, discussed in Section~\ref{class1bdy} and shown in Figures~\ref{mapA4} (left) and \ref{mapAinfty} for $n=4$ and $\infty$, where $n\;(=1, 2, \ldots, \infty)$ is the order of $r_1r_2$. When $r_0=r_2$ we obtain the regular spherical map ${\mathcal A}_n^{\circ}$ obtained from ${\mathcal A}_n$ by embedding the disc in $S^2$: see Figure~\ref{mapA4circ} for ${\mathcal A}_4^{\circ}$.

  \begin{figure}[h!]
\begin{center}
\begin{tikzpicture}[scale=0.5, inner sep=0.8mm]

\node (a) at (0,0) [shape=circle, fill=black] {};
\node (b) at (2,0) {};
\node (c) at (0,2) {};
\node (d) at (-2,0){};
\node (e) at (0,-2) {};
\draw [thick] (b) to (d);
\draw [thick] (c) to (e);

\draw (3,0) arc (0:360:3);

\end{tikzpicture}

\end{center}
\caption{The map ${\mathcal A}_4^{\circ}$.} 
\label{mapA4circ}
\end{figure}

If $r_1=1$ but $r_2\ne 1$ we obtain the regular maps ${\mathcal N}(2^*)$ and ${\mathcal N}(2^P)$ on the disc, shown in Figure~\ref{basicmaps}, as $r_0=1$ or $r_0=r_2$. If $r_2=1$ (so $r_0=1$ also) but $r_1\ne 1$ we obtain the regular map $\mathcal D$ on the disc, described in Section~\ref{class1bdy} and shown on the right in Figure~\ref{mapA4}.

This covers all the possibilities, so we have proved the following:

\begin{thm}\label{free}
The edge-transitive maps with free edges are:
\begin{itemize}
\item the maps ${\mathcal A}_n$ and ${\mathcal A}_n^{\circ}$ for $n\in{\N}\cup\{\infty\}$, together with ${\mathcal N}(1)$, ${\mathcal N}(2^*)$, ${\mathcal N}(2^P)$ and $\mathcal D$, all in class $1$;
\item the maps $\M_{2^*,0}$ and $\M_{2^P,0}$ in classes $2^*$ and $2^P$.
\end{itemize}
They are all on the closed disc, apart from the maps ${\mathcal A}_n^{\circ}$ on the sphere. Their automorphism groups are all dihedral groups, apart from ${\rm Aut}\,{\mathcal N}(1)$ which is the identity group. \hfill$\square$
\end{thm}

%{\color{blue}[Draw all the edge-transitive maps with free edges in a single diagram ($n=4$ for the infinite families)?]}

%%%%%%%%%%%%%%%%%%

\section{Medials and almost edge-transitive maps}\label{almostET}

%{\color{blue}[Postponed from earlier in Part I and largely rewritten]}

%\medskip

Here we return to the medial map construction, which was defined informally in \S\ref{medialmaps}. More formally, given a map $\M$ with a set $\Phi$ of flags, we define the set of flags of $\M^{\rm med}$ to be the set
\[\tilde\Phi=\{\phi_0, \phi_2 \mid \phi=(v,e,f)\in \Phi\}\]
of symbols $\phi_0, \phi_2$ associated with flags $\phi=(v,e,f)$ of $\M$. (Intuitively, these are flags of $\M^{\rm med}$ incident with its vertex $v_e$ and its edge joining $v_e$ to $v_{e'}$, where $e'$ is the edge of $\M$ incident with $\phi r_1\in\Phi$; they are respectively incident with the faces of $\M^{\rm med}$ containing $v$ and the centre of $f$, as in Figure~\ref{medialflags}.) Let us define the following permutations $\tilde r_i$ of $\tilde\Phi$:
\[\tilde r_0: \phi_0\mapsto (\phi r_1)_0, \; \phi_2\mapsto (\phi r_1)_2,\]
\[\tilde r_1: \phi_0\mapsto (\phi r_2)_0, \; \phi_2\mapsto (\phi r_0)_2,\]
\[\tilde r_2: \phi_0\mapsto \phi_2, \; \phi_2\mapsto \phi_0.\]

\begin{figure}[h!]
\begin{center}
\begin{tikzpicture}[scale=0.5, inner sep=0.8mm]

\node (a) at (0,0) [shape=circle, fill=black] {};
\node (b) at (10,10) [shape=circle, fill=black] {};
\node (A) at (5,5) [shape=circle, draw] {};
\node (B) at (5,-5) [shape=circle, draw] {};

\draw [thick] (a) to (A) to (b);
\draw [thick] (a) to (B) to (7,-7);
\draw [thick, dashed] (5,7) to (A) to (B);
\draw [thick, dashed] (1,5) to (A) to (9,5);

\draw (A) to (4,3) to (6,3) to (A);
\draw (B) to (4,-3) to (6,-3) to (B);
\draw (a) to (2,1) to (1,2) to (a);
\draw (a) to (2,-1) to (1.5,-1.5);
\draw (b) to (9,8) to (8.5,8.5);
\draw (3,5) to (3,4) to (A) to (7,4) to (7,5);

\node at (-0.8,0) {$v$};
\node at (4.3,5.7) {$v_e$};
\node at (4.3,-5.7) {$v_{e'}$};
\node at (2.2,2.8) {$e$};
\node at (2,-3) {$e'$};
\node at (2.5,1) {$\phi$};
\node at (10,8) {$\phi r_0$};
\node at (2.8,-1) {$\phi r_1$};
\node at (0,2.3) {$\phi r_2$};
\node at (8,0) {$f$};
\node at (4,2.5) {$\phi_0$};
\node at (6,2.5) {$\phi_2$};
\node at (8.3,4.2) {$(\phi r_0)_2$};
\node at (1.7,4.2) {$(\phi r_2)_0$};
\node at (7.5,-3) {$(\phi r_1)_2$};
\node at (3.8,-2.3) {$(\phi r_1)_0$};

\node (c) at (15,2) [shape=circle, fill=black] {};
\node (d) at (15,-2) [shape=circle, draw] {};
\draw [thick] (13,1) to (c) to (17,3);
\draw [thick, dashed] (13,-3) to (d) to (17,-1);

\node at (16,1) {$\M$};
\node at (16,-3) {$\M^{\rm med}$};

\end{tikzpicture}

\end{center}
\caption{Permutations of flags $\phi_0,\phi_2$ of a medial map $\M^{\rm med}$} 
\label{medialflags}
\end{figure}

It is straightforward to check that $\tilde r_i^2=(\tilde r_0\tilde r_2)^2=1$, that the resulting action $R_i\mapsto\tilde r_i$ of $\Gamma$ on $\tilde\Phi$ is transitive, so it defines a map, and that this map is isomorphic to the medial map of $\M$ as described above. For example, the vertices of this map are the orbits of $\langle\tilde r_1, \tilde r_2\rangle$ on $\tilde\Phi$, and these correspond to the orbits of $\langle r_0, r_2\rangle$ on $\Phi$, that is, to the edges of $\M$. Moreover, one can check that $(\tilde r_1\tilde r_2)^4=1$, so this map is $4$-valent. The isomorphism between $\M^{\rm med}$ and $D(\M)^{\rm med}$ is given by $\phi_0\mapsto \psi_2$ and $\phi_2\mapsto \psi_0$ where $\psi=(f,e,v)$ is the flag of $D(\M)$ corresponding to a flag $\phi=(v,e,f)$ of $\M$.

This formulation extends the definition of $\M^{\rm med}$ to cases where $\M$ has a non-empty boundary. In such cases, any fixed point $\phi$ of $r_1$ on $\Phi$ induces two fixed points $\phi_0, \phi_2$ of $\tilde r_0$ on $\tilde\Phi$, and any fixed point of $r_0$ or $r_2$ on $\Phi$ induces a fixed point $\phi_2$ or $\phi_0$ of $\tilde r_1$ on $\tilde\Phi$, so $\M^{\rm med}$ also has a non-empty boundary; however, $\tilde r_2$ has no fixed points so $\M^{\rm med}$ has no edges lying along the boundary. Note also that $\tilde r_0\tilde r_2$ has no fixed points, so $\M^{\rm med}$ has no internal free edges.

There is an alternative group-theoretic interpretation of the medial construction. One can regard $\Gamma$ as the extended triangle group $\Delta[\infty, 2, \infty]$. As such, it is a subgroup of index $2$ in the extended triangle group
\[\tilde\Gamma:=\Delta[4,2,\infty]
=\langle \tilde R_i\;(i=0, 1, 2)\mid \tilde R_i^2=(\tilde R_1\tilde R_2)^4=(\tilde R_2\tilde R_0)^2=1\rangle.\]

\begin{figure}[h!]
\begin{center}
\begin{tikzpicture}[scale=0.5, inner sep=0.8mm]

\node (a) at (5,0) [shape=circle, fill=black] {};
\node (b) at (0,0) [shape=circle, draw] {};

\draw [thick, dotted] (5,0) arc (0:360:5);
\draw [ultra thick] (a) to (b);
\draw (b) to (0,5);
\draw [thick, dashed] (b) to (1.45,1.45);
\draw (0,5) arc (180:270:5);

\node at (5.7,0) {$1$};
\node at (-0.7,0) {$0$};
\node at (0,5.6) {$i$};

%\draw [<->] (-0.5,1.5) to (0.5,1.5);
%\draw [<->] (1.1,1.1) to (1.8,1.8);
%\draw [<->] (1.5,0.5) to (1.5,-0.5);

%\node at (-1.3,1.5) {$R_0$};
%\node at (2.3,2.3) {$R_1$};
%\node at (1.6,-1.2) {$R_2$};
\node at (-2,-2) {$\hyp$};

\end{tikzpicture}
\end{center}
\caption{The inclusion of $\Gamma$ in $\tilde\Gamma$} 
\label{inclusion}
\end{figure}

\iffalse

\begin{figure}[h!]
\begin{center}
\begin{tikzpicture}[scale=0.4, inner sep=0.8mm]

\draw [thick] (-10,0) to (10,0);
\draw [thick] (0,0) to (0,4.14);
\draw (-0.5,0) to (-0.5,0.5) to (0.5,0.5) to (0.5,0);
\draw (-0.35,3.79) to (0,3.44) to (0.35,3.79);

\draw [thick] (-10,0) arc (270:315:14.14);
\draw [thick] (0,4.14) arc (225:270:14.14);

\end{tikzpicture}

\end{center}
\caption{Inclusion between two triangle groups} 
\label{inclusion}
\end{figure}

\fi

\noindent This is shown by Figure~\ref{inclusion}, drawn in the disc model of the hyperbolic plane. Here the triangle with vertices at $1, i$ (both ideal) and $0$ is a fundamental triangle for $\Gamma$, which is generated by reflections in the sides of this triangle as in Figure~\ref{Riacting}. The broken line, bisecting the right-angle at $0$, divides this triangle into two triangles, each of which has internal angles $\pi/4, \pi/2, 0$ and is a fundamental region for the group $\tilde\Gamma$ generated by the reflections in its sides. Just as the images under $\Gamma$ of the black vertex at $1$ and the unit interval form the vertices and edges of the universal map $\M_{\infty}$, of which every map $\M$ is a quotient by some map subgroup $M$, the images under $\tilde\Gamma$ of the white vertex at $0$ and the broken line form the vertices and edges of the map $(\M_{\infty})^{\rm med}$, of which $\M^{\rm med}$ is a quotient by $M$. Alternatively, one can simply realise $\Gamma$ as the kernel of the epimorphism $\tilde R_0, \tilde R_2\mapsto 1, \tilde R_2\mapsto -1$ from $\tilde\Gamma$ to $\{\pm 1\}\cong C_2$, with the Reidemeister--Schreier algorithm giving a presentation for this kernel in terms of the standard generators
\[R_0=\tilde R_1,\; R_2=\tilde R_0,\; R_2=\tilde R_1^{\tilde R_2}\]
for $\Gamma$. (Similar arguments show that $\Delta[n,m,n]$ is a subgroup of index $2$ in $\Delta[2m,2,n]$ for all $m, n\in\N\cup\{\infty\}$. See~\cite{Sin72} for inclusions between ordinary triangle groups; the inclusions between extended triangle groups can be deduced from these, or extracted from more general classifications of inclusions between NEC groups~\cite{CIP, EI}.) Any map $\M$ corresponds to a conjugacy class of map subgroups $M$ of $\Gamma$; by the inclusion of $\Gamma$ in $\tilde\Gamma$ these are also mutually conjugate subgroups of $\tilde\Gamma$, and as such they correspond to a $4$-valent map on the same surface but with twice as many flags as $\M$, namely $\M^{\rm med}$.

As noted earlier, ${\rm Aut}\,\M^{\rm med}={\rm Aut}\,\M$ unless $\M$ is self-dual, in which case ${\rm Aut}\M^{\rm med}$ contains ${\rm Aut}\,\M$ with index $2$.  We also noted that, as shown by \v Sir\'a\v n, Tucker and Watkins in~\cite[Lemma~2.2]{STW}, if $\M$ is edge-transitive, then $\M^{\rm med}$ is edge-transitive if and only if $\M$ is in class~$1$ or $2^{\sigma}{\rm ex}$ for some $\sigma=\emptyset, *$ or $P$.

\iffalse 
\begin{figure}[h!]
\begin{center}
\begin{tikzpicture}[scale=0.5, inner sep=0.8mm]

\node (a) at (0,3.2) {};
\node (b) at (3.2,0) {};
\node (c) at (-3.2,0) {};
\node (d) at (0,-3.2) {};
\node (e) at (0,0) [shape=circle, fill=black] {};

\draw [thick] (a) to (e);
\draw[thick]  (b) to (e);
\draw [thick] (c) to (e);
\draw [thick] (d) to (e);

\draw [dashed] (3,0) arc (0:360:3);

%%%%%

\node (E) at (13,0) [shape=circle, fill=black] {};
\draw [dashed] (13,0) arc (0:360:3);
\draw [thick] (13,0) arc (0:90:3);
\draw [thick] (13,0) arc (0:-90:3);

\end{tikzpicture}

\end{center}
\caption{The maps ${\mathcal N}(3)^{\rm med}$ and $\mathcal N(2_{\{0,2\}})$ on the closed disc} 
\label{mapsA4N0}
\end{figure}

\fi

However, it is also possible for a map which is not edge-transitive to have an edge-transitive medial map. Let us define a map $\M$ to be {\sl almost edge-transitive} if ${\rm Aut}\,\M$ acts transitively on the edges of $\M^{\rm med}$ but not on those of $\M$. Recall that in \S\ref{14classes} we defined the (non-edge-transitive) class $T=2_{\{0,2\}}$ to consist of those maps $\M$ such that $\M/{\rm Aut}\,\M$ is isomorphic to the map ${\mathcal N}(2_{\{0,2\}})$ on the closed disc shown in Figure~\ref{mapN2{02}}. It is easy to verify, both visually and by using the formal definition given above, that ${\mathcal N}(2_{\{0,2\}})^{\rm med}$ is isomorphic to the map ${\mathcal N}(3)$ in Figure~\ref{basicmaps}, with one edge.

\begin{thm}
A map is almost edge-transitive if and only if it is class~$2_{\{0,2\}}$.
\end{thm}

\noindent{\sl Proof.} If $\M$ is in class $2_{\{0,2\}}$ than ${\rm Aut}\,\M$ has two orbits on the edges of $\M$ but only one on those of $\M^{\rm med}$, so $\M$ is almost edge transitive. Conversely, suppose that $\M$ is almost edge-transitive. The edges of $\M^{\rm med}$ (the orbits of $\tilde E=\langle \tilde r_0, \tilde r_2\rangle$ on $\tilde\Phi$) can be identified with the cycles of $r_1$ on $\Phi$, and since these are permuted transitively by ${\rm Aut}\,\M$ we have $\Gamma=N_{\Gamma}(M)\langle R_1\rangle$. This implies that $|\Gamma:N_{\Gamma}(M)|\le 2$, giving just eight possibilities for $N_{\Gamma}(M)$, namely the subgroups $N(T)$ for $T=1$, $2^{\sigma}$, $2^{\sigma}{\rm ex}$ and $2_{\{0,2\}}$, where $\sigma=\emptyset, *$ or $P$, so $\M$ is in the corresponding class $T$. Only the last of these contains non-edge-transitive maps, so $T=2_{\{0,2\}}$ since $\M$ is not edge-transitive. \hfill$\square$

\medskip

Notice that since ${\mathcal N}(2_{\{0,2\}})$ has a single vertex and a single face, any almost edge-transitive map $\M$ must be vertex- and face-transitive. Indeed, ${\rm Aut}\,\M$ must act transitively on incident vertex-face pairs, since ${\mathcal N}(2_{\{0,2\}})$ has only one of these. Moreover, the edges of $\M$ must form two orbits of ${\rm Aut}\,\M$, alternating around vertices and faces of $\M$ since they are preserved by $r_0$ and $r_2$ but transposed by $r_1$. In this sense, we can say that $\M$ must be $2$-{\sl edge colourable}.

The basic subgroup $N(2_{\{0,2\}})$ for almost edge-transitive maps has a presentation
\[N(2_{\{0,2\}})=\langle R_i, S_i:=R_i^{R_1}\;(i=0, 2)\mid R_i^2=(R_0R_2)^2=S_i^2=(S_0S_2)^2=1\rangle,\]
so it is the free product $E*F$ of two Klein four-groups $E$ and $F:=E^{R_i}$. A normal subgroup $M$ of $N(2_{\{0,2\}})$ will have $N_{\Gamma}(M)=N(2_{\{0,2\}})$ if and only if it is not invariant under the automorphism of $N(2_{\{0,2\}})$ induced by conjugation by $R_1$, transposing $R_i$ and $S_i$ for each $i=0, 2$. In order to construct almost edge-transitive maps we therefore look for quotients $A=N(2_{\{0,2\}})/M$ of $N(2_{\{0,2\}})$ which do not have such a forbidden automorphism, and which can therefore be automorphism groups of such maps.

There are, of course, numerous examples of groups $A$, finite or infinite, which can be generated by the images of $E$ and $F$, and in many cases it is possible to choose the generators $r_i$ and $s_i$ corresponding to $R_i$ and $S_i$ to avoid the forbidden automorphism $r_i\longleftrightarrow s_i$. If we want $\M$ to have empty boundary and no free edges, we need $E$ and $F$ to be faithfully represented in $A$, but this is no great restriction. In this case, the face- and vertex-valencies of $\M$ are $2k$ and $2l$, where $k$ and $l$ are the orders of $r_0s_0$ and $r_2s_2$. Since there are $2|A|$ flags, there are $|A|/2k$ faces, $|A|/2l$ vertices, and $|A|/2$ edges, so in the finite case $\M$ has characteristic
\[\chi=\frac{|A|}{2}\left(\frac{1}{k}+\frac{1}{l}-\frac{1}{2}\right).\]
The map $\M$ (and hence also $\M^{\rm med}$) is orientable if and only if $M\le\Gamma^+$, or equivalently, $A$ has a subgroup of index $2$ containing none of the generators $r_i$ or $s_i$ ($i=0, 2$).

\begin{figure}[h!]
\begin{center}
\begin{tikzpicture}[scale=0.2, inner sep=0.8mm]

\node (a) at (-15, 10) [shape=circle, fill=black] {};
\node (b) at (-5, 10) [shape=circle, fill=black] {};
\node (c) at (5, 10) [shape=circle, fill=black] {};
\node (d) at (15, 10) [shape=circle, fill=black] {};
\node (a') at (-15, 0) [shape=circle, fill=black] {};
\node (b') at (-5, 0) [shape=circle, fill=black] {};
\node (c') at (5, 0) [shape=circle, fill=black] {};
\node (d') at (15, 0) [shape=circle, fill=black] {};
\node (a*) at (-15, -10) [shape=circle, fill=black] {};
\node (b*) at (-5, -10) [shape=circle, fill=black] {};
\node (c*) at (5, -10) [shape=circle, fill=black] {};
\node (d*) at (15, -10) [shape=circle, fill=black] {};

\node (A) at (-10, 10) [shape=circle, draw] {};
\node (B) at (0, 10) [shape=circle, draw] {};
\node (C) at (10, 10) [shape=circle, draw] {};
\node (D) at (-15, 5) [shape=circle, draw] {};
\node (E) at (-5, 5) [shape=circle, draw] {};
\node (F) at (5, 5) [shape=circle, draw] {};
\node (G) at (15, 5) [shape=circle, draw] {};
\node (D') at (-15, -5) [shape=circle, draw] {};
\node (E') at (-5, -5) [shape=circle, draw] {};
\node (F') at (5, -5) [shape=circle, draw] {};
\node (G') at (15, -5) [shape=circle, draw] {};
\node (A') at (-10, 0) [shape=circle, draw] {};
\node (B') at (0, 0) [shape=circle, draw] {};
\node (C') at (10, 0) [shape=circle, draw] {};
\node (A*) at (-10, -10) [shape=circle, draw] {};
\node (B*) at (0, -10) [shape=circle, draw] {};
\node (C*) at (10, -10) [shape=circle, draw] {};

\draw [thick] (a) to (A) to (b) to (B) to (c) to (C) to (d);
\draw [thick] (a') to (A') to (b') to (B') to (c') to (C') to (d');
\draw [thick] (a*) to (A*) to (b*) to (B*) to (c*) to (C*) to (d*);
\draw [thick] (a) to (D) to (a') to (D') to (a*);
\draw [thick] (b) to (E) to (b') to (E') to (b*);
\draw [thick] (c) to (F) to (c') to (F') to (c*);
\draw [thick] (d) to (G) to (d') to (G') to (d*);

\draw [thick, dashed] (A) to (D) to (A') to (E) to (A);
\draw [thick, dashed] (B) to (E) to (B') to (F) to (B);
\draw [thick, dashed] (C) to (F) to (C') to (G) to (C);
\draw [thick, dashed] (A') to (D') to (A*) to (E') to (A');
\draw [thick, dashed] (B') to (E') to (B*) to (F') to (B');
\draw [thick, dashed] (C') to (F') to (C*) to (G') to (C');

\node (x) at (25,2) [shape=circle, fill=black] {};
\node (y) at (25,-2) [shape=circle, draw] {};
\draw [thick] (23,1) to (x) to (27,3);
\draw [thick, dashed] (23,-3) to (y) to (27,-1);

\node at (31,2) {$\M$};
\node at (32.5,-2) {$\M^{\rm med}$};

\end{tikzpicture}

\end{center}
\caption{A medial map on the torus} 
\label{torusmedial}
\end{figure}

\medskip

\noindent{\bf Example.} For a very simple example, take $A=D_m\times D_n$ for distinct integers $m, n\ge 2$, choose $r_0$ and $s_2$ to be involutions generating $D_m$, and  $r_2$ and $s_0$ to be involutions generating $D_n$. The resulting map $\M$ is a torus map of type $\{4,4\}$, the quotient of the universal map of this type by the lattice generated by translations through $m$ and $n$ in the two coordinate directions. The translation subgroup $C_m\times C_n$ of $A$ acts regularly on the $mn$ vertices and faces, and $A$ has two orbits of length $mn$ on the edges, consisting of the two parallel classes; the stabiliser of each edge is generated by the reflection in that edge. The map $\M^{\rm med}$ is also of type $\{4,4\}$, with $A$ acting regularly on its $4mn$ edges; in addition the self-duality of $\M$ provides extra automorphisms of $\M^{\rm med}$, such as half-turns around the mid-points of its edges, which are not induced by $A$. (See Figure~\ref{torusmedial} for the case $m=3, n=2$, with opposite sides of the outer rectangle identified to form a torus.) The automorphism groups of $\M$ and of $\M^{\rm med}$ are quotients of the $2$-dimensional euclidean crystallographic groups $pmm$ and $cmm$.

\medskip

\noindent{\bf Example.} For a series of examples with unbounded characteristic, let $A=S_n$ for some $n\ge 4$, with
\[r_0=(1,2)(3,n)(4,n-1)\ldots \quad{\rm and}\quad s_0=(2)(1,3)(4,n)(5,n-1)\ldots,\]
so that $r_0s_0=(1,2,3,\ldots, n)$, let $r_2=(1,2)$, so that $(r_0r_2)^2=1$ and $\langle r_0s_0, r_1\rangle=A$, and let $s_2$ be any involution (other than $s_0$) commuting with $s_0$, such as $(1,3)$ if $n\ge 5$ or $(1,3)(2,4)$ if $n=4$. The resulting map $\M$ has type $\{2n, 2l\}$, where $l$ is the order of $r_2s_2$. It has characteristic
\[\chi=\frac{n!}{2}\left(\frac{1}{n}+\frac{1}{l}-1\right),\]
so that since $l\ge 2$ we have
\[-\frac{n!}{2}<\chi\le-\frac{n!}{2}\left(\frac{1}{2}-\frac{1}{n}\right)\le-\frac{n!}{8}.\]
The map $\M$ is orientable if and only if the four generators $r_i$ and $s_i$ are all odd permutations, that is, $n\equiv 3$ mod~$(4)$ and $s_2$ is chosen to be odd. It is vertex- and face-transitive, but not edge-transitive. It is not self-dual, so the edge-transitive graph $\M^{\rm med}$ also has automorphism group $A=S_n$. For instance, taking $n=4$ and $s_2=(2,4)$ gives $l=3$ and so $\chi=-5$, while taking $n=5$ and $s_2=(4,5)$ gives $l=2$ and so $\chi=-18$; the maps $\M$ and $\M^{\rm med}$ are non-orientable in both cases. The smallest orientable examples of this form arise when $n=7$ and $s_2$ commutes with $r_2$: for instance one could take $s_2=(4,7)$, giving maps of genus $451$.

%What about $\overline\M={\mathcal N}(5^P)$ with a free edge at the vertex?
\iffalse
\[r_0=(1,2),\quad r_2=(3,4)(5,6)(7,8)\ldots,\]
\[s_0=(2,3)(4,5)(6,7)\ldots,\quad s_2=(2,4)(3,5).\]
Then
\[r_i^2=s_i^2=(r_0r_2)^2=(s_0s_2)^2=1,\]
so we have a permutation representation $\theta:\Gamma_0\to G:=\langle r_0, r_2, s_0, s_2\rangle$ of $\Gamma_0$ given by $R_i\mapsto r_i, S_i\mapsto s_i$ for $i=0, 2$. Now $G$ is a transitive subgroup of $S_n$, and since $r_2s_0$ is an $(n-1)$-cycle $G$ is double transitive. Thus $G$ is primitive, and since it contains the transposition $r_0$ we have $G=S_n$. If $n\ne 6$ then all automorphisms of $G$ are inner; thus there can be no automorphism of $G$ transposing $r_i$ and $s_i$ for each $i=0, 2$, since $r_0$ and $s_0$ have different cycle structures. The same applies when $n=6$, since then these two permutations have different parities, whereas all automorphisms of $S_6$ leave $A_6$ invariant. It follows that $M:=\ker\theta$ has normaliser $\Gamma_0$, as required.

Since $r_0s_0$ and $r_2s_2$ have order $6$ and $10$ (provided $n\ge 6$), so the map $\M$ corresponding to $M$ has type $\{20, 12\}$, with automorphism group $G\cong S_n$.
\fi

\begin{figure}[h!]
\begin{center}
\begin{tikzpicture}[scale=0.5, inner sep=0.8mm]

\node (A) at (0,4) [shape=circle, fill=black] {};
\draw [thick] (4,0) arc (0:360:4);
\draw [dashed] (2,0) arc (0:360:2);
\draw [thick] (A) to (0,2);

\end{tikzpicture}

\end{center}
\caption{A map on the annulus} 
\label{Annulusmap}
\end{figure}

\medskip

Another way in which $\M^{\rm med}$ could be edge-transitive while $\M$ is not is if $\M$ is self-dual and ${\rm Aut}\,\M$ has two orbits on the edges of $\M^{\rm med}$, transposed by an automorphism of $\M^{\rm med}$ induced by the self-duality of $\M$.
%Alternatively, show that ${\rm Aut}\,\M$ and have the same orbits on edges of $\M^{\rm med}$, so one is transitive on them if and only if the other is? Do this by identifying edges of $\M^{\rm med}$ with incident vertex-face pairs of $\M$, and noting the self-duality preserves these?]
We would need $\overline\M=\M/{\rm Aut}\,\M$ to be self-dual, with two edges and two `corners' of faces, i.e.~two cycles of $r_1$ on flags, so it would have at most four flags. For example, we could take $\overline\M$ to be the map on the annulus with one boundary vertex, a loop going round that boundary component, and a free edge across to the other boundary component, so that there is a single face, as in Figure~\ref{Annulusmap}. (This the quotient of the torus map $\{4,4\}_{1,1}$ by the reflection in an edge.) Thus $|\Phi|=4$ with $r_0=(1,4)$, $r_1=(1,2)(3,4)$ and $r_2=(2,3)$.

\medskip

\noindent {\bf Problem}: Find examples of maps $\M$ for which $\overline\M$ has his form.

%Sections on Stability and Smallest in class oimitted

%%%%%%%%%%%%%%%%%%%%%
%%%%%%%%%%%%%%%%%%%%%

\bigskip

\noindent{\bf Acknowledgement} The author is grateful to Marston Conder, David Singerman, Jozef \v Sir\'a\v n and Tom Tucker for very helpful comments about maps, to David Craven, Kay Magaard and Chris Parker for many useful hints about finite simple groups, and to Dimitri Leemans for sharing details of his work with Martin Liebeck on automorphism groups of chiral polyhedra. He is also grateful to the organisers of SCDO 2016, Queenstown NZ, for the invitation to give a talk from which this paper devloped, and to the organisers and participants of the BIRS workshop on Symmetries of Discrete Structures in Geometry, Oaxaca 2017, where further progress on edge-transitive maps was made.

\end{document}